\documentclass[10pt, oneside, a4paper]{article}
\usepackage{amsmath, amsthm, amssymb, amsfonts}
\usepackage{geometry}
\usepackage{graphicx}
\usepackage{natbib}
\usepackage{setspace}
\usepackage{algorithm}
\usepackage{algorithmicx}
\usepackage{algpseudocode}
\usepackage{longtable}
\usepackage{caption}

\newtheorem{theorem}{Theorem}[section]
\newtheorem{lemma}[theorem]{Lemma}
\newtheorem{proposition}[theorem]{Proposition}

\theoremstyle{definition}
\newtheorem{definition}{Definition}[section]
\newtheorem{example}{Example}[section]
\newtheorem{assumption}{Assumption}
\newtheorem{remark}[theorem]{Remark}

\newcommand{\Rel}{\mathbb{R}}

\input{tcilatex}

\begin{document}

\title{Valid Heteroskedasticity Robust Testing\thanks{%
Financial support of the second author by the Program of Concerted Research
Actions (ARC) of the Universit\'{e} libre de Bruxelles in an early stage of
this project is gratefully acknowledged. We thank two referees, a co-editor,
and the editor for helpful comments. Address correspondence to Benedikt P%
\"{o}tscher, Department of Statistics, University of Vienna, A-1090
Oskar-Morgenstern Platz 1. E-Mail: benedikt.poetscher@univie.ac.at.}}
\author{Benedikt M. P\"{o}tscher and David Preinerstorfer \\
Department of Statistics, University of Vienna\\
SEPS-SEW, University of St.~Gallen}
\date{Preliminary draft: October 2017\\
First version: April 2021\\
Second version: November 2022\\
This version: July 2023}
\maketitle

\begin{abstract}
Tests based on heteroskedasticity robust standard errors are an important
technique in econometric practice. Choosing the right critical value,
however, is not simple at all: conventional critical values based on
asymptotics often lead to severe size distortions; and so do existing
adjustments including the bootstrap. To avoid these issues, we suggest to
use smallest size-controlling critical values, the generic existence of
which we prove in this article for the commonly used test statistics.
Furthermore, sufficient and often also necessary conditions for their
existence are given that are easy to check. Granted their existence, these
critical values are the canonical choice: larger critical values result in
unnecessary power loss, whereas smaller critical values lead to
over-rejections under the null hypothesis, make spurious discoveries more
likely, and thus are invalid. We suggest algorithms to numerically determine
the proposed critical values and provide implementations in accompanying
software. Finally, we numerically study the behavior of the proposed testing
procedures, including their power properties.
\end{abstract}

\section{Introduction\label{Intro}}

Testing hypotheses on the parameters in a regression model with potentially
heteroskedastic errors is an important problem in econometrics and
statistics; see \cite{Mackinnon2013} for a recent survey. Since the
classical $t$-statistic ($F$-statistic, respectively) is not pivotal, or
asymptotically pivotal, in such a case in general, even under Gaussianity of
the errors, so-called heteroskedasticity robust (aka heteroskedasticity
consistent) modifications of these test statistics have been proposed, which
are asymptotically standard normally (chi-square, respectively) distributed
under the null. These modifications date back to \cite{E63,E67}, see also 
\cite{H77}, and have later been popularized in econometrics by \cite{W80}
with great success (see \cite{Mackinnon2013}). Unfortunately, it turned out
that tests obtained from these heteroskedasticity robust test statistics by
relying on critical values derived from the respective asymptotic null
distributions have a tendency to overreject the null hypothesis in finite
samples (and thus are invalid), especially so if the design matrix contains
high-leverage points; see, e.g., \cite{MacW85}, \cite{DavidsonMacKinnon1985}%
, and \cite{CheshJewitt1987}. One factor contributing to this overrejection
tendency is a downward bias in the covariance matrix estimators used in
these test statistics, see \cite{CheshJewitt1987}. In an attempt to reduce
the overrejection problem, variants of the before-mentioned
heteroskedasticity robust test statistics (often denoted by HC1 through HC4,
with HC0 denoting the original proposal) have been considered; see \cite{H77}%
, \cite{MacW85}, and \cite{Crib2004}.\footnote{%
For a recent contribution geared towards high-dimensional models see \cite%
{Catt}.} These variants rescale the least-squares residuals before computing
the covariance matrix estimator employed in the construction of the test
statistic. According to simulation studies reported in, e.g., \cite%
{DavidsonMacKinnon1985} and \cite{Crib2004}, these modifications, especially
HC3 and HC4, seem to ameliorate the overrejection problem to some extent,
but do not eliminate it. Further numerical results are provided in \cite%
{CheshAust_1991}, see also \cite{Chesh_1989}. Numerical results in Section %
\ref{numerical} confirm these observations. Variants of HC0-HC3, denoted by
HC0R-HC3R, obtained by using restricted instead of unrestricted
least-squares residuals in the computation of the covariance matrix
estimators employed by the various test statistics (the restriction alluded
to being the restriction defining the null hypothesis) have been introduced
in \cite{DavidsonMacKinnon1985}. In their simulation experiments, this
typically leads to tests that do not overreject, but that may substantially
underreject; see also the simulation results in \cite{Godfrey2006}, who
additionally also considers HC4R. However, as will be shown in Section \ref%
{numerical}, also these tests are in general not immune to (sometimes
substantial) overrejection.

Note that, under the typical assumptions used in the literature, all the
modifications of HC0 discussed so far have the same asymptotic distribution
as HC0, and thus the same critical value as for HC0 (obtained from the
asymptotic null distribution) is also used for these modifications in the
before mentioned literature. Sometimes small-sample adjustments to the
asymptotic critical values are attempted by using the quantiles from a $%
t_{d} $-distribution rather than from the asymptotic normal distribution,
where the degrees of freedom $d$ are either set to $n-k$ ($n$ and $k$
denoting sample size and number of regressors, respectively), or are
obtained through proposals set down by \cite{Satterth} or \cite{BellMcCa};
see also \cite{Imbkoles2016}. While these adjustments can lead to
improvements, numerical results presented in Section \ref{numerical} show
that these adjustments are also not able to solve the overrejection problem
in general. An alternative approach is to use bootstrap methods to compute
critical values for the test statistics HC0-HC4 or HC0R-HC4R. The relevant
literature is reviewed in \cite{PPBoot}, and it is shown that such methods
are again not immune to the overrejection problem in general.\footnote{%
Another possibility is to use Edgeworth expansions to find better critical
values, see \cite{Rothenberg1988} for the case of the HC0 test statistic and 
\cite{DavidsonMacKinnon1985} for the HC0R test statistic. Simulation results
in \cite{MacW85} and \cite{DavidsonMacKinnon1985} indicate that this does
not work too well in practice. Of course, such expansions could also be
worked out for the other versions of the test statistics mentioned, but this
does not seem to have been pursued in the literature.} A referee has pointed
out the recent papers by \cite{Chuetal2021} and \cite{Hansen2021}, both of
which propose a testing procedure that can be viewed as a parametric
bootstrap method. No theoretical justification is given in those papers. In
fact, as we show in Appendix \ref{app_G}, the proposed procedures can be
considerably oversized, a feature that can already be seen to some extent in
the numerical results given in \cite{Chuetal2021} and \cite{Hansen2021}.

A result by \cite{BakiSzek2005} needs to be mentioned here which states that
-- in the special case of testing a hypothesis on the location parameter of
a heteroskedastic location model with errors that are Gaussian or scale
mixtures thereof -- the classical two-sided $t$-test (with the usual
critical value) has null rejection probability not exceeding the nominal
significance level under \emph{any} form of heteroskedasticity (for a
certain range of significance levels); see \cite{IbragMuell2010} for more
discussion. \cite{IbragMuell2016}, extending a result in \cite%
{MickeyBrown1966}, provide a related result in the case of the comparison of
two heteroskedastic populations; see also \cite{Bakirov98}. Section \ref%
{sec_BS} provides some more discussion. We note that all results mentioned
in that section are applicable only to testing \emph{certain} \emph{scalar}
linear contrasts.

Except for the \cite{BakiSzek2005} result and the variations discussed in
Section \ref{sec_BS}, which apply only to quite special situations like,
e.g., the heteroskedastic location model, none of the methods discussed so
far comes with a theoretical result implying that their associated (finite
sample) null rejection probabilities are guaranteed not to exceed the
nominal significance level whatever the form of heteroskedasticity may be.%
\footnote{%
In the special case where the number of restrictions tested equals the
number of regression parameters,\ \cite{DavidsFlach2008} have a result which
implies that certain wild bootstrap-based heteroskedasticity robust tests
have size equal to the nominal significance level (and hence do not
overreject) in finite samples. We note that this result in \cite%
{DavidsFlach2008} is not entirely correct as stated, but needs some
amendments and corrections; see \cite{PPBoot}.} In fact, it transpires from
the preceding discussion and the numerical results in Section \ref{numerical}
that for \emph{any} of these methods instances of testing problems can be
found for which the method in question overrejects substantially. Therefore,
it is imperative to be able to find size-controlling critical values for the
test statistics considered, i.e., critical values such that the resulting
worst-case rejection probability under the null hypothesis does not exceed
the nominal significance level. We shall hence pursue in this paper the
construction of size-controlling critical values for the test statistics
HC0-HC4, HC0R-HC4R, as well as for (two variants of) the classical (i.e.,
uncorrected) $F$-statistic (including the absolute value of the $t$%
-statistic as a special case).

In the present paper we consider classes of test statistics that contain the
before mentioned heteroskedasticity robust test statistics as special cases
and show under which conditions -- and how -- a critical value can be found
such that the resulting test is guaranteed to have size less than or equal
to $\alpha $, the prescribed significance level.\footnote{%
A less principled attempt at finding a valid test in a \emph{given} testing
problem (i.e., for given design matrix and restriction to be tested) could
consist in the practitioner studying the size of a handful of tests
(obtained from a few of the above mentioned test statistics in conjunction
with a few of the proposed critical values) by means of an extensive Monte
Carlo study and in hoping that one of the test procedures emerges from this
study as valid for the particular testing problem at hand. Besides being a
numerically costly procedure, it does not come with any guarantee of success.%
} It turns out that the conditions for size controllability are broadly
satisfied; in particular, for the commonly used test statistics they are
satisfied \emph{generically} in a sense made precise further below.

We want to emphasize that the existence of size-controlling critical values
for heteroskedasticity robust test statistics is \emph{not} a trivial
matter, as it has been shown in \cite{PP2016}, Section 4, that there are
cases where the size of such tests is always one, \emph{regardless} of the
choice of critical value; see also the discussion in Proposition \ref%
{rem_necessity} further below. And even in cases where size control is
possible by an appropriate choice of critical value, the standard critical
values proposed in the literature (including the small-sample adjustments
discussed above) are \emph{not} guaranteed to deliver size control; in fact,
they may fail to do so by a considerable margin (i.e., they are much too
small to control size at the desired level) as shown in Section \ref%
{numerical}. Our theoretical results also show the existence of a computable
"threshold" $C^{\ast }$, say, such that any critical value $C$ satisfying $%
C<C^{\ast }$ necessarily leads to a test with size $1$; see Proposition \ref%
{rem_C*}. Since $C^{\ast }$ is not difficult to compute, it can be used as a
simple check to weed out unsuitable proposals for critical values.

Apart from avoiding overrejection by construction, the use of smallest
size-controlling, rather than conventional, critical values offers also
advantages in terms of power in instances where conventional critical values
lead to underrejection (i.e., lead to a worst-case rejection probability
under the null hypothesis less than the nominal significance level) as is
sometimes the case; see Sections \ref{sec:trivial_het_tilde} and \ref%
{sec:powernum}. In fact, once one has decided on a test statistic to be used
for the given null hypothesis, using the smallest size-controlling critical
value (provided it exists) is obviously the optimal way to proceed.

We also discuss how the critical values that lead to size control can be
determined numerically and provide the \textsf{R}-package \textbf{hrt }(\cite%
{hrt}) for their computation. The usefulness of the proposed algorithms and
their implementation in the \textsf{R}-package are illustrated numerically
on some testing problems in Section \ref{numerical}. In particular, we
compare tests obtained from various of the above mentioned test statistics
when used with smallest size-controlling critical values in terms of their
power functions. The package \textbf{hrt} also contains a routine for
determining the size of a test obtained from a user-supplied critical value.
It is important to note that if in a particular application one uses the
observed value of the test statistic as the user-supplied critical value in
this routine, this routine actually returns a \textquotedblleft valid
p-value\textquotedblright\ in the following sense: Checking whether or not
this \textquotedblleft p-value\textquotedblright\ is smaller than the
prescribed significance level $\alpha $ is equivalent to checking whether or
not the observed value of the test statistic is larger than or equal to the
smallest size-controlling critical value. Note that the former check avoids
the need to actually compute the smallest size-controlling critical value,
which is advantageous from a computational point of view. See Section \ref%
{algor} for more details.

In the paper we work under a Gaussianity assumption. We stress, however,
that this assumption is mainly made for convenience of presentation; as
shown in Section \ref{non-gaussianity}, this assumption can be relaxed
considerably.

While a trivial remark, we would like to note that the size control results
given in this paper can easily be translated into results stating that the
minimal coverage probability of the associated confidence set obtained by
\textquotedblleft inverting\textquotedblright\ the test is not less than the
nominal confidence level.

The paper is organized as follows: After introducing notation and the most
important test statistics in Sections \ref{frame} and \ref{sec_teststatistic}%
, Section \ref{sec_intui} provides some intuition for our size-control
results which are presented in Sections \ref{sec_size_control} and \ref%
{sec_size_control_tilde}, with some further results relegated to Appendix %
\ref{app_d}. Section \ref{Generalizations} discusses ways of relaxing the
underlying assumptions. Possible extensions to other classes of test
statistics are discussed in Section \ref{sec_extensions}, while a few
comments on power are collected in Section \ref{sec:power}. Section \ref%
{numerical} provides the numerical results including a power study, with
some details relegated to Appendix \ref{details}. Section \ref{concl}
concludes. Proofs and some technical results can be found in Appendices \ref%
{app_char}-\ref{app_C}. The algorithms for computing rejection probabilities
(including size) and smallest size-controlling critical values are outlined
in Section \ref{algor}, and are presented in detail in Appendix \ref%
{app:algos}. Appendix \ref{app_G} contains a discussion of \cite{Chuetal2021}
and \cite{Hansen2021}.

\section{Framework\label{frame}}

Consider the linear regression model 
\begin{equation}
\mathbf{Y}=X\beta +\mathbf{U},  \label{lm}
\end{equation}%
where $X$ is a (real) nonstochastic regressor (design) matrix of dimension $%
n\times k$ and where $\beta \in \mathbb{R}^{k}$ denotes the unknown
regression parameter vector. We always assume $\limfunc{rank}(X)=k$ and $%
1\leq k<n$. We furthermore assume that the $n\times 1$ disturbance vector $%
\mathbf{U}=(\mathbf{u}_{1},\ldots ,\mathbf{u}_{n})^{\prime }$ has mean zero
and unknown covariance matrix $\sigma ^{2}\Sigma $, where $\Sigma $ varies
in a user-specified (nonempty) set $\mathfrak{C}$ describing the allowed
forms of heteroskedasticity, with $\mathfrak{C}$ satisfying $\mathfrak{C}%
\subseteq \mathfrak{C}_{Het}$, and where $0<\sigma ^{2}<\infty $ holds ($%
\sigma $ always denoting the positive square root).\footnote{%
Since we are concerned with finite-sample results only, the elements of $%
\mathbf{Y}$, $X$, and $\mathbf{U}$ (and even the probability space
supporting $\mathbf{Y}$ and $\mathbf{U}$) may depend on sample size $n$, but
this will not be expressed in the notation. Furthermore, the obvious
dependence of$\ \mathfrak{C}$ on $n$ will also not be shown in the notation.}
The set $\mathfrak{C}$ will be referred to as the \textquotedblleft
heteroskedasticity model\textquotedblright . Here%
\begin{equation*}
\mathfrak{C}_{Het}=\left\{ \limfunc{diag}(\tau _{1}^{2},\ldots ,\tau
_{n}^{2}):\tau _{i}^{2}>0\text{ for all }i\text{, }\sum_{i=1}^{n}\tau
_{i}^{2}=1\right\} ,
\end{equation*}%
where $\limfunc{diag}(\tau _{1}^{2},\ldots ,\tau _{n}^{2})$ denotes the $%
n\times n$ matrix with diagonal elements given by $\tau _{i}^{2}$. That is,
the errors in the regression model are uncorrelated but can be
heteroskedastic. In particular, if $\mathfrak{C}$ is chosen to be $\mathfrak{%
C}_{Het}$, one allows for heteroskedasticity of completely unknown form. The
normalization condition $\sum_{i=1}^{n}\tau _{i}^{2}=1$ is included here
only in order to guarantee identifiability of $\sigma ^{2}$ and $\Sigma $,
and could be replaced by any other normalization condition such as, e.g., $%
\max \tau _{i}^{2}=1$, or $\tau _{1}^{2}=1$, without affecting the final
results (because any of these normalizations leads to the same overall set
of covariance matrices $\sigma ^{2}\Sigma $ when $\sigma ^{2}$ varies
through the positive real line). Although a trivial observation, we stress
the fact that all conceivable forms of heteroskedasticity, including
parametric ones, can (possibly after normalization) be cast as submodels $%
\mathfrak{C}$ of $\mathfrak{C}_{Het}$.

\emph{Mainly for ease of exposition, we shall maintain in the sequel that
the disturbance vector }$\mathbf{U}$\emph{\ is normally distributed. This
assumption can be substantially relaxed as discussed in Section \ref%
{non-gaussianity}. }The linear model described in (\ref{lm}), together with
the just made Gaussianity assumption on $\mathbf{U}$ and with the given
heteroskedasticity model $\mathfrak{C}$, then induces a collection of
distributions on the Borel-sets of $\mathbb{R}^{n}$, the sample space of $%
\mathbf{Y}$. Denoting a Gaussian probability measure with mean $\mu \in 
\mathbb{R}^{n}$ and (possibly singular) covariance matrix $A$ by $P_{\mu ,A}$%
, the induced collection of distributions is then given by 
\begin{equation}
\left\{ P_{\mu ,\sigma ^{2}\Sigma }:\mu \in \mathrm{\limfunc{span}}%
(X),0<\sigma ^{2}<\infty ,\Sigma \in \mathfrak{C}\right\} ,  \label{lm2}
\end{equation}%
where $\mathrm{\limfunc{span}}(X)$ denotes the column space of $X$. Since
every $\Sigma \in \mathfrak{C}$ is positive definite by assumption, each
element of the set in the previous display is absolutely continuous with
respect to (w.r.t.) Lebesgue measure on $\mathbb{R}^{n}$.

We shall consider the problem of testing a linear (better: affine)
hypothesis on the parameter vector $\beta \in \mathbb{R}^{k}$, i.e., the
problem of testing the null $R\beta =r$ against the alternative $R\beta \neq
r$, where $R$ is a $q\times k$ matrix always of rank $q\geq 1$ and $r\in 
\mathbb{R}^{q}$. Set $\mathfrak{M}=\limfunc{span}(X)$. Define the affine
space 
\begin{equation*}
\mathfrak{M}_{0}=\left\{ \mu \in \mathfrak{M}:\mu =X\beta \text{ and }R\beta
=r\right\}
\end{equation*}%
and let 
\begin{equation*}
\mathfrak{M}_{1}=\left\{ \mu \in \mathfrak{M}:\mu =X\beta \text{ and }R\beta
\neq r\right\} .
\end{equation*}%
Adopting these definitions, this testing problem can then be written more
precisely as 
\begin{equation}
H_{0}:\mu \in \mathfrak{M}_{0},\ 0<\sigma ^{2}<\infty ,\ \Sigma \in 
\mathfrak{C}\quad \text{ vs. }\quad H_{1}:\mu \in \mathfrak{M}_{1},\
0<\sigma ^{2}<\infty ,\ \Sigma \in \mathfrak{C}.  \label{testing problem}
\end{equation}%
With $\mathfrak{M}_{0}^{lin}$ we shall denote the linear space parallel to $%
\mathfrak{M}_{0}$, i.e., $\mathfrak{M}_{0}^{lin}=\mathfrak{M}_{0}-\mu
_{0}=\left\{ X\beta :R\beta =0\right\} $ where $\mu _{0}\in \mathfrak{M}_{0}$%
. Of course, $\mathfrak{M}_{0}^{lin}$ does not depend on the choice of $\mu
_{0}\in \mathfrak{M}_{0}$.

As already mentioned, the assumption of Gaussianity is made mainly for
simplicity of presentation and can be relaxed substantially; see Section \ref%
{non-gaussianity}. The assumption of nonstochastic regressors entails little
loss of generality either, which is important to emphasize: If $X$ is random
and $\mathbf{U}$ is conditionally on $X$ distributed as $N(0,\sigma
^{2}\Sigma )$, with $\sigma ^{2}=\sigma ^{2}(X)>0$ and $\Sigma =\Sigma
(X)\in \mathfrak{C}_{Het}$, the results of the paper can be applied after
one conditions on $X$ (and a similar statement applies to the
generalizations to non-Gaussianity discussed in Section \ref{non-gaussianity}%
). See Section \ref{stoch_regr} for more discussion and details. For
arguments supporting conditional inference see, e.g., \cite{RO1979}. Note
that such a \textquotedblleft strict exogeneity\textquotedblright\
assumption is quite natural in the situation considered here.

We next collect some further terminology and notation used throughout the
paper. A (nonrandomized) \textit{test} is the indicator function of a
Borel-set $W$ in $\mathbb{R}^{n}$, with $W$ called the corresponding \textit{%
rejection region}. The \textit{size} of such a test (rejection region) is --
as usual -- defined as the supremum over all rejection probabilities under
the null hypothesis $H_{0}$ given in (\ref{testing problem}), i.e., 
\begin{equation*}
\sup_{\mu \in \mathfrak{M}_{0}}\sup_{0<\sigma ^{2}<\infty }\sup_{\Sigma \in 
\mathfrak{C}}P_{\mu ,\sigma ^{2}\Sigma }(W).
\end{equation*}%
In slight abuse of terminology, we shall sometimes refer to this quantity as
`the size of $W$ over $\mathfrak{C}$' when we want to emphasize the r\^{o}le
of $\mathfrak{C}$. Throughout the paper we let $\hat{\beta}(y)=\left(
X^{\prime }X\right) ^{-1}X^{\prime }y$, where $X$ is the design matrix
appearing in (\ref{lm}) and $y\in \mathbb{R}^{n}$. The corresponding
ordinary least-squares (OLS) residual vector is denoted by $\hat{u}(y)=y-X%
\hat{\beta}(y)$ and its elements are denoted by $\hat{u}_{t}(y)$. The
elements of $X$ are denoted by $x_{ti}$, while $x_{t\cdot }$ and $x_{\cdot
i} $ denote the $t$-th row and $i$-th column of $X$, respectively. For $%
\mathcal{A}$ an affine subspace of $\mathbb{R}^{n}$ satisfying $\mathcal{A}%
\subseteq \limfunc{span}(X)$ let $\tilde{\beta}_{\mathcal{A}}(y)$ denote the
restricted least-squares estimator, i.e., $X\tilde{\beta}_{\mathcal{A}}(y)$
solves 
\begin{equation*}
\min_{z\in \mathcal{A}}(y-z)^{\prime }(y-z).
\end{equation*}%
Lebesgue measure on the Borel-sets of $\mathbb{R}^{n}$ will be denoted by $%
\lambda _{\mathbb{R}^{n}}$, whereas Lebesgue measure on an arbitrary affine
subspace $\mathcal{A}$ of $\mathbb{R}^{n}$ (but viewed as a measure on the
Borel-sets of $\mathbb{R}^{n}$) will be denoted by $\lambda _{\mathcal{A}}$,
with zero-dimensional Lebesgue measure being interpreted as point mass. The
set of real matrices of dimension $l\times m$ is denoted by $\mathbb{R}%
^{l\times m}$ (all matrices in the paper will be real matrices) and Lebesgue
measure on this set equipped with its Borel $\sigma $-field is denoted by $%
\lambda _{\mathbb{R}^{l\times m}}$. Let $B^{\prime }$ denote the transpose
of a matrix $B\in \mathbb{R}^{l\times m}$ and let $\mathrm{\limfunc{span}}%
(B) $ denote the subspace in $\mathbb{R}^{l}$ spanned by its columns. For a
symmetric and nonnegative definite matrix $B$ we denote the unique symmetric
and nonnegative definite square root by $B^{1/2}$. For a linear subspace $%
\mathcal{L}$ of $\mathbb{R}^{n}$ we let $\mathcal{L}^{\bot }$ denote its
orthogonal complement and we let $\Pi _{\mathcal{L}}$ denote the orthogonal
projection onto $\mathcal{L}$. The Euclidean norm is denoted by $\left\Vert
\cdot \right\Vert $, but the same symbol is also used to denote a norm of a
matrix. The $j$-th standard basis vector in $\mathbb{R}^{n}$ is written as $%
e_{j}(n)$. Furthermore, we let $\mathbb{N}$ denote the set of all positive
integers. A sum (product, respectively) over an empty index set is to be
interpreted as $0$ ($1$, respectively). Finally, for $\mathcal{A}$ an affine
subspace of $\mathbb{R}^{n}$, let $G(\mathcal{A})$ denote the group of all
affine transformations $y\mapsto \delta (y-a)+a^{\ast }$ where $\delta \in 
\mathbb{R}$, $\delta \neq 0$, and $a$ as well as $a^{\ast }$ are elements of 
$\mathcal{A}$; for more information see Section 5.1 of \cite{PP2016}.

\section{Heteroskedasticity robust test statistics using unrestricted
residuals \label{sec_teststatistic}}

We now introduce two test statistics that will feature prominently in the
following. Variants thereof that use restricted residuals are discussed in
Section \ref{sec_size_control_tilde}. For results pertaining to other
classes of test statistics see Section \ref{sec_extensions}. The test
statistic we shall consider first is a standard heteroskedasticity robust
test statistic frequently encountered in the literature. It is given by%
\begin{equation}
T_{Het}\left( y\right) =\left\{ 
\begin{array}{cc}
(R\hat{\beta}\left( y\right) -r)^{\prime }\hat{\Omega}_{Het}^{-1}\left(
y\right) (R\hat{\beta}\left( y\right) -r) & \text{if }\det \hat{\Omega}%
_{Het}\left( y\right) \neq 0, \\ 
0 & \text{if }\det \hat{\Omega}_{Het}\left( y\right) =0,%
\end{array}%
\right.  \label{T_het}
\end{equation}%
where $\hat{\Omega}_{Het}=R\hat{\Psi}_{Het}R^{\prime }$ and where $\hat{\Psi}%
_{Het}$ is a heteroskedasticity robust estimator as considered in \cite%
{E63,E67}, which later on has found its way into the econometrics literature
(e.g., \cite{W80}). It is of the form%
\begin{equation*}
\hat{\Psi}_{Het}\left( y\right) =(X^{\prime }X)^{-1}X^{\prime }\limfunc{diag}%
\left( d_{1}\hat{u}_{1}^{2}\left( y\right) ,\ldots ,d_{n}\hat{u}%
_{n}^{2}\left( y\right) \right) X(X^{\prime }X)^{-1},
\end{equation*}%
where the constants $d_{i}>0$ sometimes depend on the design matrix. Typical
choices for $d_{i}$ suggested in the literature are $d_{i}=1$, $%
d_{i}=n/(n-k) $, $d_{i}=\left( 1-h_{ii}\right) ^{-1}$, or $d_{i}=\left(
1-h_{ii}\right) ^{-2}$ where $h_{ii}$ denotes the $i$-th diagonal element of
the projection matrix $X(X^{\prime }X)^{-1}X^{\prime }$, see \cite{LE2000}
for an overview. Another suggestion is $d_{i}=\left( 1-h_{ii}\right)
^{-\delta _{i}}$ for $\delta _{i}=\min (nh_{ii}/k,4)$, see \cite{Crib2004}.
For the last three choices of $d_{i}$ just given, we use the convention that
we set $d_{i}=1$ in case $h_{ii}=1$. Note that $h_{ii}=1$ implies $\hat{u}%
_{i}\left( y\right) =0$ for every $y$, and hence it is irrelevant which real
value is assigned to $d_{i}$ in case $h_{ii}=1$.\footnote{%
In fact, $h_{ii}=1$ is equivalent to $\hat{u}_{i}\left( y\right) =0$ for
every $y$, each of which in turn is equivalent to $e_{i}(n)\in $ $\limfunc{%
span}(X)$.} The five examples for the weights $d_{i}$ just given correspond
to what is often called HC0-HC4 weights in the literature.

In conjunction with the test statistic $T_{Het}$, we shall consider the
following mild assumption, which is Assumption 3 in \cite{PP2016}. As
discussed further below, this assumption is in a certain sense unavoidable
when using $T_{Het}$. It furthermore also entails that our choice of
assigning $T_{Het}\left( y\right) $ the value zero in case $\hat{\Omega}%
_{Het}\left( y\right) $ is singular has no import on the probabilistic
results of the paper (because of Lemma \ref{lem_B}(c) below and absolute
continuity of the measures $P_{\mu ,\sigma ^{2}\Sigma }$).

\begin{assumption}
\label{R_and_X}Let $1\leq i_{1}<\ldots <i_{s}\leq n$ denote all the indices
for which $e_{i_{j}}(n)\in \limfunc{span}(X)$ holds where $e_{j}(n)$ denotes
the $j$-th standard basis vector in $\mathbb{R}^{n}$. If no such index
exists, set $s=0$. Let $X^{\prime }\left( \lnot (i_{1},\ldots i_{s})\right) $
denote the matrix which is obtained from $X^{\prime }$ by deleting all
columns with indices $i_{j}$, $1\leq i_{1}<\ldots <i_{s}\leq n$ (if $s=0$ no
column is deleted). Then $\limfunc{rank}\left( R(X^{\prime }X)^{-1}X^{\prime
}\left( \lnot (i_{1},\ldots i_{s})\right) \right) =q$ holds.
\end{assumption}

Observe that this assumption only depends on $X$ and $R$ and hence can be
checked. Obviously, a simple sufficient condition for Assumption \ref%
{R_and_X} to hold is that $s=0$ (i.e., that $e_{j}(n)\notin \limfunc{span}%
(X) $ for all $j$), a generically satisfied condition. Furthermore, we
introduce the matrix%
\begin{eqnarray}
B(y) &=&R(X^{\prime }X)^{-1}X^{\prime }\limfunc{diag}\left( \hat{u}%
_{1}(y),\ldots ,\hat{u}_{n}(y)\right)  \notag \\
&=&R(X^{\prime }X)^{-1}X^{\prime }\limfunc{diag}\left( e_{1}^{\prime }(n)\Pi
_{\limfunc{span}(X)^{\bot }}y,\ldots ,e_{n}^{\prime }(n)\Pi _{\limfunc{span}%
(X)^{\bot }}y\right) .  \label{B_matrix}
\end{eqnarray}%
The facts collected in the subsequent lemma, which is taken from \cite%
{PPBoot} (but see also Lemma 4.1 in \cite{PP2016} and Lemma 5.18 in \cite%
{PP3}), will be used in the sequel.

\begin{lemma}
\label{lem_B}(a) $\hat{\Omega}_{Het}\left( y\right) $ is nonnegative
definite for every $y\in \mathbb{R}^{n}$.

(b) $\hat{\Omega}_{Het}\left( y\right) $ is singular (zero, respectively) if
and only if $\limfunc{rank}\left( B(y)\right) <q$ ($B(y)=0$, respectively).

(c) The set $\mathsf{B}$ given by $\left\{ y\in \mathbb{R}^{n}:\limfunc{rank}%
\left( B(y)\right) <q\right\} $ (or in view of (b) equivalently given by $%
\{y\in \mathbb{R}^{n}:\det (\hat{\Omega}_{Het}\left( y\right) )=0\}$) is
either a $\lambda _{\mathbb{R}^{n}}$-null set or the entire sample space $%
\mathbb{R}^{n}$. The latter occurs if and only if Assumption \ref{R_and_X}
is violated (in which case the test based on $T_{Het}$ becomes trivial, as
then $T_{Het}$ is identically zero).

(d) Under Assumption \ref{R_and_X}, the set $\mathsf{B}$ is a finite union
of proper linear subspaces of $\mathbb{R}^{n}$; in case $q=1$, $\mathsf{B}$
is even a proper linear subspace itself.\footnote{%
If Assumption \ref{R_and_X} is violated, $\mathsf{B}$ equals $\mathbb{R}^{n}$
by Part (c).}

(e) $\mathsf{B}$ is a closed set and contains $\limfunc{span}(X)$\textsf{. }%
Furthermore, $\mathsf{B}$ is $G(\mathfrak{M})$-invariant and, in particular, 
$\mathsf{B}+\limfunc{span}(X)=\mathsf{B}$ holds.
\end{lemma}

In light of Part (c) of the lemma, we see that Assumption \ref{R_and_X} is a
natural and unavoidable condition if one wants to obtain a sensible test
from $T_{Het}$.\footnote{%
If this assumption is violated then $T_{Het}$ is identically zero, an
uninteresting trivial case.} Furthermore, note that, if $\mathsf{B}=\limfunc{%
span}(X)$ is true, then Assumption \ref{R_and_X} must be satisfied (since $%
\limfunc{span}(X)$ is a $\lambda _{\mathbb{R}^{n}}$-null set due to the
maintained assumption $k<n$). As shown in Lemma A.3 in \cite{PP3}, for any
given restriction matrix $R$, the relation $\mathsf{B}=\limfunc{span}(X)$
holds generically in various universes of design matrices. For later use we
also mention that under Assumption \ref{R_and_X} the test statistic $T_{Het}$
is continuous at every $y\in \mathbb{R}^{n}\backslash \mathsf{B}$.\footnote{%
If Assumption \ref{R_and_X} is violated, then $T_{Het}$ is constant equal to
zero, and hence is trivially continuous everywhere.}

Next, we also consider the classical (i.e., uncorrected) F-test statistic,
i.e.,%
\begin{equation}
T_{uc}(y)=\left\{ 
\begin{array}{cc}
(R\hat{\beta}\left( y\right) -r)^{\prime }\left( \hat{\sigma}^{2}(y)R\left(
X^{\prime }X\right) ^{-1}R^{\prime }\right) ^{-1}(R\hat{\beta}\left(
y\right) -r) & \text{if }y\notin \limfunc{span}(X), \\ 
0 & \text{if }y\in \limfunc{span}(X),%
\end{array}%
\right.   \label{T_uncorr}
\end{equation}%
where $\hat{\sigma}^{2}(y)=\hat{u}\left( y\right) ^{\prime }\hat{u}\left(
y\right) /(n-k)\geq 0$ (which vanishes if and only if $y\in \limfunc{span}(X)
$). Our choice to set $T_{uc}(y)=0$ for $y\in \limfunc{span}(X)$ again has
no import on the probabilistic results in the paper, since $\limfunc{span}(X)
$ is a $\lambda _{\mathbb{R}^{n}}$-null set as a consequence of the
maintained assumption that $k<n$ (and since the measures $P_{\mu ,\sigma
^{2}\Sigma }$ are absolutely continuous). For reasons of comparability with (%
\ref{T_het}) we have chosen not to normalize the numerator in (\ref{T_uncorr}%
) by $q$, the number of restrictions to be tested, as is often done in the
definition of the classical F-test statistic. This also has no import on the
results as the factor $1/q$ can be absorbed into the critical value. For
later use we also mention that the test statistic $T_{uc}$ is continuous at
every $y\in \mathbb{R}^{n}\backslash \limfunc{span}(X)$.

\begin{remark}
\label{rem:GM0}(i) The test statistics $T_{Het}$ as well as $T_{uc}$ are $G(%
\mathfrak{M}_{0})$-invariant as is easily seen (with the respective
exceptional sets $\mathsf{B}$ and $\limfunc{span}(X)$ being $G(\mathfrak{M})$%
-invariant).

(ii) Both statistics actually belong to the class of nonsphericity-corrected
F-type test statistics in the sense of Section 5.4 in \cite{PP2016}
(terminology being somewhat unfortunate in case of $T_{uc}$ as no correction
for the non-sphericity is applied in this case). See Remark \ref{F-type} in
Appendix \ref{app_B} for more discussion.
\end{remark}

\begin{remark}
\label{obvious}For later use we note the following: Suppose $(R,r)$ and $(%
\bar{R},\bar{r})$ are both of dimension $q\times (k+1)$ and have $\limfunc{%
rank}(R)=\limfunc{rank}(\bar{R})=q.$ (i) Then $(R,r)$ and $(\bar{R},\bar{r})$
give rise to the same set $\mathfrak{M}_{0}$, and thus to the same testing
problem (\ref{testing problem}), if and only if $(AR,Ar)=(\bar{R},\bar{r})$
holds for a nonsingular $q\times q$ matrix $A$. (ii) The test statistics $%
T_{Het}$ and $T_{uc}$ remain the same whether they are computed using $(R,r)$
or $(\bar{R},\bar{r})$ provided $(AR,Ar)=(\bar{R},\bar{r})$ holds for a
nonsingular $q\times q$ matrix $A$. [To see this note that the respective
exceptional sets $\mathsf{B}$ and $\limfunc{span}(X)$ are the same
irrespective of whether $(R,r)$ or $(\bar{R},\bar{r})$ is used, and that $A$
cancels out in the respective quadratic forms appearing in the definitions
of the test statistics.]
\end{remark}

\section{Some intuition on why conventional critical values can lead to
overrejection\label{sec_intui}}

We begin the heuristic discussion by considering the testing problem (\ref%
{testing problem}) with heteroskedasticity model $\mathfrak{C}=\mathfrak{C}%
_{Het}$ (i.e., heteroskedasticity of unknown form). Let $T$ stand for any of
the test statistics introduced in Section \ref{sec_teststatistic}, with
rejection occurring whenever $T\geq C$, $C$ a critical value.\footnote{%
In case of $T=T_{Het}$ Assumption \ref{R_and_X} is supposed to hold.}$^{%
\text{,}}$\footnote{%
The discussion similarly applies to the test statistics introduced in
Section \ref{sec_size_control_tilde}.}. For simplicity of presentation we
assume $r=0$. As discussed in Section \ref{Intro}, basing the test on the
conventional critical value $C_{\chi ^{2}(q),0.05}$ (the 95\% quantile of a
chi-square distribution with $q$ degrees of freedom) often leads to
substantial overrejection, i.e., the size of the test (over $\mathfrak{C}%
_{Het}$) is substantially larger than the desired value $\alpha =0.05$. One
mechanism leading to such overrejection is constituted by a concentration
phenomenon discussed at some length in \cite{PP2016}: In the present
situation, the distribution $P_{0,\sigma ^{2}\Sigma }$ \textquotedblleft
concentrates\textquotedblright\ on a so-called concentration subspace (given
by $\limfunc{span}(e_{i}(n))$) when $\Sigma $ is \textquotedblleft
close\textquotedblright\ to one of the singular matrices $%
e_{i}(n)e_{i}(n)^{\prime }$.\footnote{%
There are also other concentration subspaces in the present situation which
we can ignore for the heuristic discussion.} In such a case, depending on
the design matrix $X$ and the hypothesis given by $(R,r)$, the concentration
space may fall into the rejection region $\{T\geq C_{\chi ^{2}(q),0.05}\}$,
leading to a rejection probability close to one, and thus much larger than $%
\alpha =0.05$.\footnote{%
This is an oversimplified description ignoring some technical details.} Even
if the concentration subspace $\limfunc{span}(e_{i}(n))$ is not contained in
the rejection region, but is sufficiently close to it, a considerable
portion of the mass of $P_{0,\sigma ^{2}\Sigma }$ may nevertheless fall into
the rejection region if $\Sigma $ is close to, but not too close to $%
e_{i}(n)e_{i}(n)^{\prime }$. This again leads to a relatively large
rejection probability. Overrejection will often be especially pronounced if
certain high-leverage points are present in the design matrix.\footnote{%
We note, however, that there are testing problems (e.g., testing the mean in
a heteroskedastic location model using the test statistic $T_{uc}$) for
which the text-book critical values obtained under homoskedasticity are
actually valid, see \cite{BakiSzek2005}. The reason is that the
\textquotedblleft worst case\textquotedblright\ distribution in this case
corresponds to homoskedasticity.}

In order to obtain a test that has size controlled by $\alpha $ (i.e., size $%
\leq \alpha $) in situations as just described, the rejection region $%
\{T\geq C_{\chi ^{2}(q),0.05}\}$ has to be narrowed down, i.e., $C_{\chi
^{2}(q),0.05}$ has to be replaced by a suitably larger critical value $C$.
Whether or not this can successfully be accomplished by a (finite) $C$, is a
non-trivial question, the answer depending on whether or not all possible
concentration subspaces can be made to fall outside of the rejection region $%
\{T\geq C\}$ by an appropriate choice of $C$ larger than $C_{\chi
^{2}(q),0.05}$. Sufficient conditions when this is possible are provided in
Theorems \ref{Hetero_Robust} and \ref{Hetero_Robust_tilde}. Note that, in
such a situation, the resulting size-controlling critical values $C$ are
then necessarily larger than $C_{\chi ^{2}(q),0.05}$.

In light of the preceding discussion, a natural question is whether or not
imposing a heteroskedasticity model more narrow than $\mathfrak{C}_{Het}$
such as, e.g.,%
\begin{equation*}
\mathfrak{C}_{Het,\tau _{\ast }}=\left\{ \limfunc{diag}\left( \tau
_{1}^{2},\ldots ,\tau _{n}^{2}\right) \in \mathfrak{C}_{Het}:\tau
_{i}^{2}\geq \tau _{\ast }^{2}\text{ for all }i\right\} ,
\end{equation*}%
where $\tau _{\ast }$, $0<\tau _{\ast }<n^{-1/2}$, is a pre-specified
constant set by the user, would mitigate the failure of conventional
critical values. Indeed, under the heteroskedasticity model $\mathfrak{C}%
_{Het,\tau _{\ast }}$ extreme concentration effects leading to rejection
probabilities (arbitrarily) close to one cannot occur, and it is possible to
prove that size-controlling critical values always exist when $\mathfrak{C}%
_{Het,\tau _{\ast }}$ is used, see Appendix \ref{app_d}. Unfortunately,
however, this does \emph{not} imply that conventional critical values such
as $C_{\chi ^{2}(q),0.05}$ will work. In fact, the size over $\mathfrak{C}%
_{Het,\tau _{\ast }}$ of tests using the critical value $C_{\chi
^{2}(q),0.05}$ can still be considerably larger than $\alpha $: To see this,
observe that the sets $\mathfrak{C}_{Het,\tau _{\ast }}$are an increasing
sequence of sets as $\tau _{\ast }\downarrow 0$, the union of which is $%
\mathfrak{C}_{Het}$. Consequently, if $\tau _{\ast }$ is small, the size
over $\mathfrak{C}_{Het,\tau _{\ast }}$ will be close to the size over $%
\mathfrak{C}_{Het}$, and thus the former will be much larger than $\alpha $
in case the latter is so. As a consequence, also in case of the more narrow
heteroskedasticity model $\mathfrak{C}_{Het,\tau _{\ast }}$ size-controlling
critical values larger than $C_{\chi ^{2}(q),0.05}$ will have to be used in
such a case. Furthermore, the bound $\tau _{\ast }$ has to be decided upon
prior to the data analysis and is thus part of \emph{modeling }the form of
heteroskedasticity. It is difficult to see how one would come up with a
reasonable value of $\tau _{\ast }$ in practice: If $\tau _{\ast }$ is
chosen to be small, this may result in a heteroskedasticity model under
which the test based on $C_{\chi ^{2}(q),0.05}$ is still plagued by
overrejection as just discussed, while choosing $\tau _{\ast }$ large will
typically not be defensible as it presumes considerable knowledge about the
admissible forms of heteroskedasticity.

\section{ Size control results for $T_{Het}$ and $T_{uc}$ when $\mathfrak{C}=%
\mathfrak{C}_{Het}$\label{sec_size_control}}

We introduce the following notation: For a given linear subspace $\mathcal{L}
$ of $\mathbb{R}^{n}$ we define the set of indices $I_{0}(\mathcal{L})$ via 
\begin{equation*}
I_{0}(\mathcal{L})=\left\{ i:1\leq i\leq n,e_{i}(n)\in \mathcal{L}\right\} .
\end{equation*}%
We set $I_{1}(\mathcal{L})=\left\{ 1,\ldots ,n\right\} \backslash I_{0}(%
\mathcal{L})$. Clearly, $\func{card}(I_{0}(\mathcal{L}))\leq \dim (\mathcal{L%
})$ holds. In particular, if $\dim (\mathcal{L})<n$ holds (which, in
particular, is so in the leading case $\mathcal{L}=\mathfrak{M}_{0}^{lin}$,
since $\dim (\mathfrak{M}_{0}^{lin})=k-q<n$), then $\func{card}(I_{0}(%
\mathcal{L}))<n$, and thus $\func{card}(I_{1}(\mathcal{L}))\geq 1$.

We have the following size control result for $T_{uc}$ as well as for $%
T_{Het}$ over the heteroskedasticity model $\mathfrak{C}_{Het}$ (more
precisely, over the null hypothesis $H_{0}$ described in (\ref{testing
problem}) with $\mathfrak{C}=\mathfrak{C}_{Het}$). Note that $\mathfrak{C}%
_{Het}$ is the largest possible heteroskedasticity model and reflects
complete ignorance about the form of heteroskedasticity.

\begin{theorem}
\label{Hetero_Robust} (a) For every $0<\alpha <1$ there exists a real number 
$C(\alpha )$ such that%
\begin{equation}
\sup_{\mu _{0}\in \mathfrak{M}_{0}}\sup_{0<\sigma ^{2}<\infty }\sup_{\Sigma
\in \mathfrak{C}_{Het}}P_{\mu _{0},\sigma ^{2}\Sigma }(T_{uc}\geq C(\alpha
))\leq \alpha  \label{size-control_Het_uncorr}
\end{equation}%
holds, provided that 
\begin{equation}
e_{i}(n)\notin \func{span}(X)\text{ \ \ for every \ }i\in I_{1}(\mathfrak{M}%
_{0}^{lin}).  \label{non-incl_Het_uncorr}
\end{equation}%
Furthermore, under condition (\ref{non-incl_Het_uncorr}), even equality can
be achieved in (\ref{size-control_Het_uncorr}) by a proper choice of $%
C(\alpha )$, provided $\alpha \in (0,\alpha ^{\ast }]\cap (0,1)$ holds,
where $\alpha ^{\ast }=\sup_{C\in (C^{\ast },\infty )}\sup_{\Sigma \in 
\mathfrak{C}_{Het}}P_{\mu _{0},\Sigma }(T_{uc}\geq C)$ is positive and where 
$C^{\ast }=\max \{T_{uc}(\mu _{0}+e_{i}(n)):i\in I_{1}(\mathfrak{M}%
_{0}^{lin})\}$ for $\mu _{0}\in \mathfrak{M}_{0}$ (with neither $\alpha
^{\ast }$ nor $C^{\ast }$ depending on the choice of $\mu _{0}\in \mathfrak{M%
}_{0}$).

(b) Suppose Assumption \ref{R_and_X} is satisfied.\footnote{%
Condition (\ref{non-incl_Het}) clearly implies that the set $\mathsf{B}$ is
a proper subset of $\mathbb{R}^{n}$ (as $\limfunc{card}(I_{1}(\mathfrak{M}%
_{0}^{lin}))\geq 1$) and thus implies Assumption \ref{R_and_X}. Hence, we
could have dropped this assumption from the formulation of the theorem. For
clarity of presentation we have, however, chosen to explicitly mention
Assumption \ref{R_and_X}. A similar remark applies to some of the other
results given below and will not be repeated.} Then for every $0<\alpha <1$
there exists a real number $C(\alpha )$ such that%
\begin{equation}
\sup_{\mu _{0}\in \mathfrak{M}_{0}}\sup_{0<\sigma ^{2}<\infty }\sup_{\Sigma
\in \mathfrak{C}_{Het}}P_{\mu _{0},\sigma ^{2}\Sigma }(T_{Het}\geq C(\alpha
))\leq \alpha  \label{size-control_Het}
\end{equation}%
holds, provided that 
\begin{equation}
e_{i}(n)\notin \mathsf{B}\text{ \ \ for every \ }i\in I_{1}(\mathfrak{M}%
_{0}^{lin}).  \label{non-incl_Het}
\end{equation}%
Furthermore, under condition (\ref{non-incl_Het}), even equality can be
achieved in (\ref{size-control_Het}) by a proper choice of $C(\alpha )$,
provided $\alpha \in (0,\alpha ^{\ast }]\cap (0,1)$ holds, where now $\alpha
^{\ast }=\sup_{C\in (C^{\ast },\infty )}\sup_{\Sigma \in \mathfrak{C}%
_{Het}}P_{\mu _{0},\Sigma }(T_{Het}\geq C)$ is positive and where $C^{\ast
}=\max \{T_{Het}(\mu _{0}+e_{i}(n)):i\in I_{1}(\mathfrak{M}_{0}^{lin})\}$
for $\mu _{0}\in \mathfrak{M}_{0}$ (with neither $\alpha ^{\ast }$ nor $%
C^{\ast }$ depending on the choice of $\mu _{0}\in \mathfrak{M}_{0}$).

(c) Under the assumptions of Part (a) (Part (b), respectively) implying
existence of a critical value $C(\alpha )$ satisfying (\ref%
{size-control_Het_uncorr}) ((\ref{size-control_Het}), respectively), a
smallest critical value, denoted by $C_{\Diamond }(\alpha )$, satisfying (%
\ref{size-control_Het_uncorr}) ((\ref{size-control_Het}), respectively)
exists for every $0<\alpha <1$. And $C_{\Diamond }(\alpha )$ corresponding
to Part (a) (Part (b), respectively) is also the smallest among the critical
values leading to equality in (\ref{size-control_Het_uncorr}) ((\ref%
{size-control_Het}), respectively) whenever such critical values exist.
[Although $C_{\Diamond }(\alpha )$ corresponding to Part (a) and (b),
respectively, will typically be different, we use the same symbol.]\footnote{%
Cf.~also Appendix \ref{useful}.}
\end{theorem}

We see from the theorem that the condition for size control of $T_{Het}$ ($%
T_{uc}$, respectively) over $\mathfrak{C}_{Het}$, i.e., condition (\ref%
{non-incl_Het}) ((\ref{non-incl_Het_uncorr}), respectively), only depends on 
$X$ and $R$; in particular, in case of $T_{Het}$, it does not depend on how
the weights $d_{i}$ figuring in the definition of $T_{Het}$ have been chosen
(note that the set $\mathsf{B}$ only depends on $X$ and $R$). Moreover, the
sufficient conditions for size control are generically satisfied in the
universe of all $n\times k$ design matrices $X$ (of rank $k$), see Example %
\ref{ex_general} and the attending discussion further below. Furthermore, it
is plain that the size-controlling critical values $C(\alpha )$ in Theorem %
\ref{Hetero_Robust} will depend on the choice of test statistic as well as
on the testing problem at hand. More concretely, the size-controlling
critical values in Part (b) of the theorem thus depend only on $X$, $R$, and 
$r$, as well as on the choice of weights $d_{i}$, whereas in Part (a) the
dependence is only on $X$, $R$, and $r$. We do not show these dependencies
in the notation. In fact, as discussed in Remark \ref{rem_indep_r} below, it
turns out that the size-controlling critical values in both cases actually
do \emph{not} depend on the value of $r$ at all (provided the weights $d_{i}$
are not allowed to depend on $r$ in case of $T_{Het}$). Similarly, it is
easy to see that $C^{\ast }$ and $\alpha ^{\ast }$ in Theorem \ref%
{Hetero_Robust} do not depend on $r$ (under the same provision as before in
case of $T_{Het}$).

Another observation is that any critical value delivering size control over $%
\mathfrak{C}_{Het}$ also delivers size control over \emph{any} other
heteroskedasticity model $\mathfrak{C}$ since $\mathfrak{C}\subseteq 
\mathfrak{C}_{Het}$. Of course, for such a $\mathfrak{C}$ even smaller
critical values (than needed for $\mathfrak{C}_{Het}$) may already suffice
for size control. Also note that sufficient conditions implying size control
over $\mathfrak{C}_{Het}$ may be more restrictive than sufficient conditions
implying only size control over a smaller heteroskedasticity model $%
\mathfrak{C}$. For size control results tailored to such smaller models $%
\mathfrak{C}$ see Appendix \ref{app_d}.

In light of the results of \cite{CheshJewitt1987} and \cite{Chesh_1989}, it
is useful to interpret the sufficient conditions for size control, i.e., (%
\ref{non-incl_Het_uncorr}) and (\ref{non-incl_Het}), in terms of
high-leverage points. First, note that $e_{i}(n)\in \limfunc{span}(X)$ is
equivalent to $h_{ii}=1$, which corresponds to the $i$-th observation being
an \textquotedblleft extreme high-leverage point\textquotedblright . Hence, (%
\ref{non-incl_Het_uncorr}) is equivalent to $h_{ii}<1$ for every $i\in 
\mathcal{I}_{1}(\mathfrak{M}_{0}^{lin})$. In other words, the condition for
a size-controlling critical value to exist in Part (a) of Theorem \ref%
{Hetero_Robust} requires that none of the indices in $\mathcal{I}_{1}(%
\mathfrak{M}_{0}^{lin})$ corresponds to an extreme high-leverage point. [It
is interesting to observe that all indices in $\mathcal{I}_{0}(\mathfrak{M}%
_{0}^{lin})$ (note that this set may be empty) correspond to extreme
high-leverage points.] Hence, for the condition in (\ref{non-incl_Het_uncorr}%
) \emph{not} to be satisfied, not only must extreme high-leverage points be
present, but the lever needs to be of a particular type depending on the
hypothesis given by $(R,r)$ (namely, it must have $i\in \mathcal{I}_{1}(%
\mathfrak{M}_{0}^{lin})$). Second, note that a sufficient, but not
necessary, condition for (\ref{non-incl_Het_uncorr}) is $h_{ii}<1$ for $%
i=1,\ldots ,n$. Sufficiency is obvious from the preceding discussion. That
the condition is not necessary can be seen from Example \ref{ex_b} further
below. Finally, condition (\ref{non-incl_Het}) implies condition (\ref%
{non-incl_Het_uncorr}) (since $\limfunc{span}(X)\subseteq \mathsf{B}$), and
hence implies $h_{ii}<1$ for every $i\in \mathcal{I}_{1}(\mathfrak{M}%
_{0}^{lin})$. The converse is not always true: even $h_{ii}<1$ for every $%
i=1,\ldots ,n$ does not guarantee (\ref{non-incl_Het}) to be satisfied, see
Example \ref{ex_k_pop} further below. However, \emph{generically }(\ref%
{non-incl_Het_uncorr}) and (\ref{non-incl_Het}) coincide (see Lemma A.3 in 
\cite{PP3}), in which case the discussion given above for (\ref%
{non-incl_Het_uncorr}) also applies to (\ref{non-incl_Het}).

\begin{remark}
\label{rem_indep_r}\emph{(Independence of the value of }$r$ \emph{and
implications for confidence sets)} (i) As already noted before, the
sufficient conditions for size control in both parts of Theorem \ref%
{Hetero_Robust} only depend on $X$ and $R$. In particular, they do not
depend on the value of $r$.

(ii) The size of the test based on $T_{uc}$ ($T_{Het}$, respectively) in
Theorem \ref{Hetero_Robust} as well as the size-controlling critical values $%
C(\alpha )$ (for both test statistics) do also not depend on the value of $r$
(provided the weights $d_{i}$ are not allowed to depend on $r$ in case of $%
T_{Het}$). This follows from Lemma 5.15 in \cite{PP3} combined with Remark %
\ref{F-type} in Appendix \ref{app_B}.\footnote{%
For this argument we impose Assumption \ref{R_and_X} in case of $T_{Het}$,
the case where this assumption is violated being trivial.} This observation
is of some importance, as it allows one easily to obtain confidence sets for 
$R\beta $ by \textquotedblleft inverting\textquotedblright\ the test without
the need of recomputing the critical value for every value of $r$.
\end{remark}

\begin{remark}
\label{rem:equiv}\emph{(Some equivalencies)} If the respective smallest
size-controlling critical values are used (provided they exist), the tests
obtained from $T_{Het}$ with the HC0 and the HC1 weights, respectively, are
identical, as these two test statistics differ only by a multiplicative
constant. The same reasoning applies to the test statistics based on the
HC0-HC4 weights, respectively, in case $h_{ii}$ does not depend on $i$.
\end{remark}

\begin{remark}
\label{rem_positiv}\emph{(Positivity of size-controlling critical values) }%
For every $0<\alpha <1$ any $C(\alpha )$ satisfying (\ref%
{size-control_Het_uncorr}) or (\ref{size-control_Het}) is necessarily
positive. To see this observe that $\{T_{uc}\geq C\}=\{T_{Het}\geq C\}=%
\mathbb{R}^{n}$ for $C\leq 0$, since both test statistics are nonnegative
everywhere.
\end{remark}

The next proposition complements Theorem \ref{Hetero_Robust} and provides a
useful lower bound for the size-controlling critical values (other than the
trivial bound given in the preceding remark).

\begin{proposition}
\label{rem_C*}\footnote{%
It is not difficult to show in the context of Parts (a) and (b) of the
proposition that any critical value $C>C^{\ast }$ actually leads to size
less than $1$. This follows from a reasoning similar as in Remark 5.4 of 
\cite{PP3}.}$^{\text{,}}$\footnote{%
If (\ref{non-incl_Het}) in Part (b) of the proposition does not hold, the
conclusion of Part (b) can be shown to continue to hold with $C^{\ast }$ as
defined in Theorem \ref{Hetero_Robust}(b), and also with $C^{\ast }$ as
defined in Lemma 5.11 of \ \cite{PP3} (note that under the assumptions of
Part (b) of the proposition both definitions of $C^{\ast }$ actually
coincide as shown in the proof of Theorem \ref{Hetero_Robust}). [Recall that
under violation of (\ref{non-incl_Het}) size-controlling critical values may
or may not exist.] If Assumption \ref{R_and_X} is not satisfied, then $%
T_{Het}\equiv 0$, and the conclusion of Part (b) holds trivially (as $%
C^{\ast }=0$ with both definitions). If (\ref{non-incl_Het_uncorr}) in Part
(a) of the proposition is not satisfied, then no size-controlling critical
value exists by Proposition \ref{rem_necessity}; hence, the conclusion of
Part (a) holds trivially, again regardless of which of the two definitions
of $C^{\ast }$ is adopted.}(a) Suppose that (\ref{non-incl_Het_uncorr}) is
satisfied. Then any $C(\alpha )$ satisfying (\ref{size-control_Het_uncorr})
necessarily has to satisfy $C(\alpha )\geq C^{\ast }$, where $C^{\ast }$ is
as in Part (a) of Theorem \ref{Hetero_Robust}. In fact, for any $C<C^{\ast }$
we have $\sup_{\Sigma \in \mathfrak{C}_{Het}}P_{\mu _{0},\sigma ^{2}\Sigma
}(T_{uc}\geq C)=1$ for every $\mu _{0}\in \mathfrak{M}_{0}$ and every $%
\sigma ^{2}\in (0,\infty )$.

(b) Suppose that Assumption \ref{R_and_X} and (\ref{non-incl_Het}) are
satisfied. Then any $C(\alpha )$ satisfying (\ref{size-control_Het})
necessarily has to satisfy $C(\alpha )\geq C^{\ast }$, where $C^{\ast }$ is
as in Part (b) of Theorem \ref{Hetero_Robust}. In fact, for any $C<C^{\ast }$
we have $\sup_{\Sigma \in \mathfrak{C}_{Het}}P_{\mu _{0},\sigma ^{2}\Sigma
}(T_{Het}\geq C)=1$ for every $\mu _{0}\in \mathfrak{M}_{0}$ and every $%
\sigma ^{2}\in (0,\infty )$.
\end{proposition}

The preceding observation is useful in two ways: First, critical values
suggested in the literature (such as, e.g., the $(1-\alpha )$-quantile of a
chi-square distribution with $q$ degrees of freedom or critical values
obtained from a degree of freedom adjustment) can immediately be dismissed
if they turn out to be less than $C^{\ast }$, as they then certainly will
not guarantee size control.\footnote{%
In contrast, if the critical value turns out to be larger than or equal to $%
C^{\ast }$, it does \emph{not} follow that size is less than or equal to $%
\alpha $. In fact, substantially oversized tests using a critical value $%
C>C^{\ast }$ are certainly possible; see, e.g., Table \ref{fig:size} and the
pertaining discussion.} We use this line of reasoning in the numerical
results in Section \ref{numerical}. Second, if the \emph{observed} value of
the test statistic $T_{Het}$ ($T_{uc}$, respectively) is less than $C^{\ast
} $, the decision not to reject the null hypothesis can be taken without
further need to compute size-controlling critical values. Note that $C^{\ast
}$ as given in Theorem \ref{Hetero_Robust} is quite easy to compute in any
given application.

\begin{remark}
\label{rem:larger_alphastar}Suppose the assumptions of Part (a) (Part (b),
respectively) of Theorem \ref{Hetero_Robust} are satisfied. Then we know
from that theorem that the size (over $\mathfrak{C}_{Het})$ of $\{T_{uc}\geq
C_{\Diamond }(\alpha )\}$ ($\{T_{Het}\geq C_{\Diamond }(\alpha )\}$,
respectively) equals $\alpha $ provided $\alpha \in (0,\alpha ^{\ast }]\cap
(0,1)$. If now $\alpha ^{\ast }<\alpha <1$, then the size (over $\mathfrak{C}%
_{Het})$ of $\{T_{uc}\geq C_{\Diamond }(\alpha )\}$ ($\{T_{Het}\geq
C_{\Diamond }(\alpha )\}$, respectively) equals $\alpha ^{\ast }$ (where the 
$C_{\Diamond }(\alpha )$'s pertaining to Parts (a) and (b) may be
different). This follows from $C_{\Diamond }(\alpha )\geq C^{\ast }$ (see
Proposition \ref{rem_C*} above) and Remark 5.13(i) in \cite{PP3}.\footnote{%
The assumptions for Part A of Proposition 5.12 in \cite{PP3} required in
Remark 5.13 of that paper are satisfied under the assumptions of Theorem \ref%
{Hetero_Robust} as shown in the proof of Theorem \ref%
{theorem_groupwise_hetero} in Appendix \ref{app_B}. In this proof also the
condition $\lambda _{\mathbb{R}^{n}}(T_{uc}=C^{\ast })=0$ ($\lambda _{%
\mathbb{R}^{n}}(T_{Het}=C^{\ast })=0$, respectively) required in Remark 5.13
of \cite{PP3} is verified.} This argument actually also delivers that $%
C_{\Diamond }(\alpha )=C^{\ast }$ must hold in case $\alpha ^{\ast }<\alpha
<1$.
\end{remark}

We next discuss to what extent the sufficient conditions for size control in
Theorem \ref{Hetero_Robust} are also necessary.

\begin{proposition}
\label{rem_necessity}(a) If (\ref{non-incl_Het_uncorr}) is violated, then $%
\sup_{\Sigma \in \mathfrak{C}_{Het}}P_{\mu _{0},\sigma ^{2}\Sigma
}(T_{uc}\geq C)=1$ for \emph{every} choice of critical value $C$, every $\mu
_{0}\in \mathfrak{M}_{0}$, and every $\sigma ^{2}\in (0,\infty )$ (implying
that size equals $1$ for every $C$). As a consequence, the sufficient
condition for size control (\ref{non-incl_Het_uncorr}) in Part (a) of
Theorem \ref{Hetero_Robust} is also necessary.

(b) Suppose Assumption \ref{R_and_X} is satisfied.\footnote{%
If this assumption is violated then $T_{Het}$ is identically zero, an
uninteresting trivial case.} If (\ref{non-incl_Het_uncorr}) is violated,
then $\sup_{\Sigma \in \mathfrak{C}_{Het}}P_{\mu _{0},\sigma ^{2}\Sigma
}(T_{Het}\geq C)=1$ for \emph{every} choice of critical value $C$, every $%
\mu _{0}\in \mathfrak{M}_{0}$, and every $\sigma ^{2}\in (0,\infty )$
(implying that size equals $1$ for every $C$). [In case $X$ and $R$ are such
that $\mathsf{B}=\limfunc{span}(X)$, conditions (\ref{non-incl_Het_uncorr})
and (\ref{non-incl_Het}) coincide; hence the sufficient condition for size
control (\ref{non-incl_Het}) in Part (b) of Theorem \ref{Hetero_Robust} is
then also necessary in this case.]
\end{proposition}

\begin{remark}
Suppose Assumption \ref{R_and_X} is satisfied. In case $\mathsf{B}\neq 
\limfunc{span}(X)$ and (\ref{non-incl_Het_uncorr}) hold, but (\ref%
{non-incl_Het}) is violated, neither Part (b) of Theorem \ref{Hetero_Robust}
nor Part (b) of Proposition \ref{rem_necessity} apply. We note that there
are instances of this situation (see Example \ref{ex_k_pop}) for which it
can be shown by other methods that $T_{Het}$ is size controllable despite
failure of (\ref{non-incl_Het});\footnote{%
In this example actually $e_{i}(n)\in \mathsf{B}$ holds for all $i=1,\ldots
,n$.} as a consequence, (\ref{non-incl_Het}) is not necessary for (\ref%
{size-control_Het}) in general. We conjecture that there are other instances
of the situation described here where size control is not possible, but we
have not investigated this in any detail. [What can be said in general in
this situation is that the size of the rejection region $\left\{ T_{Het}\geq
C\right\} $ over $\mathfrak{C}_{Het}$ is certainly equal to $1$ for every $%
C<\max \left\{ T_{Het}(\mu _{0}+e_{i}(n)):e_{i}(n)\notin \mathsf{B}\right\} $%
, where we use the convention that this maximum is $-\infty $ in case the
set over which the maximum is taken is empty. This follows from Lemma 4.1 in 
\cite{PP4} with $\mathbb{K}$ equal to the collection $\{\Pi _{(\mathfrak{M}%
_{0}^{lin})^{\bot }}e_{i}(n):e_{i}(n)\notin \mathsf{B}\}$.]
\end{remark}

\begin{remark}
\label{rem_strict_ineq} Let $T$ stand for either $T_{Het}$ or $T_{uc}$, and
suppose that Assumption \ref{R_and_X} is satisfied in case of $T=T_{Het}$:
By Remark \ref{F-type} in Appendix \ref{app_B} and Lemma 5.16 in \cite{PP3}
the rejection regions $\{y:T(y)\geq C\}$ and $\{y:T(y)>C\}$ differ only by a 
$\lambda _{\mathbb{R}^{n}}$-null set. Since the measures $P_{\mu ,\sigma
^{2}\Sigma }$ are absolutely continuous w.r.t.$~\lambda _{\mathbb{R}^{n}}$
when $\Sigma $ is nonsingular, $P_{\mu ,\sigma ^{2}\Sigma }(T\geq C)=P_{\mu
,\sigma ^{2}\Sigma }(T>C)$ then follows, and hence the results in this
section given for rejection probabilities $P_{\mu ,\sigma ^{2}\Sigma }(T\geq
C)$ apply to rejection probabilities $P_{\mu ,\sigma ^{2}\Sigma }(T>C)$
equally well (under the above provision in case of $T=T_{Het}$). A similar
remark applies to the results in Appendix \ref{sec_size_control_2}.
\end{remark}

\subsection{Some examples}

We illustrate Theorem \ref{Hetero_Robust} and Proposition \ref{rem_necessity}
with a few examples.

\begin{example}
\label{ex_general}(i) Suppose the design matrix satisfies $e_{i}(n)\notin 
\func{span}(X)$ for \emph{every} $1\leq i\leq n$ (which will typically be
the case). Then obviously the sufficient condition (\ref{non-incl_Het_uncorr}%
) is satisfied (in fact, for every choice of $\mathfrak{M}_{0}$, i.e., for
every choice of restriction to be tested). And the sufficient condition (\ref%
{non-incl_Het}) is also satisfied provided $\mathsf{B}=\func{span}(X)$.

(ii) Suppose the design matrix $X$ and the restriction $R$ are such that $%
e_{i}(n)\notin \mathsf{B}$ for \emph{every} $1\leq i\leq n$. Then the
sufficient condition (\ref{non-incl_Het}) is clearly satisfied.
\end{example}

This example shows, in particular, that the sufficient conditions for size
control are generically satisfied in the universe of all $n\times k$ design
matrices $X$ (of rank $k$). Given the example, this is obvious for $T_{uc}$;
and it follows for $T_{Het}$ by additionally noting that, for every given
choice of restriction to be tested, the relation $\mathsf{B}=\func{span}(X)$
holds generically in the universe of all $n\times k$ design matrices $X$ (of
rank $k$); see Lemma A.3 in \cite{PP3}. The next example discusses the case
where a standard basis vector is among the regressors.

\begin{example}
\label{ex_b}Suppose that $e_{1}(n)$ is the first column of $X$ and that $%
e_{i}(n)\notin \func{span}(X)$ for every $2\leq i\leq n$. Suppose further
that $R$ is of the form $R=(0,\tilde{R})$, where $\tilde{R}$ has dimension $%
q\times (k-1)$. That is, the restriction to be tested does not involve the
coefficient of the first regressor. Then it is easy to see that (\ref%
{non-incl_Het_uncorr}) is satisfied and size control for $T_{uc}$ is thus
possible. If also $\mathsf{B}=\func{span}(X)$ holds, then the same is true
for (\ref{non-incl_Het}) and $T_{Het}$. [In case $R$ is not as above, but
has a nonzero first coordinate, then it is easy to see that $1\in I_{1}(%
\mathfrak{M}_{0}^{lin})$, and hence (\ref{non-incl_Het_uncorr}) is violated.
It follows from Proposition \ref{rem_necessity} that the rejection region $%
\left\{ T_{uc}\geq C\right\} $ indeed has size $1$ for every choice of
critical value $C$ when $\mathfrak{C}_{Het}$ is the heteroskedasticity
model; and the same is true for $T_{Het}$, provided Assumption \ref{R_and_X}
is satisfied.\footnote{%
If Assumption \ref{R_and_X} is violated then $T_{Het}$ is identically zero,
an uninteresting trivial case.}]
\end{example}

We continue with a few more examples where $X$ has a particular structure.

\begin{example}
\label{ex_loc}\emph{(Heteroskedastic location model) }Suppose $k=1$, $%
x_{t1}=1$ for all $t$, $q=1$, $R=1$, and $r\in \mathbb{R}$. The
heteroskedasticity model is given by $\mathfrak{C}_{Het}$. Then the
conditions for size control in both parts of Theorem \ref{Hetero_Robust} are
satisfied (since it is easy to see that $\mathsf{B}$ coincides with $%
\limfunc{span}(X)$ and that Assumption \ref{R_and_X} is satisfied). Note
also that in this example $T_{Het}$ and $T_{uc}$ actually coincide in case $%
d_{i}=n/(n-1)$ for all $i$, i.e., if the HC1, HC2, or HC4 weights are used,
and differ only by a multiplicative constant if the HC0 or HC3 weights are
employed; in particular, all these test statistics give rise to one and the
same test if the respective smallest size-controlling critical values are
used (cf.~Remark \ref{rem:equiv}).\footnote{%
In fact, more is true in the location model: The test statistics $\tilde{T}%
_{Het}$ using the HC0R-HC4R weights (defined in Section \ref%
{sec_size_control_tilde} below) all coincide (cf.~Footnote \ref{FNk=q}), and
they also coincide with $\tilde{T}_{uc}$ (also defined in Section \ref%
{sec_size_control_tilde} below). Perusing the connection between $\tilde{T}%
_{uc}$ and $T_{uc}$ established in Section \ref{sec:trivial_uc_tilde}, we
can then even conclude that all the test statistics $T_{uc}$, $T_{Het}$ with
HC0-HC4 weights, $\tilde{T}_{uc}$, and $\tilde{T}_{Het}$ with HC0R-HC4R
weights give rise to (essentially) the same test, provided the respective
smallest size-controlling critical values are used.} Furthermore, note that
the here observed size controllability is in line with results in \cite%
{BakiSzek2005} stating that, for a certain range of significance levels $%
\alpha $, the usual critical values obtained from an $F_{1,n-1}$%
-distribution actually can be used as size-controlling critical values $%
C(\alpha )$ for the test statistic $T_{uc}$ (in fact, these are then the
smallest size-controlling critical values $C_{\Diamond }(\alpha )$).
\end{example}

The subsequent example is closely related to the Behrens-Fisher problem, see
Remark \ref{FB} in Appendix \ref{sec_size_control_2}.

\begin{example}
\label{ex_fish_behr}\emph{(Comparing the means of two heteroskedastic
groups) }Consider the problem of testing the equality of the means of two
independent normal populations where the variances of each item may be
different, even within a group. In our framework this corresponds to the
case $k=2$, $x_{t1}=1$ for $1\leq t\leq n_{1}$, $x_{t1}=0$ for $n_{1}<t\leq
n_{1}+n_{2}=n$, $x_{t2}=1-x_{t1}$, and $R=(1,-1)$ with $r=0$. The
heteroskedasticity model is then again $\mathfrak{C}_{Het}$. We first assume
that $n_{i}\geq 2$ holds for $i=1,2$. Note that in the present context $%
T_{uc}$ is nothing else than the square of the two-sample t-statistic that
uses a pooled variance estimator, and that $T_{Het}$ is the square of the
two-sample t-statistic that uses appropriate variance estimators from each
group (the particular form of \ the variance estimator being determined by
the choice of $d_{i}$). Now, $e_{i}(n)\notin \func{span}(X)$ for \emph{every}
$1\leq i\leq n$ holds, and hence $T_{uc}$ is size controllable (cf. Example %
\ref{ex_general}(i)). This is in line with results in \cite{Bakirov98},
cf.~also Section \ref{sec_BS}. Furthermore, it is obvious that Assumption %
\ref{R_and_X} is satisfied (as $s=0$) and a simple calculation shows that $%
B(y)=\hat{u}(y)^{\prime }A$, where $A$ is a diagonal matrix with $%
a_{ii}=n_{1}^{-1}$ for $1\leq i\leq n_{1}$ and $a_{ii}=-n_{2}^{-1}$ else.
This shows that the set $\mathsf{B}$ coincides with $\func{span}(X)$.
Consequently, also $T_{Het}$ is size controllable (again cf. Example \ref%
{ex_general}(i)). We also note here that the observed size controllability
of $T_{Het}$ is in line with results in \cite{IbragMuell2016} stating that
for a certain range of significance levels $\alpha $ and group sizes $n_{i}$
the usual critical values obtained from an $F_{1,\min (n_{1},n_{2})-1}$%
-distribution actually can be used as size-controlling critical values $%
C(\alpha )$ for the test statistic $T_{Het}$ in case $d_{i}$ is set equal to 
$\left( 1-h_{ii}\right) ^{-1}$; in fact, they are then the smallest
size-controlling critical values, cf.~the discussion preceding Theorem 1 in 
\cite{IbragMuell2016}. In the rather uninteresting case $n_{1}=1$ and $%
n_{2}\geq 2$, it is easy to see that Assumption \ref{R_and_X} is satisfied
and that the size of both tests equals $1$ for all choices of critical
values in view of Proposition \ref{rem_necessity}, since $e_{1}(n)\in \func{%
span}(X)$ and $1\in I_{1}(\mathfrak{M}_{0}^{lin})=\{1,\ldots ,n\}$. The same
is true if $n_{1}\geq 2$ and $n_{2}=1$. [The remaining and uninteresting
case $n_{1}=n_{2}=1$ falls outside of our framework since we always require $%
n>k$.]
\end{example}

The next example is an extension of the previous problem to the case of more
than two groups. An interesting phenomenon occurs here: The sufficient
conditions for size control of $T_{Het}$ given in Theorem \ref{Hetero_Robust}
are \emph{violated}, but size controllability can \emph{nevertheless} be
established by additional arguments. Hence, this example provides an
instance where the conditions in Part (b) of Theorem \ref{Hetero_Robust} are
not necessary.

\begin{example}
\label{ex_k_pop}\emph{(Comparing the means of }$k$\emph{\ heteroskedastic
groups) }We are given $k$ integers $n_{j}\geq 1$ with $\sum_{j=1}^{k}n_{j}=n$
describing group sizes where $k\geq 3$ holds. The regressors $x_{ti}$ for $%
1\leq i\leq k$ indicate group membership, i.e., they satisfy $x_{ti}=1$ for $%
\sum_{j=1}^{i-1}n_{j}<t\leq \sum_{j=1}^{i}n_{j}$ and $x_{ti}=0$ otherwise.
The heteroskedasticity model is given by $\mathfrak{C}_{Het}$. We are
interested in testing $\beta _{1}=\ldots =\beta _{k}$. We thus may choose
the $(k-1)\times k$ restriction matrix $R$ with $j$-th row $(1,0,\ldots
0,-1,\ldots ,0)$ where the entry $-1$ is at position $j+1$. Of course, $%
q=k-1 $ and $r=0$ hold. We first consider the case where $n_{j}\geq 2$ for
all $j$. Then clearly $k<n$ is satisfied. With regard to $T_{uc}$ we see
immediately that $e_{i}(n)\notin \func{span}(X)$ for every $1\leq i\leq n$
follows (since $n_{j}\geq 2$ for all $j$) and thus the sufficient condition (%
\ref{non-incl_Het_uncorr}) for size control of $T_{uc}$ is satisfied.
Turning to $T_{Het}$, it is easy to see that Assumption \ref{R_and_X} is
satisfied (since $s=0$ in view of $n_{j}\geq 2$). Furthermore, the $j$-th
row of $R(X^{\prime }X)^{-1}X^{\prime }$ is seen to be of the form%
\begin{equation*}
(n_{1}^{-1},\ldots ,n_{1}^{-1},0,\ldots ,0,-n_{j+1}^{-1},\ldots
,-n_{j+1}^{-1},0\ldots ,0),
\end{equation*}%
from which it follows that%
\begin{equation}
R(X^{\prime }X)^{-1}X^{\prime }\limfunc{diag}(d_{1}\hat{u}_{1}^{2}(y),\ldots
,d_{n}\hat{u}_{n}^{2}(y))X(X^{\prime }X)^{-1}R=S_{1}\iota \iota ^{\prime }+%
\limfunc{diag}(S_{2},\ldots ,S_{k}),  \label{B}
\end{equation}%
where $\iota $ is the $(k-1)$-dimensional vector with entries all equal to $%
1 $ and where $S_{j}=n_{j}^{-2}\tsum_{t}d_{t}\hat{u}_{t}^{2}(y)=n_{j}^{-2}%
\tsum_{t}d_{t}(y_{t}-\bar{y}_{(j)})^{2}$ with the summation index $t$
running over all elements in the $j$-th group, and where $\bar{y}_{(j)}$ is
the mean in group $j$. From (\ref{B}) it is not difficult to verify that the
set $\mathsf{B}$ is given by%
\begin{equation*}
\mathsf{B}=\tbigcup_{i,j=1,i\neq j}^{k}\left\{ y\in \mathbb{R}%
^{n}:S_{i}(y)=S_{j}(y)=0\right\} =\tbigcup_{i,j=1,i\neq j}^{k}\limfunc{span}%
\left( x_{\cdot i},x_{\cdot j},\left\{ e_{l}(n):x_{li}=x_{lj}=0\right\}
\right) .
\end{equation*}%
Note that $\mathsf{B}$ is not a linear space and is strictly larger than $%
\limfunc{span}(X)$. The set $\mathfrak{M}_{0}^{lin}$ is given by the span of
the vector $e=(1,1,\ldots ,1)^{\prime }$. Hence, $I_{1}(\mathfrak{M}%
_{0}^{lin})=\left\{ 1,\ldots ,n\right\} $. Since $e_{i}(n)\in \mathsf{B}$
holds for every $i$, we conclude that the sufficient condition (\ref%
{non-incl_Het}) for size control of $T_{Het}$ is \emph{not} satisfied and
hence Part (b) of Theorem \ref{Hetero_Robust} does \emph{not} apply.
However, it can be shown by additional arguments, see Proposition \ref%
{prop_k_pop} in Appendix\ \ref{app_B}, that $T_{Het}$ is nevertheless size
controllable, i.e., that (\ref{size-control_Het}) holds.\footnote{%
A smallest size-controlling critical value then also exists in view of
Appendix \ref{useful}.} Next, in the case where $n_{j}=1$ for some $j$, but
not for all $j$, Proposition \ref{rem_necessity} shows that the size of the
test based on $T_{uc}$ equals $1$ for all choices of critical values, since
then for some $i$ the standard basis vector $e_{i}(n)$ is one of the
regressors and thus we have $e_{i}(n)\in \func{span}(X)$ and $i\in I_{1}(%
\mathfrak{M}_{0}^{lin})=\{1,\ldots ,n\}$. For $T_{Het}$ the same is true if $%
n_{j}=1$ holds for exactly one $j$ (because of Part (b) of Proposition \ref%
{rem_necessity} and since then Assumption \ref{R_and_X} is satisfied as is
easily seen); in case $n_{j}=1$ is true for (at least) two, but not all,
values of $j$, $T_{Het}$ is identically zero (as then Assumption \ref%
{R_and_X} is violated), and thus is size-controllable in a trivial way. [The
remaining and uninteresting case $n_{j}=1$ for all $j$ falls outside of our
framework since we always require $n>k$.]
\end{example}

We close this section by one more example. Again, the sufficient conditions
in Part (b) of Theorem \ref{Hetero_Robust} fail to hold, but additional
arguments based on Example \ref{ex_loc} establish size controllability of
the test based on $T_{Het}$.

\begin{example}
\label{ex_fish_behr_2}Consider again the situation of Example \ref%
{ex_fish_behr}, except that now $R=I_{2}$, the $2\times 2$ identity matrix,
(and again $r=0$). Then $q=k=2$ holds. Consider first the case where $%
n_{i}\geq 2$ for $i=1,2$. Condition (\ref{non-incl_Het_uncorr}) is then
obviously satisfied, and hence $T_{uc}$ is size controllable. We next turn
to $T_{Het}$. Since $\mathfrak{M}_{0}^{lin}=\left\{ 0\right\} $ we have $%
I_{1}(\mathfrak{M}_{0}^{lin})=\{1,\ldots ,n\}$. Furthermore, simple
computations show that Assumption \ref{R_and_X} is satisfied and that 
\begin{equation*}
\mathsf{B}=\limfunc{span}\left( x_{\cdot 1},\left\{ e_{i}(n):i>n_{1}\right\}
\right) \cup \limfunc{span}\left( x_{\cdot 2},\left\{ e_{i}(n):i\leq
n_{1}\right\} \right) .
\end{equation*}%
Obviously, the sufficient condition (\ref{non-incl_Het}) for size control of 
$T_{Het}$ is violated. Nevertheless, $T_{Het}$ is size controllable by the
following argument:\footnote{%
A smallest size-controlling critical value then also exists in view of
Appendix \ref{useful}.} Simple computations show that $%
T_{Het}(y)=T_{1}(y)+T_{2}(y)$ for $y\notin \mathsf{B}$, where $%
T_{1}(y)=n_{1}^{2}\hat{\beta}_{1}^{2}(y)/\sum_{t=1}^{n_{1}}d_{t}\hat{u}%
_{t}^{2}(y)$ and $T_{2}(y)=n_{2}^{2}\hat{\beta}_{2}^{2}(y)/%
\sum_{t=n_{1}+1}^{n}d_{t}\hat{u}_{t}^{2}(y)$. [If the denominator in the
formula for $T_{i}(y)$ is zero for some $y\in \mathbb{R}^{n}$, we define $%
T_{i}(y)$ as zero.] Since $\mathsf{B}$ is a $\lambda _{\mathbb{R}^{n}}$-null
set, $P_{0,\sigma ^{2}\Sigma }(T_{Het}\geq C)\leq P_{0,\sigma ^{2}\Sigma
}(T_{1}\geq C/2)+P_{0,\sigma ^{2}\Sigma }(T_{2}\geq C/2)$ for $C>0$. Now, it
is easy to see that $P_{0,\sigma ^{2}\Sigma }(T_{i}\geq C/2)$ for $i=1,2$
coincides with the null rejection probability of a test for the mean in a
heteroskedastic location model (based on a test statistic of the form (\ref%
{T_het})). However, as shown in Example \ref{ex_loc}, such a test is size
controllable. [In the case $n_{1}=1$ and $n_{2}\geq 2$ (or vice versa)
condition (\ref{non-incl_Het_uncorr}) is violated and the rejection region $%
\left\{ T_{uc}\geq C\right\} $ has size $1$ for every $C$; furthermore,
Assumption \ref{R_and_X} is violated, and hence $T_{Het}$ is identically
zero. The case $n_{1}=n_{2}=1$ falls outside of our framework as then $k=n$.]
\end{example}

In Appendix \ref{app_B} we discuss yet another example where the sufficient
condition of Part (b) of Theorem \ref{Hetero_Robust} fails, but
size-controllability can nevertheless be established.

\subsection{Some variations on \protect\cite{BakiSzek2005}\label{sec_BS}}

(i) As noted in \cite{IbragMuell2010}, testing a hypothesis regarding a 
\emph{scalar} linear contrast in a heteroskedastic (Gaussian) linear
regression model more general than a location model can often be converted
to a testing problem in a heteroskedastic (Gaussian) location model by
suitably dividing the data into subgroups and by considering groupwise
least-squares estimators, thus making it amenable to the \cite{BakiSzek2005}
result mentioned in Section \ref{Intro}. However, this introduces additional
questions such as how to divide up the data. In any case, this approach is
limited to testing hypotheses on \emph{scalar} linear contrasts. It also
requires that the linear contrast subject to test is \emph{estimable} in
each subgroup.

(ii) In case the linear contrast subject to test is \emph{not} estimable in
each subgroup, but can be written as the difference of two linear contrasts
where the first contrast is estimable in the first $G_{1}$ groups whereas
the second contrast is estimable in the last $G_{2}$ groups (where we
consider a total of $G_{1}+G_{2}$ groups), \cite{IbragMuell2016} point out
that the problem can be converted into the problem of comparing two
heteroskedastic (Gaussian) populations. Now, for such a two-sample
comparison problem \cite{Bakirov98} shows for a certain two-sample $t$%
-statistic (the square of which is $T_{uc}$, cf.~Example \ref{ex_fish_behr}
above) how -- in the presence of heteroskedasticity -- size-controlling
critical values can be constructed by appropriately transforming quantiles
of a $t$-distribution; this result imposes conditions which entail that the
nominal significance level $\alpha $ must be quite small (requiring $\alpha $
not to exceed $0.01$ for many group sizes, and often to be considerably
smaller). This somewhat limits the applicability of Bakirov's result. Thus 
\cite{IbragMuell2016} go on to consider another two-sample $t$-statistic
(the square of which is $T_{Het}$ with $d_{i}=\left( 1-h_{ii}\right) ^{-1}$,
cf.~Example \ref{ex_fish_behr} above) and -- extending a result in \cite%
{MickeyBrown1966} -- provide a \cite{BakiSzek2005}-type result, i.e., they
show that the $(1-\alpha /2)$-quantile of a $t$-distribution with degrees of
freedom equal to the smaller of the two sample sizes minus $1$ provides the
smallest size-controlling critical value (for the two-sided test) even under
heteroskedasticity.\footnote{%
In the balanced case (i.e., if the two samples have the same cardinality)
the test statistic considered in \cite{Bakirov98} actually coincides with
the test statistic in \cite{IbragMuell2016}.} This result holds under
certain conditions on the sample sizes and only for small $\alpha $, but,
e.g., allows for the choice $\alpha =0.05$. [We note here that the
description of \cite{Bakirov98}'s result in \cite{IbragMuell2016} is
inaccurate in that a certain transformation of the critical value is being
ignored.]

(iii) In the problem of comparing two heteroskedastic (Gaussian) populations
based on samples of equal size (\textquotedblleft balanced
design\textquotedblright ) one can -- instead of using the two-sample $t$%
-test statistics considered in \cite{Bakirov98} and \cite{IbragMuell2016} --
employ the Bartlett test statistic, which simply is the usual $t$-test
statistic computed from the differences between the observations in the two
samples.\footnote{%
Certainly, there is some arbitrariness in how the observations are being
\textquotedblleft paired\textquotedblright .} An advantage of this approach
is that the original \cite{BakiSzek2005} result is directly applicable, and
there is no need to resort to the results described in (ii).

(iv) Another quite special case that can be brought under the realm of the 
\cite{BakiSzek2005} result is a heteroskedastic (Gaussian) regression model
with only one regressor that never takes the value zero. Dividing the $t$-th
equation in the regression model by $x_{t}$, converts this into a
heteroskedastic location problem.

(v) The results in (i)-(iv) immediately also apply if the errors in the
regression are distributed as scale mixtures of Gaussians (cf.~also Section %
\ref{non-gaussianity}).

\section{Results for heteroskedasticity robust test statistics using
restricted residuals\label{sec_size_control_tilde}}

In this section we consider two further test statistics which are versions
of $T_{Het}$ and $T_{uc}$ with the only difference that the covariance
matrix estimators used are based on restricted -- instead of unrestricted --
residuals. The first one of these test statistics has been suggested in the
literature, e.g., in \cite{DavidsonMacKinnon1985}. We thus define%
\begin{equation}
\tilde{T}_{Het}\left( y\right) =\left\{ 
\begin{array}{cc}
(R\hat{\beta}\left( y\right) -r)^{\prime }\tilde{\Omega}_{Het}^{-1}\left(
y\right) (R\hat{\beta}\left( y\right) -r) & \text{if }\det \tilde{\Omega}%
_{Het}\left( y\right) \neq 0, \\ 
0 & \text{if }\det \tilde{\Omega}_{Het}\left( y\right) =0,%
\end{array}%
\right.  \label{T_Het_tilde}
\end{equation}%
where $\tilde{\Omega}_{Het}=R\tilde{\Psi}_{Het}R^{\prime }$ and where $%
\tilde{\Psi}_{Het}$ is given by%
\begin{equation*}
\tilde{\Psi}_{Het}\left( y\right) =(X^{\prime }X)^{-1}X^{\prime }\limfunc{%
diag}\left( \tilde{d}_{1}\tilde{u}_{1}^{2}\left( y\right) ,\ldots ,\tilde{d}%
_{n}\tilde{u}_{n}^{2}\left( y\right) \right) X(X^{\prime }X)^{-1},
\end{equation*}%
where the constants $\tilde{d}_{i}>0$ sometimes depend on the design matrix
and on the restriction matrix $R$. Here $\tilde{u}\left( y\right) =y-X\tilde{%
\beta}_{\mathfrak{M}_{0}}(y)=\Pi _{(\mathfrak{M}_{0}^{lin})^{\bot }}(y-\mu
_{0})$, where the last expression does not depend on the choice of $\mu
_{0}\in \mathfrak{M}_{0}$, and where $\tilde{u}_{t}\left( y\right) $ denotes
the $t$-th component of $\tilde{u}\left( y\right) $. Typical choices for $%
\tilde{d}_{i}$ are $\tilde{d}_{i}=1$, $\tilde{d}_{i}=n/(n-(k-q))$, $\tilde{d}%
_{i}=(1-\tilde{h}_{ii})^{-1}$, or $\tilde{d}_{i}=(1-\tilde{h}_{ii})^{-2}$
where $\tilde{h}_{ii}$ denotes the $i$-th diagonal element of the projection
matrix $\Pi _{\mathfrak{M}_{0}^{lin}}$, see, e.g., \cite%
{DavidsonMacKinnon1985}. Another suggestion is $\tilde{d}_{i}=(1-\tilde{h}%
_{ii})^{-\tilde{\delta}_{i}}$ for $\tilde{\delta}_{i}=\min (n\tilde{h}%
_{ii}/(k-q),4)$ with the convention that $\tilde{\delta}_{i}=0$ if $k=q$.%
\footnote{%
Note that in case $k=q$ we have $\tilde{h}_{ii}=0$, and hence $\tilde{d}%
_{i}=1$ regardless of our convention for $\tilde{\delta}_{i}$.} For the last
three choices of $\tilde{d}_{i}$ just given we use the convention that we
set $\tilde{d}_{i}=1$ in case $\tilde{h}_{ii}=1$. Note that $\tilde{h}%
_{ii}=1 $ implies $\tilde{u}_{i}\left( y\right) =0$ for every $y$, and hence
it is irrelevant which real value is assigned to $\tilde{d}_{i}$ in case $%
\tilde{h}_{ii}=1$.\footnote{%
In fact, $\tilde{h}_{ii}=1$ is equivalent to $\tilde{u}_{i}\left( y\right)
=0 $ for every $y$, each of which in turn is equivalent to $e_{i}(n)\in $ $%
\mathfrak{M}_{0}^{lin}$.} The five examples for the weights $\tilde{d}_{i}$
just given correspond to what is often called HC0R-HC4R weights in the
literature.\footnote{\label{FNk=q}In the case $k=q$ the HC0R-HC4R weights
all coincide ($\tilde{d}_{i}=1$ for every $i$), and hence result in the same
test statistic.}

The subsequent assumption ensures that the set of $y$'s for which $\tilde{%
\Omega}_{Het}\left( y\right) $ is singular is a Lebesgue null set, implying
that our choice of assigning $\tilde{T}_{Het}\left( y\right) $ the value
zero in case $\tilde{\Omega}_{Het}\left( y\right) $ is singular has no
import on the probabilistic results of the paper (as the measures $P_{\mu
,\sigma ^{2}\Sigma }$ are absolutely continuous). Also, as discussed further
below, the assumption is in a certain sense unavoidable when using $\tilde{T}%
_{Het}$.

\begin{assumption}
\label{R_and_X_tilde}Let $1\leq i_{1}<\ldots <i_{s}\leq n$ denote all the
indices for which $e_{i_{j}}(n)\in \mathfrak{M}_{0}^{lin}$ holds where $%
e_{j}(n)$ denotes the $j$-th standard basis vector in $\mathbb{R}^{n}$. If
no such index exists, set $s=0$. Let $X^{\prime }\left( \lnot (i_{1},\ldots
i_{s})\right) $ denote the matrix which is obtained from $X^{\prime }$ by
deleting all columns with indices $i_{j}$, $1\leq i_{1}<\ldots <i_{s}\leq n$
(if $s=0$ no column is deleted). Then $\limfunc{rank}\left( R(X^{\prime
}X)^{-1}X^{\prime }\left( \lnot (i_{1},\ldots i_{s})\right) \right) =q$
holds.
\end{assumption}

Observe that this assumption only depends on $X$ and $R$ and hence can be
checked. Obviously, a simple sufficient condition for Assumption \ref%
{R_and_X_tilde} to hold is that $s=0$ (i.e., that $e_{j}(n)\notin \mathfrak{M%
}_{0}^{lin}$ for all $j$), a generically satisfied condition. Furthermore,
we introduce the matrix%
\begin{eqnarray}
\tilde{B}(y) &=&R(X^{\prime }X)^{-1}X^{\prime }\limfunc{diag}\left( \tilde{u}%
_{1}(y),\ldots ,\tilde{u}_{n}(y)\right)  \notag \\
&=&R(X^{\prime }X)^{-1}X^{\prime }\limfunc{diag}\left( e_{1}^{\prime }(n)\Pi
_{(\mathfrak{M}_{0}^{lin})^{\bot }}(y-\mu _{0}),\ldots ,e_{n}^{\prime
}(n)\Pi _{(\mathfrak{M}_{0}^{lin})^{\bot }}(y-\mu _{0})\right) .
\label{B_tilde_matrix}
\end{eqnarray}%
Note that this matrix does not depend on the choice of $\mu _{0}\in 
\mathfrak{M}_{0}$. The following lemma collects some important properties of 
$\tilde{\Omega}_{Het}$ and $\mathsf{\tilde{B}}$ (defined in that lemma) and
is reproduced from \cite{PPBoot} for ease of reference.

\begin{lemma}
\label{lem:B_tilde}(a) $\tilde{\Omega}_{Het}\left( y\right) $ is nonnegative
definite for every $y\in \mathbb{R}^{n}$.

(b) $\tilde{\Omega}_{Het}\left( y\right) $ is singular (zero, respectively)
if and only if $\limfunc{rank}(\tilde{B}(y))<q$ ($\tilde{B}(y)=0$,
respectively).

(c) The set $\mathsf{\tilde{B}}$ given by $\{y\in \mathbb{R}^{n}:\limfunc{%
rank}(\tilde{B}(y))<q\}$ (or, in view of (b), equivalently given by $\{y\in 
\mathbb{R}^{n}:\det (\tilde{\Omega}_{Het}\left( y\right) )=0\}$) is either a 
$\lambda _{\mathbb{R}^{n}}$-null set or the entire sample space $\mathbb{R}%
^{n}$. The latter occurs if and only if Assumption \ref{R_and_X_tilde} is
violated (in which case the test based on $\tilde{T}_{Het}$ becomes trivial,
as then $\tilde{T}_{Het}$ is identically zero).

(d) Suppose Assumption \ref{R_and_X_tilde} holds. Then for every $\mu
_{0}\in \mathfrak{M}_{0}$ the set $\mathsf{\tilde{B}}-\mu _{0}$ is a finite
union of proper linear subspaces; in case $q=1$, $\mathsf{\tilde{B}}-\mu
_{0} $ is even a proper linear subspace itself.\footnote{%
Consequently, $\mathsf{\tilde{B}}$ is a finite union of proper affine
subspaces, and is a proper affine subspace itself in case $q=1$.}$^{\text{,}%
} $\footnote{%
If Assumption \ref{R_and_X_tilde} is violated, then $\mathsf{\tilde{B}}-\mu
_{0}=\mathsf{\tilde{B}}=\mathbb{R}^{n}$ in view of Part (c).} [Note that $%
\mathsf{\tilde{B}}-\mu _{0}$ does not depend on the choice of $\mu _{0}\in 
\mathfrak{M}_{0}$. In particular, if $r=0$, i.e., if $\mathfrak{M}_{0}$ is
linear, we thus may set $\mu _{0}=0$.]

(e) $\mathsf{\tilde{B}}$ is a closed set and contains $\mathfrak{M}_{0}$.
Also $\mathsf{\tilde{B}}$ is $G(\mathfrak{M}_{0})$-invariant, and in
particular $\mathsf{\tilde{B}}+\mathfrak{M}_{0}^{lin}=\mathsf{\tilde{B}}$.
\end{lemma}

In light of Part (c) of the lemma, we see that Assumption \ref{R_and_X_tilde}
is a natural and unavoidable condition if one wants to obtain a sensible
test from $\tilde{T}_{Het}$.\footnote{%
If this assumption is violated then $\tilde{T}_{Het}$ is identically zero,
an uninteresting trivial case.} Furthermore, note that if $\mathsf{\tilde{B}}%
=\mathfrak{M}_{0}$ is true, then Assumption \ref{R_and_X_tilde} must be
satisfied (since $\mathfrak{M}_{0}$ is a $\lambda _{\mathbb{R}^{n}}$-null
set as $k-q<n$ is always the case). For later use we also mention that under
Assumption \ref{R_and_X_tilde} the statistic $\tilde{T}_{Het}$ is continuous
at every $y\in \mathbb{R}^{n}\backslash \mathsf{\tilde{B}}$.\footnote{%
If Assumption \ref{R_and_X_tilde} is violated, then $\tilde{T}_{Het}$ is
constant equal to zero, and hence trivially continuous everywhere.}

We finally consider for completeness, and in analogy with $T_{uc}$,%
\begin{equation}
\tilde{T}_{uc}(y)=\left\{ 
\begin{array}{cc}
(R\hat{\beta}\left( y\right) -r)^{\prime }\left( \tilde{\sigma}%
^{2}(y)R\left( X^{\prime }X\right) ^{-1}R^{\prime }\right) ^{-1}(R\hat{\beta}%
\left( y\right) -r) & \text{if }y\notin \mathfrak{M}_{0}, \\ 
0 & \text{if }y\in \mathfrak{M}_{0},%
\end{array}%
\right.  \label{T_uncorr_tilde}
\end{equation}%
where $\tilde{\sigma}^{2}(y)=\tilde{u}\left( y\right) ^{\prime }\tilde{u}%
\left( y\right) /(n-(k-q))\geq 0$ (which vanishes if and only if $y\in 
\mathfrak{M}_{0}$). Of course, our choice to set $\tilde{T}_{uc}(y)=0$ for $%
y\in \mathfrak{M}_{0}$ again has no import on the probabilistic results in
the paper, since $\mathfrak{M}_{0}$ is a $\lambda _{\mathbb{R}^{n}}$-null
set (and since the measures $P_{\mu ,\sigma ^{2}\Sigma }$ are absolutely
continuous). For later use we also mention that $\tilde{T}_{uc}$ is
continuous at every $y\in \mathbb{R}^{n}\backslash \mathfrak{M}_{0}$. As we
shall see in Section \ref{sec:trivial_uc_tilde}, there is a close connection
between $\tilde{T}_{uc}$ and $T_{uc}$.

\begin{remark}
\label{rem:tildeGM0}The test statistics $\tilde{T}_{Het}$ as well as $\tilde{%
T}_{uc}$ are $G(\mathfrak{M}_{0})$-invariant as is easily seen (with the
respective exceptional sets $\mathsf{\tilde{B}}$ and $\mathfrak{M}_{0}$ also
being $G(\mathfrak{M}_{0})$-invariant), but typically they are \emph{not}
nonsphericity-corrected F-type tests in the sense of Section 5.4 in \cite%
{PP2016}.
\end{remark}

\begin{remark}
Remark \ref{obvious} also applies to $\tilde{T}_{Het}$ and $\tilde{T}_{uc}$.
[To see this note that the respective exceptional sets $\mathsf{\tilde{B}}$
and $\mathfrak{M}_{0}$ are the same irrespective of whether $(R,r)$ or $(%
\bar{R},\bar{r})$ is used, and that $A$ cancels out in the respective
quadratic forms appearing in the definitions of the test statistics.]
\end{remark}

\subsection{Size control results for $\tilde{T}_{Het}$ and $\tilde{T}_{uc}$
when $\mathfrak{C}=\mathfrak{C}_{Het}$\label{sec_size_control_tilde_2}}

Here we discuss size control results for $\tilde{T}_{uc}$ as well as for $%
\tilde{T}_{Het}$ over the heteroskedasticity model $\mathfrak{C}_{Het}$
(more precisely, over the null hypothesis $H_{0}$ described in (\ref{testing
problem}) with $\mathfrak{C}=\mathfrak{C}_{Het}$). Some peculiar properties
of the test statistics $\tilde{T}_{uc}$ and $\tilde{T}_{Het}$ are then
discussed in the following section.

We note that the first statement in Part (a) of the subsequent theorem is
actually trivial, since $\tilde{T}_{uc}$ is bounded as shown in the next
section (which also provides a discussion when non-trivial size-controlling
critical values exist).

\begin{theorem}
\label{Hetero_Robust_tilde} (a) For every $0<\alpha <1$ there exists a real
number $C(\alpha )$ such that%
\begin{equation}
\sup_{\mu _{0}\in \mathfrak{M}_{0}}\sup_{0<\sigma ^{2}<\infty }\sup_{\Sigma
\in \mathfrak{C}_{Het}}P_{\mu _{0},\sigma ^{2}\Sigma }(\tilde{T}_{uc}\geq
C(\alpha ))\leq \alpha  \label{size-control_Het_uncorr_tilde}
\end{equation}%
holds. Furthermore, even equality can be achieved in (\ref%
{size-control_Het_uncorr_tilde}) by a proper choice of $C(\alpha )$,
provided $\alpha \in (0,\alpha ^{\ast }]\cap (0,1)$ holds, where $\alpha
^{\ast }=\sup_{C\in (C^{\ast },\infty )}\sup_{\Sigma \in \mathfrak{C}%
_{Het}}P_{\mu _{0},\Sigma }(\tilde{T}_{uc}\geq C)$ and where $C^{\ast }=\max
\{\tilde{T}_{uc}(\mu _{0}+e_{i}(n)):i\in I_{1}(\mathfrak{M}_{0}^{lin})\}$
for $\mu _{0}\in \mathfrak{M}_{0}$ (with neither $\alpha ^{\ast }$ nor $%
C^{\ast }$ depending on the choice of $\mu _{0}\in \mathfrak{M}_{0}$).

(b) Suppose Assumption \ref{R_and_X_tilde} is satisfied.\footnote{%
Condition (\ref{non-incl_Het_tilde}) clearly implies that the set $\mathsf{%
\tilde{B}}$ is a proper subset of $\mathbb{R}^{n}$ and thus implies
Assumption \ref{R_and_X_tilde}. Hence, we could have dropped this assumption
from the formulation of the theorem. A similar remark applies to some of the
other results given below and will not be repeated.} Suppose further that $%
\tilde{T}_{Het}$ is not constant on $\mathbb{R}^{n}\backslash \mathsf{\tilde{%
B}}$.\footnote{\label{FN_const}The case where $\tilde{T}_{Het}$ is constant
on $\mathbb{R}^{n}\backslash \mathsf{\tilde{B}}$ can actually occur under
Assumption \ref{R_and_X_tilde}, see Remark \ref{rem:ex_const} in Appendix %
\ref{app_C}. In such a case $\tilde{T}_{Het}$ is trivially size-controllable
(since $\mathsf{\tilde{B}}$ is a $\lambda _{\mathbb{R}^{n}}$-null set under
Assumption \ref{R_and_X_tilde} and since all probability measures in (\ref%
{lm2}) are absolutely continuous). However, neither a smallest
size-controlling critical value exists (when considering rejection regions
of the form $\{\tilde{T}_{Het}\geq C\}$) nor can exact size controllability
be achieved for $0<\alpha <1$. [If Assumption \ref{R_and_X_tilde} is
violated, $\tilde{T}_{Het}$ is identically zero and a similar remark
applies.]} Then for every $0<\alpha <1$ there exists a real number $C(\alpha
)$ such that%
\begin{equation}
\sup_{\mu _{0}\in \mathfrak{M}_{0}}\sup_{0<\sigma ^{2}<\infty }\sup_{\Sigma
\in \mathfrak{C}_{Het}}P_{\mu _{0},\sigma ^{2}\Sigma }(\tilde{T}_{Het}\geq
C(\alpha ))\leq \alpha  \label{size-control_Het_tilde}
\end{equation}%
holds, provided that for some $\mu _{0}\in \mathfrak{M}_{0}$ (and hence for
all $\mu _{0}\in \mathfrak{M}_{0}$)%
\begin{equation}
\mu _{0}+e_{i}(n)\notin \mathsf{\tilde{B}}\text{ \ \ for every \ }i\in I_{1}(%
\mathfrak{M}_{0}^{lin}).  \label{non-incl_Het_tilde}
\end{equation}%
Furthermore, under condition (\ref{non-incl_Het_tilde}), even equality can
be achieved in (\ref{size-control_Het_tilde}) by a proper choice of $%
C(\alpha )$, provided $\alpha \in (0,\alpha ^{\ast }]\cap (0,1)$ holds,
where now $\alpha ^{\ast }=\sup_{C\in (C^{\ast },\infty )}\sup_{\Sigma \in 
\mathfrak{C}_{Het}}P_{\mu _{0},\Sigma }(\tilde{T}_{Het}\geq C)$ and where $%
C^{\ast }=\max \{\tilde{T}_{Het}(\mu _{0}+e_{i}(n)):i\in I_{1}(\mathfrak{M}%
_{0}^{lin})\}$ for $\mu _{0}\in \mathfrak{M}_{0}$ (with neither $\alpha
^{\ast }$ nor $C^{\ast }$ depending on the choice of $\mu _{0}\in \mathfrak{M%
}_{0}$).

(c) Under the assumptions of Part (a) (Part (b), respectively) implying
existence of a critical value $C(\alpha )$ satisfying (\ref%
{size-control_Het_uncorr_tilde}) ((\ref{size-control_Het_tilde}),
respectively), a smallest critical value, denoted by $C_{\Diamond }(\alpha )$%
, satisfying (\ref{size-control_Het_uncorr_tilde}) ((\ref%
{size-control_Het_tilde}), respectively) exists for every $0<\alpha <1$.%
\footnote{%
Note that there are in fact no assumptions for Part (a). We have chosen this
formulation for reasons of brevity.} And $C_{\Diamond }(\alpha )$
corresponding to Part (a) (Part (b), respectively) is also the smallest
among the critical values leading to equality in (\ref%
{size-control_Het_uncorr_tilde}) ((\ref{size-control_Het_tilde}),
respectively) whenever such critical values exist. [Although $C_{\Diamond
}(\alpha )$ corresponding to Part (a) and (b), respectively, will typically
be different, we use the same symbol.]\footnote{%
Cf.~also Appendix \ref{useful}.}
\end{theorem}

We see from the theorem that $\tilde{T}_{uc}$ is always size controllable
over $\mathfrak{C}_{Het}$, but as discussed in Section \ref{sec:trivial}
below there is a caveat: Unless (\ref{non-incl_Het_uncorr}), i.e., the
necessary and sufficient condition for size-controllability of $T_{uc}$, is
satisfied, size-controlling $\tilde{T}_{uc}$ leads to trivial tests. We also
see that the condition for size control of $\tilde{T}_{Het}$ over $\mathfrak{%
C}_{Het}$, i.e., condition (\ref{non-incl_Het_tilde}) is always satisfied in
case $\mathsf{\tilde{B}}=\mathfrak{M}_{0}$ (since (\ref{non-incl_Het_tilde})
is then equivalent to $e_{i}(n)\notin \mathfrak{M}_{0}^{lin}$ for every $%
i\in I_{1}(\mathfrak{M}_{0}^{lin})$). Furthermore, condition (\ref%
{non-incl_Het_tilde}) always only depends on $X$ and $R$; in particular, it
does not depend on how the weights $\tilde{d}_{i}$ figuring in the
definition of $\tilde{T}_{Het}$ have been chosen (note that $\mu
_{0}+e_{i}(n)\notin \mathsf{\tilde{B}}$ is equivalent to $e_{i}(n)\notin 
\mathsf{\tilde{B}}-\mu _{0}$ and that the set $\mathsf{\tilde{B}}-\mu _{0}$
depends only on $X$ and $R$). Furthermore, the size-controlling critical
values $C(\alpha )$ in Part (b) of the preceding theorem depend only on $X$, 
$R$, and $r$, as well as on the choice of weights $\tilde{d}_{i}$, whereas
in Part (a) the dependence is only on $X$, $R$, and $r$. We do not show
these dependencies in the notation. In fact, as shown in Lemma \ref%
{lem_indep_r} in Appendix \ref{app_C}, it turns out that the size and the
size-controlling critical values in both cases actually do \emph{not} depend
on the value of $r$ at all (provided the weights $\tilde{d}_{i}$ are not
allowed to depend on $r$ in case of $\tilde{T}_{Het}$). Similarly, it is
easy to see that $\alpha ^{\ast }$ and $C^{\ast }$ do not depend on $r$
(under the same provision as before in case of $\tilde{T}_{Het}$).

Similarly as in Section \ref{sec_size_control}, a critical value delivering
size control over $\mathfrak{C}_{Het}$ also delivers size control over \emph{%
any} other heteroskedasticity model $\mathfrak{C}$ since $\mathfrak{C}%
\subseteq \mathfrak{C}_{Het}$. Of course, for such a $\mathfrak{C}$ even
smaller critical values (than needed for $\mathfrak{C}_{Het}$) may already
suffice for size control. Also note that sufficient conditions implying size
control over $\mathfrak{C}_{Het}$ may be more restrictive than sufficient
conditions only implying size control over a smaller heteroskedasticity
model $\mathfrak{C}$. For size control results tailored to such smaller
models $\mathfrak{C}$ see Appendix \ref{app_d}.

\begin{remark}
\label{rem:equiv_tilde}\emph{(Some equivalencies)} If the respective
smallest size-controlling critical values are used (provided they exist),
the tests obtained from $\tilde{T}_{Het}$ with the HC0R and the HC1R
weights, respectively, are identical, as these two test statistics differ
only by a multiplicative constant. The same reasoning applies to the test
statistics based on the HC0R-HC4R weights, respectively, in case $\tilde{h}%
_{ii}$ does not depend on $i$.
\end{remark}

\begin{remark}
\label{rem_positiv_tilde}\emph{(Positivity of size-controlling critical
values) }For every $0<\alpha <1$ any $C(\alpha )$ satisfying (\ref%
{size-control_Het_uncorr_tilde}) or (\ref{size-control_Het_tilde}) is
necessarily positive. To see this observe that $\{\tilde{T}_{uc}\geq C\}=\{%
\tilde{T}_{Het}\geq C\}=\mathbb{R}^{n}$ for $C\leq 0$, since both test
statistics are nonnegative everywhere.
\end{remark}

The next proposition complements Theorem \ref{Hetero_Robust_tilde} and
provides a lower bound for the size-controlling critical values (other than
the trivial bound given in the preceding remark). The lower bound is useful
for the same reasons as discussed subsequent to Proposition \ref{rem_C*}.

\begin{proposition}
\label{rem_C*_tilde}\footnote{%
It is not difficult to show in the context of Parts (a) and (b) of the
proposition that any critical value $C>C^{\ast }$ actually leads to size
less than $1$. This follows from a reasoning similar as in Remark 5.4 of 
\cite{PP3}.}$^{\text{,}}$\footnote{%
If (\ref{non-incl_Het_tilde}) in Part (b) of the proposition does not hold,
the conclusion of Part (b) can be shown to continue to hold with $C^{\ast }$
as defined in Theorem \ref{Hetero_Robust_tilde}(b), and also with $C^{\ast }$
as defined in Lemma 5.11 of \ \cite{PP3} (note that under the assumptions of
Part (b) of the proposition both definitions of $C^{\ast }$ actually
coincide as shown in the proof of Theorem \ref{Hetero_Robust_tilde}). If $%
\tilde{T}_{Het}$ is constant on $\mathbb{R}^{n}\backslash \mathsf{\tilde{B}}$
or if Assumption \ref{R_and_X_tilde} fails (the latter implying $\tilde{T}%
_{Het}\equiv 0$), the conclusion of Part (b) also holds as is easily seen
(regardless of which of the two definitions of $C^{\ast }$ is adopted).}(a)
Any $C(\alpha )$ satisfying (\ref{size-control_Het_uncorr_tilde})
necessarily has to satisfy $C(\alpha )\geq C^{\ast }$, where $C^{\ast }$ is
as in Part (a) of Theorem \ref{Hetero_Robust_tilde}. In fact, for any $%
C<C^{\ast }$ we have $\sup_{\Sigma \in \mathfrak{C}_{Het}}P_{\mu _{0},\sigma
^{2}\Sigma }(\tilde{T}_{uc}\geq C)=1$ for every $\mu _{0}\in \mathfrak{M}%
_{0} $ and every $\sigma ^{2}\in (0,\infty )$.

(b) Suppose Assumption \ref{R_and_X_tilde} and (\ref{non-incl_Het_tilde})
are satisfied, and that $\tilde{T}_{Het}$ is not constant on $\mathbb{R}%
^{n}\backslash \mathsf{\tilde{B}}$. Then any $C(\alpha )$ satisfying (\ref%
{size-control_Het_tilde}) necessarily has to satisfy $C(\alpha )\geq C^{\ast
}$, where $C^{\ast }$ is as in Part (b) of Theorem \ref{Hetero_Robust_tilde}%
. In fact, for any $C<C^{\ast }$ we have $\sup_{\Sigma \in \mathfrak{C}%
_{Het}}P_{\mu _{0},\sigma ^{2}\Sigma }(\tilde{T}_{Het}\geq C)=1$ for every $%
\mu _{0}\in \mathfrak{M}_{0}$ and every $\sigma ^{2}\in (0,\infty )$.
\end{proposition}

\begin{remark}
\label{rem:larger_alphastar_tilde}Suppose the assumptions of Part (a) (Part
(b), respectively) of Theorem \ref{Hetero_Robust_tilde} are satisfied. Then
we know from that theorem that the size (over $\mathfrak{C}_{Het})$ of $\{%
\tilde{T}_{uc}\geq C_{\Diamond }(\alpha )\}$ ($\{\tilde{T}_{Het}\geq
C_{\Diamond }(\alpha )\}$, respectively) equals $\alpha $ provided $\alpha
\in (0,\alpha ^{\ast }]\cap (0,1)$. If now $\alpha ^{\ast }<\alpha <1$, then
the size (over $\mathfrak{C}_{Het})$ of $\{\tilde{T}_{uc}\geq C_{\Diamond
}(\alpha )\}$ ($\{\tilde{T}_{Het}\geq C_{\Diamond }(\alpha )\}$,
respectively) equals $\alpha ^{\ast }$ (where the $C_{\Diamond }(\alpha )$'s
pertaining to Parts (a) and (b) may be different). This follows from $%
C_{\Diamond }(\alpha )\geq C^{\ast }$ (see Proposition \ref{rem_C*_tilde}
above) and Remark 5.13(i) in \cite{PP3}).\footnote{%
The assumptions for Part A of Proposition 5.12 in \cite{PP3} required in
Remark 5.13 of that paper are satisfied under the assumptions of Theorem \ref%
{Hetero_Robust_tilde} as shown in the proof of Theorem \ref%
{theorem_groupwise_hetero_tilde} in Appendix \ref{app_C}. In this proof also
the condition $\lambda _{\mathbb{R}^{n}}(\tilde{T}_{uc}=C^{\ast })=0$ ($%
\lambda _{\mathbb{R}^{n}}(\tilde{T}_{Het}=C^{\ast })=0$, respectively)
required in Remark 5.13 of \cite{PP3} is verified.} This argument actually
also delivers that $C_{\Diamond }(\alpha )=C^{\ast }$ must hold in case $%
\alpha ^{\ast }<\alpha <1$.
\end{remark}

\begin{remark}
\label{rem_1}In contrast to Section \ref{sec_size_control}, we have little
information on the extent to which the sufficient conditions for size
control in Part (b) of Theorem \ref{Hetero_Robust_tilde} are also necessary.
This is due to the fact that $\tilde{T}_{Het}$ is typically not a
nonsphericity-corrected F-type test as noted in Remark \ref{rem:tildeGM0}.
What can be said in general in the context of Part (b) of Theorem \ref%
{Hetero_Robust_tilde} in case (\ref{non-incl_Het_tilde}) is violated, is
that the size of the rejection region $\{\tilde{T}_{Het}\geq C\}$ over $%
\mathfrak{C}_{Het}$ is certainly equal to $1$ for every $C<\max \{\tilde{T}%
_{Het}(\mu _{0}+e_{i}(n)):\mu _{0}+e_{i}(n)\notin \mathsf{\tilde{B}}\}$,
where $\mu _{0}\in \mathfrak{M}_{0}$ is arbitrary (the maximum being
independent of the choice of $\mu _{0}\in \mathfrak{M}_{0}$) and where we
use the convention that this maximum is $-\infty $ in case the set over
which the maximum is taken is empty. This follows from Lemma 4.1 in \cite%
{PP4} with $\mathbb{K}$ equal to the collection $\{\Pi _{(\mathfrak{M}%
_{0}^{lin})^{\bot }}e_{i}(n):\mu _{0}+e_{i}(n)\notin \mathsf{\tilde{B}}\}$.
\end{remark}

\begin{remark}
\label{q=k}Suppose $q=k$. Then Assumption \ref{R_and_X_tilde} is always
satisfied (since $\mathfrak{M}_{0}$ being a singleton $\{\mu _{0}\}$ implies 
$\mathfrak{M}_{0}^{lin}=\{0\}$, and thus $s=0$ in Assumption \ref%
{R_and_X_tilde}). The subsequent claims are proved in Appendix \ref{app_C}.

(i) In case $q=k>1$, it is not difficult to see that then $\mu
_{0}+e_{i}(n)\in \mathsf{\tilde{B}}$ for every $i=1,\ldots ,n$ holds,
implying that the sufficient condition (\ref{non-incl_Het_tilde}) in Theorem %
\ref{Hetero_Robust_tilde}(b) is violated. [In contrast, in case $q=k=1$,
both examples where (\ref{non-incl_Het_tilde}) is satisfied as well as
examples where (\ref{non-incl_Het_tilde}) is not satisfied can be found.]

(ii) Despite of (i), in case $q=k\geq 1$ the test statistic $\tilde{T}_{Het}$
is always size-controllable over $\mathfrak{C}_{Het}$. This is so since in
case $q=k\geq 1$ the statistic $\tilde{T}_{Het}$ is a bounded function.

(iii) We also note that in case $q=k\geq 1$ both the case where $\tilde{T}%
_{Het}$ is constant on $\mathbb{R}^{n}\backslash \mathsf{\tilde{B}}$ as well
as the case where $\tilde{T}_{Het}$ is not constant on $\mathbb{R}%
^{n}\backslash \mathsf{\tilde{B}}$ can occur. [In the latter case a \emph{%
smallest} size-controlling critical value exists in view of Appendix \ref%
{useful}. In the former case no \emph{smallest} size-controlling critical
value exists (when considering rejection regions of the form $\{\tilde{T}%
_{Het}\geq C\}$).]
\end{remark}

\begin{remark}
\label{rem_strict_ineq_2} Let $\tilde{T}$ stand for either $\tilde{T}_{Het}$
or $\tilde{T}_{uc}$, where in case of $\tilde{T}=\tilde{T}_{Het}$ we suppose
that Assumption \ref{R_and_X_tilde} is satisfied and that $\tilde{T}_{Het}$
is not constant on $\mathbb{R}^{n}\backslash \mathsf{\tilde{B}}$: By Lemma %
\ref{lem:nullset} in Appendix \ref{app_C} the rejection regions $\{y:\tilde{T%
}(y)\geq C\}$ and $\{y:\tilde{T}(y)>C\}$ differ only by a $\lambda _{\mathbb{%
R}^{n}}$-null set. Since the measures $P_{\mu ,\sigma ^{2}\Sigma }$ are
absolutely continuous w.r.t.$~\lambda _{\mathbb{R}^{n}}$ when $\Sigma $ is
nonsingular, $P_{\mu ,\sigma ^{2}\Sigma }(\tilde{T}\geq C)=P_{\mu ,\sigma
^{2}\Sigma }(\tilde{T}>C)$ then follows, and hence the results in this and
the subsequent section given for rejection probabilities $P_{\mu ,\sigma
^{2}\Sigma }(\tilde{T}\geq C)$ apply to rejection probabilities $P_{\mu
,\sigma ^{2}\Sigma }(\tilde{T}>C)$ equally well (under the above provision
in case of $T=\tilde{T}_{Het}$). A similar remark applies to the results in
Appendix \ref{sec_size_control_2_tilde}.
\end{remark}

\subsection{Tests obtained from $\tilde{T}_{uc}$ or $\tilde{T}_{Het}$ can be
trivial\label{sec:trivial}}

For the test statistic $T_{uc}$ the rejection regions $\{T_{uc}\geq C\}$, as
well as their complements, have positive ($n$-dimensional) Lebesgue measure
for every positive real number $C$.\footnote{%
The case $C\leq 0$ is uninteresting as the rejection region of $T_{uc}$ (and
of all other test statistics considered) then are the entire space $\mathbb{R%
}^{n}$, since $T_{uc}$ (and the other test statistics considered) take on
only nonnegative values.} This follows from Parts 5\&6 of Lemma 5.15 in \cite%
{PP2016} together with Remark \ref{F-type} in Appendix \ref{app_B}. As a
consequence, all rejection probabilities -- under the null as well as under
the alternative -- are positive and less than one regardless of the choice
of $C>0$. [This is so because of our Gaussianity assumption and the fact
that all $\Sigma \in \mathfrak{C}_{Het}$ are positive definite.] For similar
reasons, the same is true for $T_{Het}$ provided Assumption \ref{R_and_X} is
satisfied.\footnote{%
If Assumption \ref{R_and_X} is not satisfied then $T_{Het}\equiv 0$, and the
resulting test (with rejection region $\{T_{Het}\geq C\}$) is trivial as it
never rejects for $C>0$, while it always rejects for $C\leq 0$.} The
situation is somewhat different for tests derived from $\tilde{T}_{uc}$ or $%
\tilde{T}_{Het}$ as we shall discuss next. In the course of this, we also
establish a connection between $T_{uc}$ and $\tilde{T}_{uc}$ that is of
independent interest. In this section the size of a test always refers to
size over $\mathfrak{C}_{Het}$.

\subsubsection{The case of $\tilde{T}_{uc}$\label{sec:trivial_uc_tilde}}

First, observe that $\tilde{T}_{uc}(y)\leq n-(k-q)$ holds for every $y\in 
\mathbb{R}^{n}$ and that this bound is sharp. To see this, note that using
standard least-squares theory%
\begin{equation}
\tilde{T}_{uc}(y)=(n-(k-q))\left( 1-\tsum_{i=1}^{n}\hat{u}%
_{i}^{2}(y)/\tsum_{i=1}^{n}\tilde{u}_{i}^{2}(y)\right) \leq n-(k-q)
\label{standard}
\end{equation}%
for $y\notin \mathfrak{M}_{0}$ and that $\tilde{T}_{uc}(y)=0$ else; the
bound is attained precisely for $y\in \limfunc{span}(X)\backslash \mathfrak{M%
}_{0}$. An immediate consequence of this observation is that any critical
value $C\geq (n-(k-q))$ leads to a test with rejection region $\{\tilde{T}%
_{uc}\geq C\}$ that is either empty (if $C>n-(k-q)$) or is a $\lambda _{%
\mathbb{R}^{n}}$-null set, namely $\limfunc{span}(X)\backslash \mathfrak{M}%
_{0}$ (if $C=n-(k-q)$). Consequently, such a test is trivial in that all
rejection probabilities (under the null as well as under the alternative)
are zero (because of our Gaussianity assumption and the fact that all $%
\Sigma \in \mathfrak{C}_{Het}$ are positive definite). As an aside we note
that any $C<n-(k-q)$ leads to a non-trivial test as is easily seen.

Of course, a critical value $C$ satisfying $C\geq n-(k-q)$ is certainly
size-controlling, but is useless since it leads to a trivial test as just
discussed. We now ask if and when the smallest size-controlling critical
value $C_{\Diamond }(\alpha )$, guaranteed to exist by Part (c) of Theorem %
\ref{Hetero_Robust_tilde}, leads to a non-trivial test. [This is certainly
so if $\alpha ^{\ast }$ in Part (a) of Theorem \ref{Hetero_Robust_tilde} is
positive, but note that the theorem is silent on this issue.] To obtain
insight, we establish a simple, but important, relationship between the test
statistics $\tilde{T}_{uc}$ and $T_{uc}$ that is of independent interest
also: Note that standard least-squares theory gives%
\begin{equation*}
T_{uc}(y)=(n-k)\left( \tsum_{i=1}^{n}\tilde{u}_{i}^{2}(y)/\tsum_{i=1}^{n}%
\hat{u}_{i}^{2}(y)-1\right)
\end{equation*}%
for $y\notin \limfunc{span}(X)$, and recall $T_{uc}(y)=0$ for $y\in \limfunc{%
span}(X)$. Hence, we obtain%
\begin{equation}
\tilde{T}_{uc}(y)=(n-(k-q))\left( T_{uc}(y)/(n-k+T_{uc}(y)\right)
=g(T_{uc}(y))  \label{identity}
\end{equation}%
for every $y\notin \limfunc{span}(X)$, where $g:[0,\infty )\rightarrow
\lbrack 0,n-(k-q))$ is continuous and strictly increasing with $%
\lim_{x\rightarrow \infty }g(x)=(n-(k-q))$. [Since $T_{uc}(y_{m})\rightarrow
\infty $ for every sequence $y_{m}\rightarrow y\in \limfunc{span}%
(X)\backslash \mathfrak{M}_{0}$, the sharpness of the bound $n-(k-q)$ can
thus also be read-off from (\ref{identity}).] As a consequence, for every
critical value $C>0$, the rejection regions $\{\tilde{T}_{uc}\geq C\}$ and $%
\{T_{uc}\geq g^{-1}(C)\}$ differ at most by $\limfunc{span}(X)$, which is a $%
\lambda _{\mathbb{R}^{n}}$-null set; in particular, the rejection
probabilities (under the null as well as under the alternative) are the same.%
\footnote{%
This is so because of our Gaussianity assumption and the fact that all $%
\Sigma \in \mathfrak{C}_{Het}$ are positive definite.} \emph{That is, the
test statistics }$\tilde{T}_{uc}$\emph{\ and }$T_{uc}$\emph{\ give rise to
(essentially) the same test, if the critical values chosen are linked by the
function }$g$\emph{\ as above. In particular, as we shall see, this is the
case if the respective smallest size-controlling critical values are used
for both test statistics (provided both these values exist).}

To see what the preceding discussion entails for the existence of
non-trivial size-controlling critical values for $\tilde{T}_{uc}$ we
distinguish two cases. In the first case we shall see that non-trivial
size-controlling critical values do not exist, whereas in the second case
they do indeed exist.

\emph{Case 1: Condition (\ref{non-incl_Het_uncorr}) is violated. }Recall
from Proposition \ref{rem_necessity} that then the size of $\{T_{uc}\geq D\}$
is $1$ for every real $D$ (in particular, implying that $T_{uc}$ is not size
controllable). It transpires from the preceding discussion, that hence the
size of $\{\tilde{T}_{uc}\geq C\}$ must equal $1$ for every $C$ satisfying $%
0<C<n-(k-q)$ (and a fortiori for $C\leq 0$), because $D:=g^{-1}(C)$ is
well-defined and real for $0<C<n-(k-q)$. As a consequence, any
size-controlling critical value $C$ for $\tilde{T}_{uc}$ must satisfy $C\geq
n-(k-q)$ (with the smallest size-controlling critical value given by $%
n-(k-q) $), thus leading to a rejection region that is trivial in that it is
empty (if $C>n-(k-q)$) or is a $\lambda _{\mathbb{R}^{n}}$-null set, namely $%
\limfunc{span}(X)\backslash \mathfrak{M}_{0}$ (if $C=n-(k-q)$). That is --
while $\tilde{T}_{uc}$ is size-controllable in the present case -- it is so
only in a trivial way.\footnote{%
The trivial size-controlling critical values $C$ for $\tilde{T}_{uc}$ sort
of correspond to using $\infty $ as a \textquotedblleft size-controlling
critical value\textquotedblright\ for $T_{uc}$.} [Another way of arriving at
the above conclusion is to use Part (a) of Proposition \ref{rem_C*_tilde}
and to observe that in Part (a) of Theorem \ref{Hetero_Robust_tilde} the
quantity $C^{\ast }$ equals $n-(k-q)$. To see the latter, note that
violation of condition (\ref{non-incl_Het_uncorr}) implies existence of an
index $i\in I_{1}(\mathfrak{M}_{0}^{lin})$ with $e_{i}(n)\in \limfunc{span}%
(X)$. In particular, $\hat{u}(\mu _{0}+e_{i}(n))=0$. Since $e_{i}(n)\notin 
\mathfrak{M}_{0}^{lin}$ must hold in view of $i\in I_{1}(\mathfrak{M}%
_{0}^{lin})$, and thus $\mu _{0}+e_{i}(n)\notin \mathfrak{M}_{0}$ for every $%
\mu _{0}\in \mathfrak{M}_{0}$ must be true, we may use (\ref{standard}) to
arrive at $\tilde{T}_{uc}(\mu _{0}+e_{i}(n))=n-(k-q)$ for this $i\in I_{1}(%
\mathfrak{M}_{0}^{lin})$. This shows $C^{\ast }\geq n-(k-q)$. Equality then
follows since $C^{\ast }\leq n-(k-q)$ trivially holds by (\ref{standard}).
As a point of interest we also note that $C^{\ast }=n-(k-q)$ implies that $%
\alpha ^{\ast }$ in Part (a) of Theorem \ref{Hetero_Robust_tilde} satisfies $%
\alpha ^{\ast }=0$.]

\emph{Case 2: Condition (\ref{non-incl_Het_uncorr}) is satisfied.} In this
case $T_{uc}$ is size controllable according to Theorem \ref{Hetero_Robust}.
In particular, for any given $\alpha \in (0,1)$ there exists a smallest real
number $D_{\Diamond }(\alpha )$ such that the size of $\{T_{uc}\geq
D_{\Diamond }(\alpha )\}$ is less than or equal to $\alpha $, with equality
holding for $\alpha \in (0,\alpha _{T_{uc}}^{\ast }]\cap (0,1)$ where $%
\alpha _{T_{uc}}^{\ast }$ refers to $\alpha ^{\ast }$ appearing in Theorem %
\ref{Hetero_Robust}(a), and recall from that theorem that $\alpha
_{T_{uc}}^{\ast }>0$; and $D_{\Diamond }(\alpha )>0$ by Remark \ref%
{rem_positiv}.\footnote{\label{fn:999}If $\alpha _{T_{uc}}^{\ast }<\alpha <1$%
, then the size, in fact, equals $\alpha _{T_{uc}}^{\ast }$; see Remark \ref%
{rem:larger_alphastar}.} Also note that the rejection region $\{T_{uc}\geq
D_{\Diamond }(\alpha )\}$ is not trivial as it has positive $\lambda _{%
\mathbb{R}^{n}}$-measure (and the same is true for its complement); see the
discussion at the very beginning of Section \ref{sec:trivial}. Setting $%
C_{\Diamond }(\alpha )=g(D_{\Diamond }(\alpha ))$ and using that $\{\tilde{T}%
_{uc}\geq C_{\Diamond }(\alpha )\}$ and $\{T_{uc}\geq g^{-1}(C_{\Diamond
}(\alpha ))\}=\{T_{uc}\geq D_{\Diamond }(\alpha )\}$ differ at most by the $%
\lambda _{\mathbb{R}^{n}}$-null set $\limfunc{span}(X)$, we see that (i) $%
0<C_{\Diamond }(\alpha )<n-(k-q)$, (ii) the size of $\{\tilde{T}_{uc}\geq
C_{\Diamond }(\alpha )\}$ is less than or equal to $\alpha $, with equality
holding for $\alpha \in (0,\alpha _{T_{uc}}^{\ast }]\cap (0,1)$, (iii) $%
C_{\Diamond }(\alpha )$ is the smallest size-controlling critical value
(recall that $g$ is strictly increasing), and (iv) the rejection region $\{%
\tilde{T}_{uc}\geq C_{\Diamond }(\alpha )\}$ is not trivial as it has
positive $\lambda _{\mathbb{R}^{n}}$-measure (and the same is true for its
complement). In particular, note that $\tilde{T}_{uc}$ and $T_{uc}$ give
rise to (essentially) the same test if the respective smallest
size-controlling critical values are used. We furthermore note that in the
present situation $C_{\tilde{T}_{uc}}^{\ast }=g(C_{T_{uc}}^{\ast })$ and $%
\alpha _{\tilde{T}_{uc}}^{\ast }=\alpha _{T_{uc}}^{\ast }$ hold, where $%
C_{T_{uc}}^{\ast }$, $\alpha _{T_{uc}}^{\ast }$ correspond to $C^{\ast }$, $%
\alpha ^{\ast }$ in Part (a) of Theorem \ref{Hetero_Robust}, whereas $C_{%
\tilde{T}_{uc}}^{\ast }$, $\alpha _{\tilde{T}_{uc}}^{\ast }$ correspond to $%
C^{\ast }$, $\alpha ^{\ast }$ in Part (a) of Theorem \ref%
{Hetero_Robust_tilde}.\footnote{%
If $\alpha _{\tilde{T}_{uc}}^{\ast }<\alpha <1$, then the size of $\{\tilde{T%
}_{uc}\geq C_{\Diamond }(\alpha )\}$ is, in fact, equal to $\alpha
_{T_{uc}}^{\ast }=\alpha _{\tilde{T}_{uc}}^{\ast }$; cf. Footnote \ref%
{fn:999} and Remark \ref{rem:larger_alphastar_tilde}.} In particular, $%
\alpha _{\tilde{T}_{uc}}^{\ast }>0$ and $0\leq C_{\tilde{T}_{uc}}^{\ast
}<n-(k-q)$ follow. These claims can be seen as follows: Under condition (\ref%
{non-incl_Het_uncorr}) we have $\mu _{0}+e_{i}(n)\notin \limfunc{span}(X)$
for every $i\in I_{1}(\mathfrak{M}_{0}^{lin})$ and every $\mu _{0}\in 
\mathfrak{M}_{0}$. Consequently, $\tilde{T}_{uc}(\mu
_{0}+e_{i}(n))=g(T_{uc}(\mu _{0}+e_{i}(n)))$, which proves $C_{\tilde{T}%
_{uc}}^{\ast }=g(C_{T_{uc}}^{\ast })$ in view of strict monotonicity of $g$.
The relation $\alpha _{\tilde{T}_{uc}}^{\ast }=\alpha _{T_{uc}}^{\ast }$
then follows from the definitions of $\alpha _{\tilde{T}_{uc}}^{\ast }$ and $%
\alpha _{T_{uc}}^{\ast }$ using that $\{\tilde{T}_{uc}\geq C\}$ and $%
\{T_{uc}\geq g^{-1}(C)\}$ differ at most by the $\lambda _{\mathbb{R}^{n}}$%
-null set $\limfunc{span}(X)$ for every $C>0$. Positivity of $\alpha _{%
\tilde{T}_{uc}}^{\ast }$ now follows from positivity of $\alpha
_{T_{uc}}^{\ast }$ discussed before, and $C_{\tilde{T}_{uc}}^{\ast }<n-(k-q)$
follows since $C_{\tilde{T}_{uc}}^{\ast }=g(C_{T_{uc}}^{\ast })$ and $%
C_{T_{uc}}^{\ast }<\infty $. [Another way of proving $\alpha _{\tilde{T}%
_{uc}}^{\ast }>0$ and $0\leq C_{\tilde{T}_{uc}}^{\ast }<n-(k-q)$ without
using relationship (\ref{identity}), is to first establish $C_{\tilde{T}%
_{uc}}^{\ast }<n-(k-q)$ (from observing that $\hat{u}(\mu _{0}+e_{i}(n))\neq
0$ (as $\mu _{0}+e_{i}(n)\notin \limfunc{span}(X)$) for every $i\in I_{1}(%
\mathfrak{M}_{0}^{lin})$, which implies $\tilde{T}_{uc}(\mu
_{0}+e_{i}(n))<n-(k-q)$ for every such $i$ in view of (\ref{standard})) and
then to proceed analogously as in the proof of Theorem \ref%
{Hetero_Robust_tilde_add} below.]

While $\tilde{T}_{uc}$ is always size-controllable, whereas $T_{uc}$ is not,
this does not represent any real advantage of $\tilde{T}_{uc}$ over $T_{uc}$%
, as we have seen that $\tilde{T}_{uc}$ admits only trivial size-controlling
critical values in the case where $T_{uc}$ is not size-controllable. Even
more importantly, and already noted above, these test statistics give rise
to (essentially) the same test if for both test statistics the respective
smallest size-controlling critical values are used (provided they both
exist).

\subsubsection{The case of $\tilde{T}_{Het}$\label{sec:trivial_het_tilde}}

For $\tilde{T}_{Het}$ we find that, not infrequently, it is also a bounded
function, although we have no proof that this is always so. We illustrate
the problems that can arise here first by an example. See also Remark \ref%
{rem:new}.

\begin{example}
\label{ex:bounded}Consider the $n\times 2$ design matrix $X$ where the first
column represents an intercept, the second column is $x:=(1,-1,0,\ldots
,0)^{\prime }$, and $n\geq 3$. Let $R=(0,1)$, $r=0$, hence $q=1$. Obviously,
the first column of $X$ spans $\mathfrak{M}_{0}^{lin}$. Since $%
e_{i}(n)\notin \mathfrak{M}_{0}^{lin}$ for every $i=1,\ldots ,n$, Assumption %
\ref{R_and_X_tilde} holds. Furthermore, $\tilde{h}_{ii}=n^{-1}$. Thus $%
\tilde{d}_{i}=\tilde{d}_{1}$ holds for every $i=1,\ldots ,n$ and for every
of the five choices HC0R-HC4R. Note that $\tilde{d}_{1}^{-1}=1$ (HC0R), $%
\tilde{d}_{1}^{-1}=1-n^{-1}$ (HC1R), $\tilde{d}_{1}^{-1}=1-n^{-1}$ (HC2R), $%
\tilde{d}_{1}^{-1}=(1-n^{-1})^{2}$ (HC3R), and $\tilde{d}_{1}^{-1}=1-n^{-1}$
(HC4R), and hence $0<\tilde{d}_{1}^{-1}\leq 1$ for all five choices.
Straightforward computations now show that $\tilde{\Omega}_{Het}(y)=\tilde{d}%
_{1}\left[ \left( y_{1}-\bar{y}\right) ^{2}+\left( y_{2}-\bar{y}\right) ^{2}%
\right] /4$ and%
\begin{equation}
\tilde{T}_{Het}(y)=\tilde{d}_{1}^{-1}\left( y_{1}-y_{2}\right) ^{2}/\left[
\left( y_{1}-\bar{y}\right) ^{2}+\left( y_{2}-\bar{y}\right) ^{2}\right]
\label{explicit}
\end{equation}%
whenever the numerator is positive, and $\tilde{T}_{Het}(y)=0$ otherwise.
Here $\bar{y}$ denotes the arithmetic mean of the observations $y_{i}$. [For
later use we also note that the set $\mathsf{\tilde{B}}$ is given by $\{y\in 
\mathbb{R}^{n}:y_{1}=y_{2}=\bar{y}\}$, and that the size control condition (%
\ref{non-incl_Het_tilde}) is satisfied, since $e_{i}(n)\notin \mathsf{\tilde{%
B}}$ for every $i=1,\ldots ,n$ (also note that $\mu _{0}$ can be chosen to
be zero because of $r=0$). Furthermore, $\tilde{T}_{Het}$ is not constant on 
$\mathbb{R}^{n}\backslash \mathsf{\tilde{B}}$, since $\tilde{T}%
_{Het}(e_{1}(n))=\tilde{T}_{Het}(e_{2}(n))=\tilde{d}%
_{1}^{-1}n^{2}/[(n-1)^{2}+1]$ and $\tilde{T}_{Het}(e_{i}(n))=0$ for $i\geq 3$
(note $n\geq 3$) and since $e_{i}(n)\notin \mathsf{\tilde{B}}$ for every $i$%
.] It is now evident from (\ref{explicit}) that $\tilde{T}_{Het}(y)\leq 2%
\tilde{d}_{1}^{-1}$ for every $y\in \mathbb{R}^{n}$ and that this bound is
attained whenever $y_{1}+y_{2}=2\bar{y}$ and $y_{1}\neq y_{2}$ (e.g., for $%
y=x$). It follows that any critical value $C\geq 2\tilde{d}_{1}^{-1}$ leads
to a test with rejection region that is empty if $C>2\tilde{d}_{1}^{-1}$,
and is a Lebesgue null-set if $C=2\tilde{d}_{1}^{-1}$ (the latter following
from Lemma \ref{lem:nullset}(d) in Appendix \ref{app_C} together with some
of the observations just noted after (\ref{explicit})); thus in both cases
all the rejection probabilities are zero under the null as well as under the
alternative (given our Gaussianity assumption and the fact that all $\Sigma
\in \mathfrak{C}_{Het}$ are positive definite); in particular, these tests
have zero power. Since $\tilde{d}_{1}^{-1}\leq 1$, this eliminates all
critical values $C\geq 2$ from practical use. In particular, this eliminates
the commonly used choice where $C$ is the $95\%$-quantile of a chi-square
distribution with $1$ degree of freedom, which is approximately equal to $%
3.8415$.
\end{example}

In the preceding example any critical value $C\geq 2\tilde{d}_{1}^{-1}$ is
trivially a size-controlling critical value for the given significance level 
$\alpha $ ($0<\alpha <1$), but it is \textquotedblleft too
large\textquotedblright\ and leads to a trivial test. Certainly, one would
prefer to use the smallest size-controlling critical value $C_{\Diamond
}(\alpha )$ instead (which in the preceding example exists by Theorem \ref%
{Hetero_Robust_tilde} and by what has been shown in the example) and one
would hope that the resulting test is not trivial. As we shall show, this is
indeed the case. To this end we first give a general result that, in
particular, is applicable to the preceding example. Recall that $C_{\Diamond
}(\alpha )$ is positive (Remark \ref{rem_positiv_tilde}), and that Theorem %
\ref{Hetero_Robust_tilde} is silent on whether $\alpha ^{\ast }>0$ or not.

\begin{theorem}
\label{Hetero_Robust_tilde_add} Suppose Assumption \ref{R_and_X_tilde} and (%
\ref{non-incl_Het_tilde}) are satisfied, and that $\tilde{T}_{Het}$ is not
constant on $\mathbb{R}^{n}\backslash \mathsf{\tilde{B}}$. Let $\alpha $
satisfy $0<\alpha <1$, and let $C^{\ast }$ and $\alpha ^{\ast }$\ be as
defined in Part (b) of Theorem \ref{Hetero_Robust_tilde}. If $C^{\ast
}<\sup_{y\in \mathbb{R}^{n}}\tilde{T}_{Het}(y)$ holds, then we have $\alpha
^{\ast }>0$, and the rejection region $\{\tilde{T}_{Het}\geq C_{\Diamond
}(\alpha )\}$ is not a $\lambda _{\mathbb{R}^{n}}$-null set, where $%
C_{\Diamond }(\alpha )$ is the smallest size-controlling critical value as
in Part (c) of Theorem \ref{Hetero_Robust_tilde}.
\end{theorem}

\begin{remark}
\label{nontrivial_rej}(i) The preceding theorem clearly implies that --
under its assumptions -- the rejection probabilities associated with the
rejection region $\{\tilde{T}_{Het}\geq C_{\Diamond }(\alpha )\}$ are
positive under the null as well as under the alternative (in view of our
Gaussianity assumption and the fact that all $\Sigma \in \mathfrak{C}_{Het}$
are positive definite). [While we already know from Theorem \ref%
{Hetero_Robust_tilde}(b) and Remark \ref{rem:larger_alphastar_tilde} that
the rejection region $\{\tilde{T}_{Het}\geq C_{\Diamond }(\alpha )\}$ has
size equal to $\alpha $ in case $\alpha \in (0,\alpha ^{\ast }]\cap (0,1)$,
and has size equal to $\alpha ^{\ast }$ if $\alpha ^{\ast }<\alpha <1$, this
by itself does \emph{not} allow one to conclude that the rejection region
has positive $\lambda _{\mathbb{R}^{n}}$-measure as the case $\alpha ^{\ast
}=0$ is not ruled out by Theorem \ref{Hetero_Robust_tilde}(b) and Remark \ref%
{rem:larger_alphastar_tilde}.]

(ii) Suppose $C^{\ast }=\sup_{y\in \mathbb{R}^{n}}\tilde{T}_{Het}(y)$, but
that the other assumptions of Theorem \ref{Hetero_Robust_tilde_add} hold.%
\footnote{%
We have not investigated whether this case can actually occur for $\tilde{T}%
_{Het}$. Recall that for $\tilde{T}_{uc}$ this case indeed can occur, see
Case 1 in Section \ref{sec:trivial_uc_tilde}.} Then the rejection region $\{%
\tilde{T}_{Het}\geq C_{\Diamond }(\alpha )\}$ is a $\lambda _{\mathbb{R}%
^{n}} $-null set; thus also the smallest (and hence any) size-controlling
critical value leads to a trivial test. To prove the claim, note that by
Proposition \ref{rem_C*_tilde} we have $C_{\Diamond }(\alpha )\geq C^{\ast }$%
, implying that the rejection regions are either empty or coincide with the
sets $\{\tilde{T}_{Het}=C^{\ast }\}$, respectively. In the latter case apply
Part (d) of Lemma \ref{lem:nullset} in Appendix \ref{app_C}. We also point
out that in the present case $\alpha ^{\ast }=0$ must hold since the
rejection regions appearing in the definition of $\alpha ^{\ast }$ are all
empty (because of $C>C^{\ast }=\sup_{y\in \mathbb{R}^{n}}\tilde{T}_{Het}(y)$
in the definition of $\alpha ^{\ast }$).

(iii) If Assumption \ref{R_and_X_tilde} holds, but $\tilde{T}_{Het}$ is
constant on $\mathbb{R}^{n}\backslash \mathsf{\tilde{B}}$, any rejection
region of the form $\{\tilde{T}_{Het}\geq C\}$ is trivial in that the
rejection region or its complement is a $\lambda _{\mathbb{R}^{n}}$-null
set. [This case can actually occur, see Remark \ref{rem:ex_const} in
Appendix \ref{app_C}.] If Assumption \ref{R_and_X_tilde} is violated, $%
\tilde{T}_{Het}$ is identically zero and a similar comment applies.
\end{remark}

\begin{example}
We continue the discussion of Example \ref{ex:bounded}. As noted prior to
Theorem \ref{Hetero_Robust_tilde_add}, any critical value $C\geq 2\tilde{d}%
_{1}^{-1}$ is size-controlling in a trivial way, but leads to trivial
rejection regions. We now show that the smallest size-controlling critical
value $C_{\Diamond }(\alpha )$ indeed leads to a non-trivial test (which, in
particular, has positive rejection probabilities in view of our Gaussianity
assumption and the fact that all $\Sigma \in \mathfrak{C}_{Het}$ are
positive definite). For this it suffices to verify the assumptions of
Theorem \ref{Hetero_Robust_tilde_add}. The first three assumptions have
already been verified above. From the calculations in Example \ref%
{ex:bounded} it is now easy to see that $C^{\ast }=\tilde{d}%
_{1}^{-1}n^{2}/[(n-1)^{2}+1]$, which is smaller than $2\tilde{d}%
_{1}^{-1}=\sup_{y\in \mathbb{R}^{n}}\tilde{T}_{Het}(y)$. This completes the
proof of the assertion. From Remark \ref{nontrivial_rej}(i) we furthermore
see that the rejection region $\{\tilde{T}_{Het}\geq C_{\Diamond }(\alpha
)\} $ has size equal to $\alpha $ if $\alpha \in (0,\alpha ^{\ast }]\cap
(0,1)$, and has size equal to $\alpha ^{\ast }$ if $\alpha ^{\ast }<\alpha
<1 $. Finally we note that size-controlling critical values that do not lead
to trivial tests must lie in the interval $[\tilde{d}%
_{1}^{-1}n^{2}/[(n-1)^{2}+1],2\tilde{d}_{1}^{-1})$ which is quite narrow as
it is contained in the interval $[\tilde{d}_{1}^{-1},2\tilde{d}_{1}^{-1})$.
\end{example}

While the situation in Example \ref{ex:bounded} is somewhat particular, the
example may perhaps contribute to a better understanding of the Monte Carlo
findings in \cite{DavidsonMacKinnon1985} and \cite{Godfrey2006}, namely that
the tests, obtained from $\tilde{T}_{Het}$ (employing HC0R-HC4R weights) in
conjunction with conventional critical values such as the $95\%$-quantile of
a chi-square distribution with appropriate degrees of freedom, can suffer
from severe underrejection under the null.

\begin{remark}
\label{rem:new} Another class of examples where $\tilde{T}_{Het}$ is bounded
is the case $q=k$ discussed in Remark \ref{q=k}. Recall from that remark
that in case $q=k>1$ condition (\ref{non-incl_Het_tilde}) is, however, never
satisfied and thus Theorem \ref{Hetero_Robust_tilde_add} is then not
applicable. We have not further investigated non-triviality of tests based
on $\tilde{T}_{Het}$ in case $q=k$ beyond the observations made in Remark %
\ref{q=k}(iii) that constancy of $\tilde{T}_{Het}$ on $\mathbb{R}%
^{n}\backslash \mathsf{\tilde{B}}$ is possible in case $q=k\geq 1$ and thus
then Remark \ref{nontrivial_rej}(iii) applies.
\end{remark}

\section{Generalizations\label{Generalizations}}

\subsection{Generalizations beyond Gaussianity\label{non-gaussianity}}

(i) All results in the preceding sections (as well as the extensions
described in Appendix \ref{app_d}) referring to properties under the null
hypothesis carry over as they stand to the situation where the error term $%
\mathbf{U}$ in (\ref{lm}) is elliptically symmetric distributed and has no
atom at zero, i.e., $\mathbf{U}$ is distributed as $\sigma \Sigma ^{1/2}%
\mathbf{z}$ where $\mathbf{z}$ has a spherically symmetric distribution on $%
\mathbb{R}^{n}$ that has no atom at zero.\footnote{%
Note that all results in the preceding sections (as well as the extensions
in Appendix \ref{app_d}), except for a few comments in Section \ref%
{sec:trivial}, are results referring to properties under the null hypothesis,%
} This is so since -- under this distributional model -- the null rejection
probabilities of any $G(\mathfrak{M}_{0})$-invariant rejection region
coincide with the corresponding null rejection probabilities under the
Gaussian model (i.e., where $\mathbf{z}$ is standard Gaussian); see the
discussion in Section 5.5 of \cite{PP2016} and Appendix E.1 of \cite{PP3}.%
\footnote{%
Note that all rejection regions considered in the preceding sections are $G(%
\mathfrak{M}_{0})$-invariant, because the test statistics considered are so.}
This implies, in particular, not only that the sufficient conditions for
size controllability under the above elliptically symmetric distributed
model as well as under the Gaussian model are the same, but that also the
numerical values of the size-controlling critical values coincide. As a
consequence, the algorithms for computing the size-controlling critical
values in the Gaussian case (used in Section \ref{numerical} and described
in Section \ref{algor} and Appendix \ref{app:algos}) can be used in the
above elliptically symmetric distributed case without any change whatsoever.
The same is actually true if $\mathbf{z}$ has a distribution in a certain
class larger than the class of spherical symmetric distributions with no
atom at zero, see Appendix E.1 of \cite{PP3}.

(ii) Furthermore, as discussed in detail in Appendix E.2 of \cite{PP3}, the
sufficient conditions for size controllability that we have derived under
Gaussianity also imply size controllability for many more forms of
distribution of $\mathbf{z}$ than those mentioned in (i); however, the
corresponding size-controlling critical values may then differ from the
size-controlling critical values that apply under Gaussianity.

(iii) Similarly as in Section 5.5 of \cite{PP2016}, the negative results
given in the preceding sections (as well as the ones described in Appendix %
\ref{app_d}) such as, e.g., size $1$ results, extend in a trivial way beyond
the Gaussian model as long as the maintained assumptions on the feasible
error distributions are weak enough to ensure that the implied (possibly
semiparametric) model, i.e., set of distributions for $\mathbf{Y}$, contains
the set given in (\ref{lm2}), but possibly contains also other distributions.

(iv) A further generalization beyond Gaussianity in the important special
case where $\mathfrak{C}=\mathfrak{C}_{Het}$ is as follows: Suppose $\mathbf{%
U}$ is distributed as $\sigma \Sigma ^{1/2}\limfunc{diag}(\mathbf{r)z}$
where $\mathbf{z}$ is standard normally distributed on $\mathbb{R}^{n}$ and
where the $n$-dimensional random vector $\mathbf{r}$ is independent of $%
\mathbf{z}$ with distribution $\rho $, where $\rho $ is a distribution on $%
(0,\infty )^{n}$. [This includes the case where the elements of $\limfunc{%
diag}(\mathbf{r)z}$ form an i.i.d. sample from a scale mixture of normals.]
Let $Q_{\mu ,\sigma ^{2}\Sigma ,\rho }$ denote the implied distribution for $%
\mathbf{Y}$ given by (\ref{lm}) where $\mu =X\beta $. Consider now instead
of (\ref{lm2}) the (semiparametric) model given by all distributions $Q_{\mu
,\sigma ^{2}\Sigma ,\rho }$ where $\mu \in \mathrm{\limfunc{span}}(X)$, $%
0<\sigma ^{2}<\infty $, $\Sigma \in \mathfrak{C}$, and $\rho $ is an
arbitrary distribution on $(0,\infty )^{n}$. Then the sufficient conditions
for size controllability derived under Gaussianity in earlier sections (and
in Appendix \ref{app_d}) also imply size controllability in this larger
model. In fact, the size-controlling critical values that apply under
Gaussianity deliver also size control under this more general model. This
follows from the following reasoning: Let $W$ be a Borel set in $\mathbb{R}%
^{n}$ such that $P_{\mu _{0},\sigma ^{2}\Sigma }(W)\leq \alpha $ for every $%
\mu _{0}\in \mathfrak{M}_{0}$, every $0<\sigma ^{2}<\infty $, and every $%
\Sigma \in \mathfrak{C}_{Het}$. Then for every such $\mu _{0}$, $\sigma ^{2}$%
, $\Sigma $, and every distribution $\rho $ on $(0,\infty )^{n}$ we have 
\begin{eqnarray*}
Q_{\mu _{0},\sigma ^{2}\Sigma ,\rho }(W) &=&\Pr (\mu _{0}+\sigma \Sigma
^{1/2}\limfunc{diag}(\mathbf{r)z}\in W)=\mathbb{E[}\Pr (\mu _{0}+\sigma
\Sigma ^{1/2}\limfunc{diag}(\mathbf{r)z}\in \left. W\right\vert \mathbf{r})]
\\
&=&\mathbb{E[}\Pr (\mu _{0}+\sigma _{\mathbf{r}}\Sigma _{\mathbf{r}}^{1/2}%
\mathbf{z}\in \left. W\right\vert \mathbf{r})]=\mathbb{E[}P_{\mu _{0},\sigma
_{\mathbf{r}}^{2}\Sigma _{\mathbf{r}}}(W)]\leq \alpha ,
\end{eqnarray*}%
where $\Sigma _{\mathbf{r}}^{1/2}:=\Sigma ^{1/2}\limfunc{diag}(\mathbf{r)/}%
s_{\mathbf{r}}$ with $s_{\mathbf{r}}$ denoting the positive square root of
the sum of the diagonal elements of $(\Sigma ^{1/2}\limfunc{diag}(\mathbf{r))%
}^{2}=\Sigma \limfunc{diag}^{2}(\mathbf{r)}$ and where $\sigma _{\mathbf{r}%
}=\sigma s_{\mathbf{r}}$. Here we have used that $P_{\mu ,\sigma _{\mathbf{r}%
}^{2}\Sigma _{\mathbf{r}}}(W)\leq \alpha $ by assumption since $\Sigma _{%
\mathbf{r}}=\Sigma \limfunc{diag}^{2}(\mathbf{r)/}s_{\mathbf{r}}^{2}\in 
\mathfrak{C}_{Het}$ and $0<\sigma _{\mathbf{r}}<\infty $ hold for every
realization of $\mathbf{r}$. In the above $\Pr $ denotes the probability
measure governing $(\mathbf{r},\mathbf{z)}$ and $\mathbb{E}$ the
corresponding expectation operator. \ As a consequence, the smallest
size-controlling critical value under Gaussianity is also the smallest
size-controlling critical value under the semiparametric model considered
here, as the latter model contains the Gaussian model as a submodel. [In the
special case where $\limfunc{diag}(\mathbf{r)}$ is a (random) multiple of
the identity matrix $I_{n}$, the assumption $\mathfrak{C}=\mathfrak{C}_{Het}$
is superfluous as then $\Sigma _{\mathbf{r}}=\Sigma $, which by assumption
belongs to the given $\mathfrak{C}$. In this case $\mathbf{U}$ satisfies the
assumptions in (i), and hence (iv) adds little new, except that -- in
contrast to (i) -- the reasoning works without use of $G(\mathfrak{M}_{0})$%
-invariance.]

(v) It is apparent from the reasoning in (iv) that Gaussianity of $\mathbf{z}
$ can be replaced by any other distributional assumption for which size
controllability has already been established. E.g., one can in (iv) choose $%
\mathbf{z}$ to have a spherically symmetric distribution without an atom at
zero or to have a distribution in the more general class mentioned in (i)
(note that all relevant rejection regions discussed in earlier sections are $%
G(\mathfrak{M}_{0})$-invariant and thus (i) applies). In a similar vein, one
can combine the results in Appendix E.2 of \cite{PP3} discussed in (ii)
above with the reasoning outlined in (iv). We abstain from presenting
details.

\subsection{Generalizations to stochastic regressors\label{stoch_regr}}

The assumption of nonstochastic regressors can be easily relaxed as follows:
Suppose $X$ is random and $\mathbf{U}$ is conditionally on $X$ distributed
as $N(0,\sigma ^{2}\Sigma )$, with $\sigma ^{2}=\sigma ^{2}(X)>0$ and $%
\Sigma =\Sigma (X)\in \mathfrak{C}_{Het}$ where $\sigma ^{2}(\cdot )$ and $%
\Sigma (\cdot )$ may vary in given classes of functions. The size control
results such as Theorems \ref{Hetero_Robust} and \ref{Hetero_Robust_tilde}
can then obviously be applied after one conditions on $X$ provided almost
all realizations of $X$ satisfy the assumptions of those theorems, which
will typically be the case (for brevity we do not provide a formal statement
here).\footnote{%
An appropriately modified statement applies to the size control results in
Appendix \ref{app_d}.} The resulting conditional size control statements
then immediately imply that the so-obtained conditional size-controlling
critical values $C=C(\alpha ,X)$ also control size unconditionally. Size $1$
results such as, e.g., Propositions \ref{rem_C*}, \ref{rem_necessity}, or %
\ref{rem_C*_tilde} also extend to conditional size $1$ results in a similar
manner provided $\sigma ^{2}(X)$ and $\Sigma (X)$ vary independently through
all of $(0,\infty )$ and $\mathfrak{C}_{Het}$, respectively, for (almost)
every realization of $X$, when the functions $\sigma ^{2}(\cdot )$ and $%
\Sigma (\cdot )$ vary in the before mentioned function classes.\footnote{%
See Footnote 40 in \cite{PPBoot} for a discussion of sufficient conditions.}
Generalizations to non-Gaussianity similarly as discussed in Section \ref%
{non-gaussianity} are also possible in the present context.

\section{Results for other classes of tests\label{sec_extensions}}

The results in Sections \ref{sec_size_control} and \ref%
{sec_size_control_tilde} (and in Appendix \ref{app_d}) have been obtained
with the help of a general theory developed in Section 5 of \cite{PP2016},
Section 5 of \cite{PP3}, and Section 3.1 of \cite{PP4} that covers a very
broad class of test statistics (and actually allows also for correlated
errors). We note that, like in Section \ref{non-gaussianity}, Gaussianity is
again not essential for a good portion of this general theory, see Section
5.5 of \cite{PP2016} as well as Appendix E of \cite{PP3}.\footnote{%
Also arguments like in (iv) and (v) of Section \ref{non-gaussianity} can be
applied to try to obtain generalizations.} We next discuss a few further
situations that can also be handled by the general theory just mentioned but
we refrain from spelling out the details:\footnote{%
Applying some of the main results of this general theory (e.g., Corollary
5.6 or Proposition 5.12 of \cite{PP3}) will require one to determine the set 
$\mathbb{J}(\mathcal{L},\mathfrak{C})$ defined in Appendix \ref{app_char}.
For the important cases $\mathfrak{C}=\mathfrak{C}_{Het}$ and $\mathfrak{C}=%
\mathfrak{C}_{(n_{1},\ldots ,n_{m})}$ (defined in Appendix \ref{app_d}),
this is already accomplished in Propositions \ref{characterization} and \ref%
{characterization_2} in Appendix \ref{app_char} below.}

(i) The test statistic considered is an OLS-based test statistic like $%
T_{Het}$, but where $\hat{\Omega}_{Het}$ is now replaced by an appropriate
estimator derived from a given (possibly misspecified) \emph{parametric}
heteroskedasticity model described by a parameter vector $\theta $.

(ii) The test statistic is a Wald-type test statistic based on a (feasible)
generalized least-squares estimator together with an appropriate covariance
matrix estimator based on a given (possibly misspecified) parametric model.
[This includes the (quasi-)maximum likelihood estimator (provided $\theta $
is unrelated to $\beta $).] Alternatively, the test statistic is the
(quasi-)likelihood ratio or (quasi-)score test statistic based on this
parametric model.

(iii) The test statistic is a Wald-type test statistic as in (ii), except
that the covariance matrix estimator is now nonparametric (in the spirit of
heteroskedasticity robust testing) as described in \cite{RomanoWolf2017}.
See also \cite{Cragg_1983, Cragg_1992}, \cite{Flachaire_2005b}, \cite%
{Wooldridge2010, Wooldridge2012}, \cite{RomanoWolf2017}, \cite{Lin_Chou_2018}%
, \cite{DiCiccio_Romao_Wolf_2019}.

\section{Some comments on power\label{sec:power}}

Under our maintained assumptions, heteroskedasticity robust tests based on $%
T_{Het}$ or $T_{uc}$ (using an arbitrary critical value $C$, including
size-controlling ones) have positive power everywhere in the alternative
(cf.~the discussion at the beginning of Section \ref{sec:trivial}). These
tests can furthermore be shown to have power that goes to one as one moves
away from the null hypothesis along sequences $(\mu _{l},\sigma
_{l}^{2},\Sigma _{l})$ where $\mu _{l}$ moves further and further away from $%
\mathfrak{M}_{0}$ (the affine space of means described by the restrictions $%
R\beta =r$) in an orthogonal direction as $l\rightarrow \infty $, where $%
\sigma _{l}^{2}$ converges to some finite and positive $\sigma ^{2}$, and $%
\Sigma _{l}$ converges to a \emph{positive definite} matrix. Despite of what
has just been said, these tests can have, in fact not infrequently will
have, \emph{infimal} power equal to zero if $\mathfrak{C}$ is sufficiently
rich, e.g., if $\mathfrak{C}=\mathfrak{C}_{Het}$; cf. Theorem 4.2 in \cite%
{PP2016}, Lemma 5.11 in \cite{PP3}, and Theorem 4.2 in \cite{PP4}. [This
does not contradict the before mentioned result as for this result sequences 
$\Sigma _{l}$ that converge to a singular matrix as $l\rightarrow \infty $
were ruled out.]

For tests based on $\tilde{T}_{Het}$ or $\tilde{T}_{uc}$ the situation is
somewhat different. As shown in Section \ref{sec:trivial}, tests based on $%
\tilde{T}_{Het}$ or $\tilde{T}_{uc}$ can be trivial for some choices of
critical values $C$ (and then will have power zero everywhere in the
alternative). However, if $C$ is chosen to be the smallest size-controlling
critical value (provided it exists), the resulting tests obtained form $%
\tilde{T}_{Het}$ or $\tilde{T}_{uc}$ will typically have positive power
(under appropriate assumptions). In particular, then the test based on $%
\tilde{T}_{uc}$ has the same power function as the test based on $T_{uc}$
that uses its smallest size-controlling critical value, provided the latter
exists, see Section \ref{sec:trivial_uc_tilde}. We have not further
investigated the power properties of the tests based on $\tilde{T}_{Het}$ in
any more detail on a theoretical level. The numerical results in Section \ref%
{sec:powernum} seem to suggest that for these tests power may not go to one
along sequences $(\mu _{l},\sigma _{l}^{2},\Sigma _{l})$ as mentioned above:
in fact, power does not rise above the significance level $\alpha $ in some
examples (on the range of alternatives considered). This feature makes tests
based on $\tilde{T}_{Het}$ rather undesirable.

\section{Computing the size and smallest size-controlling critical values 
\label{algor}}

Consider a testing problem as in Equation (\ref{testing problem}) with $%
\mathfrak{C}=\mathfrak{C}_{Het}$ and let $T$ be one of the test statistics
considered in the present article (e.g., $T_{Het}$ with some choice for the
weights $d_{i}$). Suppose we want to numerically determine the size of the
test with rejection region $\left\{ T\geq C\right\} $ for some user-supplied
critical value $C$, i.e., we want to determine%
\begin{equation}
\sup_{\mu _{0}\in \mathfrak{M}_{0}}\sup_{0<\sigma ^{2}<\infty }\sup_{\Sigma
\in \mathfrak{C}_{Het}}P_{\mu _{0},\sigma ^{2}\Sigma }(T\geq C).  \label{sz0}
\end{equation}%
Now, for all test statistics $T$ considered in the present article, this can
be simplified to%
\begin{equation}
\sup_{\Sigma \in \mathfrak{C}_{Het}}P_{\mu _{0},\Sigma }(T\geq C)  \label{sz}
\end{equation}%
where, subject to $\mu _{0}\in \mathfrak{M}_{0}$, $\mu _{0}$ can be chosen
as desired. This is due to invariance properties of $T$, cf.~Remarks \ref%
{rem:GM0} and \ref{rem:tildeGM0}. The quantity in (\ref{sz}) can now be
approximated numerically by any maximization algorithm where the
probabilities are evaluated by Monte-Carlo methods or by the algorithm
described in \cite{davies} in case $q=1$, cf.~Appendix \ref{sec:obsq1}.%
\footnote{%
Alternative to \cite{davies} other algorithms like Imhof's algorithm,
etc.~can be used, some of which are also implementented in the R-package 
\textbf{CompQuadForm} (\cite{Duchesne}).}

Suppose next that we want to numerically determine the smallest
size-controlling critical value $C_{\Diamond }(\alpha )\in \mathbb{R}$ ($%
0<\alpha <1$) when using the test statistic $T$. [We assume here that the
user knows that the smallest size-controlling critical value indeed exists,
e.g., because the user has checked that the sufficient conditions developed
in the present article hold, or because of other reasoning as, e.g., used in
Example \ref{ex_k_pop}.] Then, in view of (\ref{sz0}) and (\ref{sz}), we
need to compute $C_{\Diamond }(\alpha )$ as the smallest real number $C$ for
which 
\begin{equation}
\sup_{\Sigma \in \mathfrak{C}_{Het}}P_{\mu _{0},\Sigma }(T\geq C)\leq \alpha
\label{eqn:sizesalpha}
\end{equation}%
holds. The quantity to the left in (\ref{eqn:sizesalpha}) is non-increasing
in the critical value $C$. Hence, to determine the smallest size-controlling
critical value $C_{\Diamond }(\alpha )$, any line-search algorithm (in
combination with an algorithm to determine the sizes as described before)
can be used to compute $C_{\Diamond }(\alpha )$. We stress that it is of
foremost importance to know that the testing problem at hand actually allows
for size control before one attempts to numerically determine $C_{\Diamond
}(\alpha )$. Hence, the theoretical results of the present article are of
paramount importance also for the algorithmic aspect of the problem.

The specific algorithms we use to determine size and size-controlling
critical values in our numerical studies are based on the above observations
and are described in detail in Appendix \ref{app:algos}. They are made
available in the \textsf{R}-package \textbf{hrt }(\cite{hrt}) for the
convenience of the user. The numerical procedures we use are heuristic in
nature. Questions of efficacy of these algorithms or about theoretical
guarantees are certainly important, but are beyond the scope of the present
article.

Determining smallest size-controlling values numerically is important, e.g.,
if one wants to compare their magnitude with that of standard critical
values in some special cases, as we do inter alia in the next section, or if
one wants to obtain a confidence interval. However, a user who has observed
the data and only wants to decide whether or not to reject the null
hypothesis at significance level $\alpha $ ($0<\alpha <1$) when using $T$
combined with the smallest size-controlling critical value $C_{\Diamond
}(\alpha )$, can actually perform this test without needing to compute $%
C_{\Diamond }(\alpha )$: Let $y_{obs}$ be the observed data. Define the
\textquotedblleft maximal p-value" as%
\begin{eqnarray}
p(y_{obs}) &=&\sup_{\mu _{0}\in \mathfrak{M}_{0}}\sup_{0<\sigma ^{2}<\infty
}\sup_{\Sigma \in \mathfrak{C}_{Het}}P_{\mu _{0},\sigma ^{2}\Sigma }(\{z\in 
\mathbb{R}^{n}:T(z)\geq T(y_{obs})\})  \notag \\
&=&\sup_{\Sigma \in \mathfrak{C}_{Het}}P_{\mu _{0},\Sigma }(\{z\in \mathbb{R}%
^{n}:T(z)\geq T(y_{obs})\}),  \label{p}
\end{eqnarray}%
where the second equality in the display follows from the invariance
properties mentioned before (and $\mu _{0}\in \mathfrak{M}_{0}$ can be
chosen as desired). It is now not difficult to see that $p(y_{obs})\leq
\alpha $ is equivalent to $T(y_{obs})\geq C_{\Diamond }(\alpha )$. That is,
rejecting if and only if $p(y_{obs})\leq \alpha $ leads to exactly the same
test as rejecting if and only if $T(y_{obs})\geq C_{\Diamond }(\alpha )$,
with the former description having the advantage that the more costly
computation of $C_{\Diamond }(\alpha )$ can be avoided. What needs to be
computed is (\ref{p}), which, however, is nothing else than the size of the
test when using the \textquotedblleft critical value\textquotedblright\ $%
T(y_{obs})$. Hence, $p(y_{obs})$ can be determined by any algorithm that
determines the size (\ref{sz}) for the user-supplied \textquotedblleft
critical value\textquotedblright\ $C=T(y_{obs})$. In particular, the routine
\textquotedblleft size\textquotedblright\ provided in the \textsf{R}-package 
\textbf{hrt }(\cite{hrt}) can be used for this purpose. Note that checking
whether $p(y_{obs})\leq \alpha $ avoids the line-search part (as outlined
following (\ref{eqn:sizesalpha})), and is thus computationally \emph{more
efficient} than first determining $C_{\Diamond }(\alpha )$ (as outlined
above) and then checking whether $T(y_{obs})\geq C_{\Diamond }(\alpha )$.

Finally, we note that if (contrary to what we assume in this section) no
size-controlling critical value exists for a given significance level $%
\alpha \in (0,1)$, then the maximal p-value in (\ref{p}) is larger than $%
\alpha $ for \emph{every} possible observed value $y_{obs}$, and the
corresponding test thus never rejects and thus is uninformative. Hence,
while the explicit computation of a smallest size-controlling critical value
can be avoided for performing a single test, knowing its existence is
important as then the resulting test is guaranteed to be informative
(non-trivial) if $T_{Het}$ or $T_{uc}$ is being used; and the same is true
for $\tilde{T}_{Het}$ and $\tilde{T}_{uc}$ under the conditions discussed in
Section \ref{sec:trivial}.

We also note that in view of the discussion in Section \ref{non-gaussianity}
the algorithms for computing null rejection probabilities, size, and
smallest size-controlling critical values discussed in this section and
Appendix \ref{app:algos} remain valid for elliptically symmetric distributed
data without any need for modification. With regard to computing size and
smallest size-controlling critical values, the same is also true for the
semiparametric model described in (iv) of Section \ref{non-gaussianity}.

\section{Numerical results \label{numerical}}

In this section we pursue two goals:

\begin{enumerate}
\item In Subsection \ref{sec:chibad} we show numerically that any of the
usual heteroskedasticity robust tests can suffer from overrejection of the
null hypothesis (sometimes by a large margin) when they are based on
conventional critical values. While this adds to similar evidence already
present in the literature for the HC0-HC4 based tests (see Section \ref%
{Intro}), this seems to be a new observation for the HC0R-HC4R based tests.
In any case, this drives home the point that none of these
heteroskedasticity robust tests based on conventional critical values comes
with a guarantee that size is controlled by the nominal significance level $%
\alpha $. Consequently, instead of using conventional critical values, this
strongly suggests to use (smallest) size-controlling critical values as
investigated in this paper.

\item In Subsection \ref{sec:powernum} we then numerically compute smallest
size-controlling critical values and study the power behavior of tests based
on such size-controlling critical values in some examples.
\end{enumerate}

In this section (and in the attending Appendices \ref{app:algos} and \ref%
{details}) we shall often refer to $T_{Het}$ as HC0-HC4 when we want to
stress that the weights $d_{i}$ being used are the HC0-HC4 weights,
respectively, see Section \ref{sec_teststatistic}. Similarly, we shall refer
to $\tilde{T}_{Het}$ as HC0R-HC4R when the HC0R-HC4R weights are used, see
Section \ref{sec_size_control_tilde}. For reasons of uniformity of notation,
we shall then often denote $T_{uc}$ as UC and $\tilde{T}_{uc}$ as UCR.
Furthermore, throughout this section we consider the heteroskedastic
Gaussian linear model with $\mathfrak{C}=\mathfrak{C}_{Het}$ as introduced
in Section \ref{frame}; in particular, the notion of size in the present
section (and the attending appendices) \emph{always} refers to this model.

The algorithms for computing rejection probabilities, the size of a test,
and size-controlling critical values used in the before-mentioned numerical
computations are described in Section \ref{algor} and Appendix \ref%
{app:algos}. Implementations are available as an R-package \textbf{hrt} (%
\cite{hrt}).

\subsection{Tests based on conventional critical values\label{sec:chibad}}

We consider the important case $q=1$, and first illustrate numerically that
none of the test statistics UC, HC0-HC4, UCR, and HC0R-HC4R combined with
the critical value $C_{\chi ^{2},0.05}\approx 3.8415$ results in a test that
is \emph{guaranteed} to have size less than or equal to $\alpha =0.05$. This
is achieved by providing instances of design matrices $X$ and of hypotheses,
described by $(R,r)$, such that the respective test has size larger than the
nominal significance level $\alpha =0.05$, often by a large margin. Here $%
C_{\chi ^{2},0.05}$ denotes the $95\%$-quantile of a chi-square-distribution
with $1$ degree of freedom. [This critical value has a justification for use
with HC0-HC4 or HC0R-HC4R via asymptotic considerations, but, in general,
there is no such justification for use with UC or UCR, which we nevertheless
include here for completeness.\footnote{%
Of course, in the special case of homoskedasticity, the before mentioned
justification also applies to UC and UCR.}] That is, in the instances we
exhibit, this conventional critical value turns out to be \emph{too small}.
We next show similar results for other suggestions of critical values, e.g.,
for \textquotedblleft degree-of-freedom\textquotedblright\ adjustments to
the conventional chi-square based critical value such as the Bell-McCaffrey
adjustment (\cite{BellMcCa}, \cite{Imbkoles2016}). It is important to note
here that in all the instances mentioned our conditions for
size-controllability are satisfied, showing that size-controlling critical
values can actually be found; hence, the overrejection problems mentioned
before are \emph{not} intrinsic problems, but only reflect the fact that
conventional critical values can be a bad choice and do not guarantee size
control. [In the present context it is worth recalling that for the test
statistics HC0R-HC4R we have already shown in Example \ref{ex:bounded} in
Section \ref{sec:trivial} that other situations can be found in which
conventional critical values such as, e.g., $C_{\chi ^{2},0.05}$ are \emph{%
too large}, as the resulting tests reject with probability zero only (under
the null as well as under the alternative), rendering these tests useless.]

To uncover instances where the conventional critical value $C_{\chi
^{2},0.05}$ is too small, we make use of the following observation: In case
a given test statistic from the above list (together with a given design
matrix $X$ and hypothesis described by $(R,r)$) is such that the lower bound 
$C^{\ast }$ on size-controlling critical values obtained in Propositions \ref%
{rem_C*} (\ref{rem_C*_tilde}, respectively) exceeds $C_{\chi ^{2},0.05}$, we
are done, as we then know that the critical value $C_{\chi ^{2},0.05}$ leads
to a test that has size $1$. [As noted subsequent to Theorems \ref%
{Hetero_Robust} and \ref{Hetero_Robust_tilde}, the value of $r$ actually
plays no role here, and we may set it to zero.]

Since the lower bounds $C^{\ast }$ for size-controlling critical values in
Propositions \ref{rem_C*} (\ref{rem_C*_tilde}, respectively) depend on the
given test statistic, on $X$ and on $R$, we may -- for any given choice of
test statistic and any given $R$ -- numerically search for particularly
\textquotedblleft hostile\textquotedblright\ design matrices, i.e., for
design matrices for which the lower bound is large, to see whether matrices $%
X$ exist for which the lower bound exceeds $C_{\chi ^{2},0.05}$. We only do
this for $k=2$, $R=(0,1)$, $r=0$, and $n=25$, and restrict ourselves to
matrices $X$ with first column representing an intercept. The concrete
search used is detailed in Appendix \ref{app:numO2}, see Algorithm \ref%
{alg:O2} in particular. Table \ref{fig:lowbd} provides, for every test
statistic considered, the lower bound $C^{\ast }$ corresponding to the most
\textquotedblleft hostile\textquotedblright\ design matrix found by the
search. [As the searches are run separately for each test statistic, the
resulting \textquotedblleft hostile\textquotedblright\ design matrices will
typically differ across the runs.]\footnote{%
Since, for example, HC0 is a multiple of HC1, where the factor is $%
n/(n-k)=1.09$, we know that the \textquotedblleft hostile\textquotedblright\
design matrix obtained from the search for HC1 leads to a $C^{\ast }$-value
of $1.09\ast 1711.19=1865.20$ for HC0, larger than the value 95.56 obtained
from the search for HC0, cf.~Table \ref{fig:lowbd}. We could have reported
this larger value, but decided to present the raw results from our searches
as this is sufficient for our purposes. We also note that our search
procedure detailed in Appendix \ref{app:numO2} does not seriously attempt to
optimize the $C^{\ast }$-value (for every one of the test statistics
considered) over the set of all feasible $X$, but is only a crude search for
finding a matrix resulting in a $C^{\ast }$-value sufficiently large for our
purposes.} 
\begin{table}[th]
\parbox{.45\linewidth}{
\centering
\begin{tabular}{ll||ll}
\hline
UC &  731.60 & UCR & 23.59 \\ 
  HC0 &   95.56 & HC0R &  1.08 \\ 
  HC1 & 1711.19 & HC1R &  1.04 \\ 
  HC2 &   52.23 & HC2R &  1.04 \\ 
  HC3 &    1.00 & HC3R &  1.00 \\ 
  HC4 &    1.02 & HC4R &  1.04 \\ 
\hline
\end{tabular}
\caption{$C^*$ under respective "hostile" $X$.} 
\label{fig:lowbd}
} \hfill 
\parbox{.45\linewidth}{
\centering
\begin{tabular}{ll||ll}
\hline
UC & 0.98 & UCR & 0.98 \\ 
  HC0 & 0.99 & HC0R & 0.16 \\ 
  HC1 & 1.00 & HC1R & 0.17 \\ 
  HC2 & 0.99 & HC2R & 0.17 \\ 
  HC3 & 0.19 & HC3R & 0.14 \\ 
  HC4 & 0.11 & HC4R & 0.10 \\ 
\hline
\end{tabular}
\caption{"Worst-case" sizes using $C_{\chi^2,0.05}$.} 
\label{fig:size}
}
\end{table}
In combination with the theoretical results from Propositions \ref{rem_C*}
and \ref{rem_C*_tilde}, Table \ref{fig:lowbd} shows that for some design
matrices $X$ the critical value $C_{\chi ^{2},0.05}\approx 3.8415$ results
in a test with size equal to $1$ when combined with UC, HC0-HC2, and also
with UCR. [This is so despite the fact that, for any of the twelve test
statistics considered, the sufficient conditions for size-control in the
pertaining theorems in Sections \ref{sec_size_control} and \ref%
{sec_size_control_tilde_2} are satisfied for all relevant $X$ matrices
encountered in the numerical procedure (as we have checked), and hence it is
known that size-controlling critical values exist in all these situations!]
Table \ref{fig:lowbd} is not informative about the size of the remaining
seven tests, since the corresponding entries in that table are all less than 
$C_{\chi ^{2},0.05}$. To obtain insight into the sizes of the remaining
seven tests we do the following: for each of the tests we numerically
compute the size for various instances of design matrices (the ones that
give rise to Table \ref{fig:lowbd}) and report the largest one of these
sizes (\textquotedblleft worst case\textquotedblright\ sizes) in Table \ref%
{fig:size}.\footnote{%
Of course, considering additional design matrices $X$ would potentially lead
to even larger sizes.} We actually do this for all twelve tests considered.
The algorithm used in the size computation is the implementation of
Algorithm \ref{alg:S1} in the R-package \textbf{hrt }(\cite{hrt}), cf. the
description in Appendices \ref{app:S} and \ref{app:numO2}. Table \ref%
{fig:size} now clearly shows that for \emph{every} test statistic considered
an instance can be found, in which the size of the test (when using the
critical value $C_{\chi ^{2},0.05}$) clearly exceeds the nominal
significance level $\alpha =.05$. The lowest value in that table is attained
by HC4R, but a size of $0.10$ is still twice the nominal significance level $%
\alpha $.

We note that the numbers shown in Table \ref{fig:size} actually only
represent numerically determined lower bounds for the actual sizes, as their
computation involves (for any given $X$) a numerical search procedure (over
the set $\mathfrak{C}_{Het}$) for the worst-case null rejection probability;
that is, the numbers shown in Table \ref{fig:size} correspond to the null
rejection probability computed from a \textquotedblleft
bad\textquotedblright\ covariance matrix $\Sigma $, but potentially not for
the \textquotedblleft worst\textquotedblright\ possible one. [In this
process, for any given $\Sigma \in \mathfrak{C}_{Het}$, we have to
numerically compute the null rejection probability, which can be done quite
accurately in case $q=1$ by algorithms like the Davies algorithm, see
Appendix \ref{sec:obsq1} as well as Appendix \ref{app:S}.] In particular,
the entries in the $0.98$-$0.99$ range in Table \ref{fig:size} are
numerically determined lower bounds for the size, which, in fact, we know to
be equal to $1$ in light of Table \ref{fig:lowbd}. [We could have used this
knowledge to replace the entries in question in Table \ref{fig:size} by $1$,
but we decided otherwise in order to showcase the concrete outcome of the
numerical algorithm that has been run. Of course, one could also improve
this outcome by using a higher accuracy parameter in the optimization
procedures involved.]

Sometimes -- without much theoretical justification in general -- it is
suggested in the literature to replace $C_{\chi ^{2},0.05}$ by the $95\%$%
-quantile of an $F_{1,n-k}$-distribution, which is approximately $4.28$ in
the situation considered here ($n-k=23$). Obviously, from Table \ref%
{fig:lowbd} we see that the conclusions regarding UC, HC0-HC2, and UCR
remain the same when this critical value is used. Repeating the exercise
that has led to Table \ref{fig:size}, but with $C_{\chi ^{2},0.05}$ replaced
by the $95\%$-quantile of an $F_{1,n-k}$-distribution, gives Table \ref%
{fig:size_2}, leading essentially to the same conclusions. 
\begin{table}[th]
\centering
\begin{tabular}{ll||ll}
\hline
UC & 0.98 & UCR & 0.98 \\ 
HC0 & 0.99 & HC0R & 0.15 \\ 
HC1 & 1.00 & HC1R & 0.16 \\ 
HC2 & 0.98 & HC2R & 0.15 \\ 
HC3 & 0.18 & HC3R & 0.13 \\ 
HC4 & 0.09 & HC4R & 0.08 \\ \hline
\end{tabular}%
\caption{"Worst-case" sizes using F-critical value.}
\label{fig:size_2}
\end{table}

\textquotedblleft Degree-of-freedom\textquotedblright\ adjustments to the
conventional chi-square based critical value such as the Bell-McCaffrey
adjustment (\cite{BellMcCa}) have been discussed in the literature. In
particular, \cite{Imbkoles2016} suggested to use this adjustment with the
HC2 statistic. We have repeated the above exercise that has led to the entry
for HC2 in Table \ref{fig:size}, but with $C_{\chi ^{2},0.05}$ replaced by
the Bell-McCaffrey adjustment. For the computation of the Bell-McCaffrey
adjustment we relied on the R-package \textbf{dfadjust} (\cite{dfadjust}).
For the resulting test, the largest size that was found in our computations
was $0.24$, which is more than four times the nominal significance level. It
transpires that this adjustment does also not come with a size-guarantee.

We conclude here by stressing that the negative findings in this subsection
were obtained in a very simple model with only two regressors and where only
one of the parameters is subject to test. For more complex models and test
problems the size distortions may even be worse.

\subsection{Power comparison of tests based on size-controlling critical
values\label{sec:powernum}}

A power comparison of two tests, both conducted at a given \emph{nominal }%
significance level $\alpha $, makes sense only if both tests actually are
level $\alpha $ tests, i.e., if both tests have a size not exceeding the
given $\alpha $. For this reason, we now compare the tests obtained from the
statistics UC, HC0-HC4, UCR, HC0R-HC4R only when respective \emph{smallest
size-controlling critical values} are used. Our theoretical results
concerning the existence of size-controlling critical values, together with
the algorithms for their computation in Appendix \ref{app:algos}, allow for
such a comparison in terms of power. In all cases considered in this section 
$q=1$ will hold.

Throughout, in addition to the power functions of the before-mentioned
tests, we also show as a benchmark the power function of the \emph{infeasible%
} (i.e., oracle) GLS-based $F$-test conducted at the $5\%$-significance
level, that makes use of knowledge of $\Sigma $. For given $\Sigma \in 
\mathfrak{C}_{Het}$, the distribution of this infeasible GLS-based $F$-test
statistic is (under $P_{X\beta ,\sigma ^{2}\Sigma }$ with $\beta \in \mathbb{%
R}^{k}$, $\sigma ^{2}\in (0,\infty )$) a noncentral $F_{1,n-k}$-distribution
with noncentrality parameter $\delta ^{2}$, where%
\begin{equation*}
\delta =(R(X^{\prime }\Sigma ^{-1}X)^{-1}R^{\prime })^{-1/2}(R\beta
-r)/\sigma .
\end{equation*}%
Since the power functions of all the tests considered in our study depend on
the parameters $\beta $, $\sigma ^{2}$, and $\Sigma $ only through $(R\beta
-r)/\sigma $ and $\Sigma $ (because of $G(\mathfrak{M}_{0})$-invariance and
Proposition 5.4 in \cite{PP2016}), and thus depend only on $\delta $ and $%
\Sigma $, we shall -- for given $\Sigma $ -- present all these power
functions as a function of $\delta $. We show only results for $\delta \geq
0 $, as the power functions in fact depend on $\delta $ only through $%
\left\vert \delta \right\vert $ (for given $\Sigma $); see Proposition 5.4
in \cite{PP2016}.

\subsubsection{Comparing the means of two heteroskedastic groups\label%
{sec:powtwogroups}}

As a practically relevant example, we here compare the power of tests based
on size-controlling critical values in the context of Example \ref%
{ex_fish_behr}. That is, we treat the problem of comparing the means of two
heteroskedastic groups (e.g., a treatment and a control group), the null
hypothesis being that the difference of expected outcomes in each group is
zero. We consider the case where $n=30$ and $\alpha =0.05$. Furthermore, we
vary the size $n_{1}$ of the first group ($n_{1}\in \{3,9,15\}$),
corresponding to a \textquotedblleft strongly unbalanced\textquotedblright ,
\textquotedblleft moderately unbalanced\textquotedblright , and
\textquotedblleft balanced\textquotedblright\ design, respectively. We
compute the power for a number of covariance matrices $\Sigma _{a}$ given as
follows: For $a=1,5,9$ define%
\begin{equation*}
\Sigma _{a}=10^{-1}\limfunc{diag}\left( \frac{a}{n_{1}},\ldots ,\frac{a}{%
n_{1}},\frac{10-a}{n-n_{1}},\ldots ,\frac{10-a}{n-n_{1}}\right) \in 
\mathfrak{C}_{Het},
\end{equation*}%
where the first $n_{1}$ (and last $n-n_{1}$, respectively) diagonal entries
of each $\Sigma _{a}$ are constant. That is, we look at power functions
evaluated at covariance matrices under which the subjects in the same group
actually have the same variances. [For brevity we do not report power
functions for covariance matrices not sharing this property.] For the
balanced design, we note that $\Sigma _{1}$ and $\Sigma _{9}$ lead to the
same power of each test (but we report all results for completeness), and
that $\Sigma _{5}$ corresponds to homoskedasticity.

The critical values are chosen in each case as the smallest critical value
guaranteeing size control over $\mathfrak{C}_{Het}$ (implying, of course,
that the corresponding tests can have null rejection probabilities smaller
than $\alpha $ for the covariance matrices $\Sigma _{a}$ considered). The
existence of said critical values follows from our theory and is discussed
in detail in Example \ref{ex_fish_behr} for the test statistics UC and
HC0-HC4; in particular, all assumptions of Theorems \ref{Hetero_Robust} are
satisfied. For UCR the existence is guaranteed by Part (a) of Theorem \ref%
{Hetero_Robust_tilde}. With regard to the test statistics HC0R-HC4R, note
that Assumption \ref{R_and_X_tilde} is satisfied since $e_{i}(n)\notin 
\mathfrak{M}_{0}^{lin}=\mathfrak{M}_{0}=\limfunc{span}((1,\ldots ,1)^{\prime
})$ for every $i=1,\ldots ,n$ as $n=30>k=2$. This also shows that the
sufficient condition for size control (\ref{non-incl_Het_tilde}) is
satisfied as $\mathsf{\tilde{B}}=\mathfrak{M}_{0}$ is easily verified and
since one may set $\mu _{0}=0$. We have verified the non-constancy
assumption on the test statistics HC0R-HC4R in Theorem \ref%
{Hetero_Robust_tilde} numerically. As a consequence, all assumptions of Part
(b) of Theorem \ref{Hetero_Robust_tilde} are satisfied.

We note that some of the test statistics differ from each other only by a
known multiplicative constant and hence are equivalent in the sense that
they give rise to the same test \emph{when the respective smallest
size-controlling critical value is employed}, see Remarks \ref{rem:equiv}
and \ref{rem:equiv_tilde}: In the unbalanced case ($n_{1}\in \{3,9\}$),\ HC0
and HC1 are equivalent in this sense, as are HC0R-HC4R (the latter is so
since $\tilde{h}_{ii}=1/n$ which does not depend on $i$). In the balanced
case ($n=15$), UC and HC0-HC4 are all equivalent, and the same is true for
UCR and HC0R-HC4R as is not difficult to see. Furthermore, in the balanced
as well as in the unbalanced case, the rejection regions of the tests based
on UC and UCR coincide essentially (i.e., up to a $\lambda _{\mathbb{R}^{n}}$%
-null set) as a consequence of the relationship established in Section \ref%
{sec:trivial_uc_tilde}. In particular, it follows that in the balanced case, 
\emph{all} tests considered (essentially) coincide. We nevertheless compute
the power functions for each of the tests separately without making use of
the noted equivalencies; this provides a double-check of our numerical
results.\footnote{%
The equivalencies mentioned in this paragraph for the two-group-comparison
problem analogously hold for general $n$, $n_{1}$, and $n_{2}$ as is easily
seen.}

Numerically the critical values were determined through the implementation
of Algorithms \ref{alg:S1} and \ref{alg:CV1} in the R-package \textbf{hrt} (%
\cite{hrt}) version 1.0.0, and the power functions were computed with the
implementation of the algorithm by \cite{davies} in the R-package \textbf{%
CompQuadForm} (\cite{Duchesne}) version 1.4.3; see Appendices \ref{app:S}
and \ref{app:CVcomppow} for more details. For the sake of illustration, we
also report the critical values obtained for every test considered and every
balancedness condition in Table \ref{fig:CV}. 
\begin{table}[th]
\centering%
\begin{tabular}{l|llllll}
\hline
$n_{1}$ & UC & HC0 & HC1 & HC2 & HC3 & HC4 \\ \hline
3 & 225.97 & 26.69 & 25.63 & 17.48 & 11.86 & 5.43 \\ 
9 & 12.70 & 5.80 & 5.39 & 5.10 & 4.55 & 4.70 \\ 
15 & 4.59 & 4.91 & 4.58 & 4.59 & 4.28 & 4.58 \\ \hline\hline
$n_{1}$ & UCR & HC0R & HC1R & HC2R & HC3R & HC4R \\ 
3 & 25.82 & 3.25 & 3.14 & 3.14 & 3.02 & 3.13 \\ 
9 & 9.05 & 4.28 & 4.15 & 4.14 & 4.06 & 4.19 \\ 
15 & 4.08 & 4.23 & 3.92 & 4.08 & 3.95 & 4.09 \\ \hline
\end{tabular}%
\caption{The smallest size-controlling critical values for comparing the
means of two heteroskedastic groups.}
\label{fig:CV}
\end{table}

In relation to Table \ref{fig:CV} we note that the equivalences discussed
before predict, e.g., that the ratio between the entries in the column
labeled HC0 and the corresponding entries in the column labeled HC1 should
be equal to $n/(n-2)=30/28\approx 1.0714$. The ratios computed from the
table are $1.0414$, $1.0761$, and $1.0721$ (for $n_{1}=3,9,15)$, which is in
pretty good agreement (especially if one converts the critical values shown
in the table to critical values for the corresponding \textquotedblleft $t$%
-test\textquotedblright\ versions by computing their square roots). The
agreement between theoretical and observed ratios for the HC0R-HC4R columns
is similar. In the balanced case one can also use the additional
equivalences mentioned before and one again finds very good agreement.
Similarly, the critical values for UC and UCR in Table \ref{fig:CV} are in
excellent agreement with their theoretical relationship found in Section \ref%
{sec:trivial_uc_tilde}. The reason for the small discrepancies observed lies
in the fact that the algorithm underlying the computations for Table \ref%
{fig:CV} makes use of a random search algorithm. Concerning Table \ref%
{fig:CV}, we also mention that, in the example considered here and for the
test statistic HC2, \cite{IbragMuell2016} prove in their Theorem 1 (see also
the discussion preceding that theorem) that the smallest size-controlling
critical values are given by\emph{\ }$18.51$ ($n_{1}=3$), $5.32$ ($n_{1}=9$%
), and $4.60$ ($n_{1}=15$), respectively. The numerically determined
critical values in Table \ref{fig:CV} are reasonably close to these values
(after conversion of the critical values to corresponding \textquotedblleft $%
t$-test\textquotedblright\ critical values the maximal difference is about $%
0.1$). Of course, the accuracy of our algorithm could be increased by using
more stringent accuracy parameters in the optimization routines underlying
the computation of the critical value, but this would come with a longer
runtime.

From Table \ref{fig:CV} it is clear that for the tests based on \emph{%
unrestricted} residuals the smallest size-controlling critical values
obtained are always larger, sometimes considerably, than $C_{\chi
^{2},0.05}\approx 3.8415$, again showing that the latter critical value is
not effecting size control. For the tests based on \emph{restricted}
residuals the smallest size-controlling critical values sometimes fall below 
$C_{\chi ^{2},0.05}$ in the strongly unbalanced case (which is not
completely surprising in view of Section \ref{sec:trivial_het_tilde}); while
in this case $C_{\chi ^{2},0.05}$ effects size-control, using the smaller
size-controlling critical values given in Table \ref{fig:CV} can only be
advantageous in terms of power.

That being said, we emphasize a trivial, but important point, namely that
comparing the magnitudes\ of size-controlling critical values relating to 
\emph{different} test statistics is not very meaningful and, in particular, 
\emph{not} a valid way of comparing the quality of the resulting tests. That
is, while it may be tempting to infer from Table \ref{fig:CV} that the HC0
test should be considerably more conservative than the HC4 test, or that the
UC test should be considerably more conservative than the UCR test, such a
conclusion would be false and not warranted at all (in particular, recall
that UC and UCR in fact result in (essentially) the same test if the
critical values from Table \ref{fig:CV} are being used). While this would be
correct if the critical values were all meant to be used with the same test
statistic (which they are not), critical values belonging to different test
statistics can certainly not be compared in such a way. Instead, one has to
compare the corresponding power functions, which is what we shall do next.

\begin{figure}[tbp]
\centering
\includegraphics[width=\linewidth]{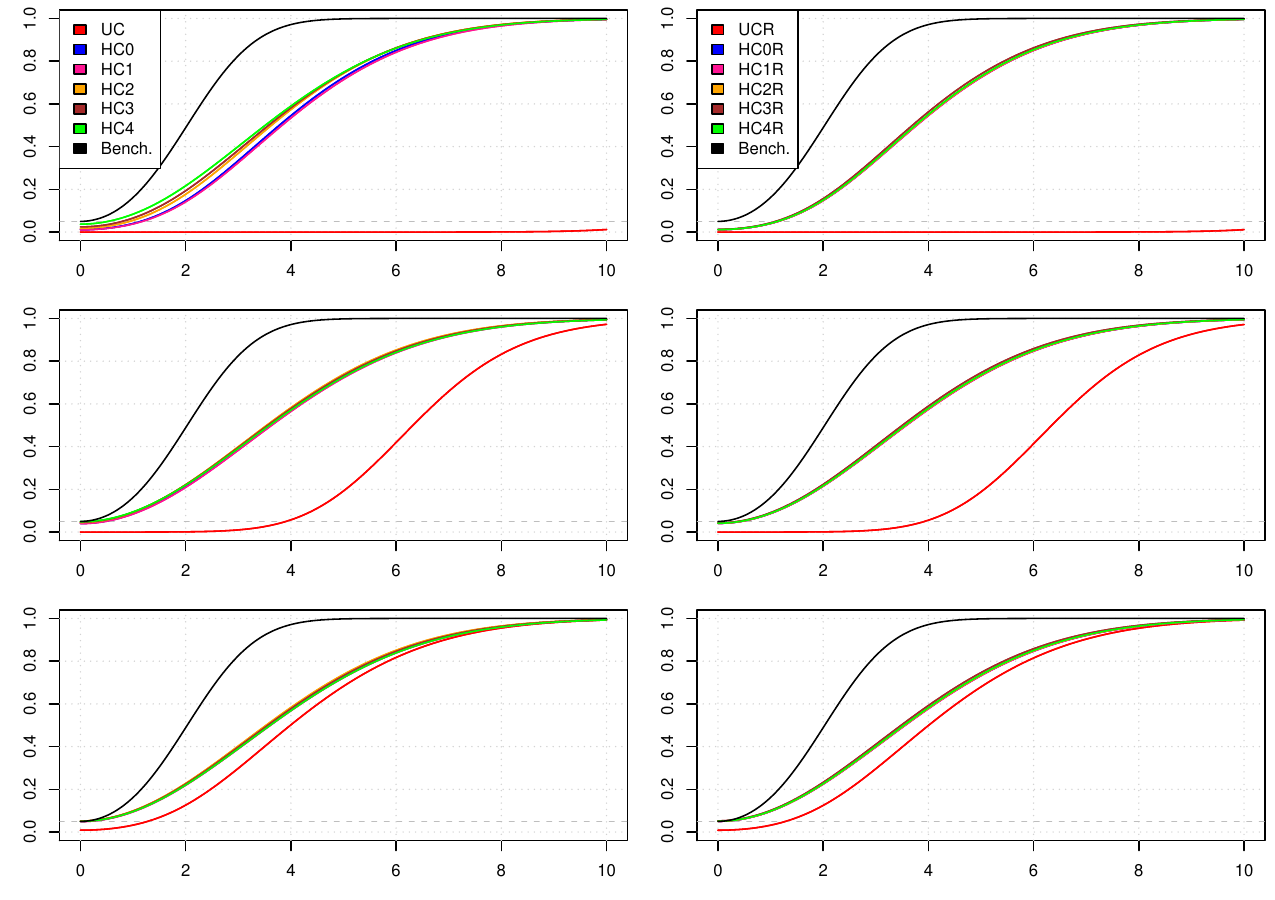}
\caption{Power functions for~$n_{1}=3$. Left column: tests based on
unrestricted residuals (cf.~legend). Right column: tests based on restricted
residuals (cf.~legend). The rows corresponds to~$\Sigma _{a}$ for~$a=1,5,9$
from top to bottom. The abscissa shows $\protect\delta $. In the left panel
the HC0-HC4-curves turn out to be barely distinguishable, with the HC1-curve
lying on top of the HC0-curve. In the right panel the HC4R-curve lies on top
of the HC0R-HC3R-curves. See the text for an explanation.}
\label{fig:bl1}
\end{figure}

\begin{figure}[tbp]
\centering
\includegraphics[width=\linewidth]{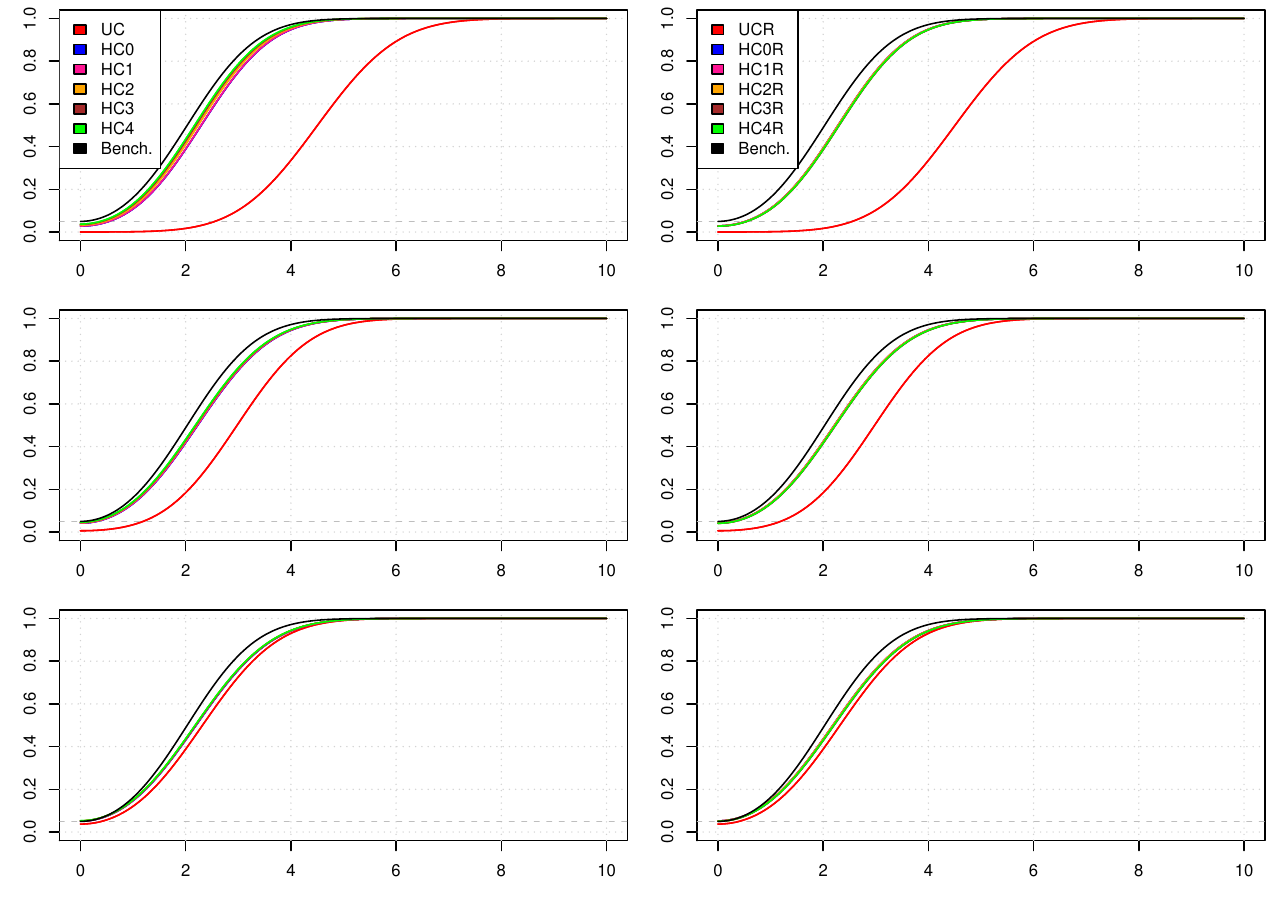}
\caption{Power functions for~$n_{1}=9$. Left column: tests based on
unrestricted residuals (cf.~legend). Right column: tests based on restricted
residuals (cf.~legend). The rows corresponds to~$\Sigma _{a}$ for~$a=1,5,9$
from top to bottom. The abscissa shows $\protect\delta $. In the left panel
the HC0-HC4-curves turn out to be barely distinguishable, with the HC1-curve
lying on top of the HC0-curve. In the right panel the HC4R-curve lies on top
of the HC0R-HC3R-curves. See the text for an explanation.}
\label{fig:bl3}
\end{figure}

The power functions are shown in Figure \ref{fig:bl1} (\textquotedblleft
strongly unbalanced\textquotedblright , $n_{1}=3$), Figure \ref{fig:bl3}
(\textquotedblleft moderately unbalanced\textquotedblright , $n_{1}=9$), and
Figure \ref{fig:bl5} (\textquotedblleft balanced\textquotedblright , $%
n_{1}=15$), where only the first two figures are shown in the main text, and
the last figure (in which the power functions of all the feasible tests lie
\textquotedblleft on top of each other\textquotedblright ) is available in
Appendix \ref{app:addfig}. Readers are referred to the online version for
colored figures. 

The power functions illustrate that the testing problem is getting easier,
(i.e., power gets closer to the oracle benchmark), for more balanced design,
which has intuitive appeal. Except for the strongly unbalanced case ($%
n_{1}=3 $), the power loss of the tests based on HC0-HC4 and HC0R-HC4R
relative to the oracle benchmark is surprisingly small (see Figure \ref%
{fig:bl3} as well as Figure \ref{fig:bl5} in Appendix \ref{app:addfig}). In
the unbalanced cases ($n_{1}\in \{3,9\}$) the HC0-HC4-based tests behave all
very similarly, with the power functions of the HC0- and HC1-based test
being virtually indistinguishable (as they should in view of the before
discussed equivalence). The UC-based test shows markedly worse power
performance. Similarly, the HC0R-HC4R-based tests have virtually
indistinguishable power functions (as they should because of the before
discussed equivalence). The UCR-based test again is inferior (and its power
function coincides with the one of UC as mentioned before). There appears
also to be little difference between basing the test statistics on
unrestricted or restricted residuals in this example. In the balanced case
we know that \emph{all} the feasible tests have exactly the same power
function in view of our earlier discussion. This is visible in Figure \ref%
{fig:bl5} in Appendix \ref{app:addfig}. Also the different forms of
heteroskedasticity considered seem not to have much effect on the power
functions (when expressed as a function of $\delta $), except for UC and UCR
in the unbalanced cases.

Hence, within the scenario considered in this section, perhaps the most
important conclusion concerning the choice of a test statistic appears to be
to avoid UC and UCR. Everything apart from that, i.e., whether one uses
unrestricted or restricted residuals to construct the test or which specific
heteroskedasticity-correction one decides to use, seems to be a comparably
irrelevant part of the problem once the right (i.e., smallest
size-controlling) critical value is used. We shall see in the next
subsection that this conclusion very much depends on the scenario considered
here and does not generalize beyond, illustrating the danger of drawing
conclusions from a limited numerical study.

\subsubsection{A high-leverage design matrix\label{sec:powhostX}}

In this section, we consider testing $\beta _{2}=0$ in a model with
intercept and a single regressor $x=(10,\cos (2),\cos (3),\ldots ,\cos
(n))^{\prime }$. Obviously, the regressor has a dominant first coordinate,
leading to diagonal elements $h_{ii}$ of $X(X^{\prime }X)^{-1}X^{\prime }$
such that the ratio of largest to smallest $h_{ii}$ is roughly $26$ ($\max
h_{ii}\simeq 0.879$, $\min h_{ii}\simeq 0.033$). Hence, the design matrix $X$
provides (on purpose) an extreme case, which leads to quite interesting
results. We consider again the case $n=30$ and $\alpha =0.05$, but now show
power functions for $\Sigma _{a}^{\ast }$, $a=0,\ldots ,4$, where 
\begin{equation*}
\Sigma _{a}^{\ast }=n^{-1}\limfunc{diag}\left( 7a+1,\frac{n-7a-1}{n-1}%
,\ldots ,\frac{n-7a-1}{n-1}\right) \in \mathfrak{C}_{Het}.
\end{equation*}%
Note that $\Sigma _{0}^{\ast }=n^{-1}I_{n}$ and that increasing $a$ from $0$
to $4$ leads to covariance matrices that approach the degenerate matrix $%
e_{1}(n)e_{1}(n)^{\prime }$. All conditions in Theorems \ref{Hetero_Robust}
and \ref{Hetero_Robust_tilde} are seen to be satisfied in this example: As
no vector $e_{i}(n)$ belongs to $\limfunc{span}(X)$ (and thus also not to $%
\mathfrak{M}_{0}^{lin}$), Assumptions \ref{R_and_X} and \ref{R_and_X_tilde}
as well as the sufficient condition for size control (\ref%
{non-incl_Het_uncorr}) are obviously satisfied. The size control conditions (%
\ref{non-incl_Het}) and (\ref{non-incl_Het_tilde}) have been checked
numerically, as has been the condition that none of the test statistics
HC0R-HC4R is constant on $\mathbb{R}^{n}\backslash \mathsf{\tilde{B}}$.

As in the preceding subsection, the critical values for each test statistic
are again chosen as \emph{the smallest critical value guaranteeing size
control} over $\mathfrak{C}_{Het}$ and they are presented in Table \ref%
{fig:CV2} below. [Existence follows from our theory since all assumptions
are satisfied as noted before.] For their computation the same algorithms
were used as in Section \ref{sec:powtwogroups}, with a similar statement
applying to the numerical routines used for computing the power functions.
Note that the critical values for the test statistics UC, HC0-HC3 are large,
reflecting the high-leverage in the design matrix; an exception is HC4, the
reason being that some of the HC4-weights are considerably larger than the
weights for HC0-HC3. Similarly as in the preceding subsection, the tests
based on HC0 and HC1 coincide (since HC0 and HC1 differ only by a
multiplicative constant and since smallest size-controlling critical values
are being used), and the same is true for the tests based on HC0R-HC4R, see
Remarks \ref{rem:equiv} and \ref{rem:equiv_tilde}. It is easily checked that
the ratios of the respective critical values provided in Table \ref{fig:CV2}
are in good agreement with the theoretical ratios predicted by theory.
Furthermore, the tests based on UC and UCR coincide (see Section \ref%
{sec:trivial_uc_tilde}), and the critical values for UC and UCR in Table \ref%
{fig:CV2} are in excellent agreement with their theoretical relationship
found in Section \ref{sec:trivial_uc_tilde}.

Table \ref{fig:CV2} shows that in this example the smallest size-controlling
critical values are -- except in one case -- always larger, sometimes
considerably larger, than $C_{\chi ^{2},0.05}\approx 3.8415$, once more
showing that the latter critical value is not effecting size control in
general. In the exceptional case, namely when the HC4 test statistic is
used, $C_{\chi ^{2},0.05}$ is considerably larger than the smallest
size-controlling critical value, which is $1.12$; while in this case $%
C_{\chi ^{2},0.05}$ effects size-control, using the smaller size-controlling
critical value $1.12$ can only be advantageous in terms of power. 
\begin{table}[th]
\centering
\begin{tabular}{llllll}
\hline
UC & HC0 & HC1 & HC2 & HC3 & HC4 \\ 
217.58 & 355.56 & 333.31 & 121.89 & 29.34 & 1.12 \\ \hline\hline
UCR & HC0R & HC1R & HC2R & HC3R & HC4R \\ 
25.69 & 5.41 & 5.45 & 5.34 & 5.29 & 5.44 \\ \hline
\end{tabular}%
\caption{Smallest size-controlling critical values for the high-leverage
design matrix.}
\label{fig:CV2}
\end{table}

\begin{figure}[tbp]
\centering
\includegraphics[width=\linewidth]{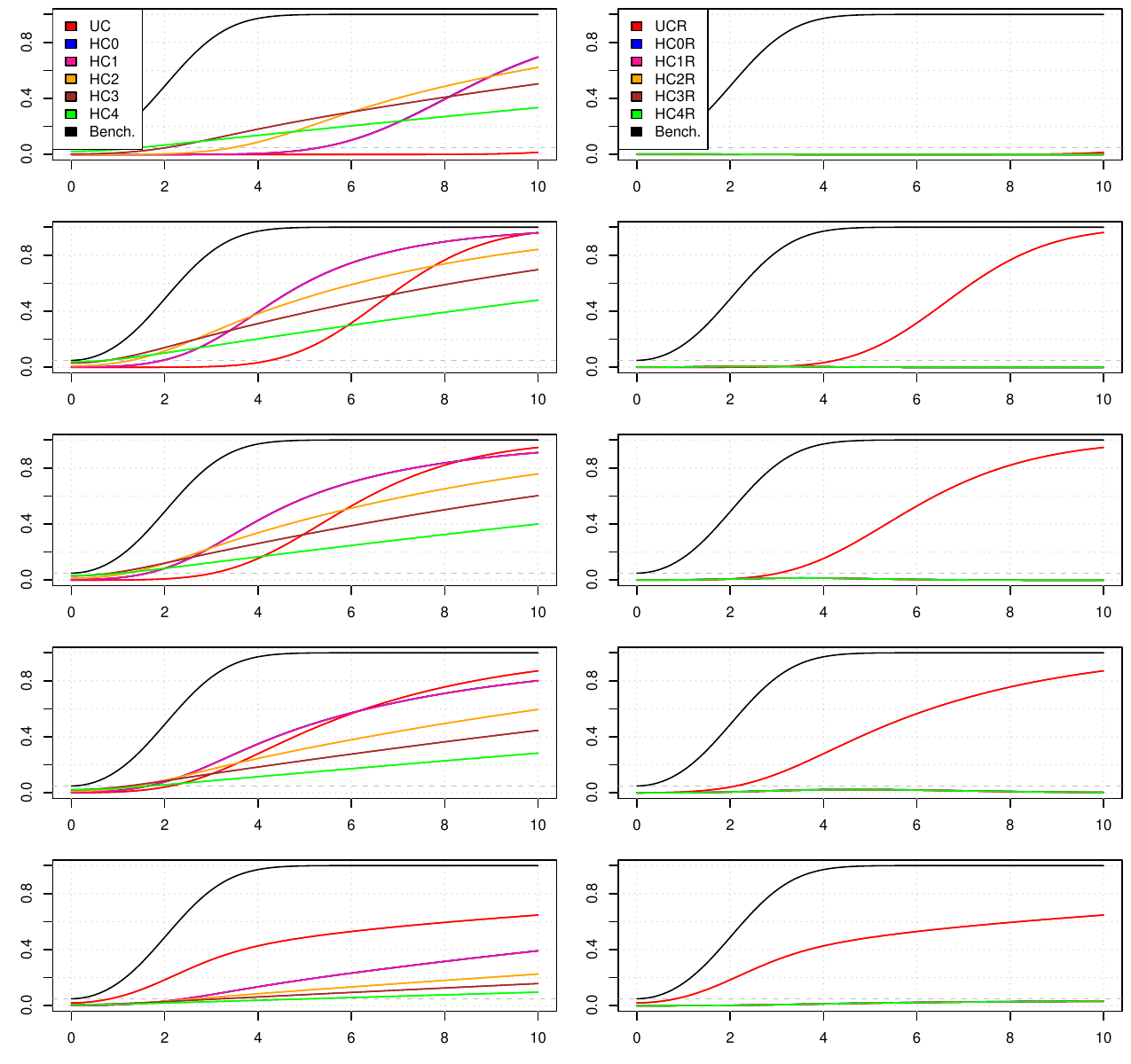}
\caption{Power functions for the design matrix considered in Section~\protect
\ref{sec:powhostX}. Left column: tests based on unrestricted residuals
(cf.~legend). Right column: tests based on restricted residuals
(cf.~legend). The rows~from top to bottom correspond to~$\Sigma _{a}^{\ast }$
for~$a=0,1,2,3,4$, the case~$a=0$ corresponding to homoskedasticity. The
abscissa shows $\protect\delta $. In the left panel the HC1-curve lies on
top of the HC0-curve. In the right panel the HC4R-curve lies on top of the
HC0R-HC3R-curves. See the text for an explanation.}
\label{fig:PW2}
\end{figure}

The power functions, when the size-controlling critical values from Table %
\ref{fig:CV2} are being used, are shown in Figure \ref{fig:PW2}. Readers are
referred to the online version for a colored figure. Again, as predicted by
theory, the power functions of the tests based on HC0 and HC1 shown in
Figure \ref{fig:PW2} coincide, as do the power functions of the tests based
on HC0R-HC4R; the same is true for the power functions of the tests based on
UC and UCR. The figure furthermore shows that in the setting considered
here, there is now a marked difference between tests based on HC0-HC4 and on
HC0R-HC4R, respectively: the power of the tests based on HC0R-HC4R is
nowhere greater than $\alpha $, their power function being even
non-monotonic, whereas the tests based on HC0-HC4 have increasing power as a
function of $\delta $. In contrast to the example considered in the
preceding subsection, the power functions of the tests based on HC0-HC4 and
on UC are now all markedly different and typically intersect, an exception
being the case of $\Sigma _{4}^{\ast }$ where the test based on UC offers
the highest power for that covariance matrix. Overall, however, there is no
clear ranking between the tests using unrestricted residuals in the example
considered here, although we note that the test based on UC (or,
equivalently, on UCR) performs very badly in the case of $\Sigma _{0}^{\ast
} $. This is not surprising as $\Sigma _{0}^{\ast }$ corresponds to
homoskedasticity and the critical value used here is much larger than the
classical critical value one would use given knowledge of this
homoskedasticity. Furthermore, and in contrast to the results in the
preceding subsection, the different forms of heteroskedasticity considered
have a noticeable effect on the power functions. The main takeaway is that
tests based on HC0R-HC4R (and probably on UC and UCR) should rather be
avoided. 

\section{Conclusion\label{concl}}

The usual heteroskedasticity robust test statistics such as $T_{Het}$ (using
HC0-HC4 weights) or $\tilde{T}_{Het}$ (using HC0R-HC4R weights), used in
conjunction with conventional critical values obtained from the asymptotic
null distribution, are often plagued by overrejection under the null. This
has been clearly documented in the literature for $T_{Het}$, and is shown
numerically for $\tilde{T}_{Het}$ (as well as for $T_{Het}$) in Section \ref%
{numerical} above. Not surprisingly, similar observations apply to the
\textquotedblleft uncorrected\textquotedblright\ test statistics $T_{uc}$
and $\tilde{T}_{uc}$. We show theoretically that all these test statistics
can be size-controlled under quite weak conditions by an appropriate choice
of critical values.

From the above discussion and the numerical results in Section \ref%
{numerical} it transpires that smallest size-controlling critical values
rather than conventional critical values should be used in order to avoid
the risk of overrejection. For the computation of smallest size-controlling
critical values we provide algorithms which have been implemented in the
R-package \textbf{hrt} (\cite{hrt}) and thus are readily available for the
user.

An additional advantage from using smallest size-controlling critical values
over conventional critical values is that this typically leads to improved
power in instances where conventional critical values lead to underrejection
(i.e., lead to worst-case rejection probability under the null less than the
nominal significance level) as is sometimes the case; see Sections \ref%
{sec:trivial_het_tilde} and \ref{sec:powernum}.

If smallest size-controlling critical values are adopted (as they should),
the numerical results in Section \ref{numerical} suggest that the test
statistic $\tilde{T}_{Het}$ (with the usual weights HC0R-HC4R) should be
avoided, as the resulting tests may have very poor power properties (see the
example in Section \ref{sec:powhostX}). The test statistic $T_{Het}$ seems
to perform better in terms of power, with no clear ranking emerging with
regards to the weights HC0-HC4 being used. The \textquotedblleft
uncorrected\textquotedblright\ test statistics $T_{uc}$ and $\tilde{T}_{uc}$
appear to be inferior to $T_{Het}$ in terms of power in almost all of the
numerical examples considered. We also point out that -- when using smallest
size-controlling critical values -- the tests based on $T_{Het}$ employing
the HC0 and the HC1 weights, respectively, in fact coincide; and the same
holds for tests based on $\tilde{T}_{Het}$ employing the HC0R and the HC1R
weights, respectively. Also, the tests based on $T_{uc}$ and $\tilde{T}_{uc}$
then (essentially) coincide. See Remarks \ref{rem:equiv}, \ref%
{rem:equiv_tilde}, and Section \ref{sec:trivial_uc_tilde} as well as the
pertaining discussion in Section \ref{numerical} for more information,
including additional equivalencies when the design matrix $X$ and the
restriction $R$ have certain special properties.

\newpage

\appendix

\section{Appendix: Size control over other heteroskedasticity models\label%
{app_d}}

As already noted earlier, if size control is possible over $\mathfrak{C}%
_{Het}$, then the same is true over \emph{any} conceivable class of
heteroskedasticity structures, since these can (possible after
normalization) be cast as a subset $\mathfrak{C}$ of $\mathfrak{C}_{Het}$;
and, in fact, any critical value delivering size control over $\mathfrak{C}%
_{Het}$ also delivers size control over any such $\mathfrak{C}$, but even
smaller critical values may already suffice for size control over $\mathfrak{%
C}$. Also, for some heteroskedasticity models $\mathfrak{C}\subseteq 
\mathfrak{C}_{Het}$, the sufficient conditions employed in Theorems \ref%
{Hetero_Robust} and \ref{Hetero_Robust_tilde} (which imply size control over 
$\mathfrak{C}_{Het}$) may be unnecessarily restrictive, if one wants to
establish size control over $\mathfrak{C}$ only. For this reason, we show in
the following how the general theory laid out in Section 5 of \cite{PP3} can
be used to derive size control results tailored to various subsets $%
\mathfrak{C}$ by exemplarily treating the cases $\mathfrak{C}=\mathfrak{C}%
_{(n_{1},\ldots ,n_{m})}$ and $\mathfrak{C}=\mathfrak{C}_{Het,\tau _{\ast }}$
defined below. Size control results over other choices of $\mathfrak{C}$ can
be derived from the results in Section 5 of \cite{PP3} in a similar manner,
see Subsection \ref{sec_size_control_2_2} further below for some discussion.
Here $\mathfrak{C}_{(n_{1},\ldots ,n_{m})}$ is defined as follows: Let $m\in 
\mathbb{N}$, and let $n_{j}\in \mathbb{N}$ for $j=1,\ldots ,m$ satisfy $%
\sum_{j=1}^{m}n_{j}=n$. Set $n_{j}^{+}=\sum_{l=1}^{j}n_{l}$ and define%
\begin{equation*}
\mathfrak{C}_{(n_{1},\ldots ,n_{m})}=\left\{ \limfunc{diag}(\tau
_{1}^{2},\ldots ,\tau _{n}^{2})\in \mathfrak{C}_{Het}:\tau
_{n_{j-1}^{+}+1}^{2}=\ldots =\tau _{n_{j}^{+}}^{2}\text{ for }j=1,\ldots
,m\right\}
\end{equation*}%
with the convention that $n_{0}^{+}=0$. This may be a natural
heteroskedasticity model when the observations come from $m$ groups and when
it is reasonable to assume homoskedasticity within groups.\footnote{%
As long as we assume that the grouping is known, there is little loss of
generality to assume that the elements belonging to the same group are
numbered contiguously, since we otherwise only need to relabel the data.}
Note that in case $n_{j}=1$ for all $j$, then $m=n$ and $\mathfrak{C}%
_{(n_{1},\ldots ,n_{m})}=\mathfrak{C}_{Het}$ hold; and in case $m=1$ we have 
$\mathfrak{C}_{(n_{1},\ldots ,n_{m})}=\{n^{-1}I_{n}\}$, i.e., we have
homoskedasticity. Furthermore, $\mathfrak{C}_{Het,\tau _{\ast }}$ is given
by 
\begin{equation*}
\mathfrak{C}_{Het,\tau _{\ast }}=\left\{ \limfunc{diag}(\tau _{1}^{2},\ldots
,\tau _{n}^{2})\in \mathfrak{C}_{Het}:\tau _{i}^{2}\geq \tau _{\ast }^{2}%
\text{ for all }i\right\} ,
\end{equation*}%
where the lower bound $\tau _{\ast }$, $0<\tau _{\ast }<n^{-1/2}$, is set by
the user.

\subsection{Size control results for $T_{Het}$ and $T_{uc}$\label%
{sec_size_control_2}}

\subsubsection{Size control over $\mathfrak{C}_{(n_{1},\ldots ,n_{m})}$}

Proofs of the results in this subsection can be found in Appendix \ref{app_B}%
.

\begin{theorem}
\label{theorem_groupwise_hetero}Let $m\in \mathbb{N}$, and let $n_{j}\in 
\mathbb{N}$ for $j=1,\ldots ,m$ satisfy $\sum_{j=1}^{m}n_{j}=n$. Then:

(a) For every $0<\alpha <1$ there exists a real number $C(\alpha )$ such that%
\begin{equation}
\sup_{\mu _{0}\in \mathfrak{M}_{0}}\sup_{0<\sigma ^{2}<\infty }\sup_{\Sigma
\in \mathfrak{C}_{(n_{1},\ldots ,n_{m})}}P_{\mu _{0},\sigma ^{2}\Sigma
}(T_{uc}\geq C(\alpha ))\leq \alpha  \label{size-control_groupHet_uncorr}
\end{equation}%
holds, provided that%
\begin{equation}
\func{span}\left( \left\{ e_{i}(n):i\in (n_{j-1}^{+},n_{j}^{+}]\right\}
\right) \nsubseteqq \func{span}(X)\text{ \ for every }j=1,\ldots ,m\text{
with }(n_{j-1}^{+},n_{j}^{+}]\cap I_{1}(\mathfrak{M}_{0}^{lin})\neq
\emptyset .  \label{non-incl_groupHet_uncorr}
\end{equation}%
Furthermore, under condition (\ref{non-incl_groupHet_uncorr}), even equality
can be achieved in (\ref{size-control_groupHet_uncorr}) by a proper choice
of $C(\alpha )$, provided $\alpha \in (0,\alpha ^{\ast }]\cap (0,1)$ holds,
where $\alpha ^{\ast }=\sup_{C\in (C^{\ast },\infty )}\sup_{\Sigma \in 
\mathfrak{C}_{(n_{1},\ldots ,n_{m})}}P_{\mu _{0},\Sigma }(T_{uc}\geq C)$ is
positive and where $C^{\ast }$ is defined as in Lemma 5.11 of \cite{PP3}
with $\mathfrak{C}=\mathfrak{C}_{(n_{1},\ldots ,n_{m})}$, $T=T_{uc}$, $%
N^{\dag }=\func{span}(X)$, and $\mathcal{L}=\mathfrak{M}_{0}^{lin}$ (with
neither $\alpha ^{\ast }$ nor $C^{\ast }$ depending on the choice of $\mu
_{0}\in \mathfrak{M}_{0}$).

(b) Suppose Assumption \ref{R_and_X} is satisfied. Then for every $0<\alpha
<1$ there exists a real number $C(\alpha )$ such that%
\begin{equation}
\sup_{\mu _{0}\in \mathfrak{M}_{0}}\sup_{0<\sigma ^{2}<\infty }\sup_{\Sigma
\in \mathfrak{C}_{(n_{1},\ldots ,n_{m})}}P_{\mu _{0},\sigma ^{2}\Sigma
}(T_{Het}\geq C(\alpha ))\leq \alpha  \label{size-control_groupHet}
\end{equation}%
holds, provided that 
\begin{equation}
\func{span}\left( \left\{ e_{i}(n):i\in (n_{j-1}^{+},n_{j}^{+}]\right\}
\right) \nsubseteqq \mathsf{B}\text{ \ for every }j=1,\ldots ,m\text{ with }%
(n_{j-1}^{+},n_{j}^{+}]\cap I_{1}(\mathfrak{M}_{0}^{lin})\neq \emptyset .
\label{non-incl_groupHet}
\end{equation}%
Furthermore, under condition (\ref{non-incl_groupHet}), even equality can be
achieved in (\ref{size-control_groupHet}) by a proper choice of $C(\alpha )$%
, provided $\alpha \in (0,\alpha ^{\ast }]\cap (0,1)$ holds, where $\alpha
^{\ast }=\sup_{C\in (C^{\ast },\infty )}\sup_{\Sigma \in \mathfrak{C}%
_{(n_{1},\ldots ,n_{m})}}P_{\mu _{0},\Sigma }(T_{Het}\geq C)$ is positive
and where $C^{\ast }$ is defined as in Lemma 5.11 of \cite{PP3} with $%
\mathfrak{C}=\mathfrak{C}_{(n_{1},\ldots ,n_{m})}$, $T=T_{Het}$, $N^{\dag }=%
\mathsf{B}$, and $\mathcal{L}=\mathfrak{M}_{0}^{lin}$ (with neither $\alpha
^{\ast }$ nor $C^{\ast }$ depending on the choice of $\mu _{0}\in \mathfrak{M%
}_{0}$).

(c) Under the assumptions of Part (a) (Part (b), respectively) implying
existence of a critical value $C(\alpha )$ satisfying (\ref%
{size-control_groupHet_uncorr}) ((\ref{size-control_groupHet}),
respectively), a smallest critical value, denoted by $C_{\Diamond }(\alpha )$%
, satisfying (\ref{size-control_groupHet_uncorr}) ((\ref%
{size-control_groupHet}), respectively) exists for every $0<\alpha <1$. And $%
C_{\Diamond }(\alpha )$ corresponding to Part (a) (Part (b), respectively)
is also the smallest among the critical values leading to equality in (\ref%
{size-control_groupHet_uncorr}) ((\ref{size-control_groupHet}),
respectively) whenever such critical values exist. [Although $C_{\Diamond
}(\alpha )$ corresponding to Part (a) and (b), respectively, will typically
be different, we use the same symbol.]\footnote{%
Cf.~also Appendix \ref{useful}.}
\end{theorem}

It is easy to see that the discussion in the first paragraph following
Theorem \ref{Hetero_Robust} applies mutatis mutandis also to the above
theorem. Similarly, Remarks \ref{rem_indep_r}, \ref{rem:equiv}, \ref%
{rem_positiv}, \ref{rem:larger_alphastar}, \ref{rem_strict_ineq}, and
Proposition \ref{rem_C*} carry over. Furthermore, we have the following
result corresponding to Proposition \ref{rem_necessity}:

\begin{proposition}
\label{rem_necessity_2}(a) If (\ref{non-incl_groupHet_uncorr}) is violated,
then $\sup_{\Sigma \in \mathfrak{C}_{(n_{1},\ldots ,n_{m})}}P_{\mu
_{0},\sigma ^{2}\Sigma }(T_{uc}\geq C)=1$ for \emph{every} choice of
critical value $C$, every $\mu _{0}\in \mathfrak{M}_{0}$, and every $\sigma
^{2}\in (0,\infty )$ (implying that size equals $1$ for every $C$). As a
consequence, the sufficient condition for size control (\ref%
{non-incl_groupHet_uncorr}) in Part (a) of Theorem \ref%
{theorem_groupwise_hetero} is also necessary.

(b) Suppose Assumption \ref{R_and_X} is satisfied.\footnote{%
If this assumption is violated then $T_{Het}$ is identically zero, an
uninteresting trivial case.} If (\ref{non-incl_groupHet_uncorr}) is
violated, then $\sup_{\Sigma \in \mathfrak{C}_{(n_{1},\ldots ,n_{m})}}P_{\mu
_{0},\sigma ^{2}\Sigma }(T_{Het}\geq C)=1$ for \emph{every} choice of
critical value $C$, every $\mu _{0}\in \mathfrak{M}_{0}$, and every $\sigma
^{2}\in (0,\infty )$ (implying that size equals $1$ for every $C$). [In case 
$X$ and $R$ are such that $\mathsf{B}=\limfunc{span}(X)$, conditions (\ref%
{non-incl_groupHet_uncorr}) and (\ref{non-incl_groupHet}) coincide; hence
the sufficient condition for size control (\ref{non-incl_groupHet}) in Part
(b) of Theorem \ref{theorem_groupwise_hetero} is then also necessary in this
case.]
\end{proposition}

\begin{remark}
\emph{(Homoskedasticity)} Theorem \ref{theorem_groupwise_hetero} allows also
for the case $m=1$, in which case $\mathfrak{C}_{(n_{1},\ldots
,n_{m})}=\left\{ n^{-1}I_{n}\right\} $, i.e., errors are homoskedastic. In
this case it is easy to see that the sufficient conditions for size control
in the theorem are trivially satisfied and size control for $T_{Het}$ (and $%
T_{uc}$) is possible.\footnote{%
A related but slightly different argument proceeds by directly noting from
its definition that $\mathbb{J}(\mathfrak{M}_{0}^{lin},\mathfrak{C}%
_{(n_{1},\ldots n_{m})})$ is empty in case $m=1$ (cf.~Apendix \ref{app_char}%
), and then to appeal to Remark 5.7 (or Proposition 5.12) in \cite{PP3}.} Of
course, this is in line with the fact that $T_{Het}$ and $T_{uc}$ are
obviously pivotal under the null if the errors are homoskedastic.
\end{remark}

\begin{remark}
\label{FB}\emph{(Behrens-Fisher problem)} Consider again the problem of
testing the equality of the means of two independent normal populations as
in Example \ref{ex_fish_behr} with the only difference that the variance
within each of the two groups is now assumed to be constant, i.e., the
heteroskedasticity model used is now given by $\mathfrak{C}_{(n_{1},n_{2})}$%
, where $n_{1}\geq 2$ and $n_{2}\geq 2$ are the group sizes. This is the
celebrated Behrens-Fisher problem. The square of the two-sample t-statistic $%
t_{FB}$, say, often used in this context coincides with $T_{Het}$ for the
choice $d_{i}=\left( 1-h_{ii}\right) ^{-1}$. The size controllability of $%
T_{Het}$ over $\mathfrak{C}_{Het}$ established in Example \ref{ex_fish_behr}
therefore a fortiori implies size controllability of $T_{Het}$ (and hence of 
$t_{FB}^{2}$) over $\mathfrak{C}_{(n_{1},n_{2})}$. Of course, this does not
add anything new to the literature on the Behrens-Fisher problem, since it
is known that under the null hypothesis $\left\vert t_{FB}\right\vert $ is
stochastically not larger than a $t$-distributed random variable with $\min
(n_{1},n_{2})-1$ degrees of freedom when $\mathfrak{C}_{(n_{1},n_{2})}$ is
the heteroskedasticity model, see \cite{MickeyBrown1966}. For more on the
Behrens-Fisher problem see \cite{KimCohen}, \cite{Ruben2002}, \cite{LR05}, 
\cite{BelloniDidier2008}, and the references cited therein.
\end{remark}

\subsubsection{Further size control results \label{sec_size_control_2_2}}

In this subsection it is understood that Assumption \ref{R_and_X} is
maintained when discussing results relating to $T_{Het}$.

(i) Given a heteroskedasticity model $\mathfrak{C}$ (i.e., $\emptyset \neq 
\mathfrak{C}\subseteq \mathfrak{C}_{Het}$), with the property that $\mathbb{J%
}(\mathfrak{M}_{0}^{lin},\mathfrak{C})$ is empty (where the collection $%
\mathbb{J}(\mathfrak{M}_{0}^{lin},\mathfrak{C})$ is defined on p.~421 of 
\cite{PP3}, see also Appendix \ref{app_char} further below), the tests based
on $T_{uc}$ and $T_{Het}$ are always size controllable over $\mathfrak{C}$.
This follows from Corollary 5.6 and Remark 5.7 in \cite{PP3}. In fact, exact
size control is then possible for every $\alpha \in (0,1)$ as a consequence
of Proposition 5.12 in the same reference upon noting that then $C^{\ast
}=-\infty $ and $\alpha ^{\ast }=1$ hold. [We note in passing that for such
a heteroskedasticity model $\mathfrak{C}$ the size of the rejection region $%
\{T_{uc}\geq C\}$ ($\{T_{Het}\geq C\}$, respectively) is less than $1$ for
every $C>0$ (this follows from Proposition 5.2 and Remark 5.4 in \cite{PP3}
as well as Part 6 of Lemma 5.15 in \cite{PP2016}).]\footnote{%
The verification of the assumptions in Corollary 5.6 and Propositions 5.2
and 5.12 of \cite{PP3} proceeds as in the proofs of Theorems \ref%
{Hetero_Robust} and \ref{theorem_groupwise_hetero}.}

(ii) A particular instance of the situation described in (i) is provided by
heteroskedasticity models $\mathfrak{C}$ that are subsets of a set of the
form $\mathfrak{C}_{Het,\tau _{\ast }}$ ($0<\tau _{\ast }<n^{-1/2}$), as in
this case $\mathbb{J}(\mathfrak{M}_{0}^{lin},\mathfrak{C})$ is easily seen
to be empty.

(iii)\ More generally, the tests based on $T_{uc}$ (on $T_{Het}$,
respectively) are size controllable over a heteroskedasticity model $%
\mathfrak{C}$, provided any $\mathcal{S}\in \mathbb{J}(\mathfrak{M}%
_{0}^{lin},\mathfrak{C})$ is not contained in $\func{span}(X)$ ($\mathsf{B}$%
, respectively). This follows easily from Corollary 5.6 and Proposition 5.12
in \cite{PP3}, the latter proposition also providing an exact size result,
which we refrain from spelling out in detail. Again there is a (partial)
converse: If an $\mathcal{S}\in \mathbb{J}(\mathfrak{M}_{0}^{lin},\mathfrak{C%
})$ exists with $\mathcal{S}\subseteq \func{span}(X)$, then the size over $%
\mathfrak{C}$ of the rejection region $\{T_{uc}\geq C\}$ ($\{T_{Het}\geq C\}$%
, respectively) is equal to $1$; see Theorem 3.1 in \cite{PP4}. Furthermore,
lower bounds for critical values that lead to size less than $1$ (in
particular, for size-controlling critical values) can be had with the help
of Corollary 5.17 in \cite{PP2016}, Lemma 5.11 and Proposition 5.12 in \cite%
{PP3}, or Lemma 4.1 in \cite{PP4}.

\subsection{Size control results for $\tilde{T}_{Het}$ and $\tilde{T}_{uc}$%
\label{sec_size_control_2_tilde}}

The proof of the subsequent theorem is given in Appendix \ref{app_C}. We
note that the first statement in Part (a) of the subsequent theorem is
actually trivial, since $\tilde{T}_{uc}$ is bounded as has been shown in
Section \ref{sec:trivial_uc_tilde}.

\begin{theorem}
\label{theorem_groupwise_hetero_tilde}Let $m\in \mathbb{N}$, and let $%
n_{j}\in \mathbb{N}$ for $j=1,\ldots ,m$ satisfy $\sum_{j=1}^{m}n_{j}=n$.
Then:

(a) For every $0<\alpha <1$ there exists a real number $C(\alpha )$ such that%
\begin{equation}
\sup_{\mu _{0}\in \mathfrak{M}_{0}}\sup_{0<\sigma ^{2}<\infty }\sup_{\Sigma
\in \mathfrak{C}_{(n_{1},\ldots ,n_{m})}}P_{\mu _{0},\sigma ^{2}\Sigma }(%
\tilde{T}_{uc}\geq C(\alpha ))\leq \alpha
\label{size-control_groupHet_uncorr_tilde}
\end{equation}%
holds. Furthermore, even equality can be achieved in (\ref%
{size-control_groupHet_uncorr_tilde}) by a proper choice of $C(\alpha )$,
provided $\alpha \in (0,\alpha ^{\ast }]\cap (0,1)$ holds, where $\alpha
^{\ast }=\sup_{C\in (C^{\ast },\infty )}\sup_{\Sigma \in \mathfrak{C}%
_{(n_{1},\ldots ,n_{m})}}P_{\mu _{0},\Sigma }(\tilde{T}_{uc}\geq C)$ and
where $C^{\ast }$ is defined as in Lemma 5.11 of \cite{PP3} with $\mathfrak{C%
}=\mathfrak{C}_{(n_{1},\ldots ,n_{m})}$, $T=\tilde{T}_{uc}$, $N^{\dag }=%
\mathfrak{M}_{0}$, and $\mathcal{L}=\mathfrak{M}_{0}^{lin}$ (with neither $%
\alpha ^{\ast }$ nor $C^{\ast }$ depending on the choice of $\mu _{0}\in 
\mathfrak{M}_{0}$).

(b) Suppose Assumption \ref{R_and_X_tilde} is satisfied. Suppose further
that $\tilde{T}_{Het}$ is not constant on $\mathbb{R}^{n}\backslash \mathsf{%
\tilde{B}}$.\footnote{%
Cf.~Footnote \ref{FN_const}.} Then for every $0<\alpha <1$ there exists a
real number $C(\alpha )$ such that%
\begin{equation}
\sup_{\mu _{0}\in \mathfrak{M}_{0}}\sup_{0<\sigma ^{2}<\infty }\sup_{\Sigma
\in \mathfrak{C}_{(n_{1},\ldots ,n_{m})}}P_{\mu _{0},\sigma ^{2}\Sigma }(%
\tilde{T}_{Het}\geq C(\alpha ))\leq \alpha
\label{size-control_groupHet_tilde}
\end{equation}%
holds, provided that for some $\mu _{0}\in \mathfrak{M}_{0}$ (and hence for
all $\mu _{0}\in \mathfrak{M}_{0}$) 
\begin{equation}
\mu _{0}+\func{span}\left( \left\{ e_{i}(n):i\in
(n_{j-1}^{+},n_{j}^{+}]\right\} \right) \nsubseteqq \mathsf{\tilde{B}}\text{
\ for every }j=1,\ldots ,m\text{ with }(n_{j-1}^{+},n_{j}^{+}]\cap I_{1}(%
\mathfrak{M}_{0}^{lin})\neq \emptyset .  \label{non-incl_groupHet_tilde}
\end{equation}%
Furthermore, under condition (\ref{non-incl_groupHet_tilde}), even equality
can be achieved in (\ref{size-control_groupHet_tilde}) by a proper choice of 
$C(\alpha )$, provided $\alpha \in (0,\alpha ^{\ast }]\cap (0,1)$ holds,
where $\alpha ^{\ast }=\sup_{C\in (C^{\ast },\infty )}\sup_{\Sigma \in 
\mathfrak{C}_{(n_{1},\ldots ,n_{m})}}P_{\mu _{0},\Sigma }(\tilde{T}%
_{Het}\geq C)$ and where $C^{\ast }$ is defined as in Lemma 5.11 of \cite%
{PP3} with $\mathfrak{C}=\mathfrak{C}_{(n_{1},\ldots ,n_{m})}$, $T=\tilde{T}%
_{Het}$, $N^{\dag }=\mathsf{\tilde{B}}$, and $\mathcal{L}=\mathfrak{M}%
_{0}^{lin}$ (with neither $\alpha ^{\ast }$ nor $C^{\ast }$ depending on the
choice of $\mu _{0}\in \mathfrak{M}_{0}$).

(c) Under the assumptions of Part (a) (Part (b), respectively) implying
existence of a critical value $C(\alpha )$ satisfying (\ref%
{size-control_groupHet_uncorr_tilde}) ((\ref{size-control_groupHet_tilde}),
respectively), a smallest critical value, denoted by $C_{\Diamond }(\alpha )$%
, satisfying (\ref{size-control_groupHet_uncorr_tilde}) ((\ref%
{size-control_groupHet_tilde}), respectively) exists for every $0<\alpha <1$.%
\footnote{%
Note that there are in fact no assumptions for Part (a). We have chosen this
formulation for reasons of brevity.} And $C_{\Diamond }(\alpha )$
corresponding to Part (a) (Part (b), respectively) is also the smallest
among the critical values leading to equality in (\ref%
{size-control_groupHet_uncorr_tilde}) ((\ref{size-control_groupHet_tilde}),
respectively) whenever such critical values exist. [Although $C_{\Diamond
}(\alpha )$ corresponding to Part (a) and (b), respectively, will typically
be different, we use the same symbol.]\footnote{%
Cf.~also Appendix \ref{useful}.}
\end{theorem}

It is easy to see that the discussion in the first paragraph following
Theorem \ref{Hetero_Robust_tilde} applies mutatis mutandis also to the above
theorem. Similarly, Remarks \ref{rem:equiv_tilde}, \ref{rem_positiv_tilde}, %
\ref{rem:larger_alphastar_tilde}, \ref{rem_strict_ineq_2}, and Proposition %
\ref{rem_C*_tilde} carry over.

A discussion of size control results for $\tilde{T}_{uc}$ and $\tilde{T}%
_{Het}$ over other choices of $\mathfrak{C}$ based on the results in Section
5 of \cite{PP3} can also be given (cf. the discussion in Subsection \ref%
{sec_size_control_2_2}), but we refrain from spelling out the details. We
only note that the test based on $\tilde{T}_{Het}$ is always size
controllable over $\mathfrak{C}_{Het,\tau _{\ast }}$ ($0<\tau _{\ast
}<n^{-1/2}$), and the same is trivially true for $\tilde{T}_{uc}$.

\subsection{A useful observation\label{useful}}

Let $\mathfrak{C}$ be an arbitrary heteroskedasticity model (i.e., $%
\emptyset \neq \mathfrak{C}\subseteq \mathfrak{C}_{Het}$), let $0<\alpha <1$%
, and let $T$ stand for $T_{uc}$ or $T_{Het}$, respectively, where in case
of $T=T_{Het}$ we assume that Assumption \ref{R_and_X} is satisfied. Suppose
that $T$ is size-controllable at significance level $\alpha $ (i.e., that $%
\sup_{\mu _{0}\in \mathfrak{M}_{0}}\sup_{0<\sigma ^{2}<\infty }\sup_{\Sigma
\in \mathfrak{C}}P_{\mu _{0},\sigma ^{2}\Sigma }(T\geq C)\leq \alpha $ holds
for some real $C$). Then a smallest size-controlling critical value $%
C_{\Diamond }(\alpha )$ always exists.\footnote{%
Note that this, e.g., covers the case discussed in Example \ref{ex_k_pop},
where size-control can be established for $T_{Het}$ despite the fact that
the sufficient conditions in Theorem \ref{Hetero_Robust} are not satisfied
(and hence Part (c) of that theorem can not be used).} And if a critical
value $C\in \mathbb{R}$ exists such that $\sup_{\mu _{0}\in \mathfrak{M}%
_{0}}\sup_{0<\sigma ^{2}<\infty }\sup_{\Sigma \in \mathfrak{C}}P_{\mu
_{0},\sigma ^{2}\Sigma }(T\geq C)=\alpha $ holds, then $C_{\Diamond }(\alpha
)$ is also the smallest among these critical values. This follows from
Remark 5.10 and Lemma 5.16 in\ \cite{PP3} combined with Remark \ref{F-type}
in Appendix \ref{app_B}. The same is true for $T=\tilde{T}_{uc}$ and $T=%
\tilde{T}_{Het}$, where in case of $T=\tilde{T}_{Het}$ we assume that
Assumption \ref{R_and_X_tilde} is satisfied and that $\tilde{T}_{Het}$ is
not constant on $\mathbb{R}^{n}\backslash \mathsf{\tilde{B}}$. This follows
again from Remark 5.10 in\ \cite{PP3} now together with Lemma \ref%
{lem:nullset} in Appendix \ref{app_C}.

\section{Appendix: Characterization of $\mathbb{J}(\mathcal{L},\mathfrak{C})$
for $\mathfrak{C}=\mathfrak{C}_{Het}$ and $\mathfrak{C}=\mathfrak{C}%
_{(n_{1},\ldots ,n_{m})}$\label{app_char}}

A key ingredient in the proof of size control results such as Theorem \ref%
{Hetero_Robust} or \ref{Hetero_Robust_tilde} is a certain collection $%
\mathbb{J}(\mathcal{L},\mathfrak{C})$ of linear subspaces of $\mathbb{R}^{n}$
introduced in \cite{PP3}. For the convenience of the reader we reproduce
this definition, specialized to the present setting, below. The leading case
in the applications will be the case $\mathcal{L}=\mathfrak{M}_{0}^{lin}$.

\begin{definition}
Let $\mathfrak{C}$ be a heteroskedasticity model (i.e., $\emptyset \neq 
\mathfrak{C}\subseteq \mathfrak{C}_{Het}$). Given a linear subspace $%
\mathcal{L}$ of $\mathbb{R}^{n}$ with $\dim (\mathcal{L})<n$ and an element $%
\Sigma \in \mathfrak{C}$, we let 
\begin{equation*}
\mathcal{L}(\Sigma )=\frac{\Pi _{\mathcal{L}^{\bot }}\Sigma \Pi _{\mathcal{L}%
^{\bot }}}{\Vert {\Pi _{\mathcal{L}^{\bot }}\Sigma \Pi _{\mathcal{L}^{\bot }}%
}\Vert }
\end{equation*}%
and $\mathcal{L}(\mathfrak{C})=\left\{ \mathcal{L}(\Sigma ):\Sigma \in 
\mathfrak{C}\right\} $. Furthermore, we define 
\begin{equation*}
\mathbb{J}(\mathcal{L},\mathfrak{C})=\left\{ \mathrm{\limfunc{span}}(\bar{%
\Sigma}):\bar{\Sigma}\in \limfunc{cl}(\mathcal{L}(\mathfrak{C})),\ \limfunc{%
rank}(\bar{\Sigma})<n-\dim (\mathcal{L})\right\} ,
\end{equation*}%
where the closure $\limfunc{cl}(\cdot )$ is to be understood w.r.t. $\mathbb{%
R}^{n\times n}$.
\end{definition}

Recalling the definition of $I_{0}(\mathcal{L})$, it is easy to see that $%
I_{0}(\mathcal{L})=\left\{ i:1\leq i\leq n,\pi _{\mathcal{L}^{\bot
},i}=0\right\} $ holds, where $\pi _{\mathcal{L}^{\bot },i}$ denotes the $i$%
-th column of $\Pi _{\mathcal{L}^{\bot }}$. Also recall that $I_{1}(\mathcal{%
L})$ is nonempty in case $\dim (\mathcal{L})<n$ holds. The characterization
of $\mathbb{J}(\mathcal{L},\mathfrak{C}_{Het})$ is now given in the next
proposition.

\begin{proposition}
\label{characterization}Suppose $\dim (\mathcal{L})<n$ holds. Then the set $%
\mathbb{J}(\mathcal{L},\mathfrak{C}_{Het})$ is given by%
\begin{equation}
\left\{ \func{span}\left( \left\{ \pi _{\mathcal{L}^{\bot },i}:i\in
I\right\} \right) :\emptyset \neq I\subseteq I_{1}(\mathcal{L})\text{, }\dim
\left( \func{span}\left( \left\{ \pi _{\mathcal{L}^{\bot },i}:i\in I\right\}
\right) \right) <n-\dim (\mathcal{L})\right\} .  \label{charac}
\end{equation}
\end{proposition}

This proposition is a special case of Proposition \ref{characterization_2}
given below since $\mathfrak{C}_{Het}$ coincides with $\mathfrak{C}%
_{(n_{1},\ldots ,n_{m})}$ in case $m=n$ and $n_{j}=1$ for all $j=1,\ldots ,m$%
.

We next turn to the characterization of $\mathbb{J}(\mathcal{L},\mathfrak{C}%
_{(n_{1},\ldots ,n_{m})})$, where $\mathfrak{C}_{(n_{1},\ldots ,n_{m})}$ has
been defined in Appendix \ref{app_d}. Here $m\in \mathbb{N}$, and $n_{j}\in 
\mathbb{N}$ for $j=1,\ldots ,m$ satisfy $\sum_{j=1}^{m}n_{j}=n$. Consider
the partition of the set $\{1,\ldots ,n\}$ into the intervals $%
(n_{0}^{+},n_{1}^{+}]$, $(n_{1}^{+},n_{2}^{+}]$,..., $%
(n_{m-1}^{+},n_{m}^{+}] $ where $n_{j}^{+}$ has been defined in Appendix \ref%
{app_d}. Let $\boldsymbol{I}_{(n_{1},\ldots ,n_{m})}$ consist of all
non-empty subsets $I$ of $\{1,\ldots ,n\}$ that can be represented as a
union of intervals of the form $(n_{j-1}^{+},n_{j}^{+}]$.

\begin{proposition}
\label{characterization_2}Suppose $\dim (\mathcal{L})<n$ holds. Let $m\in 
\mathbb{N}$, and let $n_{j}\in \mathbb{N}$ for $j=1,\ldots ,m$ satisfy $%
\sum_{j=1}^{m}n_{j}=n$. Then the set $\mathbb{J}(\mathcal{L},\mathfrak{C}%
_{(n_{1},\ldots ,n_{m})})$ is given by%
\begin{equation}
\left\{ \func{span}\left( \left\{ \pi _{\mathcal{L}^{\bot },i}:i\in
I\right\} \right) :I\in \boldsymbol{I}_{(n_{1},\ldots ,n_{m})}\text{, }%
\emptyset \neq I\cap I_{1}(\mathcal{L})\text{, }\dim \left( \func{span}%
\left( \left\{ \pi _{\mathcal{L}^{\bot },i}:i\in I\right\} \right) \right)
<n-\dim (\mathcal{L})\right\} .  \label{charac_2}
\end{equation}
\end{proposition}

Note that in (\ref{charac_2}) we have $\func{span}\left( \left\{ \pi _{%
\mathcal{L}^{\bot },i}:i\in I\right\} \right) =\func{span}\left( \left\{ \pi
_{\mathcal{L}^{\bot },i}:i\in I\cap I_{1}(\mathcal{L})\right\} \right) $.

\textbf{Proof:} Suppose $\mathcal{S}$ is an element of $\mathbb{J}(\mathcal{L%
},\mathfrak{C}_{(n_{1},\ldots ,n_{m})})$. Then there exist a sequence $%
\Sigma _{m}\in \mathfrak{C}_{(n_{1},\ldots ,n_{m})}$ such that $\Pi _{%
\mathcal{L}^{\bot }}\Sigma _{m}\Pi _{\mathcal{L}^{\bot }}/\left\Vert \Pi _{%
\mathcal{L}^{\bot }}\Sigma _{m}\Pi _{\mathcal{L}^{\bot }}\right\Vert $
converges to a limit $\bar{\Sigma}$, say, in $\mathbb{R}^{n\times n}$ with $%
\func{span}(\bar{\Sigma})=\mathcal{S}$. Now,%
\begin{eqnarray*}
\Pi _{\mathcal{L}^{\bot }}\Sigma _{m}\Pi _{\mathcal{L}^{\bot }}/\left\Vert
\Pi _{\mathcal{L}^{\bot }}\Sigma _{m}\Pi _{\mathcal{L}^{\bot }}\right\Vert
&=&\left\Vert \Pi _{\mathcal{L}^{\bot }}\Sigma _{m}\Pi _{\mathcal{L}^{\bot
}}\right\Vert ^{-1}\sum_{i=1}^{n}\tau _{i}^{2}(m)\pi _{\mathcal{L}^{\bot
},i}\pi _{\mathcal{L}^{\bot },i}^{\prime } \\
&=&\left\Vert \Pi _{\mathcal{L}^{\bot }}\Sigma _{m}\Pi _{\mathcal{L}^{\bot
}}\right\Vert ^{-1}\sum_{j=1}^{m}\sum_{i\in (n_{j-1}^{+},n_{j}^{+}]}\tau
_{i}^{2}(m)\pi _{\mathcal{L}^{\bot },i}\pi _{\mathcal{L}^{\bot },i}^{\prime }
\\
&=&\sum_{j:(n_{j-1}^{+},n_{j}^{+}]\cap I_{1}(\mathcal{L})\neq \emptyset
}\left\Vert \Pi _{\mathcal{L}^{\bot }}\Sigma _{m}\Pi _{\mathcal{L}^{\bot
}}\right\Vert ^{-1}\tau _{n_{j}^{+}}^{2}(m)\sum_{i\in
(n_{j-1}^{+},n_{j}^{+}]}\pi _{\mathcal{L}^{\bot },i}\pi _{\mathcal{L}^{\bot
},i}^{\prime },
\end{eqnarray*}%
where $\tau _{i}^{2}(m)$ denotes the $i$-th diagonal element of $\Sigma _{m}$%
. Here we have used the fact that variances are constant within groups, as
well as that $\pi _{\mathcal{L}^{\bot },i}=0$ for all $i\in
(n_{j-1}^{+},n_{j}^{+}]$ if $(n_{j-1}^{+},n_{j}^{+}]$ is disjoint from $%
I_{1}(\mathcal{L})$. Also note that the outer sum extends over a nonempty
index set since $\limfunc{card}(I_{1}(\mathcal{L}))\geq 1$ must hold in view
of $\dim (\mathcal{L})<n$. Since the l.h.s. converges to the limit $\bar{%
\Sigma}\in \mathbb{R}^{n\times n}$, since the r.h.s. is bounded from below
in the Loewner order by 
\begin{equation*}
\left\Vert \Pi _{\mathcal{L}^{\bot }}\Sigma _{m}\Pi _{\mathcal{L}^{\bot
}}\right\Vert ^{-1}\tau _{n_{j}^{+}}^{2}(m)\sum_{i\in
(n_{j-1}^{+},n_{j}^{+}]}\pi _{\mathcal{L}^{\bot },i}\pi _{\mathcal{L}^{\bot
},i}^{\prime },
\end{equation*}%
for every $j$ appearing in the range of the outer sum, and since $\pi _{%
\mathcal{L}^{\bot },i}\neq 0$ for at least one $i\in (n_{j-1}^{+},n_{j}^{+}]$
holds when $(n_{j-1}^{+},n_{j}^{+}]\cap I_{1}(\mathcal{L})\neq \emptyset $,
it follows that the sequence 
\begin{equation*}
(\left\Vert \Pi _{\mathcal{L}^{\bot }}\Sigma _{m}\Pi _{\mathcal{L}^{\bot
}}\right\Vert ^{-1}\tau _{n_{j}^{+}}^{2}(m):m\in \mathbb{N)}
\end{equation*}%
is bounded for every $j$ satisfying $(n_{j-1}^{+},n_{j}^{+}]\cap I_{1}(%
\mathcal{L})\neq \emptyset $. Possibly after passing to a subsequence, we
may thus assume that these sequences converge to nonnegative real numbers $%
\gamma _{j}$ for such $j$. It follows that 
\begin{eqnarray*}
\bar{\Sigma} &=&\sum_{j:(n_{j-1}^{+},n_{j}^{+}]\cap I_{1}(\mathcal{L})\neq
\emptyset }\gamma _{j}\sum_{i\in (n_{j-1}^{+},n_{j}^{+}]}\pi _{\mathcal{L}%
^{\bot },i}\pi _{\mathcal{L}^{\bot },i}^{\prime } \\
&=&\sum_{j:(n_{j-1}^{+},n_{j}^{+}]\cap I_{1}(\mathcal{L})\neq \emptyset
}\sum_{i\in (n_{j-1}^{+},n_{j}^{+}]}\gamma _{j}^{1/2}\pi _{\mathcal{L}^{\bot
},i}\left( \gamma _{j}^{1/2}\pi _{\mathcal{L}^{\bot },i}\right) ^{\prime }.
\end{eqnarray*}%
Let $I$ be the union of those intervals $(n_{j-1}^{+},n_{j}^{+}]$ satisfying
(i) $(n_{j-1}^{+},n_{j}^{+}]\cap I_{1}(\mathcal{L})\neq \emptyset $ and (ii) 
$\gamma _{j}>0$. Note that $I$ cannot be the empty set as this would imply $%
\bar{\Sigma}=0$, which is impossible since it is the limit of a sequence of
matrices residing in the unit sphere of $\mathbb{R}^{n\times n}$.
Furthermore, by construction, $I\in \boldsymbol{I}_{(n_{1},\ldots ,n_{m})}$
and $I\cap I_{1}(\mathcal{L})\neq \emptyset $ hold. Using the fact that $%
\func{span}(\sum_{l=1}^{L}A_{l}A_{l}^{\prime })=\func{span}(A_{1},\ldots
,A_{L})$ holds for arbitrary real matrices of the same row-dimension, we
obtain $\mathcal{S}=\func{span}(\bar{\Sigma})=\func{span}\left( \left\{ \pi
_{\mathcal{L}^{\bot },i}:i\in I\right\} \right) $ for the before constructed
set $I$. [Note that $\pi _{\mathcal{L}^{\bot },i}=0$ if $i\in
(n_{j-1}^{+},n_{j}^{+}]$ but $i\notin I_{1}(\mathcal{L})$.] Since $\mathcal{S%
}$, being an element of $\mathbb{J}(\mathcal{L},\mathfrak{C}_{(n_{1},\ldots
,n_{m})})$, satisfies $\dim (\mathcal{S})<n-\dim (\mathcal{L})$, we have
established that $\mathcal{S}$ is also an element of (\ref{charac_2}).

To prove the converse, suppose that $\mathcal{S}$ is an element of (\ref%
{charac_2}), i.e., that $\mathcal{S}=\func{span}\left( \left\{ \pi _{%
\mathcal{L}^{\bot },i}:i\in I\right\} \right) $ for some $I\in \boldsymbol{I}%
_{(n_{1},\ldots ,n_{m})}$ with $\emptyset \neq I\cap I_{1}(\mathcal{L})$ and
that $\dim \left( \mathcal{S}\right) <n-\dim (\mathcal{L})$ holds. Note that 
$\func{card}(I)<n$ holds, since otherwise $\mathcal{S}=\mathcal{L}^{\bot }$
would follow, contradicting $\dim \left( \mathcal{S}\right) <n-\dim (%
\mathcal{L})$. Also note that $\func{card}(I)\geq 1$ as $\emptyset \neq
I\cap I_{1}(\mathcal{L})$. Define diagonal $n\times n$ matrices $\Sigma _{m}$
via their diagonal elements%
\begin{equation*}
\tau _{i}^{2}(m)=\left\{ 
\begin{array}{cc}
\left( \func{card}(I)\right) ^{-1}-\delta _{m} & \text{if }i\in I \\ 
\left( \func{card}(I)/(n-\func{card}(I))\right) \delta _{m} & \text{if }%
i\notin I%
\end{array}%
\right.
\end{equation*}%
where $0<\delta _{m}<1/\func{card}(I)$ with $\delta _{m}\rightarrow 0$ for $%
m\rightarrow \infty $. Then $\tau _{i}^{2}(m)>0$ as well as $%
\sum_{i=1}^{n}\tau _{i}^{2}(m)=1$ hold, and $\tau
_{n_{j-1}^{+}+1}^{2}(m)=\ldots =\tau _{n_{j}^{+}}^{2}(m)$ holds for $%
j=1,\ldots ,m$ since $I\in \boldsymbol{I}_{(n_{1},\ldots ,n_{m})}$. That is, 
$\Sigma _{m}$ belongs to $\mathfrak{C}_{(n_{1},\ldots ,n_{m})}$. Obviously, $%
\Sigma _{m}$ converges to a diagonal matrix $\Sigma ^{\ast }$ with diagonal
elements given by%
\begin{equation*}
\tau _{i}^{\ast 2}=\left\{ 
\begin{array}{cc}
\left( \func{card}(I)\right) ^{-1} & \text{if }i\in I \\ 
0 & \text{if }i\notin I%
\end{array}%
\right. .
\end{equation*}%
Consequently, $\Pi _{\mathcal{L}^{\bot }}\Sigma _{m}\Pi _{\mathcal{L}^{\bot
}}/\left\Vert \Pi _{\mathcal{L}^{\bot }}\Sigma _{m}\Pi _{\mathcal{L}^{\bot
}}\right\Vert $ converges to $\bar{\Sigma}:=\Pi _{\mathcal{L}^{\bot }}\Sigma
^{\ast }\Pi _{\mathcal{L}^{\bot }}/\left\Vert \Pi _{\mathcal{L}^{\bot
}}\Sigma ^{\ast }\Pi _{\mathcal{L}^{\bot }}\right\Vert $, since $\Pi _{%
\mathcal{L}^{\bot }}\Sigma ^{\ast }\Pi _{\mathcal{L}^{\bot }}\neq 0$ in view
of 
\begin{equation*}
\Pi _{\mathcal{L}^{\bot }}\Sigma ^{\ast }\Pi _{\mathcal{L}^{\bot
}}=\sum_{i=1}^{n}\tau _{i}^{\ast 2}\pi _{\mathcal{L}^{\bot },i}\pi _{%
\mathcal{L}^{\bot },i}^{\prime }=\left( \func{card}(I)\right)
^{-1}\sum_{i\in I}\pi _{\mathcal{L}^{\bot },i}\pi _{\mathcal{L}^{\bot
},i}^{\prime }
\end{equation*}%
and the fact that $\emptyset \neq I\cap I_{1}(\mathcal{L})$ holds and thus $%
\pi _{\mathcal{L}^{\bot },i}\neq 0$ must hold at least for one $i\in I$.
Again using $\func{span}(\sum_{l=1}^{L}A_{l}A_{l}^{\prime })=\func{span}%
(A_{1},\ldots ,A_{L})$ we arrive at 
\begin{equation*}
\func{span}(\bar{\Sigma})=\func{span}(\Pi _{\mathcal{L}^{\bot }}\Sigma
^{\ast }\Pi _{\mathcal{L}^{\bot }})=\limfunc{span}\left( \left( \func{card}%
(I)\right) ^{-1}\sum_{i\in I}\pi _{\mathcal{L}^{\bot },i}\pi _{\mathcal{L}%
^{\bot },i}^{\prime }\right) =\func{span}\left( \left\{ \pi _{\mathcal{L}%
^{\bot },i}:i\in I\right\} \right) =\mathcal{S}.
\end{equation*}%
Because we have assumed that $\dim \left( \mathcal{S}\right) <n-\dim (%
\mathcal{L})$ holds, the preceding display shows that $\mathcal{S}\in 
\mathbb{J}(\mathcal{L},\mathfrak{C}_{(n_{1},\ldots ,n_{m})})$. $\blacksquare 
$

\begin{remark}
Note that $\mathbb{J}(\mathcal{L},\mathfrak{C}_{(n_{1},\ldots ,n_{m})})$ is
empty if $m=1$ (as can be seen directly from the definition of $\mathbb{J}(%
\mathcal{L},\mathfrak{C}_{(n_{1},\ldots ,n_{m})})$ or from (\ref{charac_2})).
\end{remark}

\begin{remark}
\label{rem_conc_spaces}It is easy to see that the concentration spaces of $%
\mathfrak{C}_{Het}$ in the sense of \cite{PP2016} are precisely given by all
spaces of the form $\func{span}\left( \left\{ e_{i}(n):i\in I\right\}
\right) $ where $I$ varies through all subsets of $\left\{ 1,\ldots
,n\right\} $ that satisfy $0<\limfunc{card}(I)<n$. More generally, the
concentration spaces of $\mathfrak{C}_{(n_{1},\ldots ,n_{m})}$ are precisely
given by all spaces of the form $\func{span}\left( \left\{ e_{i}(n):i\in
I\right\} \right) $ where $I\in \boldsymbol{I}_{(n_{1},\ldots ,n_{m})}$
satisfies $0<\limfunc{card}(I)<n$. [In view of Remark 5.1(i) in \cite{PP3}
these results correspond to the case $\dim (\mathcal{L})=0$ in the preceding
propositions.]
\end{remark}

\section{Appendix: Proofs for Section \protect\ref{sec_size_control} and
Appendix \protect\ref{sec_size_control_2}\label{app_B}}

The facts collected in the subsequent remark will be used in the proofs
further below.

\begin{remark}
\label{F-type}(i) Suppose Assumption \ref{R_and_X} holds. Then the test
statistic $T_{Het}$ is a non-sphericity corrected F-type test statistic in
the sense of Section 5.4 in \cite{PP2016}. More precisely, $T_{Het}$ is of
the form (28) in \cite{PP2016} and Assumption 5 in the same reference is
satisfied with $\check{\beta}=\hat{\beta}$, $\check{\Omega}=\hat{\Omega}%
_{Het}$, and $N=\emptyset $. Furthermore, the set $N^{\ast }$ defined in
(27) of \cite{PP2016} satisfies $N^{\ast }=\mathsf{B}$. And also Assumptions
6 and 7 of \cite{PP2016} are satisfied. All these claims follow easily in
view of Lemma 4.1 in \cite{PP2016}, see also the proof of Theorem 4.2 in
that reference.

(ii) The test statistic $T_{uc}$ is also a non-sphericity corrected F-type
test statistic in the sense of Section 5.4 in \cite{PP2016} (terminology
being somewhat unfortunate here as no correction for the non-sphericity is
being attempted). More precisely, $T_{uc}$ is of the form (28) in \cite%
{PP2016} and Assumption 5 in the same reference is satisfied with $\check{%
\beta}=\hat{\beta}$, $\check{\Omega}=\hat{\sigma}^{2}R\left( X^{\prime
}X\right) ^{-1}R^{\prime }$, and $N=\emptyset $. Furthermore, the set $%
N^{\ast }$ defined in (27) of \cite{PP2016} satisfies $N^{\ast }=\limfunc{%
span}(X)$. And also Assumptions 6 and 7 of \cite{PP2016} are satisfied. All
these claims are evident (and obviously do not rely on Assumption \ref%
{R_and_X}).

(iii) We note that any non-sphericity corrected F-type test statistic (for
testing (\ref{testing problem})) in the sense of Section 5.4 in \cite{PP2016}%
, i.e., any test statistic $T$ of the form (28) in \cite{PP2016} that also
satisfies Assumption 5 in that reference, is invariant under the group $G(%
\mathfrak{M}_{0})$. Furthermore, the associated set $N^{\ast }$ defined in
(27) of \cite{PP2016} is even invariant under the larger group $G(\mathfrak{M%
})$. See Sections 5.1 and 5.4 of \cite{PP2016} as well as Lemma 5.16 in \cite%
{PP3} for more information.
\end{remark}

\textbf{Proof of Theorem \ref{Hetero_Robust}:} We first prove Part (b). We
apply Part (b) of Theorem \ref{theorem_groupwise_hetero} with $n_{j}=1$ for $%
j=1,\ldots ,n=m$ observing that then $\mathfrak{C}_{(n_{1},\ldots ,n_{m})}=%
\mathfrak{C}_{Het}$ and that condition (\ref{non-incl_groupHet}) reduces to (%
\ref{non-incl_Het}) (exploiting that $\mathsf{B}$ is a finite union of
proper linear subspaces as discussed in Lemma \ref{lem_B}). This establishes
(\ref{size-control_Het}). The final claim in Part (b) of the theorem follows
from Part (b) of Theorem \ref{theorem_groupwise_hetero}, if we can show that 
$\alpha ^{\ast }$ and $C^{\ast }$ given there can be written as claimed in
Theorem \ref{Hetero_Robust}: To this end we proceed as follows:\footnote{%
Alternatively, one could base a proof on Lemma C.1 in \cite{PP4}.} Choose an
element $\mu _{0}$ of $\mathfrak{M}_{0}$. Observe that $I_{1}(\mathfrak{M}%
_{0}^{lin})\neq \emptyset $ (since $\dim (\mathfrak{M}_{0}^{lin})=k-q<n$),
and that for every $i\in I_{1}(\mathfrak{M}_{0}^{lin})$ the linear space $%
\mathcal{S}_{i}=\func{span}(\Pi _{\left( \mathfrak{M}_{0}^{lin}\right)
^{\bot }}e_{i}(n))$ is $1$-dimensional (since $\mathcal{S}_{i}=\{0\}$ is
impossible in view of $i\in I_{1}(\mathfrak{M}_{0}^{lin})$), and belongs to $%
\mathbb{J}(\mathfrak{M}_{0}^{lin},\mathfrak{C}_{Het})$ (since $n-k+q>1=\dim (%
\mathcal{S}_{i})$ holds) in view of Proposition \ref{characterization} in
Section \ref{app_char}. Since $T_{Het}$ is $G(\mathfrak{M}_{0})$-invariant
(Remark \ref{F-type} above), it follows that $T_{Het}$ is constant on $(\mu
_{0}+\mathcal{S}_{i})\backslash \left\{ \mu _{0}\right\} $, cf. the
beginning of the proof of Lemma 5.11 in \cite{PP3}. Hence, $\mathcal{S}_{i}$
belongs to $\mathbb{H}$ (defined in Lemma 5.11 in \cite{PP3}) and
consequently for $C^{\ast }$ as defined in that lemma%
\begin{equation*}
C^{\ast }\geq \max \left\{ T_{Het}(\mu _{0}+\Pi _{\left( \mathfrak{M}%
_{0}^{lin}\right) ^{\bot }}e_{i}(n)):i\in I_{1}(\mathfrak{M}%
_{0}^{lin})\right\}
\end{equation*}%
must hold (recall that $\Pi _{\left( \mathfrak{M}_{0}^{lin}\right) ^{\bot
}}e_{i}(n)\neq 0$). To prove the opposite inequality, let $\mathcal{S}$ be
an arbitrary element of $\mathbb{H}$, i.e., $\mathcal{S}\in \mathbb{J}(%
\mathfrak{M}_{0}^{lin},\mathfrak{C}_{Het})$ and $T_{Het}$ is $\lambda _{\mu
_{0}+\mathcal{S}}$-almost everywhere equal to a constant $C(\mathcal{S})$,
say. Then Proposition \ref{characterization} in Section \ref{app_char} shows
that $\mathcal{S}_{i}\subseteq \mathcal{S}$ holds for some $i\in I_{1}(%
\mathfrak{M}_{0}^{lin})$. Because of Condition (\ref{non-incl_Het}) we have $%
\mathcal{S}_{i}\nsubseteqq \mathsf{B}$ since $\Pi _{\left( \mathfrak{M}%
_{0}^{lin}\right) ^{\bot }}e_{i}(n)$ and $e_{i}(n)$ differ only by an
element of $\mathfrak{M}_{0}^{lin}\subseteq \limfunc{span}(X)$ and since $%
\mathsf{B}+\limfunc{span}(X)=\mathsf{B}$. Thus $\mu _{0}+\mathcal{S}%
_{i}\nsubseteqq \mathsf{B}$ by the same argument as $\mu _{0}\in \mathfrak{M}%
_{0}\subseteq \limfunc{span}(X)$. We thus can find $s\in \mathcal{S}_{i}$
such that $\mu _{0}+s\notin \mathsf{B}$. Note that $s\neq 0$ must hold,
since $\mu _{0}\in \mathfrak{M}_{0}\subseteq \limfunc{span}(X)\subseteq 
\mathsf{B}$. In particular, $T_{Het}$ is continuous at $\mu _{0}+s$, since $%
\mu _{0}+s\notin \mathsf{B}$. Now, for every open ball $A_{\varepsilon }$ in 
$\mathbb{R}^{n}$ with center $s$ and radius $\varepsilon >0$ we can find an
element $a_{\varepsilon }\in A_{\varepsilon }\cap \mathcal{S}$ such that $%
T_{Het}(\mu _{0}+a_{\varepsilon })=C(\mathcal{S})$. Since $a_{\varepsilon
}\rightarrow s$ for $\varepsilon \rightarrow 0$, it follows that $C(\mathcal{%
S})=T_{Het}(\mu _{0}+s)$. Since $s\neq 0$ and since $T_{Het}$ is constant on 
$(\mu _{0}+\mathcal{S}_{i})\backslash \left\{ \mu _{0}\right\} $ as shown
before, we can conclude that $C(\mathcal{S})=T_{Het}(\mu _{0}+s)=T_{Het}(\mu
_{0}+\Pi _{\left( \mathfrak{M}_{0}^{lin}\right) ^{\bot }}e_{i}(n))$, where
we recall that $\Pi _{\left( \mathfrak{M}_{0}^{lin}\right) ^{\bot
}}e_{i}(n)\neq 0$. But this now implies%
\begin{equation*}
C^{\ast }=\max \left\{ T_{Het}(\mu _{0}+\Pi _{\left( \mathfrak{M}%
_{0}^{lin}\right) ^{\bot }}e_{i}(n)):i\in I_{1}(\mathfrak{M}%
_{0}^{lin})\right\} .
\end{equation*}%
Using $G(\mathfrak{M}_{0})$-invariance of $T_{Het}$ we conclude that%
\begin{equation*}
C^{\ast }=\max \left\{ T_{Het}(\mu _{0}+e_{i}(n)):i\in I_{1}(\mathfrak{M}%
_{0}^{lin})\right\} .
\end{equation*}%
The expression for $\alpha ^{\ast }$ given in the theorem now follows
immediately from the expression for $\alpha ^{\ast }$ given in Part (b) of
Theorem \ref{theorem_groupwise_hetero}.

We next prove Part (a): Apply Part (a) of Theorem \ref%
{theorem_groupwise_hetero} with $n_{j}=1$ for $j=1,\ldots ,n=m$ observing
that then $\mathfrak{C}_{(n_{1},\ldots ,n_{m})}=\mathfrak{C}_{Het}$ and that
condition (\ref{non-incl_groupHet_uncorr}) reduces to (\ref%
{non-incl_Het_uncorr}) (exploiting that $\limfunc{span}(X)$ is a linear
space). This establishes (\ref{size-control_Het_uncorr}). The final claim in
Part (a) of the theorem follows similarly as the corresponding claim of Part
(b) upon replacing the set $\mathsf{B}$ by $\limfunc{span}(X)$ in the
argument and by noting that $T_{uc}$ is $G(\mathfrak{M}_{0})$-invariant.

Part (c) follows from Part (c) of Theorem \ref{theorem_groupwise_hetero}
upon setting $n_{j}=1$ for $j=1,\ldots ,n=m$ (and upon noting that then the
conditions in Theorem \ref{theorem_groupwise_hetero} reduce to the
conditions of the present theorem). $\blacksquare $

\textbf{Proof of Proposition \ref{rem_C*}:} Follows from Part A.1 of
Proposition 5.12 of \cite{PP3} and the sentence following this proposition.
Note that the assumptions of this proposition have been verified in the
proof of Theorem \ref{Hetero_Robust} (see also the proof of Theorem \ref%
{theorem_groupwise_hetero}, on which the proof of Theorem \ref{Hetero_Robust}
is based), where it is also shown that the quantity $C^{\ast }$ used in
Proposition 5.12 of \cite{PP3} coincides with $C^{\ast }$ defined in Theorem %
\ref{Hetero_Robust}. $\blacksquare $

We note that the result for $T_{Het}$ in Proposition \ref{rem_C*} can also
be obtained from Theorem 4.2 in \cite{PP2016}.

\textbf{Proof of Proposition \ref{rem_necessity}:} (a) This can be seen as
follows (cf. also the discussion on p.302 of \cite{PP2016}): By Remark \ref%
{F-type} above, $T_{uc}$ satisfies the assumptions in Corollary 5.17 in \cite%
{PP2016} (with $\check{\beta}=\hat{\beta}$, $\check{\Omega}(y)=\hat{\sigma}%
^{2}(y)R\left( X^{\prime }X\right) ^{-1}R^{\prime }$, $N=\emptyset $, and $%
N^{\ast }=\limfunc{span}(X)$). Let $e_{i}(n)$ be one of the standard basis
vectors with $i\in I_{1}(\mathfrak{M}_{0}^{lin})$ that does belong to $%
\limfunc{span}(X)$. Set $\mathcal{Z}=$ $\limfunc{span}(e_{i}(n))$ and note
that this is a concentration space of $\mathfrak{C}_{Het}$, cf.~Remark \ref%
{rem_conc_spaces} in Appendix \ref{app_char}. The nonnegative definiteness
assumption on $\check{\Omega}$ in Part 3 of Corollary 5.17 in \cite{PP2016}
is clearly satisfied. We also have $\check{\Omega}(\lambda e_{i}(n))=0$
(since $e_{i}(n)\in \limfunc{span}(X)$) for every $\lambda \in \mathbb{R}$
and $R\hat{\beta}(\lambda e_{i}(n))\neq 0$ for every $\lambda \in \mathbb{R}%
\backslash \{0\}$ (since $e_{i}(n)\in \limfunc{span}(X)$ but $e_{i}(n)\notin 
\mathfrak{M}_{0}^{lin}$ in view of $i\in I_{1}(\mathfrak{M}_{0}^{lin})$).
Part 3 of Corollary 5.17 in \cite{PP2016} then proves the claim for $C>0$. A
fortiori it then also holds for all real $C$.

(b) This follows for $C>0$ from Part 3 of Theorem 4.2 in \cite{PP2016} upon
observing that a vector $e_{i}(n)$ satisfying $e_{i}(n)\in \limfunc{span}(X)$
for some $i\in I_{1}(\mathfrak{M}_{0}^{lin})$ clearly satisfies $%
B(e_{i}(n))=0$ (as $e_{i}(n)\in \limfunc{span}(X)$) and $R\hat{\beta}%
(e_{i}(n))\neq 0$ (since $e_{i}(n)\in \limfunc{span}(X)$ but $e_{i}(n)\notin 
\mathfrak{M}_{0}^{lin}$ in view of $i\in I_{1}(\mathfrak{M}_{0}^{lin})$). A
fortiori it then also holds for all real $C$. $\blacksquare $

\textbf{Proof of Theorem \ref{theorem_groupwise_hetero}:} We first prove
Part (b). We wish to apply Part A of Proposition 5.12 of \cite{PP3} with $%
\mathfrak{C}=\mathfrak{C}_{(n_{1},\ldots ,n_{m})}$, $T=T_{Het}$, $\mathcal{L}%
=\mathfrak{M}_{0}^{lin}$, and $\mathcal{V}=\{0\}$. First, note that $\dim (%
\mathfrak{M}_{0}^{lin})=k-q<n$. Second, under Assumption \ref{R_and_X}, $%
T_{Het}$ is a non-sphericity corrected F-type test with $N^{\ast }=\mathsf{B}
$, which is a closed $\lambda _{\mathbb{R}^{n}}$-null set (see Remarks \ref%
{rem:GM0} and \ref{F-type} as well as Lemma \ref{lem_B}). Hence, the general
assumptions on $T=T_{Het}$, on $N^{\dag }=N^{\ast }=\mathsf{B}$, on $%
\mathcal{L}=\mathfrak{M}_{0}^{lin}$, as well as on $\mathcal{V}$ in
Proposition 5.12 of \cite{PP3} are satisfied in view of Part 1 of Lemma 5.16
in the same reference. [Alternatively, this can be gleaned from Lemma \ref%
{lem_B} and the attending discussion.] Next, observe that condition (\ref%
{non-incl_groupHet}) is equivalent to%
\begin{equation*}
\func{span}\left( \left\{ \Pi _{\left( \mathfrak{M}_{0}^{lin}\right) ^{\bot
}}e_{i}(n):i\in (n_{j-1}^{+},n_{j}^{+}]\right\} \right) \nsubseteqq \mathsf{B%
}
\end{equation*}%
for every $j=1,\ldots ,m$, such that $(n_{j-1}^{+},n_{j}^{+}]\cap I_{1}(%
\mathfrak{M}_{0}^{lin})\neq \emptyset $, since $\Pi _{\left( \mathfrak{M}%
_{0}^{lin}\right) ^{\bot }}e_{i}(n)$ and $e_{i}(n)$ differ only by an
element of $\mathfrak{M}_{0}^{lin}\subseteq \limfunc{span}(X)$ and since $%
\mathsf{B}+\limfunc{span}(X)=\mathsf{B}$ (as noted in Lemma \ref{lem_B}). In
view of Proposition \ref{characterization_2} in Appendix \ref{app_char},
this implies that any $\mathcal{S}\in \mathbb{J}(\mathfrak{M}_{0}^{lin},%
\mathfrak{C}_{(n_{1},\ldots ,n_{m})})$ is not contained in $\mathsf{B}$, and
thus not in $N^{\dag }$. Using $\mathfrak{M}_{0}\subseteq \limfunc{span}(X)$
and $\mathsf{B}+\limfunc{span}(X)=\mathsf{B}$, it follows that $\mu _{0}+%
\mathcal{S}\nsubseteqq \mathsf{B}=N^{\dag }$ for every $\mu _{0}\in 
\mathfrak{M}_{0}$. Since $\mu _{0}+\mathcal{S}$ is an affine space and $%
N^{\dag }=\mathsf{B}$ is a finite union of proper affine (even linear)
spaces under Assumption \ref{R_and_X} as discussed in Lemma \ref{lem_B}, we
may conclude (cf. Corollary 5.6 in \cite{PP3} and its proof) that $\lambda
_{\mu _{0}+\mathcal{S}}(N^{\dag })=0$ for every $\mathcal{S}\in \mathbb{J}(%
\mathfrak{M}_{0}^{lin},\mathfrak{C}_{(n_{1},\ldots ,n_{m})})$ and every $\mu
_{0}\in \mathfrak{M}_{0}$. This completes the verification of the
assumptions of Proposition 5.12 in \cite{PP3} that are not specific to Part
A (or Part B) of this proposition. We next verify the assumptions specific
to Part A of this proposition: Assumption (a) is satisfied (even for every $%
C\in \mathbb{R}$) as a consequence of Part 2 of Lemma 5.16 in \cite{PP3} and
of Remark \ref{F-type}(i) above. And Assumption (b) in Part A follows from
Lemma 5.19 of \cite{PP3}, since $T_{Het}$ results as a special case of the
test statistics $T_{GQ}$ defined in Section 3.4 of \cite{PP3} upon choosing $%
\mathcal{W}_{n}^{\ast }=n^{-1}\limfunc{diag}(d_{i})$. Part A of Proposition
5.12 of \cite{PP3} now immediately delivers claim (\ref%
{size-control_groupHet}), since $C^{\ast }<\infty $ as noted in that
proposition. That $C^{\ast }$ and $\alpha ^{\ast }$ do not depend on the
choice of $\mu _{0}\in \mathfrak{M}_{0}$ is an immediate consequence of $G(%
\mathfrak{M}_{0})$-invariance of $T_{Het}$. Also note that $\alpha ^{\ast }$
as defined in the theorem coincides with $\alpha ^{\ast }$ as defined in
Proposition 5.12 of \cite{PP3} in view of $G(\mathfrak{M}_{0})$-invariance
of $T_{Het}$. Positivity of $\alpha ^{\ast }$ then follows from Part 5 of
Lemma 5.15 in \cite{PP2016} in view of Remark \ref{F-type}(i), noting that $%
\lambda _{\mathbb{R}^{n}}$ and $P_{\mu _{0},\Sigma }$ are equivalent
measures (since $\Sigma \in \mathfrak{C}_{Het}$ is positive definite); cf.
Remark 5.13(vi) in \cite{PP3}. In case $\alpha <\alpha ^{\ast }$, the
remaining claim in Part (b) of the theorem, namely that equality can be
achieved in (\ref{size-control_groupHet}), follows from the definition of $%
C^{\ast }$ in Lemma 5.11 of \cite{PP3} and from Part A.2 of Proposition 5.12
of \cite{PP3} (and the observation immediately following that proposition
allowing one to drop the suprema w.r.t. $\mu _{0}$ and $\sigma ^{2}$, and to
set $\sigma ^{2}=1$); in case $\alpha =\alpha ^{\ast }<1$, it follows from
Remarks 5.13(i),(ii) in \cite{PP3} using Lemma 5.16 in the same reference.

The proof of Part (a) proceeds along the same lines with some minor
differences: Observe that $T_{uc}$ is a non-sphericity corrected F-type test
with $N^{\dag }=N^{\ast }=\func{span}(X)$, which obviously is a closed $%
\lambda _{\mathbb{R}^{n}}$-null set (see Remark \ref{F-type}(ii)), showing
similarly that the general assumptions on $T=T_{uc}$, on $N^{\dag }=N^{\ast
}=\func{span}(X)$, as well as on $\mathcal{L}=\mathfrak{M}_{0}^{lin}$ in
Proposition 5.12 of \cite{PP3} are again satisfied (with $\mathfrak{C}=%
\mathfrak{C}_{(n_{1},\ldots ,n_{m})}$). A similar, even simpler argument as
in the proof of Part (b), again shows that condition (\ref%
{non-incl_groupHet_uncorr}) implies $\lambda _{\mu _{0}+\mathcal{S}}(N^{\dag
})=0$ for every $\mathcal{S}\in \mathbb{J}(\mathfrak{M}_{0}^{lin},\mathfrak{C%
}_{(n_{1},\ldots ,n_{m})})$ and every $\mu _{0}\in \mathfrak{M}_{0}$, thus
completing the verification of the assumptions of Proposition 5.12 of \cite%
{PP3} that are not specific to Part A (or Part B) of this proposition.
Verification of Assumption (a) in Part A of Proposition 5.12 of \cite{PP3}
proceeds exactly as before. For Assumption (b) we now use Lemma 5.19(iii) of 
\cite{PP3}, since $T_{uc}$ results as a special case of the test statistics $%
T_{E,\mathsf{W}}$ defined in Section 3 of \cite{PP3} upon choosing $\mathsf{W%
}$ as $n(n-k)^{-1}I_{n}$. Part A of Proposition 5.12 of \cite{PP3} then
delivers the claim (\ref{size-control_groupHet_uncorr}), again since $%
C^{\ast }<\infty $ as noted in that proposition. Again, $G(\mathfrak{M}_{0})$%
-invariance of $T_{uc}$ implies that $C^{\ast }$ and $\alpha ^{\ast }$ do
not depend on the choice of $\mu _{0}\in \mathfrak{M}_{0}$, and that $\alpha
^{\ast }$ as defined in the theorem coincides with $\alpha ^{\ast }$ as
defined in Proposition 5.12 of \cite{PP3}. Positivity of $\alpha ^{\ast }$
follows exactly as before making now use of Remark \ref{F-type}(ii). The
remaining claim in Part (a) is proved completely analogous as the
corresponding claim in Part (b).

We finally prove Part (c):\ The claims follow from Remark 5.10 and Lemma
5.16 in \cite{PP3} combined with Remark \ref{F-type} above; cf.~also
Appendix \ref{useful}. $\blacksquare $

\textbf{Proof of Proposition \ref{rem_necessity_2}:} (a) This follows from
Part 3 of Corollary 5.17 in \cite{PP2016}: As shown in the proof of
Proposition \ref{rem_necessity}(a) $T_{uc}$ satisfies the assumptions of
this corollary (with $\check{\beta}=\hat{\beta}$, $\check{\Omega}(y)=\hat{%
\sigma}^{2}(y)R\left( X^{\prime }X\right) ^{-1}R^{\prime }$, $N=\emptyset $,
and $N^{\ast }=\limfunc{span}(X)$). Set now $\mathcal{Z}=\func{span}%
(\{e_{i}(n):i\in (n_{j-1}^{+},n_{j}^{+}]\})$, where $j$ is such that $%
(n_{j-1}^{+},n_{j}^{+}]\cap I_{1}(\mathfrak{M}_{0}^{lin})\neq \emptyset $
and $\mathcal{Z}\subseteq \func{span}(X)$ hold. Note that $\mathcal{Z}$ is
not contained in $\mathfrak{M}_{0}^{lin}$ by construction. Observe that $%
\mathcal{Z}$ is a concentration space of $\mathfrak{C}_{(n_{1},\ldots
,n_{m})}$ in view of Remark \ref{rem_conc_spaces} in Appendix \ref{app_char}
(note that $\limfunc{card}((n_{j-1}^{+},n_{j}^{+}])<n$ must hold in view of $%
\mathcal{Z}\subseteq \func{span}(X)$ and $k<n$, while $0<\limfunc{card}%
((n_{j-1}^{+},n_{j}^{+}])$ is obvious). The nonnegative definiteness
assumption on $\check{\Omega}$ in Part 3 of Corollary 5.17 in \cite{PP2016}
is clearly satisfied. Obviously $\check{\Omega}(z)=0$ holds for every $z\in 
\mathcal{Z}$ since $\mathcal{Z}\subseteq \func{span}(X)$. It remains to
establish that $R\hat{\beta}(z)\neq 0$ holds $\lambda _{\mathcal{Z}}$%
-everywhere: Clearly, $R\hat{\beta}(z)=0$ for $z\in \mathcal{Z}$ occurs
precisely for $z\in \mathcal{Z}\cap \mathfrak{M}_{0}^{lin}$ since $\mathcal{Z%
}\subseteq \func{span}(X)$. But $\mathcal{Z}\cap \mathfrak{M}_{0}^{lin}$ is
a $\lambda _{\mathcal{Z}}$-null set in view of the fact that $\mathcal{Z}$
is not contained in $\mathfrak{M}_{0}^{lin}$ as noted before (and hence $%
\mathcal{Z}\cap \mathfrak{M}_{0}^{lin}$ is a proper linear subspace of $%
\mathcal{Z}$). Part 3 of Corollary 5.17 in \cite{PP2016} then proves the
claim for $C>0$. A fortiori it then also holds for all real $C$.

(b) This follows in the same way as Part (a) by applying Part 3 of Corollary
5.17 in \cite{PP2016} now to $T_{Het}$ (with $\check{\beta}=\hat{\beta}$, $%
\check{\Omega}=\hat{\Omega}_{Het}$, $N=\emptyset $, and $N^{\ast }=\mathsf{B}
$). $\blacksquare $

We note that Propositions \ref{rem_necessity} and \ref{rem_necessity_2}
could also be proved by making use of Theorem 3.1 in \cite{PP4}.

\begin{remark}
(i) Condition (\ref{non-incl_Het_uncorr}) ((\ref{non-incl_Het}),
respectively) in Theorem \ref{Hetero_Robust} can equivalently be written as $%
\func{span}(\{\pi _{\left( \mathfrak{M}_{0}^{lin}\right) ^{\bot
},i}\})\nsubseteqq \func{span}(X)$ ($\nsubseteqq \mathsf{B}$, respectively)
for every $i\in I_{1}(\mathfrak{M}_{0}^{lin})$ as discussed in the proof.
Since the spaces $\func{span}(\{\pi _{\left( \mathfrak{M}_{0}^{lin}\right)
^{\bot },i}\})$ are one-dimensional for $i\in I_{1}(\mathfrak{M}_{0}^{lin})$
and since $1<n-k+q=n-\dim (\mathfrak{M}_{0}^{lin})$, it follows that these
spaces are necessarily elements of $\mathbb{J}(\mathfrak{M}_{0}^{lin},%
\mathfrak{C}_{Het})$; in fact, they are precisely the minimal elements of $%
\mathbb{J}(\mathfrak{M}_{0}^{lin},\mathfrak{C}_{Het})$ w.r.t. the order
induced by inclusion.

(ii)\ Condition (\ref{non-incl_groupHet_uncorr}) ((\ref{non-incl_groupHet}),
respectively) in Theorem \ref{theorem_groupwise_hetero} can equivalently be
written as 
\begin{equation*}
\func{span}(\{\Pi _{\left( \mathfrak{M}_{0}^{lin}\right) ^{\bot
}}e_{i}(n):i\in (n_{j-1}^{+},n_{j}^{+}]\})\nsubseteqq \func{span}(X)\text{ (}%
\nsubseteqq \mathsf{B}\text{, respectively)}
\end{equation*}%
for every $j=1,\ldots ,m$ with $(n_{j-1}^{+},n_{j}^{+}]\cap I_{1}(\mathfrak{M%
}_{0}^{lin})\neq \emptyset $ as discussed in the proof. However, in this
more general case, it can happen that such a space appearing on the
l.h.s.~of the non-inclusion relation has a dimension not smaller than $%
n-\dim (\mathfrak{M}_{0}^{lin})$, and hence is not a member of $\mathbb{J}(%
\mathfrak{M}_{0}^{lin},\mathfrak{C}_{(n_{1},\ldots ,n_{m})})$. In light of
the general results in \cite{PP3} (e.g., Corollary 5.6) one may wonder if
requiring the non-inclusion condition in (\ref{non-incl_groupHet_uncorr}) (%
\ref{non-incl_groupHet}, respectively) for such spaces does not add an
unnecessary restriction. However, this is not so as this non-inclusion is
easily seen to be automatically satisfied for such spaces.\footnote{%
Note that any such space is necessarily equal to $\left( \mathfrak{M}%
_{0}^{lin}\right) ^{\bot }$. If now $\left( \mathfrak{M}_{0}^{lin}\right)
^{\bot }$ were contained in $\func{span}(X)$ ($\mathsf{B}$, respectively),
then $\mathbb{R}^{n}$ would also have to be contained in $\func{span}(X)$ ($%
\mathsf{B}$, respectively), since $\mathbb{R}^{n}$ can be written as the
direct sum of $\left( \mathfrak{M}_{0}^{lin}\right) ^{\bot }$ and $\mathfrak{%
M}_{0}^{lin}$ and since $\func{span}(X)$ ($\mathsf{B}$, respectively) are
invariant under addition of elements of $\mathfrak{M}_{0}^{lin}$. However, $%
\func{span}(X)$ is a proper subspace of $\mathbb{R}^{n}$ (since we always
assume $k<n$) and $\mathsf{B}$ is a finite union of proper linear subspaces
of $\mathbb{R}^{n}$ under Assumption \ref{R_and_X}. This gives a
contradiction.} Furthermore, the collection of all spaces of the form $\func{%
span}(\{\Pi _{\left( \mathfrak{M}_{0}^{lin}\right) ^{\bot }}e_{i}(n):i\in
(n_{j-1}^{+},n_{j}^{+}]\})$ for $j=1,\ldots ,m$, such that $%
(n_{j-1}^{+},n_{j}^{+}]\cap I_{1}(\mathfrak{M}_{0}^{lin})\neq \emptyset $
and such that the dimension of these spaces is smaller than $n-\dim (%
\mathfrak{M}_{0}^{lin})$ is precisely the collection of minimal elements of $%
\mathbb{J}(\mathfrak{M}_{0}^{lin},\mathfrak{C}_{(n_{1},\ldots ,n_{m})})$
w.r.t. the order induced by inclusion. [Note that $\mathbb{J}(\mathfrak{M}%
_{0}^{lin},\mathfrak{C}_{(n_{1},\ldots ,n_{m})})$ may be empty.]
\end{remark}

\begin{proposition}
\label{prop_k_pop}Suppose we are in the setting of Example \ref{ex_k_pop}
with $n_{j}\geq 2$ for all $j$. Then $T_{Het}$ is size controllable over $%
\mathfrak{C}_{Het}$, i.e., (\ref{size-control_Het}) holds for every $%
0<\alpha <1$.
\end{proposition}

\textbf{Proof: }Note that $\mathsf{B}$ is a subset of%
\begin{equation*}
\mathsf{S}:=\left\{ y\in \mathbb{R}^{n}:\hat{u}_{i}(y)=0\text{ for some }%
i=1,\ldots ,n\right\} ,
\end{equation*}%
and that $\mathsf{S}$ is a $\lambda _{\mathbb{R}^{n}}$-null set, as it is a
finite union of $\lambda _{\mathbb{R}^{n}}$-null sets (since $e_{i}(n)\notin 
\limfunc{span}(X)$ in view of $n_{j}\geq 2$ for all $j$). Also note that $%
S_{j}>0$ holds for $y\notin \mathsf{S}$. Now, for $y\notin \mathsf{S}$, by
the Sherman-Morrison formula, the inverse of $S_{1}\iota \iota ^{\prime }+%
\limfunc{diag}(S_{2},\ldots ,S_{k})$ equals%
\begin{equation*}
\limfunc{diag}(S_{2}^{-1},\ldots ,S_{k}^{-1})-\limfunc{diag}%
(S_{2}^{-1},\ldots ,S_{k}^{-1})\iota \iota ^{\prime }\limfunc{diag}%
(S_{2}^{-1},\ldots ,S_{k}^{-1})/\sum_{j=1}^{k}1/S_{j}.
\end{equation*}%
We may thus write 
\begin{equation}
T_{Het}(y)=\sum_{j=2}^{k}\frac{(\bar{y}_{(1)}-\bar{y}_{(j)})^{2}}{S_{j}}-%
\left[ \sum_{j=2}^{k}\frac{\bar{y}_{(1)}-\bar{y}_{(j)}}{S_{j}}\right]
^{2}/\sum_{j=1}^{k}1/S_{j}\text{ \ \ for every }y\notin \mathsf{S}.
\label{eqn:anova1}
\end{equation}

As noted in Remark \ref{obvious}, for any invertible $q\times q$-dimensional
matrix $A$, the test statistic $T_{Het}$ based on $R$ and the analogous test
statistic, but computed with $AR$ instead of $R$, coincide everywhere (note $%
r=0$). We apply this observation in the following way: fix $l\in \{2,\ldots
,k\}$, and choose $A$ with $l$-th column $(-1,\ldots ,-1)^{\prime }$, $l$-th
row$(0,\ldots 0,-1,0,\ldots ,0)$, and such that after deleting the $l$-th
column and the $l$-th row we obtain $I_{q-1}$. Then%
\begin{equation*}
AR=RP_{l},
\end{equation*}%
where $P_{l}$ is the $k\times k$ permutation matrix that interchanges the
first and $l$-th coordinate (and keeps all other coordinates fixed). By a
similar computation as the one that led to the expression in (\ref%
{eqn:anova1}), but now with $RP_{l}$ in place of $R$, we can now conclude
that for every $l\in \{1,\ldots ,k\}$ we have%
\begin{equation*}
T_{Het}(y)=\sum_{j=1,j\neq l}^{k}\frac{(\bar{y}_{(l)}-\bar{y}_{(j)})^{2}}{%
S_{j}}-\left[ \sum_{j=1,j\neq l}^{k}\frac{\bar{y}_{(l)}-\bar{y}_{(j)}}{S_{j}}%
\right] ^{2}/\sum_{j=1}^{k}1/S_{j}\text{ \ \ for every }y\notin \mathsf{S}.
\end{equation*}%
For $y\notin \mathsf{S}$ we may thus upper bound $T_{Het}(y)$ by $%
\sum_{j\neq l}^{k}(\bar{y}_{l}-\bar{y}_{j})^{2}/S_{j}$, and we are free to
choose $l$. Setting $l=l(y)\in \arg \min_{j=1,\ldots ,k}S_{j}$, the upper
bound for $T_{Het}(y)$ just derived, together with $S_{j}\geq
(S_{j}+S_{l})/2>0$, gives for $y\notin \mathsf{S}$ 
\begin{eqnarray*}
T_{Het}(y) &\leq &\sum_{j=1,j\neq l}^{k}\frac{(\bar{y}_{(l)}-\bar{y}%
_{(j)})^{2}}{S_{j}}\leq 2\sum_{j=1,j\neq l}^{k}\frac{(\bar{y}_{(l)}-\bar{y}%
_{(j)})^{2}}{S_{j}+S_{l}} \\
&\leq &2\sum_{i,j=1,i\neq j}^{k}\frac{(\bar{y}_{(i)}-\bar{y}_{(j)})^{2}}{%
S_{i}+S_{j}}=2\sum_{i,j=1,i\neq j}^{k}T_{i,j}(y),
\end{eqnarray*}%
where $T_{i,j}(y)=(\bar{y}_{(i)}-\bar{y}_{(j)})^{2}/(S_{i}+S_{j})$. Note
that the quantity to the far right does not depend on our particular choice
of $l$. For $y\in \mathsf{S}$, define $T_{i,j}$ by the same formula as long
as $S_{i}+S_{j}>0$ and as $T_{i,j}=0$ else. Since $\mathsf{S}$ is a $\lambda
_{\mathbb{R}^{n}}$-null set, we have for any $C>0$%
\begin{equation*}
\sup_{\mu _{0}\in \mathfrak{M}_{0}}\sup_{0<\sigma ^{2}<\infty }\sup_{\Sigma
\in \mathfrak{C}_{Het}}P_{\mu _{0},\sigma ^{2}\Sigma }(T_{Het}\geq C)\leq
\sum_{i,j=1,i\neq j}^{k}\sup_{\mu _{0}\in \mathfrak{M}_{0}}\sup_{0<\sigma
^{2}<\infty }\sup_{\Sigma \in \mathfrak{C}_{Het}}P_{\mu _{0},\sigma
^{2}\Sigma }(T_{i,j}\geq C/2(k^{2}-k)).
\end{equation*}%
Now observe that $T_{i,j}$ depends only on the coordinates of $y$
corresponding to groups $i$ and $j$ and furthermore coincides with the test
statistic of the form (\ref{T_het}) for a two sample mean comparison as
considered in Example \ref{ex_fish_behr} (with sample size being equal to $%
n_{i}+n_{j}$). A simple argument then shows that the terms in the sum on the
r.h.s of the preceding display can be rewritten as the sizes of the test
statistic (\ref{T_het}) as considered in Example \ref{ex_fish_behr} with
sample size now being given by $n_{i}+n_{j}$. Hence, all these terms can be
made arbitrarily small by choosing $C$ large enough by what has been
established in Example \ref{ex_fish_behr}. $\blacksquare $

We provide here a further example, where the sufficient condition of Part
(b) of Theorem \ref{Hetero_Robust} fails, but size control is possible.

\begin{example}
\label{ex_unreas}Suppose we are given $k\geq 2$ integers $n_{j}$ describing
group sizes satisfying $n_{1}\geq 2$ and $n_{j}\geq 1$ for $j\geq 2$. Sample
size is $n=\sum_{j=1}^{k}n_{j}$. Clearly $k<n$ is then satisfied. The
regressors $x_{ti}$ indicate group membership, i.e., they satisfy $x_{ti}=1$
for $\sum_{j=1}^{i-1}n_{j}<t\leq \sum_{j=1}^{i}n_{j}$ and $x_{ti}=0$
otherwise. The heteroskedasticity model is again given by $\mathfrak{C}%
_{Het} $. Let $R=(1,0,\ldots ,0)$, i.e., the coefficient of the first
regressor is subject to test. Then $I_{0}(\mathfrak{M}_{0}^{lin})=\{%
\sum_{l=1}^{j}n_{l}:n_{j}=1,\;j=2,\ldots ,k\}$. With regard to $T_{uc}$ we
immediately see that $e_{i}(n)\notin \func{span}(X)$ for $i\in I_{1}(%
\mathfrak{M}_{0}^{lin})$ holds, and thus the sufficient condition (\ref%
{non-incl_Het_uncorr}) for size control of $T_{uc}$ is satisfied. Turning to 
$T_{Het}$, observe that Assumption \ref{R_and_X} is satisfied as is easily
seen. Furthermore, it is not difficult to see that $\mathsf{B}=\left\{ y\in 
\mathbb{R}^{n}:y_{1}=\ldots =y_{n_{1}}\right\} $. Note that $\limfunc{span}%
(X)\subseteq \mathsf{B}$, but $\mathsf{B}\neq \limfunc{span}(X)$, except if $%
n_{j}=1$ for all $j\geq 2$ holds. In the latter case it is then easy to see
that $e_{i}(n)\notin \func{span}(X)=\mathsf{B}$ for every $i\in I_{1}(%
\mathfrak{M}_{0}^{lin})$ holds, and thus the sufficient condition (\ref%
{non-incl_Het}) for size control of $T_{Het}$ is satisfied. But if $n_{j}>1$
for some $j\geq 2$ holds, then for any index $i$ satisfying $%
\sum_{l=1}^{j-1}n_{l}<i\leq \sum_{l=1}^{j}n_{l}$ we have $i\in I_{1}(%
\mathfrak{M}_{0}^{lin})$ as well as $e_{i}(n)\in \mathsf{B}$. Consequently,
the sufficient condition (\ref{non-incl_Het}) for size control of $T_{Het}$
is \emph{not} satisfied and hence Theorem \ref{Hetero_Robust} does not
inform us about size controllability of $T_{Het}$ in this case. However, the
following argument shows that size control for $T_{Het}$ is possible also in
this case: The test statistic $T_{Het}$\ for the given problem coincides
with a corresponding test statistic (again of the form (\ref{T_uncorr}) for
an appropriate choice of $d_{i}$'s) in the \textquotedblleft
reduced\textquotedblright\ problem that one obtains by throwing away all
data points for $t>n_{1}$ and by also deleting all regressors from the
regression model but the first one. This leads one to the heteroskedastic
location model discussed in Example \ref{ex_loc} albeit with sample size
reduced to $n_{1}$. It is now not difficult to see that the size of $T_{Het}$
in the original formulation of the problem coincides with the size of the
corresponding test statistic in the \textquotedblleft
reduced\textquotedblright\ problem, which -- in light of the discussion in
Example \ref{ex_loc} -- shows that size control for $T_{Het}$ in the
original problem is possible also in the case where $n_{j}>1$ for some $%
j\geq 2$ holds. [If $n_{1}=1$ and if $n_{j}\geq 2$ for some $j$, condition (%
\ref{non-incl_Het_uncorr}) in Theorem \ref{Hetero_Robust}(a) is violated,
implying -- in view of Proposition \ref{rem_necessity} -- that the size of
the rejection region $\left\{ T_{uc}>C\right\} $ is $1$ for every choice of $%
C$; and that the test statistic $T_{Het}$ is identically zero (since
Assumption \ref{R_and_X} is violated and, in fact, $\hat{\Omega}_{Het}$ is
identically zero). The case where all $n_{j}$ are equal to $1$ even falls
outside of our framework since we always require $n>k$.]
\end{example}

\begin{remark}
\label{alter}Alternatively to the argument given in Example \ref{ex_unreas}
for the case where $n_{j}>1$ for some $j\geq 2$ holds, size controllability
of $T_{Het}$ can also be established by the following reasoning: Keep the
sample of size $n$, but replace the regressors $x_{\cdot i}$ for $2\leq
i\leq k$ by new regressors given by the standard basis vectors $e_{j}(n)$
for $j>n_{1}$ (the number of regressors now being $k^{\ast }=n-n_{1}+1<n$
and $R=(1,0,\ldots ,0)$ now being $1\times k^{\ast }$). Then one observes
that (i) this does not affect the test statistic, (ii) makes the set $%
\mathfrak{M}_{0}$ at most larger, and (ii) in the new model the sufficient
condition (\ref{non-incl_Het}) is now satisfied (as in the new model $%
n_{j}=1 $ holds for $j>n_{1}$). Hence, size control (even over the larger $%
\mathfrak{M}_{0}$) follows. A third possibility to establish the
size-controllability result is to observe that the test statistic $T_{Het}$
as well as the set $\mathsf{B}$ in the original model are -- additional to
being $G(\mathfrak{M}_{0})$-invariant -- also invariant w.r.t. addition of
the elements $e_{i}(n)$ for $i>n_{1}$ and then to appeal to a generalization
of Theorem \ref{Hetero_Robust} that exploits this additional invariance and
provides sufficient conditions for size control that can be seen to be
satisfied in the model considered in this example. Such a generalization of
Theorem \ref{Hetero_Robust}, which we refrain from stating, can be obtained
from the general size control results presented in \cite{PP3}.
\end{remark}

\begin{remark}
Example \ref{ex_unreas} is an instance of the following observation: Suppose 
$X$ is block-diagonal of rank $k$ with blocks $X_{1}$ and $X_{2}$ where $%
X_{i}$ is $n_{i}\times k_{i}$ with $n_{1}+n_{2}=n$ and $k_{1}+k_{2}=k$.
Assume $k_{1}<n_{1}$ (which entails $k<n$). Assume that the $q\times k$
restriction matrix $R$ is of rank $q$ and has the form $R=(R_{1}:0)$ with $%
R_{1}$ of dimension $q\times k_{1}$. The heteroskedasticity model is given
by $\mathfrak{C}_{Het}$. Then, using the same reasoning as in Example \ref%
{ex_unreas}, we see that the question of size control of $T_{Het}$ is
equivalent to the question of size control of the corresponding test
statistic in the \textquotedblleft reduced\textquotedblright\ problem where
one considers the regression model with regressor matrix equal to $X_{1}$
using only observations with $t\leq n_{1}$ (and as heteroskedasticity model
the analogue of $\mathfrak{C}_{Het}$ for sample size $n_{1}$). As Example %
\ref{ex_unreas} has shown, it is possible that the sufficient conditions for
size control of $T_{Het}$ are violated in the \textquotedblleft
original\textquotedblright\ problem, while at the same time the sufficient
conditions may be satisfied in the \textquotedblleft
reduced\textquotedblright\ problem. Alternatively, one can argue similarly
as in Remark \ref{alter}.
\end{remark}

\section{Appendix: Proofs for Section \protect\ref{sec_size_control_tilde}
and Appendix \protect\ref{sec_size_control_2_tilde}\label{app_C}}

\begin{lemma}
\label{lem:nullset}(a) Let $\mathcal{S}$ be a linear subspace of $\mathbb{R}%
^{n}$ and $\mu $ an element of $\mathbb{R}^{n}$ such that $\tilde{T}_{uc}$
restricted to $\mu +\mathcal{S}$ is not equal to a constant $\lambda _{\mu +%
\mathcal{S}}$-almost everywhere. Then $\lambda _{\mu +\mathcal{S}}(\tilde{T}%
_{uc}=C)=0$ holds for every $C\in \mathbb{R}$.

(b) $\lambda _{\mathbb{R}^{n}}(\tilde{T}_{uc}=C)=0$ holds for every $C\in 
\mathbb{R}$ .

(c) Let $\mathcal{S}$ be a linear subspace of $\mathbb{R}^{n}$ and $\mu $ an
element of $\mathbb{R}^{n}$ such that $\tilde{T}_{Het}$ restricted to $\mu +%
\mathcal{S}$ is not equal to a constant $\lambda _{\mu +\mathcal{S}}$-almost
everywhere. Then $\lambda _{\mu +\mathcal{S}}(\tilde{T}_{Het}=C)=0$ holds
for every $C\in \mathbb{R}$.

(d) Suppose Assumption \ref{R_and_X_tilde} holds and $\tilde{T}_{Het}$ is
not constant on $\mathbb{R}^{n}\backslash \mathsf{\tilde{B}}$. Then $\lambda
_{\mathbb{R}^{n}}(\tilde{T}_{Het}=C)=0$ holds for every $C\in \mathbb{R}$.
\end{lemma}

\textbf{Proof:} (a) Since $\tilde{T}_{uc}$ is constant on $\mathfrak{M}_{0}$
by definition, it follows that $\mu +\mathcal{S}\not\subseteq \mathfrak{M}%
_{0}$ must hold, and hence $\mathfrak{M}_{0}$ is a $\lambda _{\mu +\mathcal{S%
}}$-null set (cf. the argument in Remark 5.9(i) in \cite{PP3}).
Consequently, $\tilde{T}_{uc}$ restricted to $(\mu +\mathcal{S})\backslash 
\mathfrak{M}_{0}$ is not constant. Suppose now there exists a $C\in \mathbb{R%
}$ so that $\lambda _{\mu +\mathcal{S}}(\{y\in \mathbb{R}^{n}:\tilde{T}%
_{uc}(y)=C\})>0$. Then, since $\mathfrak{M}_{0}$ is a $\lambda _{\mu +%
\mathcal{S}}$-null set as just shown, it follows that even $\lambda _{\mu +%
\mathcal{S}}(\{y\in \mathbb{R}^{n}\backslash \mathfrak{M}_{0}:\tilde{T}%
_{uc}(y)=C\})>0$ must hold, which can be written as $\lambda _{\mu +\mathcal{%
S}}(\{y\in \mathbb{R}^{n}\backslash \mathfrak{M}_{0}:p(y)=0\})>0$, with the
multivariate polynomial $p$ given by $p(y)=(R\hat{\beta}\left( y\right)
-r)^{\prime }\left( R(X^{\prime }X)^{-1}R^{\prime }\right) ^{-1}(R\hat{\beta}%
\left( y\right) -r)-C\tilde{\sigma}^{2}(y)$. This implies that $p$
restricted to $\mu +\mathcal{S}$ vanishes on a set of positive $\lambda
_{\mu +\mathcal{S}}$-measure. Since $p$ restricted to $\mu +\mathcal{S}$ can
clearly be expressed as a polynomial in coordinates parameterizing the
affine space $\mu +\mathcal{S}$, it follows that $p$ vanishes identically on 
$\mu +\mathcal{S}$. But this implies that $\tilde{T}_{uc}$ restricted to $%
(\mu +\mathcal{S})\backslash \mathfrak{M}_{0}$ is constant equal to $C$, a
contradiction (as $\mathfrak{M}_{0}$ is a $\lambda _{\mu +\mathcal{S}}$-null
set).

(b) Follows from Part (a) upon choosing $\mathcal{S}=\mathbb{R}^{n}$, if we
can show that $\tilde{T}_{uc}$ is not $\lambda _{\mathbb{R}^{n}}$-almost
everywhere constant. Given that $\tilde{T}_{uc}$ is continuous on $\mathbb{R}%
^{n}\backslash \mathfrak{M}_{0}$ (the complement of a proper affine
subspace), it suffices to show that $\tilde{T}_{uc}$ is not constant on $%
\mathbb{R}^{n}\backslash \mathfrak{M}_{0}$. To this end consider first $%
y=X\beta $ with $R\beta -r\neq 0$ (such a $\beta $ obviously exists).
Observe that $\tilde{\sigma}^{2}(y)\neq 0$ as $y\notin \mathfrak{M}_{0}$ and
that $R\hat{\beta}(y)-r=R\beta -r\neq 0$. Hence, $\tilde{T}_{uc}(y)\neq 0$
for this choice of $y$. Next, choose $y=X\beta +w$, where $R\beta -r=0$
(such a $\beta $ obviously exists) and where $w\neq 0$ is orthogonal to $%
\limfunc{span}(X)$ (which is possible since $k<n$ is always maintained).
Then $\hat{\beta}(y)=\beta =\tilde{\beta}(y)$, implying $R\hat{\beta}%
(y)-r=R\beta -r=0$ and $\tilde{\sigma}^{2}(y)=w^{\prime }w/(n-(k-q))\neq 0$.
Note that $y\notin \mathfrak{M}_{0}$. It follows that $\tilde{T}_{uc}(y)=0$
holds for this choice of $y$. This establishes non-constancy of $\tilde{T}%
_{uc}$ on $\mathbb{R}^{n}\backslash \mathfrak{M}_{0}$.

(c) Completely analogous to the proof of Part (a) except that $\tilde{T}%
_{uc} $ and $\mathfrak{M}_{0}$ are replaced by $\tilde{T}_{Het}$ and $%
\mathsf{\tilde{B}}$, respectively, and that $p$ now takes the form $p(y)=(R%
\hat{\beta}\left( y\right) -r)^{\prime }\limfunc{adj}(\tilde{\Omega}%
_{Het}(y))(R\hat{\beta}\left( y\right) -r)-C\det (\tilde{\Omega}_{Het}(y))$,
where $\mathrm{\limfunc{adj}}(\cdot )$ denotes the adjoint of the square
matrix indicated, with the convention that the adjoint of a $1\times 1$
dimensional matrix equals one. [We note that under the assumptions for Part
(c) the set $\mathsf{\tilde{B}}$ cannot coincide with $\mathbb{R}^{n}$
(since otherwise $\tilde{T}_{Het}$ would be constant equal to zero), and
thus Assumption \ref{R_and_X_tilde} must hold.]

(d) Follows from Part (c) upon choosing $\mathcal{S}=\mathbb{R}^{n}$, if we
can show that $\tilde{T}_{Het}$ is not $\lambda _{\mathbb{R}^{n}}$-almost
everywhere constant. Given that $\tilde{T}_{Het}$ is continuous on $\mathbb{R%
}^{n}\backslash \mathsf{\tilde{B}}$ (the complement of a finite union of
proper affine subspaces by Lemma \ref{lem:B_tilde}), this follows from the
assumed non-constancy on $\mathbb{R}^{n}\backslash \mathsf{\tilde{B}}$. $%
\blacksquare $

\begin{remark}
\label{rem:ex_const}The additional assumption that $\tilde{T}_{Het}$ is not
constant on $\mathbb{R}^{n}\backslash \mathsf{\tilde{B}}$ in Part (d) of the
preceding lemma can not be dropped as can be seen from the following
example:\ Consider the case where $k=q=1$, $R=1$, $r=0$, the regressor is
given by $e_{1}(n)$, and the constants $\tilde{d}_{i}$ satisfy $\tilde{d}%
_{i}=1$ for all $i$. Then $\mathfrak{M}_{0}=\mathfrak{M}_{0}^{lin}=\{0\}$,
Assumption \ref{R_and_X_tilde} is satisfied, and $\mathsf{\tilde{B}}=%
\limfunc{span}(e_{1}(n))^{\bot }$. Furthermore, $\tilde{T}_{Het}(y)=1$ for
every $y\in \mathbb{R}^{n}\backslash \mathsf{\tilde{B}}$. As a point of
interest we note that $\tilde{T}_{Het}$ is trivially size controllable for
every $0<\alpha <1$, but that the condition (\ref{non-incl_Het_tilde}) for
size controllability is violated since $e_{j}(n)\in \mathsf{\tilde{B}}$ for $%
j>1$. [Of course, neither a smallest size-controlling critical value exists
(when considering rejection regions of the form $\{\tilde{T}_{Het}\geq C\}$)
nor can exact size controllability be achieved for $0<\alpha <1$.] An
extension of this example to the case $q=k>1$ is discussed in the proof of
Remark \ref{q=k} given further below.
\end{remark}

\begin{lemma}
\label{lem_indep_r}Let $C$ be a given critical value. Then the rejection
probabilities $P_{\mu _{0},\sigma ^{2}\Sigma }(\tilde{T}_{uc}\geq C)$ as
well as $P_{\mu _{0},\sigma ^{2}\Sigma }(\tilde{T}_{Het}\geq C)$ for $\mu
_{0}\in \mathfrak{M}_{0}$, $\sigma ^{2}\in (0,\infty )$, $\Sigma \in 
\mathfrak{C}_{Het}$, do not depend on $r$. [It is understood here that the
constants $\tilde{d}_{i}$ appearing in the definition of $\tilde{T}_{Het}$
have been chosen independently of the value of $r$.]
\end{lemma}

\textbf{Proof:} Fix $\mu _{0}\in \mathfrak{M}_{0}$. Observe that $\tilde{T}%
_{Het}(y)=\tilde{T}_{Het}^{0}(y-\mu _{0})$, where%
\begin{equation*}
\tilde{T}_{Het}^{0}\left( z\right) =\left\{ 
\begin{array}{cc}
(R\hat{\beta}\left( z\right) )^{\prime }\left( \tilde{\Omega}%
_{Het}^{0}(z)\right) ^{-1}(R\hat{\beta}\left( z\right) ) & \text{if }%
\limfunc{rank}\tilde{B}^{0}\left( z\right) =q, \\ 
0 & \text{if }\limfunc{rank}\tilde{B}^{0}\left( z\right) <q,%
\end{array}%
\right.
\end{equation*}%
where 
\begin{equation*}
\tilde{\Omega}_{Het}^{0}(z)=R(X^{\prime }X)^{-1}X^{\prime }\limfunc{diag}%
\left( \tilde{d}_{1}(\tilde{u}_{1}^{0}\left( z\right) )^{2},\ldots ,\tilde{d}%
_{n}(\tilde{u}_{n}^{0}\left( z\right) )^{2}\right) X(X^{\prime
}X)^{-1}R^{\prime },
\end{equation*}%
where $\tilde{u}^{0}\left( z\right) =\Pi _{(\mathfrak{M}_{0}^{lin})^{\bot
}}z $, and where $\tilde{B}^{0}\left( z\right) =R(X^{\prime
}X)^{-1}X^{\prime }\limfunc{diag}(e_{1}^{\prime }(n)\Pi _{(\mathfrak{M}%
_{0}^{lin})^{\bot }}(z),\ldots ,e_{n}^{\prime }(n)\Pi _{(\mathfrak{M}%
_{0}^{lin})^{\bot }}(z))$. Here we have made use of (\ref{B_tilde_matrix})
and the fact that $\tilde{u}\left( y\right) =\tilde{u}^{0}(y-\mu _{0})$. Now%
\begin{equation*}
P_{\mu _{0},\sigma ^{2}\Sigma }(\tilde{T}_{Het}(y)\geq C)=P_{\mu _{0},\sigma
^{2}\Sigma }(\tilde{T}_{Het}^{0}(y-\mu _{0})\geq C)=P_{0,\sigma ^{2}\Sigma }(%
\tilde{T}_{Het}^{0}(z)\geq C)
\end{equation*}%
and the far right-hand side does not depend on $r$ as $\tilde{T}_{Het}^{0}$
does not depend on $r$. The proof for $\tilde{T}_{uc}$ is completely
analogous, noting that $\tilde{T}_{uc}(y)=\tilde{T}_{uc}^{0}(y-\mu _{0})$,
where%
\begin{equation*}
\tilde{T}_{uc}^{0}\left( z\right) =\left\{ 
\begin{array}{cc}
(R\hat{\beta}\left( z\right) )^{\prime }\left( (\tilde{\sigma}%
^{0}(z))^{2}R(X^{\prime }X)^{-1}R^{\prime }\right) ^{-1}(R\hat{\beta}\left(
z\right) ) & \text{if }z\notin \mathfrak{M}_{0}^{lin}, \\ 
0 & \text{if }z\in \mathfrak{M}_{0}^{lin},%
\end{array}%
\right.
\end{equation*}%
and where $(\tilde{\sigma}^{0}(z))^{2}=(\tilde{u}^{0}\left( z\right)
)^{\prime }\tilde{u}^{0}\left( z\right) /(n-(k-q))$. $\blacksquare $

\textbf{Proof of Theorem \ref{Hetero_Robust_tilde}:} We first prove Part
(b). We apply Part (b) of Theorem \ref{theorem_groupwise_hetero_tilde} with $%
n_{j}=1$ for $j=1,\ldots ,n=m$ observing that then $\mathfrak{C}%
_{(n_{1},\ldots ,n_{m})}=\mathfrak{C}_{Het}$ and that condition (\ref%
{non-incl_groupHet_tilde}) reduces to (\ref{non-incl_Het_tilde}) (exploiting
that $\mathsf{\tilde{B}}-\mu _{0}$ is a finite union of proper linear
subspaces as discussed in Lemma \ref{lem:B_tilde}). This establishes (\ref%
{size-control_Het_tilde}). The final claim in Part (b) of the theorem
follows from Part (b) of Theorem \ref{theorem_groupwise_hetero_tilde}, if we
can show that $C^{\ast }$ given there can be written as claimed in Theorem %
\ref{Hetero_Robust_tilde}: To this end we proceed as follows:\footnote{%
Alternatively, one could base a proof on Lemma C.1 in \cite{PP4}.} Choose an
element $\mu _{0}$ of $\mathfrak{M}_{0}$. Observe that $I_{1}(\mathfrak{M}%
_{0}^{lin})\neq \emptyset $ (since $\dim (\mathfrak{M}_{0}^{lin})=k-q<n$),
and that for every $i\in I_{1}(\mathfrak{M}_{0}^{lin})$ the linear space $%
\mathcal{S}_{i}=\func{span}(\Pi _{\left( \mathfrak{M}_{0}^{lin}\right)
^{\bot }}e_{i}(n))$ is $1$-dimensional (since $\mathcal{S}_{i}=\{0\}$ is
impossible in view of $i\in I_{1}(\mathfrak{M}_{0}^{lin})$), and belongs to $%
\mathbb{J}(\mathfrak{M}_{0}^{lin},\mathfrak{C}_{Het})$ (since $n-k+q>1=\dim (%
\mathcal{S}_{i})$ holds) in view of Proposition \ref{characterization} in
Section \ref{app_char}. Since $\tilde{T}_{Het}$ is $G(\mathfrak{M}_{0})$%
-invariant (Remark \ref{rem:tildeGM0}), it follows that $\tilde{T}_{Het}$ is
constant on $(\mu _{0}+\mathcal{S}_{i})\backslash \left\{ \mu _{0}\right\} $%
, cf. the beginning of the proof of Lemma 5.11 in \cite{PP3}. Hence, $%
\mathcal{S}_{i}$ belongs to $\mathbb{H}$ (defined in Lemma 5.11 in \cite{PP3}%
) and consequently for $C^{\ast }$ as defined in that lemma%
\begin{equation*}
C^{\ast }\geq \max \left\{ \tilde{T}_{Het}(\mu _{0}+\Pi _{\left( \mathfrak{M}%
_{0}^{lin}\right) ^{\bot }}e_{i}(n)):i\in I_{1}(\mathfrak{M}%
_{0}^{lin})\right\}
\end{equation*}%
must hold. To prove the opposite inequality, let $\mathcal{S}$ be an
arbitrary element of $\mathbb{H}$, i.e., $\mathcal{S}\in \mathbb{J}(%
\mathfrak{M}_{0}^{lin},\mathfrak{C}_{Het})$ and $\tilde{T}_{Het}$ is $%
\lambda _{\mu _{0}+\mathcal{S}}$-almost everywhere equal to a constant $C(%
\mathcal{S})$, say. Then Proposition \ref{characterization} in Section \ref%
{app_char} shows that $\mathcal{S}_{i}\subseteq \mathcal{S}$ holds for some $%
i\in I_{1}(\mathfrak{M}_{0}^{lin})$. Because of Condition (\ref%
{non-incl_Het_tilde}) we have $\mu _{0}+\mathcal{S}_{i}\nsubseteqq \mathsf{%
\tilde{B}}$ since $\Pi _{\left( \mathfrak{M}_{0}^{lin}\right) ^{\bot
}}e_{i}(n)$ and $e_{i}(n)$ differ only by an element of $\mathfrak{M}%
_{0}^{lin}$ and since $\mathsf{\tilde{B}}+\mathfrak{M}_{0}^{lin}=\mathsf{%
\tilde{B}}$. We thus can find $s\in \mathcal{S}_{i}$ such that $\mu
_{0}+s\notin \mathsf{\tilde{B}}$. Note that $s\neq 0$ must hold, since $\mu
_{0}\in \mathfrak{M}_{0}\subseteq \mathsf{\tilde{B}}$ (see Lemma \ref%
{lem:B_tilde}). In particular, $\tilde{T}_{Het}$ is continuous at $\mu
_{0}+s $, since $\mu _{0}+s\notin \mathsf{\tilde{B}}$. Now, for every open
ball $A_{\varepsilon }$ in $\mathbb{R}^{n}$ with center $s$ and radius $%
\varepsilon >0$ we can find an element $a_{\varepsilon }\in A_{\varepsilon
}\cap \mathcal{S}$ such that $\tilde{T}_{Het}(\mu _{0}+a_{\varepsilon })=C(%
\mathcal{S})$. Since $a_{\varepsilon }\rightarrow s$ for $\varepsilon
\rightarrow 0$, it follows that $C(\mathcal{S})=\tilde{T}_{Het}(\mu _{0}+s)$%
. Since $s\neq 0$ and since $\tilde{T}_{Het}$ is constant on $(\mu _{0}+%
\mathcal{S}_{i})\backslash \left\{ \mu _{0}\right\} $ as shown before, we
can conclude that $C(\mathcal{S})=\tilde{T}_{Het}(\mu _{0}+s)=\tilde{T}%
_{Het}(\mu _{0}+\Pi _{\left( \mathfrak{M}_{0}^{lin}\right) ^{\bot
}}e_{i}(n)) $, where we recall that $\Pi _{\left( \mathfrak{M}%
_{0}^{lin}\right) ^{\bot }}e_{i}(n)\neq 0$. But this now implies%
\begin{equation*}
C^{\ast }=\max \left\{ \tilde{T}_{Het}(\mu _{0}+\Pi _{\left( \mathfrak{M}%
_{0}^{lin}\right) ^{\bot }}e_{i}(n)):i\in I_{1}(\mathfrak{M}%
_{0}^{lin})\right\} .
\end{equation*}%
Using $G(\mathfrak{M}_{0})$-invariance of $\tilde{T}_{Het}$ we conclude that%
\begin{equation*}
C^{\ast }=\max \left\{ \tilde{T}_{Het}(\mu _{0}+e_{i}(n)):i\in I_{1}(%
\mathfrak{M}_{0}^{lin})\right\} .
\end{equation*}

We next prove Part (a): Apply Part (a) of Theorem \ref%
{theorem_groupwise_hetero_tilde} with $n_{j}=1$ for $j=1,\ldots ,n=m$,
observing that then $\mathfrak{C}_{(n_{1},\ldots ,n_{m})}=\mathfrak{C}_{Het}$%
. This establishes (\ref{size-control_Het_uncorr_tilde}).\footnote{%
This argument is actually superfluous since $\tilde{T}_{uc}$ is bounded as
noted in Section \ref{sec:trivial_uc_tilde}.} The final claim in Part (a) of
the theorem follows similarly as the corresponding claim of Part (b) upon
replacing the set $\mathsf{\tilde{B}}$ by $\mathfrak{M}_{0}$ in the
argument, by noting that $\tilde{T}_{uc}$ is $G(\mathfrak{M}_{0})$%
-invariant, and that $\mu _{0}+\mathcal{S}_{i}\nsubseteqq \mathfrak{M}_{0}$
holds because of $i\in I_{1}(\mathfrak{M}_{0}^{lin})$.

Part (c) follows from Part (c) of Theorem \ref%
{theorem_groupwise_hetero_tilde} upon setting $n_{j}=1$ for $j=1,\ldots ,n=m$
(and upon noting that then the conditions in Theorem \ref%
{theorem_groupwise_hetero_tilde} reduce to the conditions of the present
theorem). $\blacksquare $

\textbf{Proof of Proposition \ref{rem_C*_tilde}:} Follows from Part A.1 of
Proposition 5.12 of \cite{PP3} and the sentence following this proposition.
Note that the assumptions of this proposition have been verified in the
proof of Theorem \ref{Hetero_Robust_tilde} (see also the proof of Theorem %
\ref{theorem_groupwise_hetero_tilde}, on which the proof of Theorem \ref%
{Hetero_Robust_tilde} is based), where it is also shown that the quantity $%
C^{\ast }$ used in Proposition 5.12 of \cite{PP3} coincides with $C^{\ast }$
defined in Theorem \ref{Hetero_Robust_tilde}. $\blacksquare $

\textbf{Proof of Remark \ref{q=k}}: (i) From the definition of $\tilde{B}(y)$
and since here $\tilde{u}(y)=y-\mu _{0}$ we obtain for $i=1,\ldots ,n$%
\begin{equation*}
\tilde{B}(\mu _{0}+e_{i}(n))=R(X^{\prime }X)^{-1}\left( 0,\ldots
,0,x_{i\cdot }^{\prime },0,\ldots ,0\right)
\end{equation*}%
where $x_{i\cdot }^{\prime }$ appears in the $i$-th position (recall that $%
x_{i\cdot }^{\prime }$ is the $i$-th column of $X^{\prime }$). But then $%
\limfunc{rank}(\tilde{B}(\mu _{0}+e_{i}(n)))\leq 1<q$, implying that $\mu
_{0}+e_{i}(n)\in \mathsf{\tilde{B}}$. [In case $q=k=1$, $\limfunc{rank}(%
\tilde{B}(\mu _{0}+e_{i}(n)))=1=q$ for every $i=1,\ldots ,n$ whenever the
matrix $X$ has no zero entry. This then implies that (\ref%
{non-incl_Het_tilde}) is satisfied. However, if $X$ contains a zero at the $%
m $-th position, say, then $\tilde{B}(\mu _{0}+e_{m}(n))=0<1=q$, implying
that $\mu _{0}+e_{m}(n)\in \mathsf{\tilde{B}}$, thus leading to violation of
(\ref{non-incl_Het_tilde}) as $I_{1}(\mathfrak{M}_{0}^{lin})=\{1,\ldots ,n\}$%
.]

(ii) Define $\beta _{0}=(X^{\prime }X)^{-1}X^{\prime }\mu _{0}$ and note
that $R\beta _{0}=r$ holds. Observing that $R$ is nonsingular, that $\tilde{d%
}_{i}>0$ for $1,\ldots ,n$, and that $\tilde{u}(y)=y-\mu _{0}$, we obtain
for $y\notin \mathsf{\tilde{B}}$%
\begin{eqnarray*}
\tilde{T}_{Het}(y) &=&\left( y-\mu _{0}\right) ^{\prime }X\left[ X^{\prime }%
\limfunc{diag}\left( \tilde{d}_{1}\tilde{u}_{1}^{2}(y),\ldots ,\tilde{d}_{n}%
\tilde{u}_{n}^{2}(y)\right) X\right] ^{-1}X^{\prime }\left( y-\mu _{0}\right)
\\
&=&\tilde{u}(y)^{\prime }X\left[ X^{\prime }\limfunc{diag}\left( \tilde{d}%
_{1}\tilde{u}_{1}^{2}(y),\ldots ,\tilde{d}_{n}\tilde{u}_{n}^{2}(y)\right) X%
\right] ^{-1}X^{\prime }\tilde{u}(y) \\
&\leq &\left( \min_{1\leq i\leq n}\tilde{d}_{i}\right) ^{-1}e^{\prime }A(y)%
\left[ A^{\prime }(y)A(y)\right] ^{-1}A^{\prime }(y)e\leq n\left(
\min_{1\leq i\leq n}\tilde{d}_{i}\right) ^{-1}
\end{eqnarray*}%
where $e=(1,\ldots ,1)^{\prime }$ and $A(y)=\limfunc{diag}(\tilde{u}%
_{1}(y),\ldots ,\tilde{u}_{n}(y))X$. Note that $A^{\prime }(y)A(y)$ is
nonsingular for $y\notin \mathsf{\tilde{B}}$ and that the matrix in the
quadratic form is a projection matrix. For $y\in \mathsf{\tilde{B}}$ we have 
$\tilde{T}_{Het}(y)=0$. Hence, $\tilde{T}_{Het}(y)$ is bounded from above,
and is trivially bounded from below as $\tilde{T}_{Het}(y)\geq 0$ for every $%
y\in \mathbb{R}^{n}$.

(iii) In the following examples we always set $\mu _{0}=0$ (i.e., $r=0$) for
the sake of simplicity. Remark \ref{rem:ex_const} provides an example where $%
\tilde{T}_{Het}$ is constant on $\mathbb{R}^{n}\backslash \mathsf{\tilde{B}}$%
. This example has $q=k=1$. It can be easily extended to the case $q=k\geq 2$
by considering a design matrix $X$, the columns of which are given by the
first $k$ standard basis vectors, by setting $R=I_{q}$, and $\tilde{d}_{i}=1$
for every $i=1,\ldots ,n$. Then $\tilde{T}_{Het}(y)=k$ for every $y\in 
\mathbb{R}^{n}\backslash \mathsf{\tilde{B}}=\{y\in \mathbb{R}^{n}:y_{1}\neq
0,\ldots ,y_{k}\neq 0\}$. An example where $\tilde{T}_{Het}$ is not constant
on $\mathbb{R}^{n}\backslash \mathsf{\tilde{B}}$ is in case $q=k=1$ given by
the location model: Here $\tilde{T}_{Het}(y)=(\tsum_{t=1}^{n}y_{t})^{2}/%
\tsum_{t=1}^{n}y_{t}^{2}$ for every $y\in \mathbb{R}^{n}\backslash \mathsf{%
\tilde{B}}=\{y\in \mathbb{R}^{n}:y\neq 0\}$, which obviously is not constant
(as $n>k=1)$.\footnote{%
In this example condition (\ref{non-incl_Het_tilde}) is satisfied as $%
e_{i}(n)\notin \mathsf{\tilde{B}}$ for every $i=1,\ldots ,n$. To arrive at
an example where again $\tilde{T}_{Het}$ is not constant on $\mathbb{R}%
^{n}\backslash \mathsf{\tilde{B}}$ but where condition (\ref%
{non-incl_Het_tilde}) is not satisfied, consider the case where $X=(1,\ldots
,1,0)^{\prime }$ with $n\geq 2$. Observe that then $I_{1}(\mathfrak{M}%
_{0}^{lin})=\{1,\ldots ,n\}$.} This example can again be extended to the
case $q=k>1$ as follows: Let $R=I_{q}$ and let $X$ be the design matrix
where each of the columns correspond to a dummy variable describing
membership in one of $k$ disjoint groups $G_{j}$, each group of the same
cardinality $n_{1}$ with $n_{1}>1$. Consequently, $n=kn_{1}$. W.l.o.g., we
may assume that the elements $G_{1}$ have the lowest indices, followed by
the elements of $G_{2}$, and so on. It is then easy to see that%
\begin{equation}
\tilde{T}_{Het}(y)=\tsum_{j=1}^{k}\left[ (\tsum_{t\in
G_{j}}y_{t})^{2}/\tsum_{t\in G_{j}}y_{t}^{2}\right]  \label{star}
\end{equation}%
for $y\in \mathbb{R}^{n}\backslash \mathsf{\tilde{B}}=\tbigcap_{j=1}^{k}\{y%
\in \mathbb{R}^{n}:y_{t}\neq 0$ for at least one $t\in G_{j}\}$. Obviously,
the expression in (\ref{star}) is not constant: Choosing $y=e$ gives the
value $kn_{1}=n$, whereas choosing $y$ such that $%
y_{1}=y_{n_{1}+1}=y_{2n_{1}+1}=\ldots =y_{(k-1)n_{1}+1}=1$ with all the
other coordinates being zero gives a value of $k<n=kn_{1}$ since $n_{1}>1$.

\textbf{Proof of Theorem \ref{Hetero_Robust_tilde_add}:} From the definition
of $C^{\ast }$ we see that $C^{\ast }$ is nonnegative and finite. Let $C$ be
arbitrary but satisfying $C^{\ast }<C<\sup_{y\in \mathbb{R}^{n}}\tilde{T}%
_{Het}(y)$. We can then choose $y_{0}\in \mathbb{R}^{n}$ with $\tilde{T}%
_{Het}(y_{0})>C>0$. In view of the definition of $\tilde{T}_{Het}$ it
follows that $y_{0}\notin \mathsf{\tilde{B}}$, and hence $\tilde{T}_{Het}$
is continuous at $y_{0}$. We can thus find an open neighborhood $U(y_{0})$
of $y_{0}$ in $\mathbb{R}^{n}$ such that $\tilde{T}_{Het}$ is larger than $C$
on $U(y_{0})$. In particular, $P_{\mu _{0},\Sigma }(\tilde{T}_{Het}\geq
C)\geq P_{\mu _{0},\Sigma }(U(y_{0}))>0$ for every $\mu _{0}\in \mathfrak{M}%
_{0}$ and every $\Sigma \in \mathfrak{C}_{Het}$. This establishes $\alpha
^{\ast }>0$. Choose $\delta >0$ such that $\delta \leq \alpha $ and $\delta
<\alpha ^{\ast }$. Then the size of the rejection region $\{\tilde{T}%
_{Het}\geq C_{\Diamond }(\delta )\}$ is exactly equal to $\delta $ by Parts
(b) and (c) of Theorem \ref{Hetero_Robust_tilde}. Consequently, $\{\tilde{T}%
_{Het}\geq C_{\Diamond }(\delta )\}$ is not a $\lambda _{\mathbb{R}^{n}}$%
-null set. By construction, $C_{\Diamond }(\alpha )\leq C_{\Diamond }(\delta
)$ holds, and hence $\{\tilde{T}_{Het}\geq C_{\Diamond }(\alpha )\}$
contains $\{\tilde{T}_{Het}\geq C_{\Diamond }(\delta )\}$, which completes
the proof. $\blacksquare $

\textbf{Proof of Theorem \ref{theorem_groupwise_hetero_tilde}:} We first
prove Part (b). We wish to apply Part A of Proposition 5.12 of \cite{PP3}
with $\mathfrak{C}=\mathfrak{C}_{(n_{1},\ldots ,n_{m})}$, $T=\tilde{T}_{Het}$%
, $\mathcal{L}=\mathfrak{M}_{0}^{lin}$, and $\mathcal{V}=\{0\}$. First, note
that $\dim (\mathfrak{M}_{0}^{lin})=k-q<n$. Second, under Assumption \ref%
{R_and_X_tilde}, $\tilde{T}_{Het}$ is clearly Borel-measurable and is
continuous on the complement of $\mathsf{\tilde{B}}$, where $\mathsf{\tilde{B%
}}$ is a closed $\lambda _{\mathbb{R}^{n}}$-null set (see Lemma \ref%
{lem:B_tilde} and the paragraph following this lemma). Because of Remark \ref%
{rem:tildeGM0}, we hence see that the general assumptions on $T=\tilde{T}%
_{Het}$, on $N^{\dag }=\mathsf{\tilde{B}}$, on $\mathcal{L}=\mathfrak{M}%
_{0}^{lin}$, as well as on $\mathcal{V}=\{0\}$ in Proposition 5.12 of \cite%
{PP3} are satisfied. Next, observe that the validity of condition (\ref%
{non-incl_groupHet_tilde}) clearly does not depend on the choice of $\mu
_{0}\in \mathfrak{M}_{0}$ since $\mathsf{\tilde{B}}+\mathfrak{M}_{0}^{lin}=%
\mathsf{\tilde{B}}$ as shown in Lemma \ref{lem:B_tilde}. For the same reason
condition (\ref{non-incl_groupHet_tilde}) can equivalently be written as

\begin{equation*}
\mu _{0}+\func{span}\left( \left\{ \Pi _{\left( \mathfrak{M}%
_{0}^{lin}\right) ^{\bot }}e_{i}(n):i\in (n_{j-1}^{+},n_{j}^{+}]\right\}
\right) \nsubseteqq \mathsf{\tilde{B}}
\end{equation*}%
for every $j=1,\ldots ,m$, such that $(n_{j-1}^{+},n_{j}^{+}]\cap I_{1}(%
\mathfrak{M}_{0}^{lin})\neq \emptyset $, since $\Pi _{\left( \mathfrak{M}%
_{0}^{lin}\right) ^{\bot }}e_{i}(n)$ and $e_{i}(n)$ differ only by an
element of $\mathfrak{M}_{0}^{lin}$. In view of Proposition \ref%
{characterization_2} in Appendix \ref{app_char}, this implies that $\mu _{0}+%
\mathcal{S}$ for any $\mathcal{S}\in \mathbb{J}(\mathfrak{M}_{0}^{lin},%
\mathfrak{C}_{(n_{1},\ldots ,n_{m})})$ is not contained in $\mathsf{\tilde{B}%
}$, and thus not in $N^{\dag }$. Since $\mu _{0}+\mathcal{S}$ is an affine
space and $N^{\dag }=\mathsf{\tilde{B}}$ is a finite union of proper affine
spaces under Assumption \ref{R_and_X_tilde} as discussed in Lemma \ref%
{lem:B_tilde}, we may conclude (cf. Corollary 5.6 in \cite{PP3} and its
proof) that $\lambda _{\mu _{0}+\mathcal{S}}(N^{\dag })=0$ for every $%
\mathcal{S}\in \mathbb{J}(\mathfrak{M}_{0}^{lin},\mathfrak{C}_{(n_{1},\ldots
,n_{m})})$ and every $\mu _{0}\in \mathfrak{M}_{0}$. This completes the
verification of the assumptions of Proposition 5.12 in \cite{PP3} that are
not specific to Part A (or Part B) of this proposition. We next verify the
assumptions specific to Part A of this proposition: Assumption (a) is
satisfied (even for every $C\in \mathbb{R}$) as a consequence of Part (d) of
Lemma \ref{lem:nullset} (note that we have assumed that $\tilde{T}_{Het}$ is
not constant on $\mathbb{R}^{n}\backslash \mathsf{\tilde{B}}$). And
Assumption (b) in Part A follows from Part (c) of Lemma \ref{lem:nullset}.
Part A of Proposition 5.12 of \cite{PP3} now immediately delivers claim (\ref%
{size-control_groupHet_tilde}), since $C^{\ast }<\infty $ as noted in that
proposition. That $C^{\ast }$ and $\alpha ^{\ast }$ do not depend on the
choice of $\mu _{0}\in \mathfrak{M}_{0}$ is an immediate consequence of $G(%
\mathfrak{M}_{0})$-invariance of $\tilde{T}_{Het}$. Also note that $\alpha
^{\ast }$ as defined in the theorem coincides with $\alpha ^{\ast }$ as
defined in Proposition 5.12 of \cite{PP3} in view of $G(\mathfrak{M}_{0})$%
-invariance of $\tilde{T}_{Het}$. In case $\alpha <\alpha ^{\ast }$, the
remaining claim in Part (b) of the theorem, namely that equality can be
achieved in (\ref{size-control_groupHet}), follows from the definition of $%
C^{\ast }$ in Lemma 5.11 of \cite{PP3} and from Part A.2 of Proposition 5.12
of \cite{PP3} (and the observation immediately following that proposition
allowing one to drop the suprema w.r.t. $\mu _{0}$ and $\sigma ^{2}$, and to
set $\sigma ^{2}=1$); in case $\alpha =\alpha ^{\ast }<1$, it follows from
Remarks 5.13(i),(ii) in \cite{PP3} using Part (d) of Lemma \ref{lem:nullset}%
. [In case $\alpha ^{\ast }=0$, there is nothing to prove.]

The proof of Part (a) proceeds similarly, but with some differences: Noting
that $\tilde{T}_{uc}$ is clearly Borel-measurable and is continuous on the
complement of $\mathfrak{M}_{0}$, where $\mathfrak{M}_{0}$ is a closed $%
\lambda _{\mathbb{R}^{n}}$-null set, and using Remark \ref{rem:tildeGM0}, we
now see that the general assumptions on $T=\tilde{T}_{uc}$, on $N^{\dag }=%
\mathfrak{M}_{0}$, on $\mathcal{L}=\mathfrak{M}_{0}^{lin}$, as well as on $%
\mathcal{V}=\{0\}$ in Proposition 5.12 of \cite{PP3} are satisfied (again
with $\mathfrak{C}=\mathfrak{C}_{(n_{1},\ldots ,n_{m})}$). Let now $\mathcal{%
S}\in \mathbb{J}(\mathfrak{M}_{0}^{lin},\mathfrak{C}_{(n_{1},\ldots
,n_{m})}) $. In view of Proposition \ref{characterization_2} in Appendix \ref%
{app_char}, $\mathcal{S}$ must then contain an element of the form $\Pi
_{\left( \mathfrak{M}_{0}^{lin}\right) ^{\bot }}e_{i}(n)$ for some $i\in
I_{1}(\mathfrak{M}_{0}^{lin})$. Observe that $\Pi _{\left( \mathfrak{M}%
_{0}^{lin}\right) ^{\bot }}e_{i}(n)\notin \mathfrak{M}_{0}^{lin}$ must hold,
since otherwise we would have $e_{i}(n)\in \mathfrak{M}_{0}^{lin}$,
contradicting $i\in I_{1}(\mathfrak{M}_{0}^{lin})$. It follows that $%
\mathcal{S}\nsubseteqq \mathfrak{M}_{0}^{lin}$, and thus $\mu _{0}+\mathcal{S%
}\nsubseteqq \mathfrak{M}_{0}$ for every $\mu _{0}\in \mathfrak{M}_{0}$.
Since $\mu _{0}+\mathcal{S}$ is an affine space and $N^{\dag }=\mathfrak{M}%
_{0}$ is a proper affine space we may conclude (cf. Corollary 5.6 in \cite%
{PP3} and its proof) that $\lambda _{\mu _{0}+\mathcal{S}}(N^{\dag })=0$ for
every $\mathcal{S}\in \mathbb{J}(\mathfrak{M}_{0}^{lin},\mathfrak{C}%
_{(n_{1},\ldots ,n_{m})})$ and every $\mu _{0}\in \mathfrak{M}_{0}$. We have
thus now completed the verification of the assumptions of Proposition 5.12
of \cite{PP3} that are not specific to Part A (or Part B) of this
proposition. We next verify the assumptions specific to Part A of this
proposition: Verification of Assumptions (a) and (b) in Part A of
Proposition 5.12 of \cite{PP3} proceeds similar as before except for now
using Parts (b) and (a) of Lemma \ref{lem:nullset}. Part A of Proposition
5.12 of \cite{PP3} now immediately delivers claim (\ref%
{size-control_groupHet_uncorr_tilde}), again since $C^{\ast }<\infty $ as
noted in that proposition.\footnote{%
This argument is actually superfluous since $\tilde{T}_{uc}$ is bounded as
noted in Section \ref{sec:trivial_uc_tilde}. However, verification of the
assumptions of Proposition 5.12 in \cite{PP3} is essential for the proof of
the other claims in Part(a) of Theorem \ref{theorem_groupwise_hetero_tilde}.}
Again, $G(\mathfrak{M}_{0})$-invariance of $\tilde{T}_{uc}$ implies that $%
C^{\ast }$ and $\alpha ^{\ast }$ do not depend on the choice of $\mu _{0}\in 
\mathfrak{M}_{0}$, and that $\alpha ^{\ast }$ as defined in the theorem
coincides with $\alpha ^{\ast }$ as defined in Proposition 5.12 of \cite{PP3}%
. The remaining claim in Part (a) is proved completely analogous as the
corresponding claim in Part (b) except for now using Part (b) of Lemma \ref%
{lem:nullset}.

We finally prove Part (c):\ The claims follow from Remark 5.10 in \cite{PP3}
and Lemma \ref{lem:nullset}; cf.~also Appendix \ref{useful}. $\blacksquare $

\section{Appendix: Algorithms\label{app:algos}}

In this appendix, we discuss in more detail algorithms for determining (i)
rejection probabilities, (ii) the size of a test based on one of the test
statistics $T_{Het}$, $T_{uc}$, $\tilde{T}_{Het}$, or $\tilde{T}_{uc}$
together with a given candidate critical value, and (iii) size-controlling
critical values. We discuss these algorithms under the Gaussianity
assumption made in Section \ref{frame}, but recall from Section \ref%
{non-gaussianity} that the algorithms as given here can also be used to
calculate \emph{null} rejection probabilities, size, and size-controlling
critical values in the elliptically symmetric case \emph{without any changes}%
; similarly, the algorithms given here can also be used to calculate the
size and size-controlling critical values in the semiparametric model
discussed in (iv) of Section \ref{non-gaussianity} as they stand.
Furthermore, we restrict ourselves to the heteroskedasticity model $%
\mathfrak{C}_{Het}$; adapting the algorithms to subsets $\mathfrak{C}$ of $%
\mathfrak{C}_{Het}$ is rather straightforward (basically one has to
appropriately constrain the optimization routines involved, appropriately
redefine some of the quantities like $C_{low}$, and refer to the
size-control conditions pertinent to the given heteroskedasticity model $%
\mathfrak{C}$).

\subsection{Computing rejection probabilities\label{sec:obsq1}}

Suppose that a $G(\mathfrak{M}_{0})$-invariant test statistic $T:\mathbb{R}%
^{n}\rightarrow \mathbb{R}$ has the following property: for some (and hence
any) $\mu _{0}\in \mathfrak{M}_{0}$ and a critical value $C\in \mathbb{R}$,
there exists a symmetric $n\times n$ matrix $A_{C}$, such that%
\begin{equation}
T(\mu _{0}+z)\geq C\Leftrightarrow z^{\prime }A_{C}z\geq 0\text{ holds for }%
\lambda _{\mathbb{R}^{n}}\text{-almost every }z\in \mathbb{R}^{n}.
\label{eqn:Aexists}
\end{equation}%
If this property is satisfied, then for all choices of $\Sigma \in \mathfrak{%
C}_{Het}$, $\mu _{0}\in \mathfrak{M}_{0}$, $\mu \in \mathfrak{M}$, and $%
\sigma ^{2}\in (0,\infty )$, setting $\nu :=\sigma ^{-1}\Sigma ^{-1/2}(\mu
-\mu _{0})$, we may write 
\begin{equation}
P_{\mu ,\sigma ^{2}\Sigma }(\left\{ z\in \mathbb{R}^{n}:T(z)\geq C\right\}
)=P_{\nu ,I_{n}}(\{\zeta \in \mathbb{R}^{n}:\zeta ^{\prime }\Sigma
^{1/2}A_{C}\Sigma ^{1/2}\zeta \geq 0\});  \label{eqn:qformgen}
\end{equation}%
in case $\mu \in \mathfrak{M}_{0}$, we may set $\mu _{0}=\mu $ to further
simplify the right-hand-side in (\ref{eqn:qformgen}) to%
\begin{equation}
P_{0,I_{n}}(\{\zeta \in \mathbb{R}^{n}:\zeta ^{\prime }\Sigma
^{1/2}A_{C}\Sigma ^{1/2}\zeta \geq 0\}).  \label{eqn:nullQF}
\end{equation}

The probability that a Gaussian quadratic form is not less than $0$ (such as
(\ref{eqn:qformgen}) or (\ref{eqn:nullQF})) can numerically be determined by
standard algorithms such as \cite{davies}. Relation (\ref{eqn:Aexists}) can
thus be exploited for efficiently computing rejection probabilities (for a
given critical value), and thus plays an instrumental r\^{o}le in
numerically determining the size of a test, size-controlling critical
values, or the power function of a test.

For the important case $q=1$ we now show that the above approach can indeed
be used. It follows from the subsequent lemma that for any critical value $C$
the property in (\ref{eqn:Aexists}) holds for the following test statistics:
(i) $T_{Het}$ provided Assumption \ref{R_and_X} holds; (ii) $T_{uc}$; (iii) $%
\tilde{T}_{Het}$ provided Assumption \ref{R_and_X_tilde} holds; (iv) $\tilde{%
T}_{uc}$. Recall from Lemmata \ref{lem_B} and \ref{lem:B_tilde} that under
Assumption \ref{R_and_X} (Assumption \ref{R_and_X_tilde}, respectively), the
set $\mathsf{B}$ ($\mathsf{\tilde{B}}$, respectively) is a $\lambda _{%
\mathbb{R}^{n}}$-null set. Note that $v$ defined in the lemma satisfies $%
v\neq 0$.

\begin{lemma}
\label{lem:qform} Suppose $q=1$. Let $v=v_{R,X}:=X(X^{\prime
}X)^{-1}R^{\prime }$. Then, for every $C\in \mathbb{R}$ and every $\mu
_{0}\in \mathfrak{M}_{0}$, we have:

(a) If $\mu _{0}+z\notin \mathsf{B}$, then $T_{Het}(\mu _{0}+z)\geq C$ ($%
\leq C$) is equivalent to $z^{\prime }A_{Het,C}z\geq 0$ ($\leq 0$), where%
\begin{equation}
A_{Het,C}:=vv^{\prime }-C\Pi _{\limfunc{span}(X)^{\bot }}\limfunc{diag}%
\left( v_{1}^{2}d_{1},\ldots ,v_{n}^{2}d_{n}\right) \Pi _{\limfunc{span}%
(X)^{\bot }}.  \label{eqn:MHetdef}
\end{equation}%
(b) If $\mu _{0}+z\notin \limfunc{span}(X)$, then $T_{uc}(\mu _{0}+z)\geq C$
($\leq C$) is equivalent to $z^{\prime }A_{uc,C}z\geq 0$ ($\leq 0$), where%
\begin{equation}
A_{uc,C}:=vv^{\prime }-C\frac{v^{\prime }v}{n-k}\Pi _{\limfunc{span}%
(X)^{\bot }}.  \label{eqn:Mucdef}
\end{equation}%
(c) If $\mu _{0}+z\notin \mathsf{\tilde{B}}$, then $\tilde{T}_{Het}(\mu
_{0}+z)\geq C$ ($\leq C$) is equivalent to $z^{\prime }\tilde{A}%
_{Het,C}z\geq 0$ ($\leq 0$), where%
\begin{equation}
\tilde{A}_{Het,C}:=vv^{\prime }-C\Pi _{(\mathfrak{M}_{0}^{lin})^{\bot }}%
\limfunc{diag}\left( v_{1}^{2}\tilde{d}_{1},\ldots ,v_{n}^{2}\tilde{d}%
_{n}\right) \Pi _{(\mathfrak{M}_{0}^{lin})^{\bot }}.
\label{eqn:tildeMHetdef}
\end{equation}%
(d) If $\mu _{0}+z\notin \mathfrak{M}_{0}$, then $\tilde{T}_{uc}(\mu
_{0}+z)\geq C$ ($\leq C$) is equivalent to $z^{\prime }\tilde{A}_{uc,C}z\geq
0$ ($\leq 0$), where%
\begin{equation}
\tilde{A}_{uc,C}:=vv^{\prime }-C\frac{v^{\prime }v}{n-(k-1)}\Pi _{(\mathfrak{%
M}_{0}^{lin})^{\bot }}.  \label{eqn:tildeMucdef}
\end{equation}
\end{lemma}

\textbf{Proof:} We first observe that there is nothing to prove in Part (a)
(Part (c), respectively) if Assumption \ref{R_and_X} (Assumption \ref%
{R_and_X_tilde}, respectively) is violated, since then $\mathsf{B}=\mathbb{R}%
^{n}$ ($\mathsf{\tilde{B}}=\mathbb{R}^{n}$, respectively) by Lemma \ref%
{lem_B} (Lemma \ref{lem:B_tilde}, respectively). In the following we hence
may assume for Part (a) (Part (c), respectively) that Assumption \ref%
{R_and_X} (Assumption \ref{R_and_X_tilde}, respectively) hold, in which case 
$\mathsf{B}$ ($\mathsf{\tilde{B}}$, respectively) is a $\lambda _{\mathbb{R}%
^{n}}$-null set. The expressions in (\ref{eqn:MHetdef}), (\ref{eqn:Mucdef}),
(\ref{eqn:tildeMHetdef}), and (\ref{eqn:tildeMucdef}) now follow directly
from the definitions of the test statistics since $q=1$, recalling in
particular that $\hat{u}(\mu _{0}+z)=\Pi _{\limfunc{span}(X)^{\bot }}(\mu
_{0}+z)=\Pi _{\limfunc{span}(X)^{\bot }}z$, and $\tilde{u}(\mu _{0}+z)=\Pi
_{(\mathfrak{M}_{0}^{lin})^{\bot }}((\mu _{0}+z)-\mu _{0})=\Pi _{(\mathfrak{M%
}_{0}^{lin})^{\bot }}z$, and noting that for $q=1$%
\begin{align*}
& R\hat{\beta}(\mu _{0}+z)=r+v^{\prime }z, \\
& \hat{\Omega}_{Het}(\mu _{0}+z)=z^{\prime }\Pi _{\limfunc{span}(X)^{\bot }}%
\limfunc{diag}(v_{1}^{2}d_{1},\ldots ,d_{n}v_{n}^{2})\Pi _{\limfunc{span}%
(X)^{\bot }}z, \\
& \tilde{\Omega}_{Het}(\mu _{0}+z)=z^{\prime }\Pi _{(\mathfrak{M}%
_{0}^{lin})^{\bot }}\limfunc{diag}(v_{1}^{2}\tilde{d}_{1},\ldots ,\tilde{d}%
_{n}v_{n}^{2})\Pi _{(\mathfrak{M}_{0}^{lin})^{\bot }}z, \\
& \hat{\sigma}^{2}(\mu _{0}+z)=\frac{z^{\prime }\Pi _{\limfunc{span}%
(X)^{\bot }}z}{n-k},~\tilde{\sigma}^{2}(z)=\frac{z^{\prime }\Pi _{(\mathfrak{%
M}_{0}^{lin})^{\bot }}z}{n-(k-1)}
\end{align*}%
hold.$\ \blacksquare $

\begin{remark}
The algorithm in \cite{davies} applied to (\ref{eqn:qformgen}) requires that
the matrix $A_{C}$ is not the zero matrix. In (i)-(iii) below we always have 
$q=1$.

(i) It is easy to see that $A_{Het,C}$, $A_{uc,C}$, and $\tilde{A}_{uc,C}$
are never equal to the zero matrix: Note that $v^{\prime
}A_{Het,C}v=(v^{\prime }v)^{2}>0$, since $v\in \limfunc{span}(X)$ and $v\neq
0$. The same argument applies to $A_{uc,C}$. Furthermore, for $C=0$ the
matrix $\tilde{A}_{uc,C}$ is obviously not the zero matrix; for $C\neq 0$
let $w\in (\mathfrak{M}_{0}^{lin})^{\bot }$, $w\neq 0$, $w$ orthogonal to $v$%
, then $w^{\prime }\tilde{A}_{uc,C}w=-w^{\prime }wCv^{\prime
}v/(n-(k-1))\neq 0$ (note that such a $w$ exists, since $v\in (\mathfrak{M}%
_{0}^{lin})^{\bot }$ and $\dim ((\mathfrak{M}_{0}^{lin})^{\bot
})=n-(k-q)>n-k\geq 1$ hold).

(ii) For $\tilde{A}_{Het,C}$ we have the following: Since $v\in (\mathfrak{M}%
_{0}^{lin})^{\bot }$ holds, $v^{\prime }\tilde{A}_{Het,C}v=(v^{\prime
}v)^{2}-Cv^{\prime }\limfunc{diag}(v_{1}^{2}\tilde{d}_{1},\ldots ,v_{n}^{2}%
\tilde{d}_{n})v$, which is zero only for $C=C_{0}$ where $%
C_{0}=\sum_{i=1}^{n}v_{i}^{2}/\sum_{i=1}^{n}v_{i}^{4}\tilde{d}_{i}$ (note
that the ratio is well-defined since all the $\tilde{d}_{i}$ are positive
and since $v\neq 0$). Hence, $\tilde{A}_{Het,C}$ is not the zero matrix,
except possibly for $C=C_{0}$. We now show that -- in case Assumption \ref%
{R_and_X_tilde} is satisfied -- $\tilde{A}_{Het,C_{0}}=0$ is equivalent to $%
\tilde{T}_{Het}(y)$ being constant for $y\in $ $\mathbb{R}^{n}\backslash 
\mathsf{\tilde{B}}$: Suppose $\tilde{A}_{Het,C_{0}}=0$. Since $C_{0}>0$, we
obtain $\Pi _{(\mathfrak{M}_{0}^{lin})^{\bot }}\limfunc{diag}(v_{1}^{2}%
\tilde{d}_{1},\ldots ,v_{n}^{2}\tilde{d}_{n})\Pi _{(\mathfrak{M}%
_{0}^{lin})^{\bot }}=vv^{\prime }/C_{0}$ and thus $\tilde{A}%
_{Het,C}=vv^{\prime }(1-C/C_{0})$. Fix $\mu _{0}\in \mathfrak{M}_{0}$
arbitrary. For every $C>C_{0}$ we have $z^{\prime }\tilde{A}_{Het,C}z\leq 0$
for every $z$, and hence for every $z$ with $\mu _{0}+z\notin \mathsf{\tilde{%
B}}$ (note that $\mathbb{R}^{n}\backslash \mathsf{\tilde{B}}$ is nonempty
under Assumption \ref{R_and_X_tilde}). By Lemma \ref{lem:qform} we can
conclude that $\tilde{T}_{Het}(\mu _{0}+z)\leq C$ for every $\mu
_{0}+z\notin \mathsf{\tilde{B}}$. By the same token, we obtain that $\tilde{T%
}_{Het}(\mu _{0}+z)\geq C$ for every $\mu _{0}+z\notin \mathsf{\tilde{B}}$
when $C<C_{0}$ holds. We conclude that $\tilde{T}_{Het}(\mu _{0}+z)=C_{0}$
for every $\mu _{0}+z\notin \mathsf{\tilde{B}}$, i.e., $\tilde{T}%
_{Het}(y)=C_{0}$ for every $y\in $ $\mathbb{R}^{n}\backslash \mathsf{\tilde{B%
}}$. To prove the converse, assume $\tilde{T}_{Het}(y)=C_{1}$ for every $%
y\notin \mathsf{\tilde{B}}$. Fix $\mu _{0}\in \mathfrak{M}_{0}$ arbitrary.
Then $\tilde{T}_{Het}(\mu _{0}+z)=C_{1}$ for every $z$ with $\mu
_{0}+z\notin \mathsf{\tilde{B}}$. By Lemma \ref{lem:qform} we get $z^{\prime
}\tilde{A}_{Het,C}z\geq 0$ ($\leq 0$, respectively) for $C\leq C_{1}$ ($%
C\geq C_{1}$, respectively) for every $z\notin \mathsf{\tilde{B}}-\mu _{0}$.
Under Assumption \ref{R_and_X_tilde}\ the set $\mathsf{\tilde{B}}-\mu _{0}$
is a $\lambda _{\mathbb{R}^{n}}$-null set, hence its complement is dense in $%
\mathbb{R}^{n}$. By continuity of the quadratic forms, we get $z^{\prime }%
\tilde{A}_{Het,C}z\geq 0$ ($\leq 0$, respectively) for $C\leq C_{1}$ ($C\geq
C_{1}$, respectively) for all $z\in \mathbb{R}^{n}$. We thus obtain $%
z^{\prime }\tilde{A}_{Het,C_{1}}z=0$ for every $z\in \mathbb{R}^{n}$. Since $%
\tilde{A}_{Het,C_{1}}$ is symmetric, $\tilde{A}_{Het,C_{1}}=0$ follows and $%
C_{1}=C_{0}$ must hold.

(iii) Before applying the algorithm in \cite{davies} to (\ref{eqn:qformgen})
with $T=T_{Het}$ and $A_{C}=A_{Het,C}$ we first check that Assumption \ref%
{R_and_X} holds since otherwise Part (a) of the preceding lemma does not
apply. In case of $T=\tilde{T}_{Het}$ and $A_{C}=\tilde{A}_{Het,C}$ we check
that Assumption \ref{R_and_X_tilde} holds for similar reasons; and, in case
this assumption is satisfied, we then always also compute $C_{0}$ and check
numerically that $\tilde{A}_{Het,C_{0}}$ (and hence any $\tilde{A}_{Het,C}$)
is not the zero matrix.
\end{remark}

In case $q>1$, the algorithm in \cite{davies} could also be used to compute
rejection probabilities for the tests based on $T_{uc}$ and $\tilde{T}_{uc}$
as is easy to see. Since this is not so for $T_{Het}$ and $\tilde{T}_{Het}$,
we do not proceed in this way for reasons of comparability. In case $q>1$ we
thus compute the required rejection probabilities by Monte Carlo.

\subsection{Determining the size of a test\label{app:S}}

\emph{For simplicity throughout this subsection $T$ denotes any one of the
test statistics UC, HC0-HC4, UCR, HC0R-HC4R. In case of HC0-HC4 we assume in
our discussion that the design matrix $X$ and $R$ are such that Assumption %
\ref{R_and_X} is satisfied, and in case of HC0R-HC4R we assume that
Assumption \ref{R_and_X_tilde} holds and that the test statistic is not
constant on }$\mathbb{R}^{n}\backslash \mathsf{\tilde{B}}$.\footnote{%
This rules out trivial cases only.} These conditions should be checked
either theoretically or numerically before using the algorithms described
below. Such numerical checks are implemented in the R-package \textbf{hrt} (%
\cite{hrt}) realizing these algorithms.

We now discuss algorithms for determining the size (over $\mathfrak{C}_{Het}$%
) of the test that rejects if $T\geq C$ for a given critical value $C>0$
(note that any $C\leq 0$ leads to a trivial test that always rejects). By $G(%
\mathfrak{M}_{0})$-invariance of $T$, for any given $\mu _{0}\in \mathfrak{M}%
_{0}$, the size of this test simplifies to%
\begin{equation}
\sup_{\Sigma \in \mathfrak{C}_{Het}}P_{\mu _{0},\Sigma }(T\geq C),
\label{eqn:sizeHet}
\end{equation}%
which is what the algorithms described below compute numerically.

Before trying to determine the size numerically, it is advisable to check
whether $C$ is not less than the pertinent lower bound $C^{\ast }$ for
size-controlling critical values obtained in our theoretical results in
Propositions \ref{rem_C*} and \ref{rem_C*_tilde} (and the attending
footnotes), since otherwise one already knows that the size of the test is
equal to $1$, and hence there is no need to run the algorithm. The
implementations of the algorithms in the R-package \textbf{hrt} (\cite{hrt})
have an option that provides such a check and outputs $1$ if the check fails
without running the algorithm.

Of course, the design matrix $X$, the restriction $(R,r)$, and the
particular choice of test statistic from the above list, are inputs to all
the algorithms that are discussed in this and the subsequent section \ref%
{app:SCCV}, but we do not show these inputs explicitly in the descriptions
of the algorithms given further down.

\subsubsection{Case $q=1$}

In the important special case $q=1$ we can use (\ref{eqn:nullQF}) and Lemma %
\ref{lem:qform} to compute the rejection probabilities $P_{\mu _{0},\Sigma
}(T\geq C)$ appearing in (\ref{eqn:sizeHet}) efficiently via, e.g., \cite%
{davies} (referred to as DA in what follows). A generic algorithm based on
this observation is summarized in Algorithm \ref{alg:S1}. 
\begin{algorithm}
\caption{Computing the size in case~$q = 1$.} \label{alg:S1} 
\begin{algorithmic}[1]
\State \textbf{Input} A real number $C > 0$ and positive integers $M_0 \geq M_1 \geq M_2$.
\State \textbf{Stage 0: Initial value search}
\For{$j = 1$ \texttt{to} $j = M_0$} 
\State Choose a candidate~$\Sigma_j \in \mathfrak{C}_{Het}$. 
\State Obtain~$\tilde{p}_j := P_{\mu_0, \Sigma_j}(T \geq C)$ using DA.
\EndFor
\State Rank the candidates $\Sigma_j$ according to the value (from largest to smallest) of the corresponding quantities $\tilde{p}_j $ to obtain $\Sigma_{1:M_0}, \hdots, \Sigma_{M_1:M_0}$, the initial values for the next stage.
\State \textbf{Stage 1: Coarse localized optimizations}
\For{$j = 1$ \texttt{to} $j = M_1$} 
\State Obtain $\Sigma^*_{j}$ by running a numerical algorithm for the optimization problem~\eqref{eqn:sizeHet} initialized at $\Sigma_{j:M_0}$ and obtain $\bar{p}_{j, \Sigma^*_{j}} := P_{\mu_0, \Sigma^*_j}(T \geq C)$ (using DA to evaluate probabilities).
\EndFor
\State Rank the obtained matrices $\Sigma^*_{j}$ according to the value (from largest to smallest) of the corresponding $\bar{p}_{j, \Sigma^*_{j}}$ to obtain $\Sigma^*_{1:M_1}, \hdots, \Sigma^*_{M_2:M_1}$, the initial values for the next stage.
\State \textbf{Stage 2: Refined localized optimization}
\For{$j = 1$ \texttt{to} $j = M_2$} 
\State Obtain $\Sigma^{**}_{j}$ by running a (refined) numerical algorithm for the optimization problem~\eqref{eqn:sizeHet} initialized at $\Sigma^*_{j:M_1}$ and obtain $\bar{\bar{p}}_{j, \Sigma^{**}_{j}} := P_{\mu_0, \Sigma^{**}_j}(T \geq C)$ (using DA to evaluate probabilities).
\EndFor
\State \textbf{Return} $\max_{j = 1, \hdots, M_2} \bar{\bar{p}}_{j, \Sigma^{**}_{j}}$
\end{algorithmic}
\end{algorithm}

\begin{remark}
\label{rem:starging} The initial values $\Sigma _{j}$ in Stage 0 of
Algorithm \ref{alg:S1} can, for example, be obtained randomly (e.g., by
sampling the diagonal elements of $\Sigma _{j}$ from a uniform distribution
on the unit simplex in $\mathbb{R}^{n}$). Such random choices may then be
supplemented by \textquotedblleft special\textquotedblright\ elements of $%
\mathfrak{C}_{Het}$, e.g., matrices that are close to $e_{i}(n)e_{i}(n)^{%
\prime }$, $i=1,\ldots ,n$, or the matrix $n^{-1}I_{n}$, or a matrix $\Sigma 
$ that maximizes the expectation of the quadratic form $y\mapsto y^{\prime
}\Sigma ^{1/2}A_{C}\Sigma ^{1/2}y$ under $P_{0,I_{n}}$ (where $A_{C}$ is
obtained via Lemma \ref{lem:qform}, cf.~also the discussion preceding that
lemma), the latter choice being motivated by (\ref{eqn:nullQF}). For the
particular choice of initial values used in the R-package \textbf{hrt} and
in our numerical calculations see \cite{hrt} and Appendix \ref{details}.
\end{remark}

\begin{remark}
\label{rem:initval_2} If Algorithm \ref{alg:S1} is to be applied to a
relatively large critical value $C$ (say $C$ larger than $5$ times the $%
(1-\alpha )$-quantile of the cdf of $P_{0,I_{N}}\circ T$), then one may run
Algorithm \ref{alg:S1} on a smaller critical value first (e.g., the just
mentioned quantile), and use the covariance matrix realizing the maximal
rejection probability for this smaller critical value (in line 17 of
Algorithm \ref{alg:S1}) as an additional initial value when running
Algorithm \ref{alg:S1} for determining the size corresponding to the
originally given $C$. This can help to ameliorate numerical difficulties due
to the rejection probabilities being close to zero over large portions of $%
\mathfrak{C}_{Het}$. The just described procedure is available as an option
in the R-package \textbf{hrt}.
\end{remark}

\begin{remark}
\label{rem:maxalg} The concrete choice of the numerical optimization
algorithm used in Stages 1 and 2 of Algorithm \ref{alg:S1} is left
unspecified here, but may, for example, be a constrained \cite{nelder}
algorithm (as provided in R's \textquotedblleft
constrOptim\textquotedblright\ function), where in Stage 2 the parameters in
this algorithm (and in principle also in DA) should be chosen to guarantee a
higher accuracy. For the particular choice of optimization routines used in
the R-package \textbf{hrt} and in our numerical calculations see \cite{hrt}
and Appendix \ref{details}.
\end{remark}

Remarks \ref{rem:starging} and \ref{rem:maxalg} also apply to other
algorithms introduced further down, and will not be repeated.

\subsubsection{General case}

An algorithm that is similar to Algorithm \ref{alg:S1}, but uses Monte-Carlo
simulation instead of DA to compute the rejection probabilities $P_{\mu
_{0},\Sigma }(T\geq C)$ is discussed in Algorithm \ref{alg:S2}; this
algorithm is a modification of Algorithm 2 in \cite{PP3}.\footnote{%
This algorithm involves evaluating the test statistic $T$. Since the
definition of $T$ depends on invertibility of a covariance matrix estimator,
an invertibility check is required. We use the same invertibility check as
discussed in the second paragraph in Appendix E.3 of \cite{PPBoot}, with a
tolerance parameter that can be specified by the user.} In Algorithm \ref%
{alg:S2} the number of replications used in the Monte-Carlo simulations (and
thus their accuracy but also their runtime) is increased in each stage,
leading to an improved accuracy in the rejection probabilities computed.
While this algorithm is also applicable in case $q=1$, Algorithm \ref{alg:S1}
is to be preferred (and is automatically applied by the R-package \textbf{hrt%
} in this case), as it is based on a preferable way of computing the
rejection probabilities. 
\begin{algorithm}
\caption{Computing the size for general $q$.} \label{alg:S2} 
\begin{algorithmic}[1]
\State \textbf{Input} A real number $C > 0$ and positive integers $M_0 \geq M_1 \geq M_2$, $N_0 \leq  N_1 \leq N_2$.
\State \textbf{Stage 0: Initial value search}
\For{$j = 1$ \texttt{to} $j = M_0$} 
\State Generate a pseudorandom sample $Z_{1}, \hdots, Z_{N_0}$ from $P_{0, I_n}$.
\State Obtain a candidate $\Sigma_j \in \mathfrak{C}_{Het}$. 
\State Compute $\tilde{p}_j = N_0^{-1} \sum_{i = 1}^{N_0} \mathbf{1}_{[C, \infty)}(T(\mu_0 + \Sigma_j^{1/2} Z_{i}))$.
\EndFor
\State Rank the candidates $\Sigma_j$ according to the value (from largest to smallest) of the corresponding quantities $\tilde{p}_j $ to obtain $\Sigma_{1:M_0}, \hdots, \Sigma_{M_1:M_0}$, the initial values for the next stage.
\State \textbf{Stage 1: Coarse localized optimizations}
\For{$j = 1$ \texttt{to} $j = M_1$} 
\State Generate a pseudorandom sample $Z_{1}, \hdots, Z_{N_1}$ from $P_{0, I_n}$. 
\State Define $\bar{p}_{j, \Sigma} = N_1^{-1} \sum_{i = 1}^{N_1} \mathbf{1}_{[C, \infty)}(T(\mu_0 + \Sigma^{1/2} Z_{i}))$ for $\Sigma \in \mathfrak{C}_{Het}$.
\State Obtain $\Sigma^*_{j}$ by running a numerical optimization algorithm for the problem $\sup_{\Sigma \in \mathfrak{C}_{Het}} \bar{p}_{j, \Sigma }$ initialized at $\Sigma_{j:M_0}$.
\EndFor
\State Rank the obtained numbers $\Sigma^*_{j}$ according to the value (from largest to smallest) of the corresponding $\bar{p}_{j, \Sigma^*_{j}}$ to obtain $\Sigma^*_{1:M_1}, \hdots, \Sigma^*_{M_2:M_1}$, the initial values for the next stage.
\State \textbf{Stage 2: Refined localized optimization}
\For{$j = 1$ \texttt{to} $j = M_2$}
\State Generate a pseudorandom sample $Z_{1}, \hdots, Z_{N_2}$ from $P_{0, I_n}$. 
\State Define $\bar{\bar{p}}_{j, \Sigma} = N_2^{-1} \sum_{i = 1}^{N_2} \mathbf{1}_{[C, \infty)}(T(\mu_0 + \Sigma^{1/2}Z_{i}))$ for $\Sigma \in \mathfrak{C}_{Het}$.
\State Obtain $\Sigma^{**}_{j}$ by running a numerical optimization algorithm for the problem $\sup_{\Sigma \in \mathfrak{C}_{Het}} \bar{\bar{p}}_{j, \Sigma}$ initialized at $\Sigma^*_{j:M_1}$. 
\EndFor
\State \textbf{Return} $\max_{j = 1, \hdots, M_2} \bar{\bar{p}}_{j, \Sigma^{**}_{j}} $.
\end{algorithmic}
\end{algorithm}

\subsection{Determining smallest size-controlling critical values\label%
{app:SCCV}}

\emph{Again, in this subsection $T$ denotes any one of the test statistics
UC, HC0-HC4, UCR, HC0R-HC4R. In case of HC0-HC4 we assume in our discussion
that the design matrix $X$ and $R$ are such that Assumption \ref{R_and_X} is
satisfied, and in case of HC0R-HC4R we assume that Assumption \ref%
{R_and_X_tilde} holds and that the test statistic is not constant on }$%
\mathbb{R}^{n}\backslash \mathsf{\tilde{B}}$.\footnote{%
This rules out trivial cases only.} \emph{Furthermore, we assume that
size-controlling critical values exist. }These conditions should be checked
either theoretically or numerically before using the algorithms described
below. The last mentioned existence can be guaranteed by checking
(theoretically or numerically) the respective sufficient conditions for size
control in Theorems \ref{Hetero_Robust} and \ref{Hetero_Robust_tilde}.%
\footnote{%
In case the respective sufficient conditions are violated, but
size-controlling critical values nevertheless exist (as, e.g., in Example %
\ref{ex_k_pop} or in Remark \ref{q=k}), the algorithm still works.} We note
that the implementations of the algorithms presented below in the R-package 
\textbf{hrt (}\cite{hrt}) include such numerical checks.

We now proceed to discussing several algorithms for determining the \emph{%
smallest} critical value $C_{\Diamond }(\alpha )$ such that the size of the
test, which rejects if $T\geq C_{\Diamond }(\alpha )$, does not exceed $%
\alpha $ ($0<\alpha <1$).\footnote{%
Such a \emph{smallest} size-controlling critical value indeed exists under
the assumptions of this subsection (which includes existence of a
size-controlling critical value) in view of Appendix \ref{useful}. [Under
the sufficient conditions for size control in the respective theorems, this
can also be read off directly from these theorems.]} [In fact, for $%
C_{\Diamond }(\alpha )$ the size then equals $\alpha $ \emph{provided} a
critical value that results in size \emph{equal} to $\alpha $ actually
exists.] Note that $C_{\Diamond }(\alpha )>0$ must hold, in view of Remarks %
\ref{rem_positiv} and \ref{rem_positiv_tilde} since $\alpha <1$. By $G(%
\mathfrak{M}_{0})$-invariance, for some fixed $\mu _{0}\in \mathfrak{M}_{0}$%
, the algorithms numerically compute the smallest critical value that
satisfies%
\begin{equation}
\sup_{\Sigma \in \mathfrak{C}_{Het}}P_{\mu _{0},\Sigma }(T\geq C)\leq \alpha
,  \label{eqn:sizeHet_2}
\end{equation}%
cf. the discussion surrounding (\ref{eqn:sizeHet}). For later use we denote
by $F_{\Sigma }$ the cdf of $P_{\mu _{0},\Sigma }\circ T$, which by $G(%
\mathfrak{M}_{0})$-invariance does not depend on the particular choice for $%
\mu _{0}\in \mathfrak{M}_{0}$.

\subsubsection{Computing smallest size-controlling critical values via line
search based on algorithms in Section~\protect\ref{app:S} \label{app:line}}

Given an algorithm $\mathsf{A}:(0,\infty )\rightarrow \lbrack 0,1]$ that for 
$C>0$ returns the size of the test that rejects if $T\geq C$, one can use a
line-search algorithm to determine the smallest critical value $%
C=C_{\Diamond }(\alpha )$ satisfying $\mathsf{A}(C)\leq \alpha $. To this
end, one starts at the lower bound $C_{low}=\max (C^{\ast },C_{\hom })$,
where $C^{\ast }$ is given in the pertinent parts of Theorems \ref%
{Hetero_Robust} and \ref{Hetero_Robust_tilde}, respectively (cf.~also
Propositions \ref{rem_C*} and \ref{rem_C*_tilde}, respectively, and the
attending footnotes), and $C_{\hom }$ denotes the \emph{smallest} $1-\alpha $
quantile of $F_{I_{n}}$, i.e., of the cdf of the test statistic under
homoskedasticity. Note that then $P_{\mu _{0},I_{n}}(T\geq C_{\hom })=\alpha 
$ and that $P_{\mu _{0},I_{n}}(T\geq C)>\alpha $ for $C<C_{\hom }$ (to see
this note that $F_{I_{n}}$ is continuous as $\{T=C\}$ is a $\lambda _{%
\mathbb{R}^{n}}$-null set for all real $C$, cf. Lemma 5.16 in\ \cite{PP3}
and Lemma \ref{lem:nullset} in Appendix \ref{app_C}). Furthermore, $C_{\hom
}>0$ (since $T\geq 0$ and $\{T=0\}$ is a $\lambda _{\mathbb{R}^{n}}$-null
set), and consequently $C_{low}>0$ holds. Starting from $C_{low}$, one then
keeps increasing the critical value \textquotedblleft in a reasonable
way\textquotedblright\ until one obtains, for the first time, a $C$ such
that $\mathsf{A}(C)\leq \alpha $ holds. This procedure is summarized in
Algorithm \ref{alg:CV1}, in which the particular algorithm $\mathsf{A}$ used
is an input to Algorithm \ref{alg:CV1}. For $\mathsf{A}$ one may either use
Algorithm \ref{alg:S1} if $q=1$, or Algorithm \ref{alg:S2} for general $q$.
Note that one may need to terminate the while-loop after a maximal number of
iterations. 
\begin{algorithm}
\caption{Numerical approximation of the smallest size-controlling critical value via a line search algorithm.} \label{alg:CV1} 
\begin{algorithmic}[1]
\State \textbf{Input} ~$\alpha \in (0, 1)$, $\mathsf{A}$,~$C_{low}$,~$\epsilon \in [0,1-\alpha )$  ($\epsilon$ a small tolerance parameter).
\State $C \leftarrow C_{low}$
\While{$\mathsf{A}(C) > \alpha + \epsilon$}
\State Let~$\Sigma^{**}$ be such that~$P_{\mu_0, \Sigma^{**}}(T \geq C) \approx \mathsf{A}(C)$.
\State Determine, by an upward line search initialized at~$C$, the smallest value~$C_+$ such that~$P_{\mu_0, \Sigma^{**}}(T \geq C_+) \leq \alpha$ .
\State $C \leftarrow C_+$.
\EndWhile
\State 	 	\Return $C$
\end{algorithmic}
\end{algorithm}

\begin{remark}
(i) Note that a matrix $\Sigma ^{\ast \ast }$ as required for the while-loop
in Algorithm \ref{alg:CV1} can easily be obtained by implementing Algorithm %
\ref{alg:S1} or \ref{alg:S2} in such a fashion as to also return the
covariance matrix for which the maximal rejection probability is attained in
the respective Stage 2.

(ii)\ A smallest $C_{+}$ as required in line 5 of Algorithm \ref{alg:CV1}
indeed exists since $\{T=C\}$ is a $\lambda _{\mathbb{R}^{n}}$-null set for
all real $C$ as noted before.

(iii) For details regarding the computation of $C_{low}$ in the R-package 
\textbf{hrt} see \cite{hrt} and Appendix \ref{app:CVcomppow}.
\end{remark}

\subsubsection{Computing smallest size-controlling critical values via
quantile maximization}

For completeness and comparison with \cite{PP3}, we briefly describe an
algorithm that is a modification of Algorithm 1 in \cite{PP3}. In contrast
to the algorithm discussed in the previous section, it does not make use of
size-computations, but determines the smallest size-controlling critical
value as 
\begin{equation}
\sup_{\Sigma \in \mathfrak{C}_{Het}}F_{\Sigma }^{-1}(1-\alpha )
\label{quant-funct}
\end{equation}%
where $F_{\Sigma }^{-1}$ denotes the quantile function of the cdf $F_{\Sigma
}$. That (\ref{quant-funct}) indeed gives the smallest size-controlling
critical value is not difficult to see keeping in mind that $P_{\mu
_{0},\Sigma }(T=C)=0$ for every real $C$, every $\mu _{0}\in \mathfrak{M}%
_{0} $, and every $\Sigma \in \mathfrak{C}_{Het}$ (in view of $\lambda _{%
\mathbb{R}^{n}}(\{T=C\})=0$ as noted before). The algorithm is summarized in
Algorithm \ref{alg:CV2}. 
\begin{algorithm}
\caption{Numerical approximation of the smallest size-controlling critical value via quantiles.} \label{alg:CV2} 
\begin{algorithmic}[1]
\State \textbf{Input} Positive integers $M_0 \geq M_1 \geq M_2$, $N_0 \leq  N_1 \leq N_2$.
\State \textbf{Stage 0: Initial value search}
\For{$j = 1$ \texttt{to} $j = M_0$} 
\State Generate a pseudorandom sample $Z_{1}, \hdots, Z_{N_0}$ from $P_{0, I_n}$. 
\State Obtain a candidate $\Sigma_j \in \mathfrak{C}_{Het}$. 
\State Compute $\tilde{F}_{j}^{-1}(1-\alpha)$ where $\tilde{F}_j(x) = N_0^{-1} \sum_{i = 1}^{N_0} \mathbf{1}_{(-\infty, x]}(T(\mu_0 + \Sigma_j^{1/2} Z_{i}))$ for $x \in \mathbb{R}$.
\EndFor
\State Rank the candidates $\Sigma_j$ according to the value (from largest to smallest) of the corresponding quantities $\tilde{F}_{j}^{-1}(1-\alpha)$ to obtain $\Sigma_{1:M_0}, \hdots, \Sigma_{M_1:M_0}$, the initial values for the next stage.
\State \textbf{Stage 1: Coarse localized optimizations}
\For{$j = 1$ \texttt{to} $j = M_1$} 
\State Generate a pseudorandom sample $Z_{1}, \hdots, Z_{N_1}$ from $P_{0, I_n}$. 
\State Define $\bar{F}_{j, \Sigma}(x) = N_1^{-1} \sum_{i = 1}^{N_1} \mathbf{1}_{(-\infty, x]}(T(\mu_0 + \Sigma^{1/2} Z_{i}))$ for $x \in \mathbb{R}$ and $\Sigma \in \mathfrak{C}_{Het}$.
\State Obtain $\Sigma^*_{j}$ by running a numerical optimization algorithm for the problem $\sup_{\Sigma \in \mathfrak{C}_{Het}} \bar{F}^{-1}_{j, \Sigma}(1-\alpha)$ initialized at $\Sigma_{j:M_0}$.
\EndFor
\State Rank the obtained~$\Sigma^*_{j}$ according to the value (from largest to smallest) of the corresponding $\bar{F}^{-1}_{j, \Sigma^*_{j}}(1-\alpha)$ to obtain $\Sigma^*_{1:M_1}, \hdots, \Sigma^*_{M_2:M_1}$, the initial values for the next stage.
\State \textbf{Stage 2: Refined localized optimization}
\For{$j = 1$ \texttt{to} $j = M_2$}
\State Generate a pseudorandom sample $Z_{1}, \hdots, Z_{N_2}$ from $P_{0, I_n}$. 
\State Define $\bar{\bar{F}}_{j, \Sigma}(x) = N_2^{-1} \sum_{i = 1}^{N_2} \mathbf{1}_{(-\infty, x]}(T(\mu_0 + \Sigma^{1/2}Z_{i}))$ for $x \in \mathbb{R}$ and $\Sigma \in \mathfrak{C}_{Het}$.
\State Obtain $\Sigma^{**}_{j}$ by running a numerical optimization algorithm for the problem $\sup_{\Sigma \in \mathfrak{C}_{Het}} \bar{\bar{F}}^{-1}_{j, \Sigma}(1-\alpha)$ initialized at $\Sigma^*_{j:M_1}$. 
\EndFor
\State \textbf{Return} $\max_{j = 1, \hdots, M_2} \bar{\bar{F}}^{-1}_{j, \Sigma^{**}_{j}}(1-\alpha)$. 
\end{algorithmic}
\end{algorithm}

\section{Appendix: Details concerning numerical computations in Section 
\protect\ref{numerical}\label{details}}

\subsection{Details concerning Section \protect\ref{sec:chibad} \label%
{app:numO2}}

To obtain Tables \ref{fig:lowbd} and \ref{fig:size}, for each of the test
statistics UC, HC0-HC4, UCR, HC0R-HC4R, we repeated the procedure summarized
in Algorithm \ref{alg:O2} below 15 times (recall that $n=25$, $R=(0,1)$, and 
$r=0$). Each time this algorithm returned a design matrix, the corresponding
size of the rejection region $\{T\geq C_{\chi ^{2},0.05}\}$ was obtained for
the specific test statistic used, as well as a corresponding lower bound for
the smallest size-controlling critical value. Then, we computed the maximum
out of the 15 lower bounds, which (for each test statistic) is reported in
Table \ref{fig:lowbd}. We also computed the maximum out of the 15 sizes,
which (for each test statistic) is reported in Table \ref{fig:size}. We also
did the same with the critical value $C_{\chi ^{2},0.05}$ replaced by the $%
95\%$-quantile of an $F_{1,n-k}$-distribution ($n-k=23$), the corresponding
results being reported in Table \ref{fig:size_2}.

In the description of Algorithm \ref{alg:O2}, the function $f(x)$ is an
abbreviation for $C^{\ast }=\max \{T(\mu _{0}+e_{i}(n)):i\in I_{1}(\mathfrak{%
M}_{0}^{lin})\}$, the lower bound for the size-controlling critical values
(cf.~Propositions \ref{rem_C*}, \ref{rem_C*_tilde}, and the attending
footnotes), with the $n\times 2$ design matrix $X$ given by an intercept $e$%
, say, as the first column and a regressor $x$ as the second one. Note that
computing $C^{\ast }$ necessitates the evaluation of the test statistic on a
finite set of elements of $\mathbb{R}^{n}$, and then determining the maximum
among the values obtained.\footnote{%
In the present context $\mathfrak{M}_{0}^{lin}$ is spanned by the intercept.
Thus, $I_{1}(\mathfrak{M}_{0}^{lin})=\{1,\ldots ,n\}$ holds since $n\geq 2$.
In general, to determine $I_{1}(\mathfrak{M}_{0}^{lin})$ numerically, the
algorithm implemented in the R-package \textbf{hrt }(\cite{hrt}) first
obtains a basis for $\mathfrak{M}_{0}^{lin}$, and then checks for every $%
i=1,\ldots ,n$ whether or not the rank of the matrix obtained by appending
the basis with $e_{i}(n)$ increases. This is done by a rank computation
analogous to the one described in the last-but-one paragraph of Appendix E.3
of~\cite{PPBoot}, using the same function \textquotedblleft
rank\textquotedblright\ referred to there, and with tolerance parameter~$%
10^{-8}$.} Concerning the evaluation of test statistics, the definition of
which depends on the invertibility of a covariance matrix estimator, we used
the same invertibility check as discussed in the second paragraph in
Appendix E.3 of \cite{PPBoot} with a tolerance parameter of $10^{-8}$. For $%
R=(0,1)$ and for each matrix $X$ returned by Algorithm \ref{alg:O2} all
relevant assumptions (i.e., the assumptions in the pertinent parts of
Theorems \ref{Hetero_Robust} and \ref{Hetero_Robust_tilde}, respectively)
have been checked numerically. 
\begin{algorithm}
\caption{Search procedure used for generating Tables~\ref{fig:lowbd},~\ref{fig:size}, and~\ref{fig:size_2}.} 
\label{alg:O2} 
\begin{algorithmic}[1]
\State Initialize~$x \leftarrow 0 \in \Rel^n$.
\For{$i = 1$ \texttt{to} $i = 5$}
\State Generate an~$n$-dimensional pseudo-random vector~$z$ of independent coordinates each from a log-standard normal distribution.
\State Run a \cite{nelder} algorithm initialized at~$z$ to maximize~$f$ over~$\Rel^n$ (with a maximal number of iterations of~$50$, and otherwise the default parameters in R's ``optim'' function) to obtain~$z^*$, say.
\If{$i = 1$, or $i \geq 2$ and $f(z^*) > f(x)$}
\State $x \leftarrow z^*$.
\EndIf
\If{$f(x) > 4$}
\State Go to line~\ref{alg:ret}.
\EndIf
\EndFor
\State Use Algorithm~\ref{alg:S1} to determine the size of the test for the test statistic under consideration for the design matrix~$(e, x)$ and based on either of the following two critical values: (i) $C_{\chi^2, 0.05}$ and (ii) the~$95\%$ quantile of an~$F_{1,n-k}$ distribution. \label{alg:ret}
\State \Return~$x$, $f(x)$, and the two sizes determined in the previous step.
\end{algorithmic}
\end{algorithm}

Algorithm \ref{alg:O2} uses Algorithm \ref{alg:S1} in determining the size
of a given test. We made the following choices concerning the parameters
required in Algorithm \ref{alg:S1} (and used default settings if not
mentioned otherwise):

\begin{enumerate}
\item The candidates in Stage 0 of Algorithm \ref{alg:S1} were determined by
combining the suggestions in Remarks \ref{rem:starging} and \ref%
{rem:initval_2}. That is, denoting $M_{p}=200~000$, we combined: (i)
sampling $M_{p}/4-1$ points from the unit simplex in $\mathbb{R}^{n}$, each
corresponding to the diagonal of a matrix in $\mathfrak{C}_{Het}$, and
sampling $3M_{p}/4+1$ points $\xi =(\xi _{1},\ldots ,\xi _{n})$, say,
analogously, each point $\xi $ giving rise to a diagonal of a matrix in $%
\mathfrak{C}_{Het}$ via $(\xi _{1}^{2},\ldots ,\xi
_{n}^{2})/\sum_{i=1}^{n}\xi _{i}^{2}$; (ii) trying all diagonal matrices
with a single dominant coordinate $0.9999$ and the other coordinates all
equal to $0.0001/(n-1)$, so that the trace equals $1$; (iii) $n^{-1}I_{n}$;
(iv) using a maximizer of the quadratic form described in Remark \ref%
{rem:starging}; and (v) using an additional initial value in case of a
\textquotedblleft large\textquotedblright\ critical value $C$ as described
in Remark \ref{rem:initval_2}, making use of the conventions discussed in
parentheses in that remark. This results in $M_{0}=M_{p}+n+2$ and possible
one more (in case $C$ is large) candidates for initial values.

\item $M_{1}$ was chosen as $500$, the optimization algorithm run in Stage~1
was a constrained \cite{nelder} algorithm (the default in R's
\textquotedblleft constrOptim\textquotedblright\ function), which was run
with a relative tolerance parameter of $10^{-2}$ and a maximal number of
iterations of $20n$.

\item $M_{2}$ was chosen as $1$, the optimization algorithm run in Stage~1
was a constrained \cite{nelder} algorithm (the default in R's
\textquotedblleft constrOptim\textquotedblright\ function), which was run
with a relative tolerance parameter of $10^{-3}$ and a maximal number of
iterations of $30n$.

\item DA (used by Algorithm \ref{alg:S1}) was run with the parameters
\textquotedblleft $\text{acc}=10^{-3}$\textquotedblright\ and
\textquotedblleft $\text{lim}=30000$\textquotedblright\ using the function
\textquotedblleft davies\textquotedblright\ of the package \textbf{%
CompQuadForm}.
\end{enumerate}

\subsection{Details concerning Section \protect\ref{sec:powernum} \label%
{app:CVcomppow}}

The smallest size-controlling critical values reported in Tables \ref{fig:CV}
and \ref{fig:CV2} in Section \ref{sec:powernum} were obtained by running
Algorithm \ref{alg:CV1} (with algorithm $\mathsf{A}$ given by Algorithm \ref%
{alg:S1} and a maximal number of $25$ iterations in the while loop) as
implemented in the R-package \textbf{hrt} (\cite{hrt}) version 1.0.0.
Concerning $\mathsf{A}$, the same input parameters as described in the
enumeration at the end of Appendix \ref{app:numO2} were used but with $M_p =
500~000$ (and with $n = 30$). Concerning Algorithm \ref{alg:CV1} we made the
following choices for the required inputs:

\begin{enumerate}
\item $C_{low}=\max (C^{\ast },C_{\hom })$ is determined as follows: $%
C_{\hom }$ is determined by a line-search algorithm (using R's uniroot
function and monotonicity of the rejection probabilities in the critical
value) with the rejection probabilities obtained from DA (in case $q=1$) or
via Monte Carlo, whereas $C^{\ast }$ is determined as described in Appendix %
\ref{app:numO2}. For more detail see \cite{hrt}.

\item $\epsilon $ was set to $10^{-3}$.
\end{enumerate}

For computing the power functions in Section \ref{sec:powernum}, we made use
of (\ref{eqn:qformgen}) with the matrices $A_{C}$ given in Lemma \ref%
{lem:qform} together with the implementation of the algorithm by \cite%
{davies} in the R-package \textbf{CompQuadForm} (\cite{Duchesne}) version
1.4.3 and with default parameters.

\subsection{Additional figures for Section~\protect\ref{sec:powernum}\label%
{app:addfig}}

\begin{figure}[tbp]
\centering
\includegraphics[width=\linewidth]{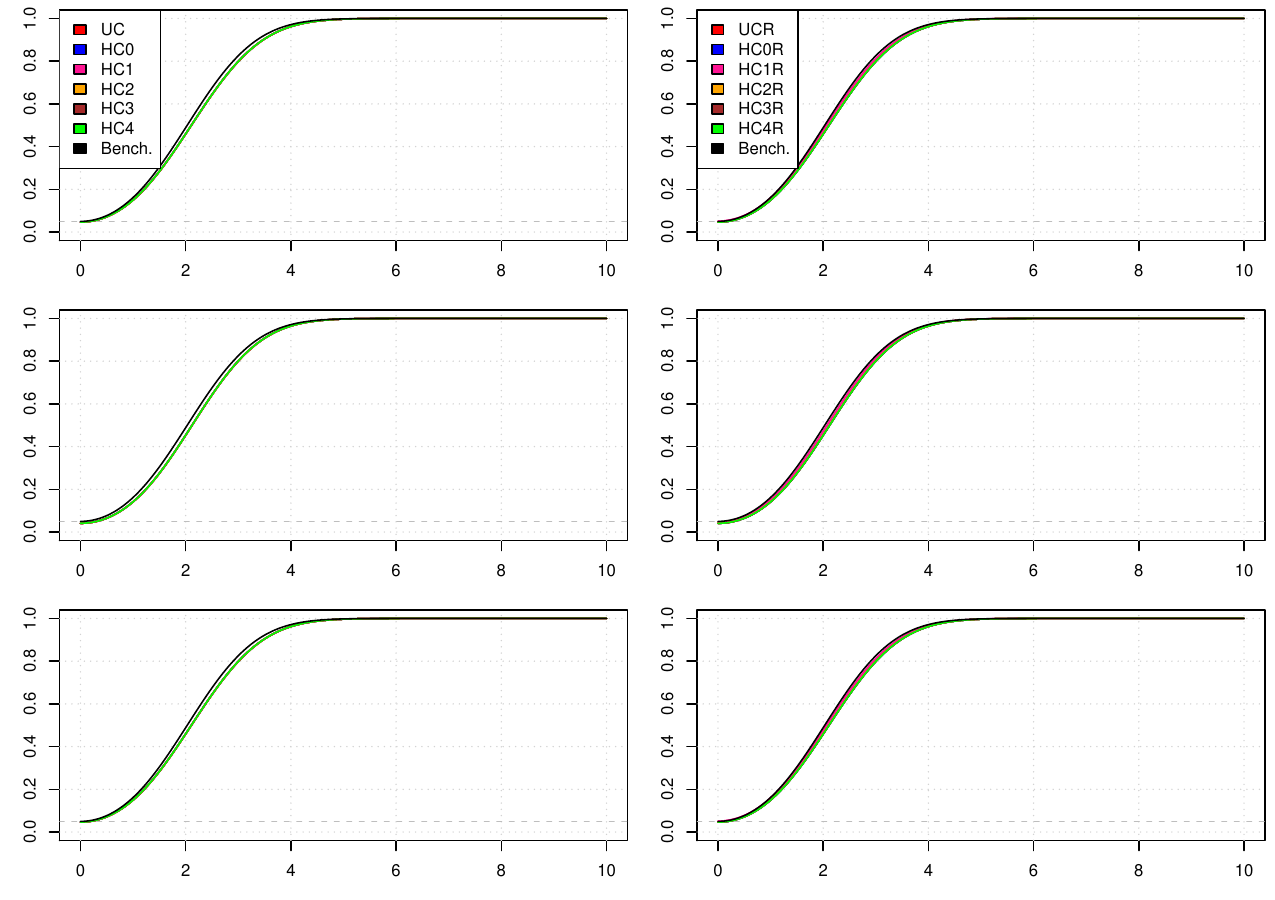}
\caption{Power functions for~$n_{1}=15$. Left column: tests based on
unrestricted residuals (cf.~legend). Right column: tests based on restricted
residuals (cf.~legend). The rows corresponds to~$\Sigma _{a}$ for~$a=1,5,9$
from top to bottom. The abscissa shows $\protect\delta $. In the left panel
the HC4-curve lies on top of the HC0--HC3-curves and the UC-curve. In the
right panel the HC4R-curve lies on top of the HC0R-HC3R-curves and the
UCR-curve. See the text for an explanation.}
\label{fig:bl5}
\end{figure}

The power functions for $n_{1}=15$ are given in Figure \ref{fig:bl5}. 

\section{Appendix: Comments on \protect\cite{Chuetal2021} and \protect\cite%
{Hansen2021}\label{app_G}}

In the special case of testing only one restriction (i.e., $q=1$), \cite%
{Chuetal2021} and \cite{Hansen2021} recently considered an interesting
alternative approach to obtain tests based on the test statistics $T_{Het}$
(for the commonly used choices of the weights $d_{i}$). Their suggestions
are based on the observation (cf. also Section \ref{sec:obsq1} above) that,
assuming Gaussianity of the errors, the null rejection probability of the
test that rejects if $T_{Het}$ exceeds a given critical value $C$ can be
rewritten as the probability that a quadratic form in Gaussian variables is
nonnegative, which can efficiently be determined numerically for any given $%
\Sigma \in \mathfrak{C}_{Het}$ by a number of methods.\footnote{%
A reader has pointed out that the results in \cite{phill93} could also be
developed into numerical approximations similar to the ones in \cite%
{Hansen2021}.} That is, if $\Sigma $ were known, one could use this
observation to numerically determine a critical value (an observation that
is also exploited by our algorithms in case $q=1$) or a p-value. Because $%
\Sigma $ is, however, not known, this approach is infeasible. One solution,
put forward in the present paper, is to work instead with a
\textquotedblleft worst-case\textquotedblright\ critical value, i.e., the
smallest critical value that controls size (if such a critical value
exists). In contrast, the idea in \cite{Chuetal2021} and \cite{Hansen2021}
to obtain a feasible test is a parametric bootstrap idea (cf.~their papers
for details):\footnote{%
A similar approach has already been put forward earlier by \cite{welch1938,
Welch1951} and \cite{Satterth}.~} (i) replace $\Sigma $ by an estimate $\hat{%
\Sigma}$; e.g., 
\begin{equation}
\limfunc{diag}\left( d_{1}\hat{u}_{1}^{2}\left( y\right) ,\ldots ,d_{n}\hat{u%
}_{n}^{2}\left( y\right) \right)  \label{eqn:covest}
\end{equation}%
based on typical choices of $d_{i}$; (ii) numerically determine a critical
value (or p-value) from the cdf of the test statistic acting as if $\Sigma =%
\hat{\Sigma}$ (e.g., as outlined above); and (iii) reject the null
hypothesis if the observed test statistic exceeds the so-computed critical
value (or, equivalently, if the corresponding p-value obtained is less than
the desired significance level). Note that the critical value in (ii)
depends on the data $Y$ through $\hat{\Sigma}$ (and is thus data-dependent
in this sense).

No theoretical guarantees concerning the size of the tests proposed in \cite%
{Chuetal2021} and \cite{Hansen2021} are given in these papers. Numerical
results in both papers suggest that these parametric bootstrap tests can
work well for certain design matrices $X$ and hypotheses $(R,r)$, but the
authors also document some situations where the tests are considerably
oversized. Hence, these tests are not valid, in general, which is in
contrast to the procedure we suggest in the present paper. That a parametric
bootstrap approach does not deliver size control is in line with results in 
\cite{Loh1985} (see also \cite{LP2017}) showing that under appropriate
conditions parametric bootstrap procedures are oversized. It is also in line
with a large body of literature on size distortions of (other)
bootstrap-based tests for the testing problem under consideration, cf.
Section \ref{Intro} and \cite{PPBoot}. As an aside we note that any \emph{%
valid }data-dependent critical value, i.e., one that leads to a test with
correct size (which is \emph{not} the case for the proposals in \cite%
{Chuetal2021} and \cite{Hansen2021}), must exceed the smallest
size-controlling critical value with positive probability (or must be equal
to the smallest size-controlling critical value with probability $1$).
Hence, a \emph{valid} data-dependent critical value cannot always be smaller
than the smallest size-controlling critical values, an observation that
seems to have gone unnoticed in the discussion of the present article given
in the introduction of \cite{Hansen2021} (a discussion that also overlooks
that one needs to take the square root of our critical values and lower
bounds when discussing them in the context of the corresponding $t$%
-statistics).

To demonstrate further that the parametric bootstrap tests in \cite%
{Chuetal2021} and \cite{Hansen2021} can be considerably oversized, we now
report some numerical results for these tests. In particular, we report null
rejection probabilities for a selection of points in the null hypothesis
(i.e., for a selection of $\Sigma $'s) and demonstrate that procedures
suggested by \cite{Chuetal2021} and \cite{Hansen2021} are not valid in the
sense that these null rejection probabilities are considerably larger than
the nominal significance level $\alpha =0.05$ that is being used. Note that
what we report are lower bounds for the size of the procedures investigated,
which can even be larger, i.e., the overrejection problem can, in fact, be
even more serious than what is seen in the tables below. Throughout, we
study the following procedures

\begin{enumerate}
\item[C:] the procedure in \cite{Chuetal2021}, when $T_{Het}$ based on
HC0-HC4 weights, respectively, is combined with the estimator $\hat{\Sigma}$
in (\ref{eqn:covest}) using the same weights as in the construction of the
test statistic;

\item[C3:] the procedure as above, but where the estimator $\hat{\Sigma}$ in
(\ref{eqn:covest}) always makes use of the HC3 weights;

\item[H:] the procedure in \cite{Hansen2021} when $T_{Het}$ is based on
HC0-HC4 weights, respectively, and where for $\hat{\Sigma}$ the estimator
suggested in Section 7 of \cite{Hansen2021} is used.
\end{enumerate}

We note here that procedure C3 is not considered in \cite{Chuetal2021}; we
include it, because $\hat{\Sigma}$ based on HC3 weights can be expected to
perform better than if, e.g., HC0 weights are used. We also point out that 
\cite{Hansen2021} only considers $T_{Het}$ based on HC0-HC3 weights, but 
\emph{not} on HC4; we also report rejection probabilities for the latter
choice, because, in the examples we consider, it actually works better in
terms of size than the choices considered in \cite{Hansen2021}.

Our implementations of the procedures in \cite{Chuetal2021} and \cite%
{Hansen2021} rely on the algorithm in \cite{davies} (cf.~Section \ref%
{sec:obsq1}) to decide whether or not to reject (i.e., in Step (ii) of the
description of that approach given further above in this section). To
compute the rejection probabilities for the tests we used a Monte Carlo
sample of size 100.000 for each of them. The nominal significance level used
is $\alpha =0.05$ throughout.

We consider three testing problems: The testing problems considered in
Sections \ref{sec:powtwogroups} and \ref{sec:powhostX}, as well as an
additional one. Note that in all these examples the test statistics are size
controllable and thus our test procedures based on smallest size-controlling
critical values are applicable. [For Examples \ref{ex:G1} and \ref{ex:G2}
this has already been discussed in Sections \ref{sec:powtwogroups} and \ref%
{sec:powhostX}, respectively. For Example \ref{ex:G3} validity of Assumption %
\ref{R_and_X} is obvious while condition (\ref{non-incl_Het}) we have
verified numerically.]

\begin{example}
\label{ex:G1}\emph{(Comparing the means of two groups) }We here consider the
same testing problem and setting (same $n$, $n_{1}$, $n_{2}$, $\alpha $) as
in Section \ref{sec:powtwogroups}. Table \ref{tab:sizeCH2G} below shows the
null rejection probabilities for the procedures C, C3, and H for the case $%
n_{1}=3$ and $a=3$ (i.e., for $\Sigma =\Sigma _{3}$ defined in Section \ref%
{sec:powtwogroups}). 
\begin{table}[th]
\centering
\begin{tabular}{l||lllll}
\hline
& HC0 & HC1 & HC2 & HC3 & HC4 \\ \hline
C & 0.14 & 0.14 & 0.12 & 0.11 & 0.10 \\ 
C3 & 0.12 & 0.12 & 0.12 & 0.11 & 0.10 \\ 
H & 0.08 & 0.08 & 0.08 & 0.07 & 0.06 \\ \hline
\end{tabular}%
\caption{Null-rejection probabilities of the procedures C, C3, and H for
comparing the means of two groups when~$n_{1}=3$ and $a=3$.}
\label{tab:sizeCH2G}
\end{table}
We see from that table that the procedures suggested in \cite{Chuetal2021},
i.e., procedures C, as well as the modification C3 are all considerably
oversized (i.e., show rejection probabilities greater or equal to $2\alpha $%
). The methods using the idea in \cite{Hansen2021} (including the case using
HC4 weights not considered in \cite{Hansen2021}) are slightly oversized in
this example.
\end{example}

\begin{example}
\label{ex:G2}\emph{(High-leverage design) }We here consider the same testing
problem and setting (same $n$, $\alpha $, $X$) as in Section \ref%
{sec:powhostX}. Table \ref{tab:sizeCHHL} shows the null rejection
probabilities for the procedures C, C3, and H for the case $a=1$ (i.e., for $%
\Sigma =\Sigma _{1}^{\ast }$ defined in Section \ref{sec:powhostX}). 
\begin{table}[th]
\centering
\begin{tabular}{l||lllll}
\hline
& HC0 & HC1 & HC2 & HC3 & HC4 \\ \hline
C & 0.65 & 0.65 & 0.30 & 0.14 & 0.09 \\ 
C3 & 0.18 & 0.18 & 0.17 & 0.14 & 0.10 \\ 
H & 0.16 & 0.16 & 0.15 & 0.12 & 0.08 \\ \hline
\end{tabular}%
\caption{Null-rejection probabilities of the procedures C, C3, and H for the
high-leverage design matrix when $a=1$.}
\label{tab:sizeCHHL}
\end{table}
The methods based on the approach in \cite{Chuetal2021} i.e., procedures C,
as well as the modification C3 are all considerably oversized also in this
example. The methods using the idea in \cite{Hansen2021} are now also
considerably oversized. The test using the HC4 estimator (which was not
considered in \cite{Hansen2021}) performs somewhat better and has a null
rejection probability that exceeds the nominal significance level $\alpha
=0.05$ by a factor of $1.6$.
\end{example}

The tables in the two preceding examples already show that the tests
proposed by \cite{Chuetal2021} and \cite{Hansen2021} can be considerably
oversized. Note that the overrejection problem potentially is even more
serious than what is seen from the tables as we have not searched over the
space of $\Sigma $ matrices, i.e., we have not reported size but only the
null rejection probability at a particular value of $\Sigma $. Also, we have
not made any attempt to search for design matrices $X$ where overrejection
is even more pronounced, but have only used design matrices from Sections %
\ref{sec:powtwogroups} and \ref{sec:powhostX}.

We have seen in the preceding examples that pairing the method in \cite%
{Hansen2021} with a HC4 based $T_{Het}$ statistic performs more reasonably
in these settings (it also is oversized, but less so). The question then
arises whether there is some hope that this generalizes to other settings.
The next example shows that this is unfortunately not the case.

\begin{example}
\label{ex:G3}We consider the same model and null hypothesis as in Section %
\ref{sec:powhostX} except that the regressor $x$ ($x\in \mathbb{R}^{n}$, $%
n=30$) is different. Its entries $x_{i}$ can be found plotted (against the
index $i$) in Figure \ref{fig:Xsov}. For this scenario one can prove (using
similar arguments as in \cite{PPBoot}) that the size of the test obtained
from pairing \cite{Hansen2021}'s method with a HC4 based $T_{Het}$ statistic
actually equals $1$. We do not give the details here but rather compute the
null rejection probability of this test for $\Sigma $ equal to the diagonal
matrix with $0.999$ at the $21^{st}$ entry and the other diagonal entries
constant so that the diagonal sums up to one. We used a Monte Carlo
simulation (with 100.000 replications) and obtained a null rejection
probability of $0.28$, which is more than the five-fold nominal significance
level. 
\begin{figure}[tbp]
\centering\includegraphics[width=0.5\linewidth]{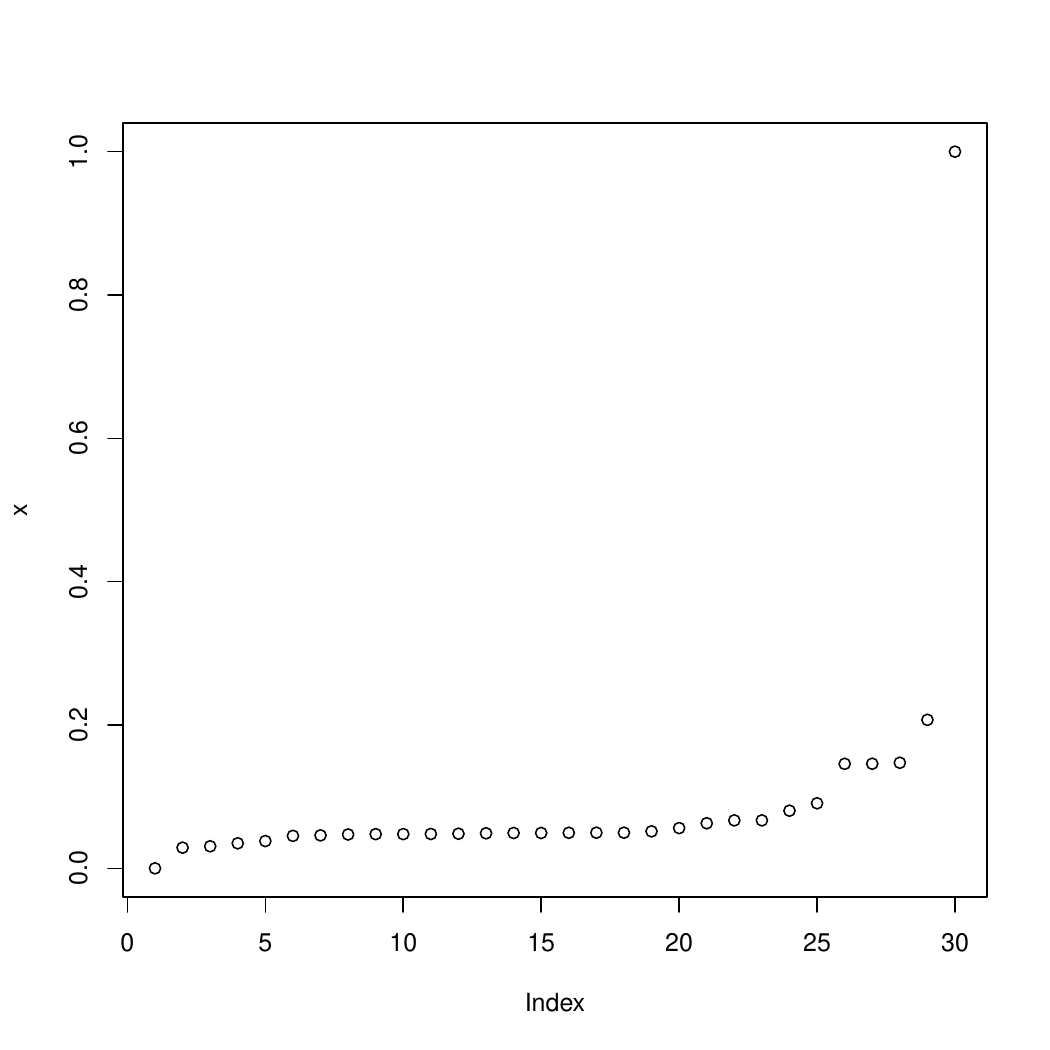}
\caption{Regressor used in Example \protect\ref{ex:G3}.}
\label{fig:Xsov}
\end{figure}
\end{example}

\bibliographystyle{ims}
\bibliography{refs}

\def\cprime{$'$}
\begin{thebibliography}{53}
\expandafter\ifx\csname natexlab\endcsname\relax\def\natexlab#1{#1}\fi
\expandafter\ifx\csname url\endcsname\relax
  \def\url#1{\texttt{#1}}\fi
\expandafter\ifx\csname urlprefix\endcsname\relax\def\urlprefix{URL }\fi
\providecommand{\eprint}[2][]{\url{#2}}

\bibitem[{Bakirov and Sz{\'e}kely(2005)}]{BakiSzek2005}
\textsc{Bakirov, N.} and \textsc{Sz{\'e}kely, G.} (2005).
\newblock Student's t-test for {G}aussian scale mixtures.
\newblock \textit{Zapiski Nauchnyh Seminarov POMI}, \textbf{328} 5--19.

\bibitem[{Bakirov(1998)}]{Bakirov98}
\textsc{Bakirov, N.~K.} (1998).
\newblock Nonhomogeneous samples in the {B}ehrens-{F}isher problem.
\newblock \textit{Journal of Mathematical Sciences (New York)}, \textbf{89}
  1460--1467.

\bibitem[{Bell and McCaffrey(2002)}]{BellMcCa}
\textsc{Bell, R.~M.} and \textsc{McCaffrey, D.} (2002).
\newblock Bias reduction in standard errors for linear regression with
  multi-stage samples.
\newblock \textit{Survey Methodology}, \textbf{28} 169--181.

\bibitem[{Belloni and Didier(2008)}]{BelloniDidier2008}
\textsc{Belloni, A.} and \textsc{Didier, G.} (2008).
\newblock On the {B}ehrens-{F}isher problem: a globally convergent algorithm
  and a finite-sample study of the {W}ald, {LR} and {LM} tests.
\newblock \textit{The Annals of Statistics}, \textbf{36} 2377--2408.

\bibitem[{Cattaneo et~al.(2018)Cattaneo, Jansson and Newey}]{Catt}
\textsc{Cattaneo, M.~D.}, \textsc{Jansson, M.} and \textsc{Newey, W.~K.}
  (2018).
\newblock Inference in linear regression models with many covariates and
  heteroscedasticity.
\newblock \textit{Journal of the American Statistical Association},
  \textbf{113} 1350--1361.

\bibitem[{Chesher and Jewitt(1987)}]{CheshJewitt1987}
\textsc{Chesher, A.} and \textsc{Jewitt, I.} (1987).
\newblock The bias of a heteroskedasticity consistent covariance matrix
  estimator.
\newblock \textit{Econometrica}, \textbf{55} 1217--1222.

\bibitem[{Chesher(1989)}]{Chesh_1989}
\textsc{Chesher, A.~D.} (1989).
\newblock H\'{a}jek inequalities, measures of leverage, and the size of
  heteroskedasticity robust {W}ald tests.
\newblock \textit{Econometrica}, \textbf{57} 971--977.

\bibitem[{Chesher and Austin(1991)}]{CheshAust_1991}
\textsc{Chesher, A.~D.} and \textsc{Austin, G.} (1991).
\newblock The finite-sample distributions of heteroskedasticity robust {W}ald
  statistics.
\newblock \textit{Journal of Econometrics}, \textbf{47} 153--173.

\bibitem[{Chu et~al.(2021)Chu, Lee, Ullah and Xu}]{Chuetal2021}
\textsc{Chu, J.}, \textsc{Lee, T.-H.}, \textsc{Ullah, A.} and \textsc{Xu, H.}
  (2021).
\newblock Exact distribution of the {F}-statistic under heteroskedasticity of
  unknown form for improved inference.
\newblock \textit{Journal of Statistical Computation and Simulation},
  \textbf{91} 1782--1801.

\bibitem[{Cragg(1983)}]{Cragg_1983}
\textsc{Cragg, J.~G.} (1983).
\newblock More efficient estimation in the presence of heteroscedasticity of
  unknown form.
\newblock \textit{Econometrica}, \textbf{51} 751--763.

\bibitem[{Cragg(1992)}]{Cragg_1992}
\textsc{Cragg, J.~G.} (1992).
\newblock Quasi-{A}itken estimation for heteroscedasticity of unknown form.
\newblock \textit{Journal of Econometrics}, \textbf{54} 179--201.

\bibitem[{Cribari-Neto(2004)}]{Crib2004}
\textsc{Cribari-Neto, F.} (2004).
\newblock Asymptotic inference under heteroskedasticity of unknown form.
\newblock \textit{Computational Statistics \& Data Analysis}, \textbf{45} 215
  -- 233.

\bibitem[{Davidson and Flachaire(2008)}]{DavidsFlach2008}
\textsc{Davidson, R.} and \textsc{Flachaire, E.} (2008).
\newblock The wild bootstrap, tamed at last.
\newblock \textit{Journal of Econometrics}, \textbf{146} 162--169.

\bibitem[{Davidson and MacKinnon(1985)}]{DavidsonMacKinnon1985}
\textsc{Davidson, R.} and \textsc{MacKinnon, J.~G.} (1985).
\newblock Heteroskedasticity-robust tests in regressions directions.
\newblock \textit{Minist\`ere de l'\'{E}conomie et des Finances. Institut
  National de la Statistique et des \'{E}tudes \'{E}conomiques. Annales}
  183--218.

\bibitem[{Davies(1980)}]{davies}
\textsc{Davies, R.~B.} (1980).
\newblock Algorithm {A}{S} 155: The distribution of a linear combination of
  $\chi^2$ random variables.
\newblock \textit{Journal of the Royal Statistical Society. Series C (Applied
  Statistics)}, \textbf{29} 323--333.

\bibitem[{DiCiccio et~al.(2019)DiCiccio, Romano and
  Wolf}]{DiCiccio_Romao_Wolf_2019}
\textsc{DiCiccio, C.~J.}, \textsc{Romano, J.~P.} and \textsc{Wolf, M.} (2019).
\newblock Improving weighted least squares inference.
\newblock \textit{Econometrics and Statistics}, \textbf{10} 96--119.

\bibitem[{Duchesne and {de Micheaux}(2010)}]{Duchesne}
\textsc{Duchesne, P.} and \textsc{{de Micheaux}, P.~L.} (2010).
\newblock Computing the distribution of quadratic forms: Further comparisons
  between the {L}iu-{T}ang-{Z}hang approximation and exact methods.
\newblock \textit{Computational Statistics and Data Analysis}, \textbf{54}
  858--862.

\bibitem[{Eicker(1963)}]{E63}
\textsc{Eicker, F.} (1963).
\newblock Asymptotic normality and consistency of the least squares estimators
  for families of linear regressions.
\newblock \textit{Annals of Mathematical Statistics}, \textbf{34} 447--456.

\bibitem[{Eicker(1967)}]{E67}
\textsc{Eicker, F.} (1967).
\newblock Limit theorems for regressions with unequal and dependent errors.
\newblock In \textit{Proc. {F}ifth {B}erkeley {S}ympos. {M}ath. {S}tatist. and
  {P}robability ({B}erkeley, {C}alif., 1965/66)}. Univ. California Press,
  Berkeley, Calif., Vol. I: Statistics, pp. 59--82.

\bibitem[{Flachaire(2005)}]{Flachaire_2005b}
\textsc{Flachaire, E.} (2005).
\newblock More efficient tests robust to heteroskedasticity of unknown form.
\newblock \textit{Econometric Reviews}, \textbf{24} 219--241.

\bibitem[{Godfrey(2006)}]{Godfrey2006}
\textsc{Godfrey, L.~G.} (2006).
\newblock Tests for regression models with heteroskedasticity of unknown form.
\newblock \textit{Computational Statistics \& Data Analysis}, \textbf{50}
  2715--2733.

\bibitem[{Hansen(2021)}]{Hansen2021}
\textsc{Hansen, B.} (2021).
\newblock The exact distribution of the {W}hite t-ratio.

\bibitem[{Hinkley(1977)}]{H77}
\textsc{Hinkley, D.~V.} (1977).
\newblock Jackknifing in unbalanced situations.
\newblock \textit{Technometrics}, \textbf{19} 285--292.

\bibitem[{Ibragimov and M{\"u}ller(2010)}]{IbragMuell2010}
\textsc{Ibragimov, R.} and \textsc{M{\"u}ller, U.~K.} (2010).
\newblock t-statistic based correlation and heterogeneity robust inference.
\newblock \textit{Journal of Business and Economic Statistics}, \textbf{28}
  453--468.

\bibitem[{Ibragimov and M{\"u}ller(2016)}]{IbragMuell2016}
\textsc{Ibragimov, R.} and \textsc{M{\"u}ller, U.~K.} (2016).
\newblock Inference with few heterogeneous clusters.
\newblock \textit{The Review of Economics and Statistics}, \textbf{98} 83--96.

\bibitem[{Imbens and Koles{\'a}r(2016)}]{Imbkoles2016}
\textsc{Imbens, G.~W.} and \textsc{Koles{\'a}r, M.} (2016).
\newblock Robust standard errors in small samples: Some practical advice.
\newblock \textit{The Review of Economics and Statistics}, \textbf{98}
  701--712.

\bibitem[{Kim and Cohen(1998)}]{KimCohen}
\textsc{Kim, S.-H.} and \textsc{Cohen, A.~S.} (1998).
\newblock On the {B}ehrens-{F}isher problem: A review.
\newblock \textit{Journal of Educational and Behavioral Statistics},
  \textbf{23} 356--377.

\bibitem[{Koles{\'a}r(2019)}]{dfadjust}
\textsc{Koles{\'a}r, M.} (2019).
\newblock \textit{dfadjust: Degrees of Freedom Adjustment for Robust Standard
  Errors}.
\newblock R package version 1.0.1,
  \urlprefix\url{https://CRAN.R-project.org/package=dfadjust}.

\bibitem[{Leeb and P\"{o}tscher(2017)}]{LP2017}
\textsc{Leeb, H.} and \textsc{P\"{o}tscher, B.~M.} (2017).
\newblock Testing in the presence of nuisance parameters: some comments on
  tests post-model-selection and random critical values.
\newblock In \textit{Big and complex data analysis}. Contrib. Stat., Springer,
  Cham, 69--82.

\bibitem[{Lehmann and Romano(2005)}]{LR05}
\textsc{Lehmann, E.~L.} and \textsc{Romano, J.~P.} (2005).
\newblock \textit{Testing Statistical Hypotheses}.
\newblock 3rd ed. Springer Texts in Statistics, Springer, New York.

\bibitem[{Lin and Chou(2018)}]{Lin_Chou_2018}
\textsc{Lin, E.~S.} and \textsc{Chou, T.-S.} (2018).
\newblock Finite-sample refinement of {GMM} approach to nonlinear models under
  heteroskedasticity of unknown form.
\newblock \textit{Econometric Reviews}, \textbf{37} 1--28.

\bibitem[{Loh(1985)}]{Loh1985}
\textsc{Loh, W.-Y.} (1985).
\newblock A new method for testing separate families of hypotheses.
\newblock \textit{Journal of the American Statistical Association}, \textbf{80}
  362--368.

\bibitem[{Long and Ervin(2000)}]{LE2000}
\textsc{Long, J.~S.} and \textsc{Ervin, L.~H.} (2000).
\newblock Using heteroscedasticity consistent standard errors in the linear
  regression model.
\newblock \textit{The American Statistician}, \textbf{54} 217--224.

\bibitem[{MacKinnon(2013)}]{Mackinnon2013}
\textsc{MacKinnon, J.~G.} (2013).
\newblock Thirty years of heteroskedasticity-robust inference.
\newblock In \textit{Recent Advances and Future Directions in Causality,
  Prediction, and Specification Analysis} (X.~Chen and N.~R.~E. Swanson, eds.).
  Springer, 437--462.

\bibitem[{MacKinnon and White(1985)}]{MacW85}
\textsc{MacKinnon, J.~G.} and \textsc{White, H.} (1985).
\newblock Some heteroskedasticity-consistent covariance matrix estimators with
  improved finite sample properties.
\newblock \textit{Journal of Econometrics}, \textbf{29} 305 -- 325.

\bibitem[{Mickey and Brown(1966)}]{MickeyBrown1966}
\textsc{Mickey, M.~R.} and \textsc{Brown, M.~B.} (1966).
\newblock Bounds on the distribution functions of the {B}ehrens-{F}isher
  statistic.
\newblock \textit{Annals of Mathematical Statistics}, \textbf{37} 639--642.

\bibitem[{Nelder and Mead(1965)}]{nelder}
\textsc{Nelder, J.~A.} and \textsc{Mead, R.} (1965).
\newblock A simplex method for function minimization.
\newblock \textit{The Computer Journal}, \textbf{7} 308--313.

\bibitem[{Phillips(1993)}]{phill93}
\textsc{Phillips, P.~C.} (1993).
\newblock Operational algebra and regression t-tests.
\newblock In \textit{Models, Methods and Applications of Econometrics: Essays
  in Honor of A.R. Bergstrom} (P.~C. Phillips, ed.). Oxford: Basil Blackwell,
  140--152.

\bibitem[{P{\"o}tscher and Preinerstorfer(2018)}]{PP3}
\textsc{P{\"o}tscher, B.~M.} and \textsc{Preinerstorfer, D.} (2018).
\newblock Controlling the size of autocorrelation robust tests.
\newblock \textit{Journal of Econometrics}, \textbf{207} 406--431.

\bibitem[{P{\"o}tscher and Preinerstorfer(2019)}]{PP4}
\textsc{P{\"o}tscher, B.~M.} and \textsc{Preinerstorfer, D.} (2019).
\newblock Further results on size and power of heteroskedasticity and
  autocorrelation robust tests, with an application to trend testing.
\newblock \textit{Electronic Journal of Statistics}, \textbf{13} 3893--3942.

\bibitem[{P{\"o}tscher and Preinerstorfer(2022)}]{PPBoot}
\textsc{P{\"o}tscher, B.~M.} and \textsc{Preinerstorfer, D.} (2022).
\newblock How reliable are bootstrap-based heteroskedasticity robust tests?
\newblock \textit{Econometric Theory,} first published online 22 April 2022.
  DOI.10.1017/S0266466622000184.

\bibitem[{Preinerstorfer(2021)}]{hrt}
\textsc{Preinerstorfer, D.} (2021).
\newblock \textit{hrt: Heteroskedasticity Robust Testing}.
\newblock R package version 1.0.0.

\bibitem[{Preinerstorfer and P{\"o}tscher(2016)}]{PP2016}
\textsc{Preinerstorfer, D.} and \textsc{P{\"o}tscher, B.~M.} (2016).
\newblock On size and power of heteroskedasticity and autocorrelation robust
  tests.
\newblock \textit{Econometric Theory}, \textbf{32} 261--358.

\bibitem[{Robinson(1979)}]{RO1979}
\textsc{Robinson, G.} (1979).
\newblock Conditional properties of statistical procedures.
\newblock \textit{Annals of Statistics}, \textbf{7} 742--755.

\bibitem[{Romano and Wolf(2017)}]{RomanoWolf2017}
\textsc{Romano, J.~P.} and \textsc{Wolf, M.} (2017).
\newblock Resurrecting weighted least squares.
\newblock \textit{Journal of Econometrics}, \textbf{197} 1--19.

\bibitem[{Rothenberg(1988)}]{Rothenberg1988}
\textsc{Rothenberg, T.~J.} (1988).
\newblock Approximate power functions for some robust tests of regression
  coefficients.
\newblock \textit{Econometrica}, \textbf{56} 997--1019.

\bibitem[{Ruben(2002)}]{Ruben2002}
\textsc{Ruben, H.} (2002).
\newblock A simple conservative and robust solution of the {B}ehrens-{F}isher
  problem.
\newblock \textit{Sankhy\={a} Ser. A}, \textbf{64} 139--155.

\bibitem[{Satterthwaite(1946)}]{Satterth}
\textsc{Satterthwaite, F.~E.} (1946).
\newblock An approximate distribution of estimates of variance components.
\newblock \textit{Biometrics Bulletin}, \textbf{2} 110--114.

\bibitem[{Welch(1938)}]{welch1938}
\textsc{Welch, B.~L.} (1938).
\newblock The significance of the difference between two means when the
  population variances are unequal.
\newblock \textit{Biometrika}, \textbf{29} 350--362.

\bibitem[{Welch(1951)}]{Welch1951}
\textsc{Welch, B.~L.} (1951).
\newblock On the comparison of several mean values: an alternative approach.
\newblock \textit{Biometrika}, \textbf{38} 330--336.

\bibitem[{White(1980)}]{W80}
\textsc{White, H.} (1980).
\newblock A heteroskedasticity-consistent covariance matrix estimator and a
  direct test for heteroskedasticity.
\newblock \textit{Econometrica}, \textbf{48} 817--838.

\bibitem[{Wooldridge(2010)}]{Wooldridge2010}
\textsc{Wooldridge, J.~M.} (2010).
\newblock \textit{Econometric Analysis of Cross Section and Panel Data}.
\newblock 2nd ed. MIT Press, Cambridge, MA.

\bibitem[{Wooldridge(2012)}]{Wooldridge2012}
\textsc{Wooldridge, J.~M.} (2012).
\newblock \textit{Introductory Econometrics}.
\newblock 5th ed. South-Western, OH.

\end{thebibliography}

\end{document}